\documentclass[a4paper,11pt,titlepage]{article}

\usepackage{amsmath,amssymb}
\usepackage{amsthm}
\usepackage{palatino}
\usepackage{amscd}
\usepackage{graphicx}
\usepackage{epsfig}
\usepackage{verbatim}

\title{\textbf{Irreducibility of the Lawrence-Krammer\\
representation\\of the BMW algebra of type $\mathbf{A_{n-1}}$}}

\author{Claire Levaillant}%\\\\updated version}
\date{September 2008}

\newcommand{\m}{\mu}
\newcommand{\n}{\nu}
\newcommand{\Q}{\mathbb{Q}}
\newcommand{\R}{\mathbb{R}}
\newcommand{\be}{\beta}
\newcommand{\al}{\alpha}
\newcommand{\xb}{x_{\beta}}
\newcommand{\e}{\epsilon}
\newcommand{\lra}{\longrightarrow}
\newcommand{\ra}{\rightarrow}
\newcommand{\Ra}{\Rightarrow}
\newcommand{\la}{\lambda}

\newcommand{\unsurr}{\frac{1}{r}}
\newcommand{\unsur}{\frac{1}}
\newcommand{\x}{\Upsilon}
\newcommand{\ulra}{\underleftrightarrow}

\newcommand{\p}{\prec}
\newcommand{\q}{\succ}
\newcommand{\ps}{(\be|\al_i)}
\newcommand{\nts}{\negthickspace}
\newcommand{\ali}{\al_i}
\newcommand{\xali}{x_{\al_i}}
\newcommand{\xalj}{x_{\al_j}}
\newcommand{\balip}{x_{\be+\al_i}}
\newcommand{\balim}{x_{\be-\al_i}}
\newcommand{\baljm}{x_{\be-\al_j}}
\newcommand{\baljp}{x_{\be+\al_j}}
\newcommand{\baljip}{x_{\be+\al_i+\al_j}}

\newcommand{\xalone}{x_{\al_1}}

\newcommand{\T}{\mathcal{T}}
\newcommand{\U}{\mathcal{U}}
\newcommand{\V}{\mathcal{V}}
\newcommand{\IH}{\mathcal{H}}

\newcommand{\nui}{\n_i}
\newcommand{\lb}{\lbrace}
\newcommand{\rb}{\rbrace}
\newcommand{\noin}{\noindent}
\newcommand{\bigp}{\Big(}
\newcommand{\bigpd}{\Big)}
\newcommand{\unsurl}{\frac{1}{l}}
\newcommand{\iha}{\IH_{F,r^2}(n)}
\newcommand{\ih}{\IH_{F,r^2}}
\newcommand{\D}{\mathcal{D}}
\newcommand{\ti}{\times}
\newcommand{\W}{\mathcal{W}}
\newcommand{\di}{\text{dim}}
\newcommand{\nn}{\n^{(n)}}

\newcommand{\da}{\downarrow}
\newcommand{\Lra}{\Leftrightarrow}
\newcommand{\G}{\mathcal{G}}
\newcommand{\dw}{\text{dim}\,\W}
\newcommand{\dna}{\downarrow}

\newcommand{\B}{\mathcal{B}}
\newcommand{\X}{\mathcal{X}}
\newcommand{\Y}{\mathcal{Y}}
\newcommand{\divon}{\divideontimes}
\newcommand{\del}{\delta}
\newcommand{\tr}{\triangle}

\newcommand{\eg}{&=&}
\newcommand{\g}{\gamma}
\newcommand{\ds}{\bigoplus}
\newcommand{\cil}{\frac{n(n-3)}{2}}
\newcommand{\dbw}{\frac{(n-1)(n-2)}{2}}
\newcommand{\chl}{\frac{n(n-1)}{2}}
\newcommand{\eq}{(\divon)_}

\newcommand{\pal}{\V_{1,\unsur{r^{2n-3}}}}
\newcommand{\mfl}{\unsur{r^{2n-3}}}
\newcommand{\jca}{\unsur{r^{n-3}}}
\newcommand{\ze}{\mathbf{0}&}
\newcommand{\C}{\mathcal{C}}
\newcommand{\Rrond}{\mathcal{R}}
\newcommand{\sm}{\text{submatrix}}
\newcommand{\dete}{\text{det}}
\newcommand{\sq}{\square}
\newcommand{\overlra}{\overleftrightarrow}
\newcommand{\matp}{\mathcal{P}}

\newcommand{\ccc}{\circledcirc}
\newcommand{\Com}{\mathbb{C}}
\newcommand{\Z}{\mathbb{Z}}

\begin{document}
\maketitle \pagestyle{empty}
\noin{\huge{\textbf{Acknowledgements}}}\\\\\\
%\doublespacing
\indent I am extremely happy to thank my advisor David Wales for his
incommensurably generous guidance throughout my graduate studies. I
am indebted to him for committing to the difficult task of having me
as a graduate student and for suggesting a beautiful problem that
kept me excited and filled with wonder from the beginning until the
end of my graduate studies at Caltech. I express my deep gratitude
to him for meeting with me on a weekly basis.\\
\indent My thanks also go to Richard Wilson, Dinakar Ramakrishnan
and Cheng Yeah Ku for serving on my defense committee and to Robert
Guralnick for his helpful referral to James'paper.\\
\indent I would like to specially thank Daniel Katz and Emmanuel
Bresson for their support before I started my research at Caltech
and my friend and mathematician Magda Sebestean for her constant
friendship since
I have known her.\\
\indent I thank Pietro Perona for making me come to Caltech through
the Summer Undergraduate Research Fellowship Program and for making
me pursue captivating research in neuroscience under his guidance. I
also thank him for his constant caring eye and thoughtful advice.\\
This thesis could not have been undertaken without Albert and
Marguerite Ramond's funding through the "Albert and Marguerite
Ramond Scholarship". I am indebted to them and very much wish I had
known them. My thanks also go to the graduate office for always
being friendly and supportive, and especially to Dean Hoffmann for
being a humane dean, and to Natalie Gilmore and Rosa Carrasco for
their always gracious and very nice support to the students. I thank
the mathematics department for hosting me during these five years
and the mathematics department staff for their assistance. My thanks
go to Seraj for dayly support, to Pamela Fong for much
friendliness, concern and help and to Professor Tom Apostol and his wife Jane for their kindness.\\
\indent I thank Eric Leichtnam for guiding my feet into the
mathematical world and Marc Hindry for introducing me to algebra
with much kindness and beauty at the Ecole Normale of Cachan. I also
thank Marie-France Vigneras and Joseph Oesterle for teaching me
basic undergraduate algebra and giving me a strong background in
algebra.\\
\indent I am immensely grateful to the whole mathematics department
of the Ecole Normale Superieure de Cachan for being a family and
particularly to Jean-Michel Morel, Sylvie Fabre, Laurent
Desvillettes and Yves Meyer. They taught me much more than
mathematics.\\ \indent I thank the mathematicians Michael Maller and
Jennifer Whitehead for giving me self-confidence and introducing me
to research in mathematics at CUNY with much generosity.\\
\indent During my time at Caltech, many maths graduate students in
the years above mine have impressed me by their very exceptional
human qualities: David Whitehouse, Daniel Katz, Jennifer Johnson,
Kimball Martin, Bahattin Yildiz. It is a chance and a pleasure to
know them.\\ \indent I thank Visiting Professors Maria Jose Cantero
and Paula Trettkov for their kind support during their respective
visits at Caltech. To meet with both of them was a great pleasure. I
am grateful to Joseph Maher, Ada Chan and Cheng Yeaw Ku for their
friendliness to the graduate students during their time as postdocs
at the mathematics department.

\newpage

\pagenumbering{arabic} \pagestyle{plain}

\noin{\huge{\textbf{Contents}}}\\\\\\

%\noin\textbf{{\large{Acknowledgements}}}\hfill$\;\mathbf{iv}$\\\\
%\textbf{{\large{Abstract}}}\hfill$\;\mathbf{v}$\\\\
\noin $\mathbf{1}\;$ \textbf{{\large{Introduction}}}\hfill$\;3$\\

The Problem and the Main Results\dotfill $\;3$\\

The Approach\dotfill$\;5$\\

Some History and Recent Developments\dotfill $\;7$\\

\noin $\mathbf{2}\;$ \textbf{{\large{Background and Notations}}}\hfill$\mathbf{\;8}$\\\\
$\mathbf{3}\;$ \textbf{{\large{The Case $n=3$}}}\hfill$\mathbf{\;13}$\\\\
$\mathbf{4}\;$ \textbf{{\large{The Case $n=4$}}}\hfill$\mathbf{\;15}$\\\\
$\mathbf{5}\;$ \textbf{{\large{The Case $n=5$}}}\hfill$\mathbf{\;25}$\\\\
$\mathbf{6}\;$ \textbf{{\large{Generalization}}}\hfill$\mathbf{\;34}$\\\\
$\mathbf{7}\;$ \textbf{{\large{The Representation itself}}}\hfill$\mathbf{\;39}$\\\\
$\mathbf{8}\;$ \textbf{{\large{Proof of the Main Theorem}}}\hfill$\mathbf{\;47}$\\

$8.1\;$ Properties of the Representation and the Case
$n=6$\dotfill$47$\\

$8.2\;$ The Case $l=\unsur{r^{2n-3}}$\dotfill$48$\\

$8.3\;$ The Cases $l=\unsur{r^{n-3}}$ and $l=-\unsur{r^{n-3}}$\dotfill$51$\\

$8.4\;$ Proof of the Main Theorem\dotfill$86$\\

\noin $\mathbf{9}\;$ \textbf{{\large{More Properties of the L-K
Representation}} when $\mathbf{l\in\lb r,-r^3\rb}$}\hfill$\mathbf{\;99}$\\

$9.1\;$ Properties of Irreducibility of the BMW modules
$K(n)$\dotfill$99$\\

$9.2\;$ Properties of Dimension of the $K(n)$'s\dotfill$111$\\

$9.3\;$ A Proof of Theorems $C$ and $D$\dotfill$133$\\\\
%$9.4$ A proof of Theorem $E$
%$9.5$ The Branching Rule Point of View and Another Proof of Theorem
%$C$ and $D$\\
%$9.5$

\noin $\mathbf{10}\;$ \textbf{{\large{The Uniqueness Theorems and a Complete Description of \\the Invariant subspace of $\mathbf{\V^{(n)}}$ when $l=r$}}}\hfill$\;\mathbf{142}$\\\\
\noin$\mathbf{11}\;$ \textbf{{\large{A Proof of the Main Theorem
without
Maple}}}\hfill$\;\mathbf{149}$\\\\
\noin \textbf{{\large{Appendix $\mathbf{A}$. The program}}}\hfill$\mathbf{154}$\\\\
\noin \textbf{{\large{Appendix $\mathbf{B}$. Table for $\mathbf{Sym(8)}$}}}\hfill$\mathbf{160}$\\\\
\noin \textbf{{\large{Appendix $\mathbf{C}$. How the $\mathbf{X_{ij}}$'s act}}}\hfill$\mathbf{162}$\\\\
\noin \textbf{{\large{Appendix $\mathbf{D}$.  Determinants of submatrices of size $\mathbf{5}$ of $\mathbf{T(5)}$}}}\hfill$\mathbf{164}$\\\\
\noin \textbf{{\large{Appendix $\mathbf{E}\,$}. Reducibility of the L-K representation for $\mathbf{n=3}$}}\hfill$\mathbf{166}$\\\\
\noin \textbf{{\large{Appendix $\mathbf{F}\,$}. Reducibility of the L-K representation for $\mathbf{n=4}$}}\hfill$\mathbf{167}$\\\\
\noin \textbf{{\large{Appendix $\mathbf{G}$}. Reducibility of the L-K representation for $\mathbf{n=5}$}}\hfill$\mathbf{168}$\\\\
\noin \textbf{{\large{Appendix $\mathbf{H}$}. Reducibility of the L-K representation for $\mathbf{n=6}$}}\hfill$\mathbf{169}$\\\\

\noin {\large{\textbf{Bibliography}}}\hfill
$\mathbf{170}$\\\\

\noin {\large{\textbf{Index}}}\hfill $\mathbf{171}$
\newpage

\section{Introduction}
\subsection*{The Problem and the Main Results}
In this thesis, we build and examine a representation of degree
$\chl$ of the BMW algebra of type $A_{n-1}$ over the field $\Q(l,m)$
that is equivalent to the Lawrence-Krammer representation introduced
by Lawrence and Krammer in \cite{K4} and \cite{RL}. A result in
$\cite{CGW}$ states that this representation is generically
irreducible. Here we recover this result by a different method, and
show further that when $l$ and $m$ are specified in the field of
complex numbers, the Lawrence-Krammer representation is always
irreducible except when $l$ and $m$ are related in a certain way. We
let $r$ and $-\unsurr$ be the roots of the quadratics $X^2+mX-1=0$.
We define $\ih(n)$ as the Iwahori-Hecke algebra of the symmetric
group $Sym(n)$ over the field $F=\Q(l,r)$, with parameter $r^2$. Its
generators $g_1,\dots,g_{n-1}$ satisfy the braid relations and the
relation $g_i^2+m\,g_i=1$ for all $i$. This definition is the same
as the one given in \cite{IHA}, after the generators have been
rescaled by a factor $\unsurr$. We assume that $\ih(n)$ is
semisimple. It suffices to assume that $(r^2)^k\neq 1$ for every
$k\in\lb 1,\,\dots,\,n\rb$. We prove the following theorem:
\newtheorem*{maintheo}{Main Theorem}
\begin{maintheo}\hfill\\
Let $m=\unsurr-r$ and $l$ be two complex numbers. Let $n$ be an
integer with $n\geq 3$. Assume that $(r^2)^k\neq 1$ for every
$k\in\lb 1,\dots,n\rb$ and so $\ih(n)$ is semisimple. \\When $n\geq
4$, the Lawrence-Krammer representation of the $BMW$ algebra
$B(A_{n-1})$ with parameters $l$ and $m=\unsurr-r$ over the field
$\Q(l,r)$ is irreducible, except when $l\in\lb
r,-r^3,\unsur{r^{n-3}},-\unsur{r^{n-3}},\unsur{r^{2n-3}}\rb$ when it
is reducible. \\When $n=3$, the Lawrence-krammer representation of
the $BMW$ algebra $B(A_2)$ with parameters $l$ and $m=\unsurr-r$
over the field $\Q(l,r)$ is irreducible except when $l\in\lb
-r^3,\unsur{r^3},1,-1\rb$ when it is reducible.
\end{maintheo}
Moreover, for each of the values of $l$ and $r$ when the
representation is reducible, we show that there exists a unique
nonzero proper invariant subspace of the Lawrence-Krammer space,
except when $l=-r^3$ and $r^{2n}=-1$ and we are able to give its
dimension. When $l=-r^3$ and $r^{2n}=-1$, we prove that the
Lawrence-Krammer representation contains a direct sum of a
representation of degree one and an irreducible representation of
degree $\dbw$ and that there are no other irreducibles inside the
Lawrence-Krammer space. Suppose $\V^{(n)}$ is the Lawrence-Krammer
vector space over $\Q(l,r)$ of dimension $\chl$. We prove the
following five theorems and part of Conjecture $A$ below. In the
Theorems, the integer $n$ is taken to be greater than or equal to
$3$. When the lower bound is not $3$ in the theorems below, it means
these cases for those small integers are special and need to be
formulated apart. In fact, we prove that Conjecture $A$ below, that
we believe to be true for a general $n\geq 3$, holds in these cases.
In what follows we still assume that $\ih(n)$ is semisimple, or
which is equivalent $(r^2)^k\neq 1$ for every positive integer $k$
with $1\leq k\leq n$.
% and they
%require to be formulated apart.
\newtheorem*{theo1}{Theorem $A$}
\begin{theo1}
Let $n$ be an integer with $n\geq 4$. There exists a one-dimensional
invariant subspace of $\V^{(n)}$ if and only if
$l=\unsur{r^{2n-3}}$. If so it is unique.
\end{theo1}

\newtheorem*{theo2}{Theorem $B$}
\begin{theo2}
Let $n$ be an integer with $n\geq 3$ and $n\neq 4$. There exists an
irreducible $(n-1)$-dimensional invariant subspace of $\V^{(n)}$ if
and only if $l\in\lb\unsur{r^{n-3}},-\unsur{r^{n-3}}\rb$. If so, it
is unique.
\end{theo2}

\newtheorem*{theo3}{Theorem $C$}
\begin{theo3}
Let $n$ be an integer with $n\geq 4$. There exists an irreducible
$\cil$-dimensional invariant subspace of $\V^{(n)}$ if and only if
$l=r$. If so, it is unique.
\end{theo3}

\newtheorem*{theo4}{Theorem $D$}
\begin{theo4}
Let $n$ be an integer with $n\geq 5$. There exists an irreducible
$\dbw$-dimensional invariant subspace of $\V^{(n)}$ if and only if
$l=-r^3$. If so, it is unique.
\end{theo4}

\newtheorem*{theo5}{Theorem $E$}
\begin{theo5}
When the Lawrence-Krammer representation is reducible, $\V^{(n)}$
has a unique nontrivial proper invariant subspace except when
$l=-r^3$ and $r^{2n}=-1$ when $\V^{(n)}$ contains a direct sum of a
one-dimensional invariant subspace and an irreducible
$\dbw$-dimensional one.
% Moreover, these are the unique irreducibles
%of $\V^{(n)}$.
\end{theo5}

Moreover, in the case $l=\unsur{r^{2n-3}}$, we can give a spanning
vector for the one-dimensional invariant subspace and in the case
$l\in\lb\unsur{r^{n-3}},-\unsur{r^{n-3}}\rb$, we are also able to
describe the two respective irreducible $(n-1)$-dimensional
subspaces of $\V^{(n)}$. Also, when $l=r$ we can describe the
irreducible $\cil$-dimensional invariant subspace of $\V^{(n)}$.
% and so can we with the
% $\dbw$-dimensional invariant subspace when $l=-r^3$.

Furthermore, the cases $n=3$ and $n=4$ are described in the
following conjecture. The conjecture holds in these cases, by
Theorems $4$ and $5$ when $n=3$ and by Theorems $4$ and $6$ and
Corollary $4$ when $n=4$. It is proven that the Conjecture is also
true for $n=5$, namely by Theorems $4,5,8$ and Result $1$ in the
thesis. Part of the arguments lies in the proof of these theorems
when we show that the conjugate Specht modules $S^{(1^n)}$ (cases
$n=4,5$ in the proof of Theorem $4$), $S^{(2,1^{n-2})}$ (case
$n\in\lb 3,5\rb$ in the proof of Theorem $5$), $S^{(2,2,1^{n-4})}$
(case $n=5$ in the proof of Result $1$) cannot occur in the L-K
space $\V^{(n)}$. The results for these cases are respectively
gathered in Appendices $E$, $F$ and $G$ at the end of the thesis. We
will sometimes abbreviate L-K representation instead of
Lawrence-Krammer representation and L-K space instead of
Lawrence-Krammer space.

\newtheorem*{conjA}{Conjecture $A$}
\begin{conjA}
Let $n$ be an integer with $n\geq 3$. Assume $\ih(n)$ is semisimple
and $r^{2n}\neq -1$ when $l=-r^3$. \\
When the L-K representation is reducible, its unique nontrivial
proper invariant subspace is isomorphic to one of the Specht modules
$$S^{(n)},\,S^{(n-1,1)},\,S^{(n-2,2)},\,S^{(n-1,1,1)},$$
which respectively arise when and only when
$$l=\unsur{r^{2n-3}},\,l\in\bigg\lb\unsur{r^{n-3}},-\unsur{r^{n-3}}\bigg\rb,\,l=r,\,l=-r^3$$
When $l=-r^3$ and $r^{2n}=-1$, there are exactly three proper
nontrivial invariant subspaces of $\V^{(n)}$ and they are
respectively isomorphic to $S^{(n)}$, $S^{(n-2,1,1)}$ and
$S^{(n)}\oplus S^{(n-2,1,1)}$
\end{conjA}

The BMW algebra $B(A_{n-1})$ is generated by invertible elements
$g_i$'s, that are analogue to those from the Braid group and that we
will call generators of the first type. There is a second set of
$(n-1)$ elements in $B(A_{n-1})$ that generate an ideal, namely the
$e_i$'s (see section $2$ entitled "Background and Notations" of the
thesis). They are called generators of the second type. We will see
in the thesis that the generators of the second type and some of
their conjugates by products of generators of the first type play a
critical role. We think that the following result is true. If so it
will yield a different approach to solving Lemma $10$ of the thesis,
which is a key lemma on the road to Theorem $D$ (and $C$ as we did
in the thesis without the use of the Branching Rule).

\newtheorem*{conjB}{Conjecture $B$}
\begin{conjB}
Let $T$ denote the left action on the Lawrence-Krammer space of the
sum of all the conjugates of the $e_i$'s by products of the $g_i$'s.
Fix $r\in\Com$. Then, $det\,T\in\Q(l)$, the degree of $T$ is the
degree $\chl$ of the Lawrence-Krammer representation and
$$det(T)=\frac{(l-r)^{\cil}(l-r^3)^{\dbw}(l-\unsur{r^{n-3}})^{n-1}(l+\unsur{r^{n-3}})^{n-1}(l-\unsur{r^{2n-3}})}{\la_r\,l^{\chl}}$$
where $\la_r$ is a constant that depends on the value of the
parameter $r$.
\end{conjB}

The conjecture is shown through Maple for $n\in\lb 3,4,5,6\rb$ in
the Appendix $A$ where the procedure NOTIRR provides a necessary
condition on $l$ and $r$ so that the L-K space $\V^{(n)}$ is
reducible. When the representation is reducible, $l$ must be a root
of some polynomial of $\Q(r)[X]$, with multiplicity given in the
output. Visibly, this multiplicity is the dimension of the
corresponding invariant subspace.

%\noin The cases $n=3$ and $n=4$ are respectively described further
%in Theorem $4$ and Theorem $6$ of the thesis. \\
%Moreover, in the case $l=\unsur{r^{2n-3}}$, we can give a spanning
%vector for the one-dimensional invariant subspace and in the case
%$l\in\lb\unsur{r^{n-3}},-\unsur{r^{n-3}}\rb$, we are also able to
%explicitly describe the two respective irreducible
%$(n-1)$-dimensional subspaces of $\V$. \\
%Also, in the case $l=-r^3$,
%we can describe the unique irreducible $\dbw$-invariant subspace of
%$\V$.

\subsection*{The Approach}
Our approach is to use knot theory as a tool to build the
Lawrence-Krammer representation. We use a deformation of the Brauer
centralizer algebra, namely the Kauffman tangle algebra. This
algebra $MT_n$ constructed by Morton, Traczyk and Kauffman was shown
by Morton and Wassermann in $1989$ to be generically isomorphic to
the $BMW$ algebra of type $A_{n-1}$ (see \cite{MW}). One difference
between the Brauer centralizer algebra and the Morton-Traczyk
algebra is that the group algebra of the symmetric group has been
replaced with the Iwahori-Hecke algebra of the symmetric group. This
introduces a new type of braids where over-crossings have to be
distinguished from under-crossings. A tangle with an over-crossing
is related to a tangle with an under-crossing by the Kauffman skein
relation that uses the parameter of the Iwahori-Hecke algebra of the
symmetric group. The geometric approach of the Kauffman tangle
algebra leads us to visualize (and later on prove for the algebraic
version of the representation) that if $\W$ is a proper invariant
subspace of the Lawrence-Krammer vector space $\V$, then the action
on $\W$ must be an Iwahori-Hecke action. This was incidentally shown
by Arjeh Cohen, Die Gijsbers and David Wales in a more general
setting when they also deal with representations of the BMW algebra
of types $D$ and $E$. By using the Lawrence-Krammer representation,
we can find the two inequivalent irreducible matrix representations
of $\ih(n)$ of degree $(n-1)$. We then use those irreducible
representations to show that if there exists an irreducible
$(n-1)$-dimensional invariant subspace of $\V$, it forces
$l\in\lb\unsur{r^{n-3}},-\unsur{r^{n-3}}\rb$. Conversely, if $l$
takes one of these fractional values in $r$, we show that there
exists a unique irreducible $(n-1)$-dimensional invariant subspace
of $\V$. For the small values $3,4,5$ of $n$ we are able to find
with Maple from the Lawrence-Krammer representation the inequivalent
irreducible matrix representations of $\ih(n)$; we then deduce from
these matrix representations wether or not the Lawrence-Krammer
representation can be reducible and if so for which values of the
parameters $l$ and $r$. The table in Appendix $E$, (resp $F$, $G$,
$H$) lists which irreducible representations of $\ih(n)$, $n=3$
(resp $4,5,6$) can possibly occur in the Lawrence-Krammer
representation and for which corresponding values of $l$ and $r$.
All the cases described in the tables are proven in the thesis. The
case $n=6$ could be completed by using the Branching Rule. In the
table for $n$, each Ferrers diagram associated with a partition of
the integer $n$ corresponds to a class of irreducible
$\ih(n)$-module. Only the irreducible representations of $\ih(n)$ of
degree less than the degree $\chl$ of the Lawrence-Krammer
representation are represented in the table. Let the $\al_i$'s for
$i=1,\dots,n-1$ denote the $(n-1)$ simple roots of the root system
of type $A_{n-1}$ and let $\phi^{+}$ denote the set of all the
positive roots. There are $\chl$ of these. The Lawrence-Krammer
vector space $\V^{(n)}$ is the vector space over $\Q(l,r)$ with the
generating set $\lb x_{\be}|\be\in\phi^{+}\rb$, indexed by all the
positive roots. Thus it has dimension $\chl$. If
$\be=\al_i+\dots+\al_{j-1}$ is a positive root where $i<j$, then to
simplify the notations, we let $w_{ij}$ denote the vector
$x_{\al_i+\dots+\al_{j-1}}$. Thus, $\V^{(n)}$ is spanned over
$\Q(l,r)$ by the $\chl$ vectors $w_{ij}$ with $1\leq i<j\leq n$.
These are the vectors that appear in the tables of the Appendices.

\subsection*{Some History and Recent Developments}

The BMW algebras of type $A$ were introduced by J. Murakami in
\cite{MUR} and separately by J.S. Birman and H. Wenzl in \cite{BW}
in order to try to find faithful representations of the Braid group.
%[what they are. Cubic relation]
%The BMW algebras happen to be interesting in many regards. For
%instance, when their parameters are specialized over the field of
%complex numbers like in the thesis, they can be realized as
%endomorphism algebras over a quantum group.
Birman and Wenzl wanted to use these algebras to construct
representations of the Braid group $B_n$ over $n-1$ generators in
order to investigate wether this group is linear. Linearity of a
group means that there exists a faithful representation into
$GL(m,\R)$ for some positive integer $m$. Burau discovered an
$n$-dimensional representation of $B_n$ that was faithful for $n=3$,
but shown not to be faithful for $n\geq 9$ (Moody, \cite{MOO}). For
a long time, if it was known that $B_n$ was linear for $n\in\lb
2,3\rb$, for $n\geq 4$, the problem remained open. Krammer later
introduced a representation of $B_n$ and showed it to be faithful
for $n=4$ (cf \cite{K4}). Since the same representation occurs in an
earlier work of Lawrence, it is called the Lawrence-Krammer
representation. In \cite{BIG}, using topology, Bigelow proved the
Lawrence-Krammer representation to be faithful for all $n$, thus
stating that all the Artin groups of type $A$ are linear. Shortly
after, Krammer also showed with specific real values of the
parameters that his representation is a faithful representation of
the Braid group, using algebraic methods this time (see \cite{Kn}).
In $2002$, A.M. Cohen and D.B. Wales generalized the
Lawrence-Krammer representation to Artin groups of finite types, as
well as their proof of linearity, thus showing in \cite{CW} that
every Artin group of finite type is linear. Following the
observation that the Lawrence-Krammer representation of the Artin
group of type $A_{n-1}$ factors through the BMW algebra of the same
type, they then generalized the definition of BMW algebras to other
types $D$ and $E$ and built representations of these algebras that
they showed to be generically irreducible.
%They worked over the field $\Q(l,m)$ and
%built for each irreducible representation over $\Q(r)$ of degree say
%$k$, of the Hecke algebra of type $C$ over $Q(l,m)$ (where $C$ is
%the set of nodes whose associated simple roots are orthogonal to the
%highest root) a representation over $\Q(l,r)$ of degree
%$k\times\chl$ of the generalized BMW algebra of simply laced type.
%When the BMW algebra is of type $A_{n-1}$, $C$ is of type $A_{n-3}$
%and the

%A consequence of this fact is that if the Lawrence-Krammer space
%$\V^{(n)}$ decomposes as a sum of irreducibles
%$$\V^{(n)}=\I_1\oplus\dots\oplus\I_s,$$
%with all the $\I_j$'s proper, then by the considerations above, each
%of them is annihilated by $e_1$, so that their sum $\V^{(n)}$ is
%also annihilated by $e_1$. This yields a contradiction, forcing in
%fact $\V^{(n)}$ irreducible. Consequently, the L-K space is
%completely irreducible if and only if it is irreducible. That is how
%the theorem on the front page

As part of this work, we show that if the Lawrence-Krammer space has
a proper invariant subspace $\W$, the $g_i$ conjugates of the
$e_i$'s all annihilate $\W$. That is how the Main theorem implies
that for these specific values of the parameters $l$ and $r$, the
BMW algebra $B$ of type $A_{n-1}$ with parameters $l$ and
$m=\unsurr-r$ is not semisimple. Hans Wenzl showed a similar result
when he states in \cite{HW1} that these algebras are semisimple
except possibly if $r$ is a root of unity or $l=r^k$ for some
$k\in\Z$, where he also considers complex parameters. However, he
does not prove or mention for which integers $k$ or which roots of
unity the algebras fail to be semisimple. Hebing Rui and Mei Si
recently obtained the same result as ours in showing that for the
values of $l$ and $r$ that we found, $B$ is not semisimple. They use
the representation theory of cellular algebras (see \cite{RUI}).
%Their result is however stronger in
%that they show that $B_n$, $n\geq 3$ is semisimple if and only if

%$$o(r)>2n\;\;\&\;\;l\not\in\bigcup_{k=3}^{n}\lb r^{3-2k},\pm
%r^{3-k},-r^{2k-3},\pm\,r^{k-3}\rb\qquad(\diamond)$$

%\noin When doing so they assume that $l\not\in\lb\unsurr,-r\rb$. In
%the case when $l\in\lb\unsurr,-r\rb$, they show that $B_n$ is never
%semisimple except possibly for $n=3$ when $o(r)>6$ and $r^4\neq -1$
%and for $n=5$ when $o(q)>10$, $r^6\neq -1$ and $r^8\neq -1$. In this
%thesis, we only show that if $B_n$ is semisimple, we must have:
%$$\begin{array}{l} r^2\neq 1,\,(r^2)^2\neq 1,\dots,(r^2)^n\neq 1\\\qquad\&\\
%l\not\in\lb r^{3-2n},\pm r^{3-n},-r^{2n-3},\pm
%r^{n-3},r,-\unsurr,-r^3,\unsur{r^3}\rb
%\end{array}$$

\section{Background and Notations} We consider the BMW algebra
$B(A_{n-1})$, with nonzero complex parameters $m$ and $l$, as
defined in \cite{CGW}. In this thesis we build a representation of
the algebra $B(A_{n-1})$ and we find necessary and sufficient
conditions on $l$ and $m$ so that it is irreducible. Throughout the
thesis we will use the change of parameter $m=1/r-r$, so that the
two relevant parameters will be $l$ and $r$. We recall below the
relations defining $B(A_{n-1})$, where we assumed $r^2\neq 1$ and so
$m\neq 0$:
\begin{eqnarray}
g_ig_j & = & g_jg_i \;\;\;\;\;\;\;\;\;\;\;\;\;\;\;\;\;\;\;\;\;\;\;\;\;\;\;\,\;\;\;\;\; \forall  1\leq i,j\leq n-1,  \; |i-j|\geq 2\;\;\\
g_ig_{i+1}g_i &=& g_{i+1}g_ig_{i+1} \;\;\;\;\;\;\;\;\;\;\;\;\;\;\;\;\;\;\;\;\;\;\; \forall  1\leq i\leq n-2 \;\;\\
e_i & = & \frac{l}{m}*(g_i^2+m\;g_i-1) \;\;\;\!\;\;\;\;\; \forall 1\leq i\leq n-1 \\
g_ie_i &=& l^{-1}e_i  \;\;\;\;\;\;\;\;\;\;\;\;\;\;\;\;\;\;\;\;\;\;\;\;\;\;\;\;\;\;\;\forall 1\leq i \leq n-1 \;\;\\
e_ig_{i+1}e_i & = & le_i \;\;\;\;\;\;\;\;\;\;\;\;\;\;\;\;\;\;\;\;\;\;\;\;\;\;\;\;\;\;\;\;\;\;\;\forall 1\leq i \leq n-2 \;\;\\
e_ig_{i-1}e_i & = &
le_i\;\;\;\;\;\;\;\;\;\;\;\;\;\;\;\;\;\;\;\;\;\;\;\;\;\;\;\;\;\;\;\;\;\;\;
\forall 2\leq i\leq n-1\;\;
\end{eqnarray}
as well as some important and direct consequences of these defining
relations:
\begin{eqnarray}
e_ig_i & = & l^{-1}e_i \;\;\;\;\;\;\;\;\;\;\;\;\;\;\;\;\;\;\;\;\;\;\;\;\;\;\;\;\;\;\; \forall 1\leq i\leq n-1\;\;\;\;\;\;\;\;\;\;\;\;\;\;\;\;\\
g_i^2 & = & 1-mg_i+ml^{-1}e_i\;\;\;\;\;\;\;\;\;\;\;\forall 1\leq i\leq n-1\\
g_i^{-1} & = & g_i +
m-m\;e_i\;\;\;\;\;\;\;\;\;\;\;\;\;\;\;\;\;\forall 1\leq i\leq n-1
\end{eqnarray}
We will also use the "mixed Braid relations":
\begin{eqnarray}
g_ig_{i+1}e_i & = & e_{i+1}e_i\;\;\;\;\;\;\;\;\;\;\;\;\;\;\;\;\;\;\;\;\;\;\;\;\;\;\;\;\;\forall 1\leq i\leq n-2 \\
g_ig_{i-1}e_i & = & e_{i-1}e_i\;\;\;\;\;\;\;\;\;\;\;\;\;\;\;\;\;\;\;\;\;\;\;\;\;\;\;\;\;\forall 2\leq i\leq n-1\\
g_ie_{i+1}e_i & = & g_{i+1}e_i + m(e_i-e_{i+1}e_i)\;\;\:\forall 1\leq i\leq n-2\;\;\;\;\;\;\;\;\;\;\;\;\;\;\;\;\;\;\;\;\;\;\\
g_ie_{i-1}e_i & = & g_{i-1}e_i + m(e_i-e_{i-1}e_i)\;\;\:\forall
2\leq i\leq n-1
\end{eqnarray}

\noindent Finally, for any node $i=1,\dots n-1$, we will make use of
the idempotent relation $e_i^2=x\,e_i$ where the parameter $x$ is
defined by $x=1-\frac{l-1/l}{1/r-r}$\\\\
We will work over the field $F=\Q (r,l)$. Let $I_1$ and $I_2$ be the
two-sided ideals of $B(A_{n-1})$ respectively generated by the
$e_i$'s and the $e_{i}e_j$'s with $i\nsim j$ as in \cite{CGW}. The
representation will be built inside $I_1/I_2$. It is shown in
\cite{MW} that the BMW algebra $B(A_{n-1})$ is isomorphic to the
Kauffman algebra $TM_n$ and from now on we will think of BMW
elements in terms of tangles. The algebra $B(A_{n-1})$ acts by the
left on the subspace of $I_1/I_2$ spanned by the elements with a
fixed bottom horizontal line joining any two nodes, giving to this
subspace a structure of $B(A_{n-1})$-module. We fix such a bottom
horizontal line and we denote by $V$ the corresponding subspace.
Without loss of generality we take this horizontal line to be the
one joining nodes $1$ and $2$ on the bottom line. From now on, all
the elements that we consider are in $V$ and so they are linear
combinations over $F$ of tangles all having their nodes $1$ and $2$
joined on the bottom line. When two tangles have the same top
horizontal line, we will say that they are "similar" as in the
definition below:

\newtheorem{definition}{Definition}
\begin{definition}\hfill\\
Two tangles in $V$ are "similar" if their respective top horizontal
line are the same, with this top horizontal line overcrossing the
eventual vertical strings.
\end{definition}

\noindent We now recall some basic facts about the Artin group of
type $A_{n-1}$ and its associated root system and we introduce a set
of tangles in $V$ that contains as many elements as there are
positive roots and which is in correspondence with the set of
positive roots. These tangles are obtained by picking any random top
horizontal edge. Moreover, we agree that this top horizontal edge
always overcrosses any of the vertical strings that it
intersects.\\
Let $M=(m_{ij})$ be the Coxeter matrix of the Artin group $A$ of
type $A_{n-1}$ with generators $s_1, ..., s_{n-1}$. Let $(\epsilon
_{i})_{i=1...n-1}$ be the canonical basis of $\R^{n-1}$ and let
$B_{M}$ be the canonical symmetric bilinear form over $\R^{n-1}$
associated to $M$.
$$B_M(\e_i,\e_j) = - cos(\pi/m_{ij}) \;\;\text{where}\;\; m_{ij}=\begin{cases} 2&\text{if $|i-j|>1$},\\3&\text{if $|i-j|=1$}\\1&\text{if $i=j$}\end{cases}$$
For $i=1,...,n-1$, let $s_{H_i}$ denote the reflection with respect
to the hyperplane $H_i:= Ker B_M(.,\e_i)$. By the theory in
\cite{BOUR}, we know that the Coxeter group $\mathcal{W}$ of the
Artin group of type $A_{n-1}$ is isomorphic to the reflection group
spanned by the $s_{H_i}$'s. We have a root system and $\mathcal{W}$
is finite as it permutes the roots. Hence $B_M$ is an inner product
and throughout the paper we will denote it by ( , ). Finally,
instead of $\e_1,\e_2,...,\e_{n-1}$, we denote the simple roots by
$\al_1,...,\al_{n-1}$ and the associated simple reflections
$s_{H1},...,s_{H_{n-1}}$ by $r_1,...,r_{n-1}$. We associate to each
positive root $\be$ a tangle in $V$ in the following way: $\al_1$ is
$e_1$. Each positive root is then built from $\al_1$. For
$\be\in\phi^{+}$, where $\phi^{+}$ is the set of positive roots, let
$w_{\be 1}$ be the unique element in the Weyl group of minimal
length mapping $\al_1$ to $\be$, as in \cite{CGW}. For instance, we
have $w_{\al_21}=r_1r_2$ and $w_{\al_2+\al_3+\al_4,1}=r_4r_3r_1r_2$.
Taking the same notations as in \cite{CGW}, there is a map:
$$\left(\begin{array}{ccc}
\mathcal{W} & \stackrel{\psi}{\longrightarrow} & A\\
r_{i_1}...r_{i_l} & \longmapsto & s_{i_1}...s_{i_l}\\
\end{array}\right)$$\\
where $r_{i_1}...r_{i_l}$ is a reduced decomposition. To have this
map well defined, we need to show that if $r_{i_1}...r_{i_l}$ and
$r_{j_1}...r_{j_l}$ are two reduced decomposition in $W$ such that
$r_{i_1}...r_{i_l}=r_{j_1}...r_{j_l}$ then $s_{i_1}...s_{i_l}$ can
be transformed into $s_{j_1}...s_{j_l}$ by using the braid relations
\begin{eqnarray}
s_is_j &=& s_js_i \;\;\;\;\;\;\;\;\;\;\;\;\;\!\forall 1\leq i,j\leq n-1, |i-j|>2\\
s_is_{i+1}s_i &=& s_{i+1}s_is_{i+1} \;\;\;\forall 1\leq i\leq n-2
\end{eqnarray}
This is Matsumoto's theorem and it is stated (and proved) in \cite{IHA}.\\
And there is a morphism of groups:
$$\left(\begin{array}{ccc}
A & \stackrel{b}{\rightarrow} & B(A_{n-1})^{\times}\\
s_i & \mapsto & g_i\\
\end{array}\right)$$
where $ B(A_{n-1})^{\times}$ is the group of units of $B(A_{n-1})$.\\
We denote by $\varphi$ the composition of these two maps and we now
associate to $\be$ the BMW element $\varphi (w_{\be 1})e_1$. Thus
$\al_2$ is associated with $g_1g_2e_1$ and using the relations
defining $B(A_{n-1})$, this element is in fact $e_2e_1$. Similarly,
$\al_2+\al_3+\al_4$ is associated with $g_4g_3g_1g_2e_1$, which can
be rewritten to $g_4g_3e_2e_1$ thanks to relation $(10)$. We note
that more generally, the $\al_i$'s are associated with the
$e_ie_{i-1}...e_1$'s and are obtained from $\al_1$ by shifting the
top horizontal line to the right. Then from $\al_i$, get a positive
root $\al_i+...+\al_{i+k}$ by acting by the left on $e_i...e_1$ with
$g_{i+k}...g_{i+1}$, the top horizontal line crossing now $k$
vertical lines. More generally the top horizontal line of the tangle
representing a root of height $h$ always crosses $h-1$ vertical
lines. We denote by $W$ the subspace of $V$ spanned by these tangles
and by $\mathcal{B}$ the basis of $W$ composed of these tangles. By
construction, the cardinality of $\mathcal{B}$ equals the number of
positive roots.

We also introduce elements $X_{ij}$ in $TM_n$ which will be useful
throughout the paper. By definition, $X_{ij}$ is the element with
two horizontal lines, one at the top and the other one at the
bottom, both joining nodes $i$ and $j$, all the other nodes being
joined by straight vertical lines. Moreover, in this definition we
will take the horizontal lines overcrossing the vertical lines. We
note that there are as many elements $X_{ij}$'s as there are
positive roots since for each positive root we built a unique tangle
in $V$ by taking $\al_1$ to be $e_1$ and by letting the horizontal
edge on the top line vary for the other ones, allowing all the
possible combinations. And in fact,
$$\big|\{X_{ij};1\leq i<j\leq
n\}\big|=|\phi^{+}|=|\mathcal{B}|=\binom{n}{2}$$
%Finally, we let $S$
%be the sum of all the $X_{ij}$'s:
%$$S:=\sum_{1\leq i<j\leq n} X_{ij}$$

\noindent The following lemma relies on the fact that the $X_{ij}$'s
are conjugate to the $e_i$'s:
%\newtheorem{prop}{Proposition}
%\begin{prop}
%$Ker (X_{k-1k+1})\bigcap Ker\,e_k\bigcap Ker(X_{kk+2})$ is a
%$B(A_{n-1})$-module for all $2\leq k\leq n-2$.
%\end{prop}
%\noindent\textsc{Proof:} the proof relies on the fact that the
%$X_{ij}$'s are conjugate to the $e_i$'s. First we show the following
%lemma:
\newtheorem{lemma}{Lemma}
\begin{lemma}
$\forall x,\Big(x\in\underset{1\leq i<j\leq n}{\bigcap}Ker X_{ij}
\Rightarrow\forall k,\forall l,\;g_l\;x\in Ker\;e_k\Big)$
\end{lemma}
\noindent\textsc{Proof of the lemma:} let $x$ be in the intersection
of the kernels of the $X_{ij}$'s. By definition, we have $e_kx=0$.
Then, by $(7)$, it comes $e_kg_kx=0$. So we get: $g_kx\in Ker\;e_k$.
Next, from $X_{k\, k+2}x=0$, we derive $g_{k+1}e_kg_{k+1}^{-1}x=0$,
which implies $e_kg_{k+1}^{-1}x=0$, by multiplying by the left the
last equality with $g_{k+1}^{-1}$. Hence, $g_{k+1}^{-1}x\in Ker\;
e_k$. By using $(9)$, it follows that $g_{k+1}x\in Ker\;e_k$. Also,
from $X_{k-1k+1}x=0$, we derive $g_{k-1}^{-1}e_kg_{k-1}x=0$ which
implies that $g_{k-1}x\in Ker\;e_k$. At this point we have for any
two nodes $k$ and $l$:
$$\text{if}\; k=l\; \text{or if}\; k\sim l \;\text{then}\; g_lx\in
Ker\; e_k$$ Further, if $k\nsim l$, by using the braid relation
$(1)$ and the polynomial expression $(3)$, $e_k$ and $g_l$ commute
and it immediately follows that $g_lx\in Ker\; e_k$. Thus we have
proved the lemma.
\newtheorem{Cor}{Corollary}
\begin{Cor}
$\forall x\in\underset{1\leq i<j\leq n}{\cap} Ker\; X_{ij},\;\forall
 1\leq l\leq n-1,\; g_l^{-1}x\in\underset{1\leq k\leq n-1}{\cap}
 Ker\; e_k$
 \end{Cor}
\noindent\textsc{Proof:} direct by using the lemma and relation
$(9)$. \\\\
We now prove the stronger lemma:
\begin{lemma}
$\forall x,\Big(x\in\underset{1\leq i<j\leq n}{\bigcap}Ker X_{ij}
\Rightarrow\forall l,\;g_l\;x\in\underset{1\leq i<j\leq
n}{\bigcap}Ker X_{ij}\Big)$
\end{lemma}
\noindent\textsc{Proof:} Given an x in the intersection of the
kernels of the $X_{ij}$'s and given $j\geq i+2$, we want to show
that $g_l\,x\in Ker\,X_{ij}$ for any node $l$. The case $j=i+1$ has
already been seen in lemma $1$. The idea is to consider different
conjugation formulas for $X_{ij}$. And indeed, let's first write:
$$ X_{ij}=g_{j-1}\dots g_{i+1}\,e_i\,g_{i+1}^{-1}\dots
g_{j-1}^{-1}$$ From there we deduce two things:
\begin{list}{\texttt{-}}{}
\item If $k\nsim\{ i,\dots,j-1\}$, then $g_k\,x\in Ker\,X_{ij}$
as $g_k$ commutes with all the elements composing the word above.
\item $g_{j-1}x\in Ker\,X_{ij}$. Indeed, we have in the case $j>i+2$:
\begin{eqnarray*}
X_{ij}g_{j-1}x &=& (g_{j-1}\dots g_{i+1}\,e_i\,g_{i+1}^{-1}\dots
g_{j-1}^{-1})g_{j-1}x\\
&=& g_{j-1}(g_{j-2}\dots g_{i+1}\,e_i\,g_{i+1}^{-1}\dots g_{j-2}^{-1})x\\
&=& 0 \qquad\text{since $x\in Ker\,X_{i\,j-1}$}
\end{eqnarray*}
If $j=i+2$, we have $X_{i\,i+2}=g_{i+1}\,e_i\,g_{i+1}^{-1}$ and it
follows that\\ $X_{i\,i+2}\,g_{i+1}x=g_{i+1}\,e_i\,x=0$ since $x\in
Ker\,e_i$.
\end{list}
\noindent We now deal with the nodes $i,\,i+1,\dots j-2$. Instead of
starting from the left with $e_i$ like in the previous decomposition
for $X_{ij}$, we start from the right with $e_{j-1}$ and get the
following expression: \begin{eqnarray}X_{ij}=g_{i}^{-1}\dots
g_{j-2}^{-1}\,e_{j-1}\,g_{j-2}\dots g_i\end{eqnarray} Again we
deduce two things from this writing:
\begin{list}{\texttt{-}}{}
\item  Let $k$ be any value in $\{2,\dots ,j-i-1\}$. We look at the
action of $X_{ij}$ on $g_{j-k}\,x$. We get successively:
\begin{eqnarray*}
X_{ij}\,g_{j-k}\,x &=& g_{i}^{-1}\dots
g_{j-2}^{-1}\,e_{j-1}\,g_{j-2}\dots g_i\,g_{j-k}\,x\\
&=& g_{i}^{-1}\dots g_{j-2}^{-1}\,e_{j-1}\,g_{j-2}\dots
g_{j-k}g_{j-k-1}g_{j-k}g_{j-k-2}\dots g_i\,x\\
&=& g_{i}^{-1}\dots g_{j-2}^{-1}\,e_{j-1}\,g_{j-2}\dots g_{j-k+1}
g_{j-k-1}g_{j-k}g_{j-k-1}\dots g_ix \\
&=& g_{i}^{-1}\dots g_{j-2}^{-1}g_{j-k-1}\,e_{j-1}\,g_{j-2}\dots
g_i\,x\\
&=& g_{i}^{-1}\dots
g_{j-2}^{-1}g_{j-k-1}^{-1}\,e_{j-1}\,g_{j-2}\dots g_i\,x
-m\,X_{ij}\,x\\
&=& g_i^{-1}\dots g_{j-k-1}^{-1}g_{j-k}^{-1}g_{j-k-1}^{-1}\dots
g_{j-2}^{-1}e_{j-1}g_{j-2}\dots g_{i}\,x\\
&=& g_i^{-1}\dots
g_{j-k-2}^{-1}g_{j-k}^{-1}g_{j-k-1}^{-1}g_{j-k}^{-1}\dots
g_{j-2}^{-1}e_{j-1}g_{j-2}\dots g_{i}\,x\\
&=& g_{j-k}^{-1}\,g_{i}^{-1}\dots g_{j-2}^{-1}e_{j-1}g_{j-2}\dots
g_i\,x\\
&=& g_{j-k}^{-1}\, X_{ij}\,x\\
&=& 0
\end{eqnarray*}

To quickly manipulate the equalities above, it will be useful to
notice that moving from the right to the left, the indices increase
on the right of $e_{j-1}$ and they decrease on the left of
$e_{j-1}$. We have the first equality by definition of $X_{ij}$. The
second equality is obtained by using the braid relation $(1)$,
commuting $g_{j-k}$ till we reach $g_{j-k-1}$ and the third one is
obtained by using the braid relation $(2)$ on the nodes $j-k$ and
$j-k-1$: $g_{j-k}g_{j-k-1}g_{j-k}=g_{j-k-1}g_{j-k}g_{j-k-1}$. One of
the factors $g_{j-k-1}$ lies now next to $g_{j-k+1}$ and we may
again apply the braid relation $(1)$ repeatedly to move it to the
left of $e_{j-1}$. This is the object of the equality number four.
Next since all the terms on the left of $e_{j-1}$ are inverses, it
is natural to express $g_{j-1}$ in term of its inverse by using
$(9)$. The fifth equality is obtained by using $(9)$ and the fact
that $e_{j-k-1}e_{j-1}=0$ (since $j-k-1\in\{i,\dots,j-3\}$). When
the factor $g_{j-k-1}$ "disappears" in $(9)$, it yields $X_{ij}\,x$,
which is zero by choice of $x$ in the intersection of all the
kernels. Thus, in the product, $g_{j-k-1}$ may in fact be replaced
by $g_{j-k-1}^{-1}$ without modifying the word. Next, in the
equality number $6$, we commute $g_{j-k-1}^{-1}$ to the left till we
get stuck. Then we apply again the braid relation $(2)$ on the
inverses to get the equality number $7$. We may now move
$g_{j-k}^{-1}$ to the extreme left by using again the commuting
relation $(1)$ on the inverses. It finally allows us to isolate
$X_{ij}$ which acting by the left on $x$ yields $0$. Thus, we have
proved the following:
$$\forall k\in\{i+1,\dots j-2\},\,g_k\,x\in Ker\, X_{ij}$$
\item The second thing that we derive from the expression $(16)$ is direct: $g_{i}^{-1}\dots g_{j-2}^{-1}\,e_{j-1}\,g_{j-2}\dots
g_{i-1} g_i\,g_i\,x =0$ by expanding the square $g_i^2$ with $(8)$
and using the fact that $X_{i-1\,j}x=X_{ij}x=e_i\,x=0$. Thus, we
have $g_i\,x\in Ker\,X_{ij}$.
\end{list}
\noindent It remains to show that $g_{i-1}\,x\in Ker\,X_{ij}$ and
$g_j\,x\in Ker\, X_{ij}$ Let's write:
$$X_{ij}=g_j^{-1}X_{i\,j+1}g_j$$
The trick here is to go "one node too far" to make a $g_j$ appear on
the right extremity of the word. Then, using the formula
$g_j^2=1-m\,g_j+ml^{-1}\,e_j$ and the fact that
$X_{i\,j+1}\,x=X_{ij}\,x=e_j\,x=0$ by our hypothesis on $x$, it
immediately follows that $X_{ij}\,g_j\,x=0$. Hence we have
$g_j\,x\in Ker\, X_{ij}$.\\
Finally we write:
$$X_{ij}=g_{i-1}X_{i-1\,j}g_{i-1}^{-1}$$
Here, the trick is to go one node backwards to make a "$g_{i-1}$"
appear. It follows that $X_{ij}g_{i-1}\,x=g_{i-1}X_{i-1\,j}\,x=0$
since $x\in Ker\,X_{i-1\,j}$. Hence $g_{i-1}\,x\in Ker\,X_{ij}$ and
we have now proved that
$$\forall l,\,g_l\,x\in Ker\,X_{ij}$$
We have the immediate corollary:
\begin{Cor}
$\forall n$, $\underset{1\leq i<j\leq n}{\bigcap}Ker X_{ij}$ is a
$B(A_{n-1})$-module
\end{Cor}
%it suffices to notice that the $X_{ij}$'s are conjugate of the
%$e_i$'s in the following way:
%$$\forall 1\leq i<j\leq n,\;
%X_{ij+1}=g_{j}.\;X_{ij}.\;g_j^{-1}$$ From there, it is easy to see
%that if $x\in\cap_{1\leq i<j\leq n} Ker(X_{ij})$ then all the
%$g_{i}^{-1},\; i=1...n-1$ are in this intersection

\section{The case n=3} In this section we explicitly build a
representation of $B(A_2)$ in $V$.
\\There are three positive roots in $A_2$: $\al_1$, $\al_2$,
$\al_1+\al_2$ and $V$ is spanned by \\$e_1$, $e_2e_1$ and $g_2e_1$.
Acting by $g_1$ on these elements yields:
\\$g_1e_1=l^{-1}e_1$, $g_1e_2e_1=g_2e_1 + m e_1 - m e_2e_1$ and
$g_1g_2e_1=e_2e_1$.
\\Thus, the matrix of the action of $g_1$ is given by:

$$\begin{array}{ll}
\;\;\;\;\;\;\;\;\;\;\;\;\;\;\;\;\;\;\;
e_1\;\;\;\; e_2e_1\;\;\;g_2e_1\\
\begin{array}{lll}
e_1 \\ e_2e_1 \\ g_2e_1\\
\end{array}
\left( \begin{array}{ccccccc}
& \l^{-1} &  m  & & 0 & & \\
& 0       & \!\!\!-m  & & 1 & & \\
& 0       &  1  & & 0 & & \\
\end{array} \right)\\
\end{array}=:G_1(3)
$$
\\
\noindent Similarly, we have $g_2e_2e_1=l^{-1}e_2e_1$ and
$g_2g_2e_1=e_1+ml^{-1}e_2e_1-mg_2e_1$. Thus, the matrix of the
action of $g_2$ is:

$$\begin{array}{ll}
\;\;\;\;\;\;\;\;\;\;\;\;\;\;\;\;\;\;
e_1 \;\;\;\; e_2e_1\;\;\;\;\;\; g_2e_1\\
\begin{array}{lll}
e_1 \\ e_2e_1 \\ g_2e_1\\
\end{array}
\left( \begin{array}{cccccc}
& 0 &  0  &  1 & \\
& 0 & \;\;\;l^{-1} & \;\;\;ml^{-1} & \\
& 1 &  0  &  \!\!\!-m & \\
\end{array} \right)\\
\end{array}=:G_2(3)
$$
\\
\noindent We verify that the braid relation $G_1G_2G_1=G_2G_1G_2$ is
satisfied. Next, for $i\in\{1,2\}$, we define $E_i$ by the matrix
relation:
\begin{eqnarray}
E_i:=l/m*(G_i^2+mG_i-Id)
\end{eqnarray}
and we check that for $i\in\{1,2\}$ and $j\in\{1,2\}$ the matrix
relations
\begin{eqnarray}
G_iE_i &=& l^{-1}E_i \\
E_iG_jE_i&=&lE_i
\end{eqnarray}
are satisfied. Hence we have a left representation: $B(A_2)\ra
End_{F}(V)$. \\Suppose now that $V$ is not irreducible. Then it has
a proper nonzero invariant subspace $U$. Let $u$ be a nonzero
element of $U$. Let's decompose $u$ over the basis
$\mathcal{B}=(e_1, e_2e_1, g_2e_1)$, say $u=\la_1\; e_1 + \la_2\;
e_2e_1 + \la_3\; g_2e_1$ with $(\la_1, \la_2, \la_3)\neq (0, 0, 0)$.
We look at the action of $X_{12}=e_1$. We have:
 $$ \left\lbrace \begin{array}{ccc}
 e_1e_1=e_1^2 & = &xe_1\\
 e_1(e_2e_1) & = & e_1\\
 e_1(g_2e_1) & = & le_1\\
 \end{array}\right.$$
Thus, we get $X_{12}u= (\la_1x+ \la_2+\la_3l)\;e_1$. Assume that
$X_{12}u\neq 0$ i.e. $\la_1x+ \la_2+\la_3l\neq 0$. It follows that
$e_1$ is in $U$. Then $e_2e_1$ and $g_2e_1$ are also in $U$ as $U$
is a $B(A_2)$-submodule of $V$. Then $U$ would not be proper which
is a contradiction. Hence we have $X_{12}u=0$. This last equality
being true for any nonzero element of $U$, it comes $X_{12}U=0$.
Similarly, acting by $X_{23}=e_2$ on an element $u$ of $U$ gives a
multiple of $e_2e_1$ and acting by $X_{13}$ gives a multiple of
$g_2e_1$. Since
$$\;\left\lbrace \begin{array}{ccc}
g_1^{-1}(e_2e_1) & = & g_2e_1 \\
e_1(e_2e_1) & = & e_1 \end{array} \right. \text{and   }\;
\;\left\lbrace
\begin{array}{ccc}
g_2^{-1}(g_2e_1) & = & e_1 \\
(le_2)(g_2e_1) & = & e_2e_1 \end{array} \right. , $$ the first set
of equalities implies that $X_{23}U=0$ and the second set of
equalities implies that $X_{13}U=0$ (otherwise $U$ would be the
whole space $V$, which contradicts $U$ proper). Let
$S:=X_{12}+X_{23}+X_{13}$. Thus, we have shown that if $V$ is not
irreducible and if $U$ is a nontrivial proper invariant subspace of
$V$ then $SU=0$. This equality implies in turn that $det(S)=0$ as
$U$ is nonzero. Let's compute the matrix of $S$ in the basis
$\mathcal{B}$. The action of $X_{12}$ gives the first row of the
matrix, the action of $X_{23}$ the second row of the matrix and the
action of $X_{13}$ the third row of the matrix as we make it appear
on the matrix below:
$$Mat_{\mathcal{B}}S\;= \begin{array}{ll}
\;\;\;\;\;\;\;\;\;\;\;\;\;\;
e_1\;\;\;\; e_2e_1\;\;\; g_2e_1\\
\begin{array}{lll}
X_{12} \\ X_{23} \\ X_{13}\\
\end{array}
\left( \begin{array}{ccccc}
 x & & 1 & & l  \\
 1 & & x & & \frac{1}{l} \\
 \frac{1}{l} & & l & & x \\
\end{array} \;\;\right)\\
\end{array}$$
Using the defining relation $x=1-\frac{l-1/l}{1/r-r}$ and solving
the equation $det(S)=0$ in Maple, we get: $det(S)=0 \Leftrightarrow
l\in \lbrace -r^3, -1, 1, 1/r^3\rbrace $. Thus if l does not belong
to any of these values, then $V$ is irreducible. Conversely, we show
that if $V$ is irreducible then $l\not\in \lbrace -r^3, -1, 1,
1/r^3\rbrace $. Indeed, for each of these values of $l$ we claim
that $\cap_{1\leq i<j\leq 3} (Ker\;X_{ij}\cap V)$ is a nontrivial
proper invariant subspace of $V$. It suffices to show that
$\cap_{1\leq i<j\leq 3} Ker\; X_{ij}$ is a $B(A_2)$-module and that
for each of the values of $l$ above, there exists a non-zero element
in $V$ which is annihilated by all the $X_{ij}$'s. The first point
comes from corollary $2$ applied with $n=3$.
%To show the first point,
%by using the lemma $1$, we only need to check that given
%$x\in\cap_{1\leq i<j\leq 3} Ker\;X_{ij}$, we have $\forall
%l\in\left\{ 1,2\right\},\, g_lx\in KerX_{13}$. It suffices to notice
%that $X_{13}$ is a conjugate of $e_1$ and $e_2$ by the formulas
%$g_2e_1g_2^{-1}=X_{13}=g_1^{-1}e_2g_1$. From there, we have $g_2x\in
%KerX_{13}$ by using the first equality and the fact that $e_1x=0$.
%And we have $g_1x\in KerX_{13}$ by using the second equality, the
%relation $(8)$ and the fact that $x$ is annihilated by $X_{13}, e_2$
%and $e_1$.
To see the second point, we prove the following lemma:
\begin{lemma}
\begin{alignat}{4}
%\text{If} \;l & \in \{1,-1\}, & t(3):=e_1-e_2e_1 &
%\in\underset{1\leq i<j\leq
%3}{\bigcap}X_{ij}\\
\text{If}\; l &= -r^3, & y(3):=-re_1-\frac{1}{r}e_2e_1+g_2e_1 \;&
\in \underset{1\leq i<j\leq
3}{\bigcap}Ker\,X_{ij}\\
\text{If} \;l & \in \{1,-1\}, & z(3):=e_1-e_2e_1 &
\in\underset{1\leq i<j\leq
3}{\bigcap}Ker\,X_{ij}\\
\text{If}\; l &= \frac{1}{r^3}, &
t(3):=\frac{1}{r}e_1+re_2e_1+g_2e_1 & \in \underset{1\leq i<j\leq
3}{\bigcap}Ker\,X_{ij}
\end{alignat}
\end{lemma}
\noindent\textsc{Proof:} for $l=1$ or $l=-1$, it is easy to see that
$e_1-e_2e_1\in Ker\, e_1$ by using $x=1$ and the relations
$e_1^2=xe_1$ and $e_1e_2e_1=e_1$. Also, it is immediate that
$e_1-e_2e_1\in Ker\, e_2$. It remains to show that
$X_{13}(e_1-e_2e_1)=0$. On one hand we have
$X_{13}e_2=g_1^{-1}e_2g_1e_2=g_1^{-1}le_2$ by $(6)$ and so
$X_{13}e_2e_1=g_1^{-1}le_2e_1$. On the other hand we have
$X_{13}e_1=g_1^{-1}e_2g_1e_1=g_1^{-1}e_2l^{-1}e_1$. Replacing $l$ by
its value yields the desired result. If $l=-r^3$, we compute
$x=-r^2-1/r^2$. The left action by $e_1$ on $y(3)$ gives after the
use of the classical relations $(-rx-1/r+l)e_1$, which after
simplification is zero. Similarly, we obtain
$e_2.y(3)=(-r-x/r+1/l)e_2e_1$ and the coefficient in the parenthesis
is zero after replacing $x$ and $l$ by their respective values.
Let's now compute $X_{13}.y(3)$. We have:
$X_{13}g_2e_1=(g_2e_1g_2^{-1})g_2e_1=g_2e_1^2=x\,g_2e_1$. Next, from
$e_2e_1=g_1g_2e_1$, we derive $g_1^{-1}e_2e_1=g_2e_1$. Using the
formulas above that give the actions of $X_{13}$ on $e_2e_1$ and
$e_1$ respectively and replacing now $g_1^{-1}e_2e_1$ with $g_2e_1$,
we get: $X_{13}.y(3)=(-r/l-l/r+x)g_2e_1=0$. This finishes the proof
in the case $l=-r^3$. Finally, $(22)$ is obtained the same way as
$(20)$ and is left to the reader. We may now state the following
theorem:
\newtheorem{thm}{Theorem}
\begin{thm}\hfill\\
$Span_{F}(e_1,\,e_2e_1,\, g_2e_1)$ is an irreducible $B(A_2)$-module
\textit{iff} $l\not\in \lbrace -r^3, -1, 1, 1/r^3\rbrace $
\end{thm}

\section{The case n=4} Here $W$ is spanned over $F$ by
$$e_1,e_2e_1,g_2e_1,e_3e_2e_1,g_3e_2e_1,g_3g_2e_1$$ and
the six positive roots are (in the same order):
$$\al_1,\al_2,\al_1+\al_2,\al_3,\al_2+\al_3,\al_1+\al_2+\al_3$$
Generally speaking, in $B(A_{n-1})$ we will always order the roots
in the following way:
$$\al_1,\al_2,\al_2+\al_1,\al_3,\al_3+\al_2,\al_3+\al_2+\al_1,...,\al_{n-1},\al_{n-1}+\al_{n-2},...,\al_{n-1}+...+\al_1$$
Recall that one top horizontal line is enough to determine a
positive root. Geometrically, moving from the left to the right on
the top line, start with a node $\geq 2$ and join it to its left
neighbors, the closest nodes being considered first. We note and we
recall that the height of a root increases with the number of
crossings.
 The main interest of ordering the roots this way is that the action of
 $g_1,g_2,...,g_{n-2}$ on a positive root $\be$ whose support does not contain
the node $n-1$ has already been computed in the case of
$B(A_{n-2})$, as the tangle representing $\be$ in $B(A_{n-1})$ is
obtained from the tangle representing $\be$ in $B(A_{n-2})$ by
adding a vertical line on the right side, which is left invariant
with the action of the elements $g_1,g_2,...,g_{n-2}$ (such actions
don't impact the last node). As we will see later on, this ordering
of the basis will allow us to define the matrices of
$g_1,...,g_{n-2}$ by blocks, inductively on $n$.\\
We try the same method as in the case $n=3$. The difference is now
that there are two vertical braids instead of one and they can
either cross or not. Only in this part, we will denote by $W^c$ the
subspace of $V$ spanned by the tangles representing the positive
roots as in $W$ with the difference that their two vertical strings
are overcrossing. $\mathcal{B}^c$ denotes the basis in $W^c$ formed
by these elements. We look at the action of $g_1,\,g_2,\,g_3$ on
$\mathcal{B}$ and $\mathcal{B}^c$ respectively. We first deal with
the non crossed tangles. When acting by the $g_i$'s, crossings
appear. However, we still use matrices as a way of representing the
actions. We will indicate that there is a crossing by adding a c
when the crossing is over and a c' when the crossing is under, as an
exponent on the right hand side of the coefficient (those special
coefficients are indicated by boxes in the matrices below). For
instance, with our conventions,
$$g_1(g_3e_2e_1)=\mathbf{m^{c}e_1}+g_3g_2e_1-mg_3e_2e_1$$ means
$$g_1(g_3e_2e_1)=\mathbf{m\, e_1g_3}+g_3g_2e_1-mg_3e_2e_1$$ and
$$g_3(g_3g_2e_1)=g_2e_1-mg_3g_2e_1+\mathbf{(ml^{-1})^{c'}e_3e_2e_1}$$ means
$$g_3(g_3g_2e_1)=g_2e_1-mg_3g_2e_1+\mathbf{ml^{-1}\,e_3e_2e_1g_3^{-1}}$$ We gather the
results of the action of $g_1$ in the following matrix:
$$\begin{matrix}
 & \,\,e_1 \!\;\;\;\;\! e_2e_1 \;\;\;\;\!\! g_2e_1 \; e_3e_2e_1 \; g_3e_2e_1 \;
 g_3g_2e_1\\
 & \\
 & \vspace{-0.3in}\\
\begin{matrix}
e_1\\
e_2e_1\\
g_2e_1\\
e_3e_2e_1\\
g_3e_2e_1\\
g_3g_2e_1
\end{matrix} &
\begin{pmatrix} l^{-1} & m & & 0 & & 0 &  \boxed{m^c} & &  0 & \\
 0 & \!\!\!\!-m & &  1 & & 0 &  0 & & 0 & \\
 0& 1 & & 0 & & 0 &  0 & &  0 & \\
 0 & 0 & &  0 & & \boxed{1^c} &  0 & & 0 & \\
 0 &0 & & 0 & & 0 & -m & &  1 & \\
 0 & 0 & & 0 & & 0 &  1 & &  0 &
\end{pmatrix}
\end{matrix}$$

\noindent And similarly, we get the matrix corresponding to the
action of $g_2$. For instance, we have:
\begin{eqnarray*}
g_2(g_3g_2e_1) &=& g_3g_2g_3e_1\;\text{by the braid relation (1)}\\
                &=& (g_3g_2e_1)g_3 \;\text{by commutativity of}\;
                e_1\;\text{and}\;g_3,
\end{eqnarray*}
which yields the coefficient of the bottom right hand side.

$$\begin{matrix}
 &  \;\;\; e_1 \!\;\;\;\; e_2e_1 \;\;\;\; g_2e_1 \;\; e_3e_2e_1 \; g_3e_2e_1 \;
 g_3g_2e_1\\
 & \\
 & \vspace{-0.3in}\\
\begin{matrix}
e_1\\
e_2e_1\\
g_2e_1\\
e_3e_2e_1\\
g_3e_2e_1\\
g_3g_2e_1
\end{matrix} &
\begin{pmatrix} & 0 & 0 & 1  & 0 & & 0 & &  0 & \\
 & 0 & \;\,\, l^{-1} & ml^{-1}  & m & & 0 & & 0 & \\
 & 1& 0  & -m &  0 & & 0 & &  0 &  \\
 & 0 & 0 &   0 & -m & & 1 & & 0 &  \\
 & 0 &0 &  0 &  1 & & 0 & &  0 &  \\
 & 0 & 0 &  0  & 0 & & 0 & &  \boxed{1^c} &
\end{pmatrix}
\end{matrix}$$

\noindent Finally the matrix corresponding to the action of $g_3$
is:

$$\begin{matrix}
 &  \!\!\!\!\!\;\,\,\; e_1 \;\;\;\, e_2e_1 \;\; g_2e_1 \;\;\; e_3e_2e_1 \;\;\;\,\;\; g_3e_2e_1
 \;\;\;\;\;\;\,
 g_3g_2e_1\\
 & \\
 & \vspace{-0.3in}\\
\begin{matrix}
e_1\\
e_2e_1\\
g_2e_1\\
e_3e_2e_1\\
g_3e_2e_1\\
g_3g_2e_1
\end{matrix} &
\begin{pmatrix} & \boxed{1^c} & 0 & \;\,\,\; 0 & & 0 & & 0 &   0 & \\
 & 0 & 0 & \;\,\,\; 0  & &  0 & & 1 &  0 & \\
 & 0 & 0  & \;\,\,\; 0 & & 0 & & 0 &   1 & \\
 & 0 & 0 & \;\,\,\; 0 & & \;\;l^{-1}  & & \;\; ml^{-1} &  \boxed{(ml^{-1})^{c'}} & \\
 & 0 & 1 & \;\,\,\; 0 & & 0 & & \!\!\!\! -m &   0 & \\
 & 0 & 0 & \;\,\,\; 1 & & 0 & & 0 &  -m &
\end{pmatrix}
\end{matrix}$$

\noindent Note that these matrices can be obtained by using the
classical relations of the $B(A_3)$ algebra looking at the action of
the $g_i$'s on the elements of the basis of $W$, or can be obtained
by manipulating the tangles in $TM_4$. Both ways lead of course to
the same result and it is to the reader to determine which way (or
which combination of ways) he or she decides to use. While it is
easy to see the action of $g_1$ on
$e_1,e_2e_1,g_2e_1,g_3e_2e_1,g_3g_2e_1$ by using the rules of the
first page, it is direct to see the overcrossing for the tangle
$e_3e_2e_1$ after the left action by $g_1$. Because of the terms in
the boxes, we see that $W$ is not a $B(A_3)$-module. An idea is to
look at the action of the $g_i$'s on the crossed roots
$e_1g_3,\,e_2e_1g_3,\,g_2e_1g_3,\,e_3e_2e_1g_3,\,g_3e_2e_1g_3,\,g_3g_2e_1g_3$
which span $W^c$ to see in turn what we get and exhibit from there
an invariant subspace. Obviously, the coefficients which are not in
the boxes are unchanged since all we do is make the words bigger by
adding to them a factor $g_3$ on the right (which makes the crossing
appear) and the action by the $g_i$'s on these words takes place on
the left.
%Direct computations show that the terms in the boxes must be
%replaced by                 in the order that they appear on the
%matrices.
Next, we have:
\begin{eqnarray}
g_1(e_3e_2e_1g_3) &=& e_3(g_1e_2e_1)g_3\;\;\text{by commutativity of $e_1$ and $g_3$}\notag \\
                &=& e_3(g_2e_1+m(e_1-e_2e_1))g_3 \;\;\text{by
                $(12)$}\notag\\
                &=& e_3g_2e_1g_3-m\,e_3e_2e_1g_3\;\;\text{since
                $e_3e_1=0$ in $I_1/I_2$}\notag\\
                &=& e_3g_2g_3e_1-m\,e_3e_2e_1g_3\;\;\text{by
                commutativity of $e_1$ and $g_3$}\notag\\
                &=& \mathbf{e_3e_2e_1-m\,e_3e_2e_1g_3}
\end{eqnarray}

\noindent The last equality is obtained by using $(11)$ and the
anti-involution on products of generators of $B(A_4)$:
$$g_{i_1}...g_{i_q}\longmapsto g_{i_q}...g_{i_q}$$
We also have:
\begin{eqnarray}
g_1(g_3e_2e_1g_3) &=& \mathbf{m(e_1-me_1g_3)}-mg_3e_2e_1g_3+g_3g_2e_1g_3\\
g_2(g_3g_2e_1g_3) &=& \mathbf{g_3g_2e_1-mg_3g_2e_1g_3}\\
g_3(e_1g_3) &=& \mathbf{e_1-me_1g_3}
\end{eqnarray}

\noindent We notice that:\begin{align*} 1^c \;\; & \text{is replaced
by}\;
"noncrossed-m\,crossed"\;\text{as in}\; (23),(25)\;\text{and}\; (26)\\
m^c \;\; & \text{is replaced by}\; "m(noncrossed-m\, crossed)"
\;\text{as in}\; (24)
\end{align*}
\noindent And indeed, we look at the left action of the $g_i$'s on
words of the form $wg_3$. For such prefix $w$ we have established
that: \begin{eqnarray} \boxed{g_iw=\la \;w^{'}g_3+s} \end{eqnarray}
where $i$ is adequately picked, with $\la $ the appropriate
coefficient and with $s$ a sum of noncrossed terms (eventually
zero). Now we have:
\begin{eqnarray*}
g_i(wg_3)&=& (g_iw)g_3 \qquad\qquad\qquad\qquad\qquad\;\;\;\:\qquad\:\:\: \text{by associativity}\\
         &=& \la \;w^{'}g_3^2\,+\,s\,g_3 \qquad\qquad\qquad\:\:\;\qquad\qquad\; \text{by replacing with (27)}\\
         &=& \la (w^{'}-m\,w^{'}g_3 +\,ml^{-1}w^{'}e_3) +
         s\,g_3\;\;\;\;\;\:\:\,
         \text{by application of}\; (8)
\end{eqnarray*}
\noindent But $w^{'}$ belongs to $W$ and has its nodes $1$ and $2$
joined on the bottom line. Then, the tangle obtained by multiplying
$w^{'}$ with $e_3$ has nodes $1$ and $2$, $3$ and $4$ respectively
joined on the bottom line, hence is zero in $I_1/I_2$. So we finally
get: \begin{eqnarray}\boxed{g_i(wg_3)= \la (w^{'}-m\,w^{'}g_3)
+\,s\,g_3}\end{eqnarray} \noindent Thus, $$ \la^c\;\; \text{is
replaced by}\;\; \la\,(noncrossed-m\; crossed) $$ \noindent It
remains to compute
\begin{eqnarray*}
g_3(g_3g_2e_1g_3) &=&
g_2e_1g_3-mg_3g_2e_1g_3+ml^{-1}e_3g_2e_1g_3\;\;\text{by the rule}
\;(8)\\
&=& g_2e_1g_3-mg_3g_2e_1g_3+ml^{-1}\;e_3e_2e_1
\end{eqnarray*}
\noindent to see how the last box is being modified. Here the last
last term of the last equation is obtained from $e_3g_2e_1g_3$ by
commuting $e_1$ and $g_3$, then applying the rule $(11)$ with the
anti-involution described above. Note that another way to see it is
to notice that the term in the box is nothing else but
$ml^{-1}e_3e_2e_1g_3^{-1}$, where $g_3^{-1}$ makes the undercrossing
appear. Then acting by $g_3$ on
the left yields the non crossed term $ml^{-1}e_3e_2e_1$.\\
We summarize our results in the matrices below. Like in the non
crossed case, these matrices are not mathematical objects, but are
used
as a convenient way of representing the actions of the $g_i$'s.\\\\
For the action of $g_1$ on the elements of $\mathcal{B}^c$, we get:

$$\begin{matrix}
 & \;\;\;\;\;\, e_1g_3 \;\;\, e_2e_1g_3 \,\,g_2e_1g_3 \,\,e_3e_2e_1g_3 \;\;\;\;\;\; g_3e_2e_1g_3
 \;\;
 g_3g_2e_1g_3\\
 & \\
 & \vspace{-0.3in}\\
\begin{matrix}
e_1g_3\\
e_2e_1g_3\\
g_2e_1g_3\\
e_3e_2e_1g_3\\
g_3e_2e_1g_3\\
g_3g_2e_1g_3
\end{matrix} &
\begin{pmatrix} & l^{-1} & & m & &  0 &  0 &  \boxed{m(1^{nc}-m)} &   0 & \\
 & 0  & & \!\!\!\!-m & & 1 &  0 &  0 &  0 & \\
 & 0& & 1 & & 0 &  0 &  0 &   0 & \\
 & 0  & & 0 & & 0 & \boxed{1^{nc}-m} &  0 &  0 & \\
 & 0 & & 0 & & 0 &  0 & -m &  1 & \\
 & 0  & & 0 & & 0 &  0 &  1 &  0 &
\end{pmatrix}
\end{matrix}$$

\noindent And similarly we represent the action of $g_2$ on
$\mathcal{B}^c$ by the matrix

$$\begin{matrix}
 & \,e_1g_3\,\, e_2e_1g_3 \,\,\;\;\; g_2e_1g_3 \,\, e_3e_2e_1g_3 \,\, g_3e_2e_1g_3
 \,\,
 g_3g_2e_1g_3\\
 & \\
 & \vspace{-0.3in}\\
\begin{matrix}
e_1g_3\\
e_2e_1g_3\\
g_2e_1g_3\\
e_3e_2e_1g_3\\
g_3e_2e_1g_3\\
g_3g_2e_1g_3
\end{matrix} &
\begin{pmatrix}  0 & \;\;\;0 & & 1 & & 0 & &  0 & & 0  \\
  0 & \;\;\;\;\,\, l^{-1}& & ml^{-1} & & m & &  0 & & 0  \\
  1& \;\;\; 0 & & -m & & 0 & &  0 & & 0  \\
 0 & \;\;\; 0 & &  0 & & -m & &  1 & & 0  \\
  0 & \;\;\; 0 & & 0 & & 1 & &  0 & & 0  \\
  0 & \;\;\; 0 & & 0 & & 0 & &  0 & & \boxed{1^{nc}-m}
\end{pmatrix}
\end{matrix}$$

\noindent and we represent the action of $g_3$ on $\mathcal{B}^c$ by
the matrix

$$\begin{matrix}
 & \;\;\;\;\;\; e_1g_3 \,\;\; e_2e_1g_3 \;\; g_2e_1g_3 \;\; e_3e_2e_1g_3 \;\;
 g_3e_2e_1g_3
 \;\;
 g_3g_2e_1g_3\\
 & \\
 & \vspace{-0.3in}\\
\begin{matrix}
e_1g_3\\
e_2e_1g_3\\
g_2e_1g_3\\
e_3e_2e_1g_3\\
g_3e_2e_1g_3\\
g_3g_2e_1g_3
\end{matrix} &
\begin{pmatrix} & \boxed{1^{nc}-m} & 0 & & \;\,\,\; 0 & & 0 & & 0 &   0 \\
 & 0 & 0 & & \;\,\,\; 0  & & 0 & & 1 &  0 \\
 & 0 & 0  & & \;\,\,\; 0 & & 0 & & 0 &   1  \\
 & 0 & 0 & & \;\,\,\; 0 & & \;\;l^{-1} & & \;\; ml^{-1} &  \boxed{(ml^{-1})^{nc}}  \\
 & 0 & 1 & & \;\,\,\; 0 & & 0 & & \!\!\!\! -m &   0 \\
 & 0 & 0 & & \;\,\,\; 1 & & 0 & & 0 &  -m
\end{pmatrix}
\end{matrix}$$

\noindent Let's now consider the subspace $W_r$ of $W$ defined as
follows:
\begin{equation*}
\begin{split}
W_r:=\text{Span}_F(e_1 & +r\,e_1g_3, \, e_2e_1+re_2e_1g_3,\,
g_2e_1+r\,g_2e_1g_3,\\ & e_3e_2e_1+r\,e_3e_2e_1g_3, \,
g_3e_2e_1+r\,g_3e_2e_1g_3,\, g_3g_2e_1+r\,g_3g_2e_1g_3)
\end{split}
\end{equation*}
$W_r$ is spanned over $F$ by linear combinations of non crossed
tangles spanning $W$ and their corresponding crossed tangles
spanning $W^c$. The spanning elements above form a basis of $W_r$
that we denote by $\mathcal{B}_r$. We will show that $W_r$ is an
invariant subspace. To that aim, let's consider the linear
combination of non crossed tangle and crossed tangle $w+r\,wg_3$.
Using $(27)$ and $(28)$ above, we compute:
\begin{eqnarray}
g_i(w+rwg_3) &=& \la (r\, w^{'}+(1-mr)w^{'}g_3)+(s\,+\,r\,sg_3)\notag \\
             &=& \la r(\underbrace{w^{'}+r\,w^{'}g_3}_{\in W_r})+(\underbrace{s\,+\,r\,sg_3}_{\in
             W_r})
\end{eqnarray} Further, we have:
\begin{equation}\begin{split}
g_3(g_3g_2e_1+r\,g_3g_2e_1g_3)=(g_2e_1
+r&\,g_2e_1g_3)-m(g_3g_2e_1+r\,g_3g_2e_1g_3)\\
& +ml^{-1}(e_3e_2e_1g_3^{-1}+re_3e_2e_1)
\end{split}\end{equation}
But using the tangle formula:
\begin{eqnarray} T^{-}=T^{+}+m(T^{\circ}-T^{\infty})\end{eqnarray} which is one of the
defining relations of the Kauffman's tangle algebra (see \cite{MW}),
we can express $e_3e_2e_1g_3^{-1}$ as a sum of two similar tangles
as follows:
$$e_3e_2e_1g_3^{-1}=e_3e_2e_1g_3+m\,e_3e_2e_1$$ \noindent $T^{\infty}$ being
zero for this tangle as we work in $I_1/I_2$. It follows that $(30
)$ can be rewritten in the following way:
\begin{equation*}\begin{split}
g_3(g_3g_2e_1+r\,g_3g_2e_1g_3)=(g_2e_1
+r&\,g_2e_1g_3)-m(g_3g_2e_1+r\,g_3g_2e_1g_3)\\
& +ml^{-1}((m+r)e_3e_2e_1+ e_3e_2e_1g_3)
\end{split}\end{equation*}
Recall that by definition, $m=r^{-1}-r$, so that $m+r=r^{-1}$.
Hence, we finally get:
\begin{equation}\begin{split}
g_3(g_3g_2e_1+r\,g_3g_2e_1g_3)=(&\overbrace{g_2e_1
+r\,g_2e_1g_3}^{\in W_r})-m(\overbrace{g_3g_2e_1+r\,g_3g_2e_1g_3}^{\in W_r})\\
& +\frac{ml^{-1}}{r}(\underbrace{e_3e_2e_1+ r\,e_3e_2e_1g_3}_{\in
W_r})
\end{split}\end{equation}
Since the columns that don't contain any box in the matrices are the
same in both cases crossed and non crossed, we conclude that $W_r$
is an invariant subspace of $V$ as announced. Let's give the
matrices of the left actions by $g_1$, $g_2$ and $g_3$ in the basis
$\mathcal{B}_r$. We denote them by $G_1(4)$, $G_2(4)$ and $G_3(4)$
respectively. We need to replace the coefficients in the boxes by
appropriate new ones. And in fact, we see that
\begin{align}
&\text{by}\; (29), \;\boxed{\la^c\;\;\text{must be replaced by}\;\;
\la\,r}\\ \text{and}\notag
\\ & \text{by}\;(32),\;\boxed{(ml^{-1})^{c'}\;\;\text{must be replaced
by}\;\;\frac{ml^{-1}}{r}}
\end{align}
We obtain the following matrices:
$$G_1(4):=
\begin{pmatrix} l^{-1} & m &  0 & 0 & m\,r & 0 \\
 0 & \!\!\!\!-m & 1 & 0 & 0 & 0 \\
 0 & 1 & 0 & 0 & 0 & 0 \\
 0 & 0 & 0 & r & 0 & 0 \\
 0 & 0 & 0 & 0 & -m & 1 \\
 0 & 0 & 0 & 0 & 1 & 0
\end{pmatrix}$$

$$G_2(4):= \begin{pmatrix}  0 & 0 & 1 & 0 & 0 & 0 \\
 0 & \;\,\, l^{-1} & \;\;ml^{-1}  & m & 0 & 0 \\
 1& 0  & \!\!\!-m &  0 & 0 & 0  \\
 0 & 0 &   0 & -m & 1 & 0 \\
 0 & 0 &  0 &  1 & 0 & 0  \\
 0 & 0 &  0 & 0 & 0 & r
\end{pmatrix}$$

$$G_3(4):= \begin{pmatrix} r & 0 &  0 & 0 & 0 & 0  \\
  0 & 0 & 0 & 0 & 1 & 0 \\
  0 & 0 & 0 & 0 & 0 & 1 \\
  0 & 0 & 0 & \;\;l^{-1}  & \;\; ml^{-1} & \frac{ml^{-1}}{r} \\
  0 & 1 & 0 & 0 & \!\!\!\! -m &   0 \\
  0 & 0 & 1 & 0 & 0 & -m
\end{pmatrix}$$

\noindent We verify that the braid relations $(1)$ and $(2)$ of the
front page are satisfied on these matrices and using these matrices,
we define new matrices $E_i(4)$ for $i\in\{1,2,3\}$. They are given
by the formula:
\begin{eqnarray}
E_i(4)=\frac{l}{m}(G_i^2(4)+mG_i(4)-I_6)
\end{eqnarray}
where $I_6$ is the identity matrix of size $6$. We check in turn
that the relations $(4)$, $(5)$ and $(6)$ are satisfied on the
$G_i(4)$'s and $E_i(4)$'s. Hence we get a representation of $B(A_3)$
in $W_r$ which is defined on the generators $e_1$, $e_2$, $e_3$,
$g_1$, $g_2$, $g_3$ of the BMW algebra $B(A_3)$ by the map:
$$\begin{array}{ccc}
B(A_3)& \lra & \mathcal{M}(6,F)\\
g_1,\,g_2,\,g_3 & \longmapsto & G_1(4),\,G_2(4),\,G_3(4)\\
e_1,\,e_2,\e_3 & \longmapsto & E_1(4),\,E_2(4),\,E_3(4)
\end{array}$$
Suppose that $W_r$ is not irreducible. Let $U$ be a non trivial
proper invariant subspace of $W_r$. We claim that the $X_{ij}$'s,
$1\leq i<j\leq 4$ act trivially on $U$. Indeed, suppose that there
exists a non zero element $u$ in $U$ such that $X_{ij}u\neq 0$.
Since the tangle resulting from the product $X_{ij}u$ is in $W_r$
and is a linear combination of tangles having their nodes $i$ and
$j$ joined at the top, with the edge $(ij)$ overcrossing the
eventual vertical lines that it intersects as in $X_{ij}$, and their
nodes $1$ and $2$ joined at the bottom, it must be proportional to
the element of $\mathcal{B}_r$ having the same horizontal lines.
Let's rename the six elements of the basis $\mathcal{B}_r$:

$$\mathcal{B}_r=(x_{\al_1},\,x_{\al_2},\,x_{\al_2+\al_1},\,x_{\al_3},\,x_{\al_3+\al_2},\,x_{\al_3+\al_2+\al_1})$$

\noindent So $X_{ij}u$ is proportional to an $x_{\be}$ where $\be$
is one of the six positive roots. Recall that by part $2$, an
expression for $x_{\be}$ is:
\begin{eqnarray}
x_{\be}&=&\varphi(w_{\be 1})e_1+r\varphi(w_{\be 1})e_1g_3\notag \\
       &=&\varphi(w_{\be 1})(e_1+re_1g_3)\notag \\
       &=&\varphi(w_{\be 1})x_{\al_1}
\end{eqnarray}
We see with $(36)$ that any $x_{\be}$ can be obtained from
$x_{\al_1}$ by multiplying $x_{\al_1}$ by the left with an element
in $B(A_3)$. Moreover, $\varphi(w_{\be 1})$ is invertible as it is a
word containing only $g_i$'s by construction. We conclude that if
$x_{\be}$ is in $U$, then $x_{\al_1}$ is in $U$ as $U$ is a
$B(A_3)$-module and then all the other $x_{\gamma}$ with
$\gamma\in\phi^{+}\setminus\{\al_1,\,\be\}$ are also in $U$ by
above. Now if $X_{ij}u$ is non zero and proportional to $x_{\be}$,
then $x_{\be}$ is in $U$ and by what preceeds all the elements of
the basis $\mathcal{B}_r$ are also in $U$. It implies that $U=W_r$,
which contradicts the fact that $U$ is proper. Thus, $X_{ij}u=0$ for
all $u\in U$ and since the argument can be applied to any of the
$X_{ij}$'s, we have actually proved that all the $X_{ij}$'s
annihilate $U$. Their sum also annihilates $U$. Let's consider their
sum:
$$S:=X_{12}+X_{13}+X_{14}+X_{23}+X_{24}+X_{34}$$ It comes: $$S\,U=0$$
Since $U$ is non zero, we must have $det(S)=0$. We will compute the
matrix of the left action of $S$ in the basis $\mathcal{B}_r$ and
then compute its determinant. Note that each row in this matrix
corresponds to the action of one of the $X_{ij}$'s. The matrix can
be obtained by calculating the actions of the $X_{ij}$'s directly on
the tangles
$x_{\al_1},\,x_{\al_2},\,x_{\al_2+\al_1},\,x_{\al_3},\,x_{\al_3+\al_2},\,x_{\al_3+\al_2+\al_1}$
or by using Maple, the expressions of $G_1(4),\,G_2(4),\,G_3(4)$,
the formula $(35)$, the formula $(9)$ and the conjugation formulas
for $j>i+1$:
$$X_{ij}=g_{j-1}\dots g_{i+1}\,e_i\,g_{i+1}^{-1}\dots g_{j-1}^{-1}$$
We get the matrix:
$$Mat_{\mathcal{B}_r}(S)\;= \begin{array}{ll}
\;\;\;\;\;\;\;\;\;\;\;\;\;\;
x_{\al_1}\, x_{\al_2}\;\;\;\;\;\;\;\;\; x_{\al_1+\al_2}\;\;\;\;\;\;\;x_{\al_3}\;\;\;\;\;\;x_{\al_3+\al_2}\;x_{\al_3+\al_2+\al_1}\\
\begin{array}{llllll}
X_{12} \\ X_{23} \\ X_{13}\\X_{34}\\X_{24}\\X_{14}
\end{array}
\left( \begin{array}{cccccc}
 x & 1 & l & 0 & r & lr\\
 1 & x & \frac{1}{l} & 1 & l & 0\\
 \frac{1}{l} & l & x & r & (\frac{1}{r}-r)(r-l) & l\\
 0 & 1 & \frac{1}{r} & x & \frac{1}{l} & \frac{1}{lr}\\
 \frac{1}{r} & \frac{1}{l} & (r-\frac{1}{r})(\frac{1}{r}-\frac{1}{l}) & l & x & \frac{1}{l}\\
 \frac{1}{lr} & 0 & \frac{1}{l} & lr & l & x
\end{array} \;\;\right)\\
\end{array}$$

\noindent Using Maple to compute the determinant of this $6\times 6$
matrix, we get that $det(S)$ is zero if and only if
$l\in\{r,\,-r^3,\,\frac{1}{r},\,-\frac{1}{r},\,\frac{1}{r^5}\}$ and
we conclude that if $W_r$ is not irreducible then $l$ must be one of
these values. Conversely, we show that if $l$ takes one of these
values, then $W_r$ is not irreducible. It suffices to exhibit a
nontrivial proper invariant subspace of $W_r$ for each of these
values of $l$. In what follows, we will denote by $T_4$ the set
$\{r,\,-r^3,\,\frac{1}{r},\,-\frac{1}{r},\,\frac{1}{r^5}\}$
\newtheorem{Prop}{Proposition}
\begin{Prop}
If $l$ belongs to the set $T_4$ then $\underset{1\leq i<j\leq
4}{\bigcap} (Ker\,X_{ij}\cap W_r)$ is a non trivial proper invariant
subspace of $W_r$.
\end{Prop}
\noindent\textsc{Proof:} by corollary $(2)$ applied with $n=4$, we
know that $\underset{1\leq i<j\leq 4}{\bigcap} Ker\,X_{ij}$ is a
$B(A_3)$-module. Hence $\underset{1\leq i<j\leq
4}{\bigcap}(Ker\,X_{ij}\cap W_r)$ is an invariant subspace of $W_r$.
It is obviously not $W_r$ itself since for instance
$e_1+r\,e_1g_3\in W_r$, but $e_2(e_1+r\,e_1g_3)=\neq 0$, so that
$e_1+r\,e_1g_3\not\in\,Ker\,e_2$, which implies
$e_1+r\,e_1g_3\not\in\underset{1\leq i<j\leq 4}{\bigcap} (Ker
\,X_{ij}\cap W_r)$. It remains to show that $\underset{1\leq i<j\leq
4}{\bigcap}(Ker\,X_{ij}\cap W_r)$ is non trivial. To this aim, we
show the following lemma that exhibits for each value of $l$ in
$T_4$ an element in $W_r$ that is annihilated by all the $X_{ij}$'s:

\begin{lemma}
\hfill\\\\
Let $\,x(4)\:\!\,:=r^2\,x_{\al_1+\al_2}+x_{\al_2+\al_3}-r\,x_{\al_1+\al_2+\al_3}-r\,x_{\al_2}$\\\\
and $y(4)\,:=-r^2\,x_{\al_1}-\frac{1}{r^2}x_{\al_3}+x_{\al_1+\al_2+\al_3}-x_{\al_2}$\\\\
and $z(4)\,:=x_{\al_1}+x_{\al_3}+x_{\al_1+\al_2+\al_3}-x_{\al_2}$\\\\
and $z^{'}(4)\!:=
x_{\al_1+\al_2+\al_3}-x_{\al_2}$\\\\
and $t(4)\,:=
\frac{1}{r^2}x_{\al_1}+r^2x_{\al_3}+\frac{1}{r}x_{\al_1+\al_2}+r\,x_{\al_2+\al_3}+x_{\al_1+\al_2+\al_3}+x_{\al_2}$\\
\\We have for the distinct values of $l$:
\begin{alignat}{5}
\text{If}\;& l=r& \;\text{then}\; & x(4)& \in \underset{1\leq
i<j\leq
4}{\bigcap} Ker\,X_{ij}\\
\text{If}\; & l=-r^3 & \;\text{then}\; & y(4)& \in \underset{1\leq
i<j\leq
4}{\bigcap} Ker\,X_{ij}\\
\text{If}\;& l=-\frac{1}{r}& \;\text{then}\; & z(4)& \in
\underset{1\leq i<j\leq
4}{\bigcap} Ker\,X_{ij}\\
\text{If}\;& l=\frac{1}{r}&\; \text{then}\; & z'(4)& \in
\underset{1\leq i<j\leq
4}{\bigcap} Ker\,X_{ij}\\
\text{If}\;& l=\frac{1}{r^5}&\; \text{then}\; & t(4)& \in
\underset{1\leq i<j\leq 4}{\bigcap} Ker\,X_{ij}
\end{alignat}
\end{lemma}
\noindent\textsc{Proof of the lemma:} Since each row of the matrix,
say $M(l)$, of the left action of $S$ on the basis $\mathcal{B}_r$
corresponds to the action of exactly one of the $X_{ij}$'s, it
suffices to verify that $M(r)X=0$, $M(-r^3)Y=0$, $M(-1/r)Z=0$,
$M(1/r)Z'=0$ and $M(1/r^5)T=0$ where $X$, $Y$, $Z$, $Z'$, $T$ are
the column vectors respectively corresponding to $x(4)$, $y(4)$,
$z(4)$, $z'(4)$, $t(4)$ in the basis $\mathcal{B}_r$.

$$\begin{array}{cc}
\begin{pmatrix}
2 & 1 & r & 0 & r & r^2\\
1 & 2 & \frac{1}{r} & 1 & r & 0\\
\frac{1}{r} & r & 2 & r & 0 & r\\
0 & 1 & \frac{1}{r} & 2 & \frac{1}{r} & \frac{1}{r^2}\\
\frac{1}{r} & \frac{1}{r} & 0 & r & 2 & \unsurr \\
\frac{1}{r^2} & 0 & \unsurr & r^2 & r & 2\\
\end{pmatrix} &\!\!\!\!\!\!\!
\begin{pmatrix}
0\\ \!\!-r\\ \,r^2\\ 0\\1\\-r
\end{pmatrix}
\end{array}=0$$

$$\begin{array}{cc}
\begin{pmatrix}
-\frac{1}{r^2}-r^2 & 1 & -r^3 & 0 & r & -r^4\\
1 & -\frac{1}{r^2}-r^2 & -\frac{1}{r^3} & 1 & -r^3 & 0\\
-\frac{1}{r^3} & -r^3 & -\frac{1}{r^2}-r^2 & r & 1-r^4 & -r^3\\
0 & 1 & \unsurr & -\frac{1}{r^2}-r^2 & -\frac{1}{r^3} &
-\frac{1}{r^4}\\
\unsurr & -\frac{1}{r^3} & 1-\frac{1}{r^4} & -r^3
& -\frac{1}{r^2}-r^2 & -\frac{1}{r^3}\\
-\frac{1}{r^4} & 0 & -\frac{1}{r^3} & -r^4 & -r^3 &
-\frac{1}{r^2}-r^2
\end{pmatrix} &\!\!\!\!\!\!\!
\begin{pmatrix}
-r^2\\ \!\!-1\\0\\-\frac{1}{r^2}\\0\\1
\end{pmatrix}
\end{array}=0$$

$$\begin{array}{cc}
\begin{pmatrix}
2 & 1 & -\unsurr & 0 &\;\,\,\, r & -1\\
1 & 2 & -r & 1 & -\unsurr & \;0\\
-r & -\unsurr & 2 & r & \frac{1}{r^2}-r^2 & -\unsurr \\
0 & 1 & \frac{1}{r} & 2 & -r & -1\\
\frac{1}{r} & -r & r^2-\frac{1}{r^2} & -\unsurr & \;2 & -r \\
-1 & 0 & -r & -1 & \!\!-\unsurr & 2\\
\end{pmatrix} &\!\!\!\!\!\!\!
\begin{pmatrix}
1\\ \!\!-1\\0\\1\\0\\1
\end{pmatrix}
\end{array}=0$$

$$\begin{array}{cc}
\begin{pmatrix}
0 & 1 & \unsurr & 0 & r & 1\\
1 & 0 & r & 1 & \unsurr & \;0\\
r & \unsurr & 0 & r & 2-\frac{1}{r^2}-r^2 & \unsurr \\
0 & 1 & \frac{1}{r} & 0 & r & 1\\
\frac{1}{r} & r & 0 & \unsurr & 0 & r \\
1 & 0 & r & 1 & \unsurr & 0\\
\end{pmatrix} &\!\!\!\!\!\!\!
\begin{pmatrix}
0\\ \!\!-1\\0\\0\\0\\1
\end{pmatrix}
\end{array}=0$$

$$\begin{array}{cc}
\begin{pmatrix}
z & 1 & \frac{1}{r^5} & 0 & r & \frac{1}{r^4}\\
1 & z & r^5 & 1 & \unsur{r^5} & 0\\
r^5 & \unsur{r^5} & z & r & 1-r^2+\unsur{r^4}-\unsur{r^6} & \unsur{r^5}\\
0 & 1 & \unsurr & z & r^5 &
r^4\\
\unsurr & r^5 & 1-\frac{1}{r^2}+r^4-r^6 & \unsur{r^5}
& z & r^5\\
r^4 & 0 & r^5 & \unsur{r^4} & \unsur{r^5} & z
\end{pmatrix} &\!\!\!\!\!\!\!
\begin{pmatrix}
\unsur{r^2}\\ 1 \\ \unsur{r}\\r^2\\r\\1
\end{pmatrix}
\end{array}=0$$

\noindent\\ where in the last matrix
$z:=-\big(r^4+\unsur{r^4}+r^2+\unsur{r^2}\big)$. These equalities
end the proof of the lemma. We are now able to conclude. By
proposition $1$, if $l$ is in $T_4$ then $W_r$ is not irreducible as
it contains a non trivial proper invariant subspace. Conversely, we
have seen above that if $W_r$ is not irreducible then $l$ must
belong to $T_4$. We summarize this result in a theorem:

\begin{thm}
$Span_F(x_{\al_1},x_{\al_2},x_{\al_2+\al_1},x_{\al_3},x_{\al_3+\al_2},x_{\al_3+\al_2+\al_1})$
is an irreducible $B(A_3)$-module if and only if
$l\not\in\{r,-r^3,-\unsurr,\unsurr,\unsur{r^5}\}$
\end{thm}

\section{The case n=5} There are ten positive roots:
\begin{equation*}
\begin{split}
\phi^{+}=\{\al_1,\, \al_2,\,& \al_2+\al_1,\, \al_3,\, \al_3+\al_2,\,
\al_3+\al_2+\al_1,\\& \al_4,\, \al_4+\al_3, \,\al_4+\al_3+\al_2,\,
\al_4+\al_3+\al_2+\al_1\}
\end{split}
\end{equation*}
corresponding to the following tangles in $V$:
\begin{equation*}
\begin{split}
e_1,\,e_2e_1,\,& g_2e_1,\,e_3e_2e_1,\,g_3e_2e_1,g_3g_2e_1,\\
& e_4e_3e_2e_1,\,g_4e_3e_2e_1,\,g_4g_3e_2e_1,g_4g_3g_2e_1
\end{split}
\end{equation*}

\noindent The first six tangles are obtained from the same six
tangles in $TM_4$ by adding a straight vertical line on the right
hand side. Hence the actions of $g_1,g_2,g_3$ (which affect only the
first four nodes) are the same as in the case $n=4$. We only need to
compute the actions of $g_1,g_2,g_3$ on the four remaining tangles.
But we need to do all the computations for the action of $g_4$.
Crossings eventually appear between the three vertical braids. We
represent them by means of permutations that are either transposes
or cycles of lengths three. We will again denote those crossed
tangles by their non crossed analogues, but we will specify the
permutation corresponding to the vertical braids by an exponent on
the right hand side of the algebra element. When the permutation is
the identity, the exponent will be omitted. Moreover and unless
otherwise mentioned, all the crossings are over. When there is one
undercrossing (on the vertical braids), we will indicate it with a
prime on the transpose and if there are two undercrossings, we will
indicate them with a double prime on the $3$-cycle. Let's compute
the action of $g_1$ on the four tangles involving $\al_4$:

\begin{eqnarray}
g_1(e_4e_3e_2e_1)&=& e_4e_3g_1e_2e_1\;\;\;\;\;\;\;\;\;\;\;\;\;\;\;\;\;\;\;\text{by $(1)$}\notag\\
&=& e_4e_3g_2e_1-me_4e_3e_2e_1\,\,\,\text{by $(12)$ and the fact
that
$e_3e_1=0$}\notag\\
&=&
e_4e_3e_2e_1^{(12)}\;\;\;\;\;\;\;\;\;\;\;\;\;\;\;\;\;\;\;\!\,\text{by
$(31)$}\\
g_1(g_4e_3e_2e_1)&=&g_4e_3e_2e_1^{(12)}\\
%g_1(g_4e_3e_2e_1)&=& g_4e_3g_1e_2e_1\;\;\;\;\;\;\;\;\;\;\;\;\;\;\;\;\;\;\;\text{by $(1)$}\notag\\
%&=& g_4e_3g_2e_1-mg_4e_3e_2e_1\,\,\,\text{by $(12)$ and the fact
%that
%$e_3e_1=0$}\notag\\
%&=&
%g_4e_3e_2e_1^{(12)}\;\;\;\;\;\;\;\;\;\;\;\;\;\;\;\;\;\;\;\!\,\text{by
%$(31)$}\\
g_1(g_4g_3e_2e_1)&=& g_4g_3g_1e_2e_1\;\;\;\;\;\;\qquad\qquad\qquad\qquad\;\;\,\,\text{by $(1)$}\notag\\
&=& g_4g_3g_2e_1 + m\,g_4g_3e_1-m\,g_4g_3e_2e_1\;\;\text{by $(12)$}\notag\\
&=& g_4g_3g_2e_1 + m\,e_1^{(123)}-m\,g_4g_3e_2e_1\;\;\,\,\text{(rewriting)}\\
g_1(g_4g_3g_2e_1)&=& g_4g_3g_1g_2e_1\;\;\;\;\text{by $(1)$}\notag\\
&=& g_4g_3e_2e_1\qquad\;\;\!\!\!\text{by $(10)$}
\end{eqnarray}
\\\noindent $(43)$ was not detailed because it is obtained the same way as
$(42)$. Let's do the same with $g_2$:
\begin{eqnarray}
g_2(e_4e_3e_2e_1)&=& e_4g_2e_3e_2e_1\qquad\qquad\;\;\;\;\,\text{by $(1)$}\notag\\
&=& e_4g_3e_2e_1-m\,e_4e_3e_2e_1\;\text{by $(12)$ and the fact that
$e_4e_2=0$}\notag\\
&=& e_4e_3e_2e_1^{(23)}\qquad\qquad\;\;\;\,\,\text{by $(31)$}\\
g_2(g_4e_3e_2e_1)&=& g_4g_2e_3e_2e_1\qquad\qquad\qquad\qquad\qquad\;\text{by $(1)$}\notag\\
&=& g_4g_3e_2e_1+m\,g_4e_2e_1-m\,g_4e_3e_2e_1\;\text{by
$(12)$}\notag\\
&=& g_4g_3e_2e_1+m\,e_2e_1^{(23)}-m\,g_4e_3e_2e_1\;\text{(rewriting)}\\
g_2(g_4g_3e_2e_1)&=& g_4g_2g_3e_2e_1\qquad\text{by $(1)$}\notag\\
&=& g_4e_3e_2e_1\qquad\;\;\;\;\!\!\,\text{by $(10)$}\\
g_2(g_4g_3g_2e_1)&=& g_4g_2g_3g_2e_1\qquad\text{by $(1)$}\notag\\
&=& g_4g_3g_2g_3e_1\qquad\text{by $(2)$}\notag\\
&=& g_4g_3g_2e_1^{(12)}\qquad\!\text{(rewriting)}
\end{eqnarray}
\noindent We have for $g_3$:
\begin{eqnarray}
g_3(e_4e_3e_2e_1)&=& g_4e_3e_2e_1+m\,e_3e_2e_1-m\,e_4e_3e_2e_1\;\;\text{by $(12)$}\\
g_3(g_4e_3e_2e_1)&=& e_4e_3e_2e_1\qquad\;\;\;\;\!\text{by $(10)$}\\
g_3(g_4g_3e_2e_1)&=& g_4g_3g_4e_2e_1\qquad\text{by $(2)$}\notag\\
&=& g_4g_3e_2e_1g_4\qquad\text{by $(1)$}\notag\\
&=& g_4g_3e_2e_1^{(23)}\qquad\!\text{(rewriting)}\\
g_3(g_4g_3g_2e_1)&=&g_4g_3g_4g_2e_1\qquad\text{by $(2)$}\notag\\
&=& g_4g_3g_2e_1g_4\qquad\text{by $(1)$}\notag\\
&=& g_4g_3g_2e_1^{(23)}\qquad\!\text{(rewriting)}
\end{eqnarray}
\noindent Let's now compute the action of $g_4$ on all the elements
of the basis $\mathcal{B}$ of $W$. We have:
\begin{eqnarray}
g_4e_1&=& e_1^{(23)}\;\;\;\;\;\;\;\qquad\text{(rewriting)}\\
g_4(e_2e_1)&=&e_2e_1^{(23)}\;\;\;\;\qquad\text{(rewriting)}\\
g_4(g_2e_1)&=&g_2e_1^{(23)}\;\;\;\;\qquad\text{(rewriting)}\\
g_4(e_3e_2e_1)&=& g_4e_3e_2e_1\;\qquad\text{(unchanged)}\\
g_4(g_3e_2e_1)&=& g_4g_3e_2e_1\;\qquad\text{(unchanged)}\\
g_4(g_3g_2e_1)&=& g_4g_3g_2e_1\;\qquad\text{(unchanged)}\\
g_4(e_4e_3e_2e_1)&=& l^{-1}e_4e_3e_2e_1\;\;\;\;\text{by $(4)$}\\
g_4(g_4e_3e_2e_1)&=& e_3e_2e_1-m\,g_4e_3e_2e_1+ml^{-1}e_4e_3e_2e_1\;\;\text{by $(8)$}\\
g_4(g_4g_3e_2e_1)&=&
g_3e_2e_1-m\,g_4g_3e_2e_1+ml^{-1}e_4g_3e_2e_1\;\;\text{by
$(8)$}\notag\\
&=&
g_3e_2e_1-m\,g_4g_3e_2e_1+ml^{-1}e_4e_3e_2e_1^{(23)^{'}}\\
g_4(g_4g_3g_2e_1)&=&
g_3g_2e_1-m\,g_4g_3g_2e_1+ml^{-1}e_4g_3g_2e_1\;\;\text{by
$(8)$}\notag\\
&=& g_3g_2e_1-m\,g_4g_3g_2e_1+ml^{-1}e_4e_3e_2e_1^{(132)^{''}}
\end{eqnarray}
\noindent In $(62)$, there is one undercrossing indicated with a
prime. In $(63)$, there are two undercrossings that are indicated
with a double prime. We summarize these actions in matrices: the
action of $g_1$ is represented by the matrix:

$$\begin{pmatrix}
l^{-1} & m & 0 & 0 & \boxed{m^{(12)}} & 0 & 0 & 0 & \boxed{m^{(123)}} & 0\\
0 & \!\!\!\!-m & 1 & 0 &  0 & 0 & 0 & 0 & 0 & 0\\
0 & 1 & 0 & 0 & 0 & 0 & 0 & 0 & 0 & 0\\
0 & 0 & 0 & \boxed{1^{(12)}} &  0 & 0 & 0 & 0 & 0 & 0\\
0 & 0 & 0 & 0 & -m & 1 & 0 & 0 & 0 & 0\\
0 & 0 & 0 & 0 &  1 & 0 & 0 & 0 & 0 & 0\\
0 & 0 & 0 & 0 &  0 & 0 & \boxed{1^{(12)}} & 0 & 0 & 0\\
0 & 0 & 0 & 0 &  0 & 0 & 0 & \boxed{1^{(12)}} & 0 & 0\\
0 & 0 & 0 & 0 &  0 & 0 & 0 & 0 & -m & 1\\
0 & 0 & 0 & 0 &  0 & 0 & 0 & 0 & 1 & 0\\
\end{pmatrix}$$

\noindent The matrix giving the action of $g_2$ is:

$$\begin{pmatrix}
0 & 0 & 1 & 0 & 0 & 0 & 0 & 0 & 0 & 0\\
0 & \;\,\, l^{-1} & ml^{-1} & m & 0 & 0 & 0 & \boxed{m^{(23)}} & 0 & 0\\
1 & 0  & -m &  0 & 0 & 0 & 0 & 0 & 0 & 0\\
0 & 0 &  0 & -m & 1 & 0 & 0 & 0 & 0 & 0\\
0 & 0 &  0 &  1 & 0 & 0 & 0 & 0 & 0 & 0\\
0 & 0 &  0 & 0 & 0 & \boxed{1^{(12)}} & 0 & 0 & 0 & 0\\
0 & 0 & 0 & 0 & 0 & 0 & \boxed{1^{(23)}} & 0 & 0 & 0\\
0 & 0 & 0 & 0 & 0 & 0 & 0 & -m & 1 & 0\\
0 & 0 & 0 & 0 & 0 & 0 & 0 & 1 & 0 & 0\\
0 & 0 & 0 & 0 & 0 & 0 & 0 & 0 & 0 & \boxed{1^{(12)}}
\end{pmatrix}$$

And the one giving the action of $g_3$ is:

$$\begin{pmatrix} \boxed{1^{(12)}} & 0 & 0 & 0 & 0 & 0 & 0 & 0 & 0 & 0 \\
 0 & 0 & 0 & 0 & 1 &  0 & 0 & 0 & 0 & 0 \\
 0 & 0 & 0 & 0 & 0 & 1 & 0 & 0 & 0 & 0 \\
 0 & 0 & 0 & l^{-1} & ml^{-1} & \boxed{(ml^{-1})^{(12)^{'}}} & m & 0 & 0 & 0 \\
 0 & 1 & 0 & 0 & \!\!\!\! -m & 0 & 0 & 0 & 0 & 0 \\
 0 & 0 & 1 & 0 & 0 & -m & 0 & 0 & 0 & 0 \\
 0 & 0 & 0 & 0 & 0 & 0 & -m & 1 & 0 & 0 \\
 0 & 0 & 0 & 0 & 0 & 0 & 1 & 0 & 0 & 0 \\
 0 & 0 & 0 & 0 & 0 & 0 & 0 & 0 & \boxed{1^{(23)}} & 0 \\
 0 & 0 & 0 & 0 & 0 & 0 & 0 & 0 & 0 & \boxed{1^{(23)}}
\end{pmatrix}$$

\noindent \\As already mentioned above, in those three matrices
respectively representing the actions of $g_1$, $g_2$ and $g_3$ on
the tangles spanning $W$, the upper left $6\times 6$ block is the
same as in the case $n=4$ where the $c$ (resp $c^{'}$) indicating
the over (resp under) crossing on the two first vertical braids has
been replaced with our new notations by the transpose $(12)$ (resp
$(12)^{'}$). Finally, let's give the matrix representing the action
of $g_4$:

$$\!\!\!\!\!\!\!
\begin{pmatrix}
\boxed{1^{(23)}} & 0 & 0 & 0 & 0 & 0 & 0 & 0 & 0 & 0 \\
 0 & \boxed{1^{(23)}} & 0 & 0 & 0 &  0 & 0 & 0 & 0 & 0 \\
 0 & 0 & \boxed{1^{(23)}}& 0 & 0 & 0 & 0 & 0 & 0 & 0 \\
 0 & 0 & 0 & 0 & 0 & 0 & 0 & 1 & 0 & 0 \\
 0 & 0 & 0 & 0 & 0 & 0 & 0 & 0 & 1 & 0 \\
 0 & 0 & 0 & 0 & 0 & 0 & 0 & 0 & 0 & 1 \\
 0 & 0 & 0 & 0 & 0 & 0 & \;\;\, l^{-1} & \;\;\, ml^{-1} & \boxed{(ml^{-1})^{(23)^{'}}} & \boxed{(ml^{-1})^{(132)^{''}}} \\
 0 & 0 & 0 & 1 & 0 & 0 & 0 & \!\!\!\!\! -m & 0 & 0 \\
 0 & 0 & 0 & 0 & 1 & 0 & 0 & 0 & -m & 0 \\
 0 & 0 & 0 & 0 & 0 & 1 & 0 & 0 & 0 & -m
\end{pmatrix}$$

\noindent In the case $n=5$, it appears to be harder to explicit an
invariant subspace of $W$ like we did previously. In the case $n=4$,
we got a representation by multiplying the coefficients in the boxes
by $r$ when the crossing was over (as in $(33)$) and by dividing
them by $r$ when the crossing was under (as in $(34)$). Here we try
the same method, the power of $r$ increasing with the number of
crossings: if there are two over crossings among the vertical
strands, we multiply the corresponding coefficient in the box by
$r^2$. If there are two under crossings among the vertical strands,
we will divide the corresponding coefficient in the box by $r^2$. In
other words, a transpose is replaced by a multiplication by $r$ and
a cycle of length three is replaced by a multiplication by $r^2$
when the vertical crossings are over. A transpose is replaced by a
division by $r$ and a cycle of length three is replaced by a
division by $r^2$ when all the vertical crossings are under. It
yields the matrices:

$$G_1(5)=\begin{pmatrix}
l^{-1} & m & 0 & 0 & mr & 0 & 0 & 0 & mr^2 & 0\\
0 & \!\!\!\!-m & 1 & 0 &  0 & 0 & 0 & 0 & 0 & 0\\
0 & 1 & 0 & 0 & 0 & 0 & 0 & 0 & 0 & 0\\
0 & 0 & 0 & r &  0 & 0 & 0 & 0 & 0 & 0\\
0 & 0 & 0 & 0 & -m & 1 & 0 & 0 & 0 & 0\\
0 & 0 & 0 & 0 &  1 & 0 & 0 & 0 & 0 & 0\\
0 & 0 & 0 & 0 &  0 & 0 & r & 0 & 0 & 0\\
0 & 0 & 0 & 0 &  0 & 0 & 0 & r & 0 & 0\\
0 & 0 & 0 & 0 &  0 & 0 & 0 & 0 & -m & 1\\
0 & 0 & 0 & 0 &  0 & 0 & 0 & 0 & 1 & 0\\
\end{pmatrix}$$

$$G_2(5)=\begin{pmatrix}
0 & 0 & 1 & 0 & 0 & 0 & 0 & 0 & 0 & 0\\
0 & \;\,\, l^{-1} & ml^{-1} & m & 0 & 0 & 0 & mr & 0 & 0\\
1 & 0  & -m &  0 & 0 & 0 & 0 & 0 & 0 & 0\\
0 & 0 &  0 & -m & 1 & 0 & 0 & 0 & 0 & 0\\
0 & 0 &  0 &  1 & 0 & 0 & 0 & 0 & 0 & 0\\
0 & 0 &  0 & 0 & 0 & r & 0 & 0 & 0 & 0\\
0 & 0 & 0 & 0 & 0 & 0 & r & 0 & 0 & 0\\
0 & 0 & 0 & 0 & 0 & 0 & 0 & -m & 1 & 0\\
0 & 0 & 0 & 0 & 0 & 0 & 0 & 1 & 0 & 0\\
0 & 0 & 0 & 0 & 0 & 0 & 0 & 0 & 0 & r
\end{pmatrix}$$

$$G_3(5)=\begin{pmatrix} r & 0 & 0 & 0 & 0 & 0 & 0 & 0 & 0 & 0 \\
 0 & 0 & 0 & 0 & 1 &  0 & 0 & 0 & 0 & 0 \\
 0 & 0 & 0 & 0 & 0 & 1 & 0 & 0 & 0 & 0 \\
 0 & 0 & 0 & l^{-1} & ml^{-1} & \frac{ml^{-1}}{r} & m & 0 & 0 & 0 \\
 0 & 1 & 0 & 0 & \!\!\!\! -m & 0 & 0 & 0 & 0 & 0 \\
 0 & 0 & 1 & 0 & 0 & -m & 0 & 0 & 0 & 0 \\
 0 & 0 & 0 & 0 & 0 & 0 & -m & 1 & 0 & 0 \\
 0 & 0 & 0 & 0 & 0 & 0 & 1 & 0 & 0 & 0 \\
 0 & 0 & 0 & 0 & 0 & 0 & 0 & 0 & r & 0 \\
 0 & 0 & 0 & 0 & 0 & 0 & 0 & 0 & 0 & r
\end{pmatrix}$$

$$G_4(5)=\begin{pmatrix}
 r & 0 & 0 & 0 & 0 & 0 & 0 & 0 & 0 & 0 \\
 0 & r & 0 & 0 & 0 &  0 & 0 & 0 & 0 & 0 \\
 0 & 0 & r & 0 & 0 & 0 & 0 & 0 & 0 & 0 \\
 0 & 0 & 0 & 0 & 0 & 0 & 0 & 1 & 0 & 0 \\
 0 & 0 & 0 & 0 & 0 & 0 & 0 & 0 & 1 & 0 \\
 0 & 0 & 0 & 0 & 0 & 0 & 0 & 0 & 0 & 1 \\
 0 & 0 & 0 & 0 & 0 & 0 & \;\;\, l^{-1} & \;\;\, ml^{-1} & \frac{ml^{-1}}{r} & \frac{ml^{-1}}{r^2} \\
 0 & 0 & 0 & 1 & 0 & 0 & 0 & \!\!\!\!\! -m & 0 & 0 \\
 0 & 0 & 0 & 0 & 1 & 0 & 0 & 0 & -m & 0 \\
 0 & 0 & 0 & 0 & 0 & 1 & 0 & 0 & 0 & -m
\end{pmatrix}$$

\noindent All these matrices are invertible and have the same
determinant $\frac{-r^3}{l}$. As in $(17)$ and $(35)$, we define for
all $i\in\{1,\,2,\,3,\,4\}$:
$$E_i(5)=\frac{l}{m}(G_i(5)^2+mG_i(5)-I_{10})$$ where $I_{10}$ is the identity matrix of size $10$.
We verify with Maple that the braid relations $(1)$ and $(2)$ are
satisfied on the $G_i$'s and that the relations $(4)$, $(5)$ and
$(6)$ combining the $e_i$'s and the $g_i$'s are also satisfied on
the matrices $E_i$'s and $G_i$'s. We define a morphism of groups:
\begin{eqnarray*} B(A_4)^{\times}& \lra & GL_{10}(F)\\
g_i & \longmapsto & G_i(5)\end{eqnarray*} \noindent where
$B(A_4)^{\times}$ is the group of units of $B(A_4)$ that we extend
to a morphism of algebras defined on the generators $g_i$'s and
$e_i$'s of $B(A_4)$:
\begin{eqnarray*}
B(A_4) & \lra & \mathcal{M}(10,F)\\
g_i & \longmapsto & G_i(5)\\
e_i & \longmapsto & E_i(5)
\end{eqnarray*}
\noindent Thus, we have a representation $\x$ of $B(A_4)$ in
$F^{10}$:
\begin{eqnarray*}
B(A_4)&\stackrel{\x}{\lra}& End_{F}(F^{10})
\end{eqnarray*}
\noindent defined by $\x(g_i)(X)=G_i(5)X$ for any $X\in F^{10}$ and
any $i\in\{1,\,\dots,\,4\}$. Let us denote by
$\mathcal{E}=(\e_1,\,\e_2\dots,\e_{10})$ the canonical basis of
$F^{10}$. These vectors have one on their $ith$ coordinate and zeros
elsewhere. $G_i(5)$ is the matrix of $\x(g_i)$ in the basis
$\mathcal{E}$. Suppose that $\x$ is not irreducible. Then there
exists a nonzero proper $F$-vector subspace $H$ of $F^{10}$ such
that $\x (w)(H)\subseteq H$ for any word $w$ in $B(A_4)$. In
particular, we must have $\x(X_{ij})(H)\subseteq H$ for all $1\leq
i<j\leq 5$. Computations in Maple show that the matrices
representing the $X_{ij}$'s each have exactly one nonzero row: the
one corresponding to the positive root which has a top horizontal
line joining the nodes $i$ and $j$. For instance, the tangle
associated to the positive root $\al_2+\al_3$ has its top horizontal
line joining the nodes $2$ and $4$ and is the fifth element of the
basis $\mathcal{B}$. And the fifth row is the only nonzero row in
the matrix of the endomorphism $\x (X_{24})$ in the basis
$\mathcal{E}$. Let's give an explicit expression for the matrix $S$
representing the sum of the $X_{ij}$'s. We indicated in front of
each row the corresponding $X_{ij}$:

$$\!\!\!\!\!\!\!\!\!\!\!\!\!\!\!\!\!\!\!\!\!\!\!\!\!\!\!\!\!\!\!\!\!\!\!\!\!\!
\begin{pmatrix}
\vspace{0.07in} \!\!\!\!\!\!\!\!\!\!\!\!\!\!\!\!\!\!\!\!\!
X_{12} \;\;\;\;\;\;\; x & 1 & l & 0 & r & lr & 0 & 0 & r^2 & lr^2\\
\vspace{0.07in} \!\!\!\!\!\!\!\!\!\!\!\!\!\!\!\!\!\!\!\!\!
X_{23} \;\;\;\;\;\;\; 1 & x & \frac{1}{l} & 1 & l & 0 & 0 & r & lr & 0\\
\vspace{0.07in} \!\!\!\!\!\!\!\!\!\!\!\!\!\!\!\!\!\!\!\!\! X_{13}
\;\;\;\;\;\;\; \frac{1}{l} & l & x & r & (r-\frac{1}{r})(l-r) & l &
0 & r^2 & (r^2-1)(l-r) & lr\\ \vspace{0.07in}
\!\!\!\!\!\!\!\!\!\!\!\!\!\!\!\!\!\!\!\!\!\! X_{34} \;\;\;\;\;\;\; 0
& 1 & \frac{1}{r} & x & \frac{1}{l} & \frac{1}{lr} & 1 &
l & 0 & 0\\
\vspace{0.07in} \!\!\!\!\!\!\!\!\!\!\!\!\!\!\!\!\!\!\!\!\!\! X_{24}
\;\;\;\;\;\;\; \frac{1}{r} & \frac{1}{l} & (\unsurr
-r)(\unsur{l}-\unsurr)
& l & x & \unsur{l} & r & (r-\unsurr)(l-r) & l & 0\\
\vspace{0.07in} \!\!\!\!\!\!\!\!\!\!\!\!\!\!\!\!\!\!\!\!\! X_{14}
\;\;\;\;\;\;\; \unsur{lr} & 0 & \unsur{l} & lr & l & x & r^2 &
(r^2-1)(l-r) & (r-\unsurr)(l-r) & l\\
\vspace{0.07in} \!\!\!\!\!\!\!\!\!\!\!\!\!\!\!\!\!\!\!\!\!\!\!
X_{45} \;\;\;\;\;\;\; 0 & 0 & 0 & 1 & \unsurr & \unsur{r^2} & x &
\unsur{l} &
\unsur{lr} & \unsur{lr^2}\\
\vspace{0.07in} \!\!\!\!\!\!\!\!\!\!\!\!\!\!\!\!\!\!\!\!\!\!\!
X_{35} \;\;\;\;\;\;\; 0 & \unsurr & \unsur{r^2} & \unsur{l} &
(\unsurr -r)(\unsur{l}-\unsurr) & (\unsur{r^2}-1)(\unsur{l}-\unsurr)
& l & x
& \unsur{l} & \unsur{lr}\\
\vspace{0.07in} \!\!\!\!\!\!\!\!\!\!\!\!\!\!\!\!\!\!\!\!\! X_{25}
\;\;\;\;\;\;\; \unsur{r^2} & \unsur{lr} &
(\unsur{r^2}-1)(\unsur{l}-\unsur{r}) & 0 & \unsur{l} & (\unsurr
-r)(\unsur{l}-\unsurr) & lr & l & x & \unsur{l}\\
\vspace{0.07in} \!\!\!\!\!\!\!\!\!\!\!\!\!\!\!\!\!\!\!\! X_{15}
\;\;\;\;\;\;\; \unsur{lr^2} & 0 & \unsur{lr} & 0 & 0 & \unsur{l} &
lr^2 & lr & l & x \end{pmatrix}=S$$

\noindent Suppose now that $\x (X_{24})(H)$ is nonzero. Let $h$ be a
vector in $H$ such that $\x (X_{24})(h)\neq 0$. $\x (X_{24})(h)$ is
a vector in $H$ whose only nonzero coordinate is the fifth one. It
follows that $\e_5$ belongs to $H$. Next, we read on the fifth
column of the matrix $S$:
$$\begin{array}{ccccc}
\vspace{0.05in}\x(X_{12})(\e_5)&=& r\,\e_1 & \Ra & \e_1\in H\\
\vspace{0.05in}\x(X_{34})(\e_5)&=& \frac{1}{l}\,\e_4 & \Ra & \e_4\in H\\
\vspace{0.05in}\x(X_{45})(\e_5)&=& \frac{1}{r}\,\e_7 & \Ra & \e_7\in H\\
\vspace{0.05in}\x(X_{25})(\e_5)&=& \unsur{l}\,\e_9 & \Ra & \e_9\in H\\
\end{array}$$
\noindent Then we read on the first column of the matrix $S$:
$$\begin{array}{ccccc}
\vspace{0.05in}\x(X_{13})(\e_1)&=& \unsur{l}\,\e_3 &\Ra& \e_3\in H\\
\vspace{0.05in}\x(X_{14})(\e_1)&=& \unsur{lr}\,\e_6 &\Ra& \e_6\in H\\
\end{array}$$
\noindent Finally, we read on the third column of $S$:
$$\begin{array}{ccccc}
\vspace{0.05in}\x(X_{35})(\e_3)&=& \unsur{r^2}\,\e_8 &\Ra& \e_8\in H\\
\vspace{0.05in}\x(X_{15})(\e_3)&=& \unsur{lr}\,\e_{10}&\Ra& \e_{10}\in H\\
\vspace{0.05in}\x(X_{23})(\e_3)&=& \unsur{l}\,\e_2 &\Ra& \e_2\in H\\
\end{array}$$
\noindent We conclude that the whole basis $\mathcal{E}$ is
contained in $H$. Thus, $H=F^{10}$, which contradicts $H$ proper. So
$\x(X_{24})(H)=0$. Similarly, the positive root having nodes $i$ and
$j$ joined on the top line in the tangles is
$\al_{j-1}+\dots+\al_i$. The associated tangle is the $\{[1+2+\dots
+(j-2)]+(j-i)\}$-th vector of the basis $\mathcal{B}$ of $W$. Hence
we have:
$$\big( \x(X_{ij})(h)\neq 0,\,\text{some}\; h\in
H\big)\Longrightarrow\Big(\e_{\frac{(j-1)(j-2)}{2}+j-i}\in H\Big)$$
Next, we observe that in each column of $S$, there are at least six
nonzero and non diagonal coefficients. In particular, there are at
least six nonzero and non diagonal coefficients in the column $k$ of
$S$, where $$k:=\frac{(j-1)(j-2)}{2}+j-i$$ From there, reasoning
like above, we deduce that there are six other elements of the basis
$\mathcal{B}$, say $\e_{i_2},\dots\e_{i_7}$, that are in $H$. We
need to verify that the three remaining ones, say $\e_{i_8}$,
$\e_{i_9}$ and $\e_{i_{10}}$, also belong to $H$. We notice again
that in each row of $S$, there are at least six nonzero and non
diagonal coefficients. Let $s_{i_8j_1},\dots,s_{i_8j_6}$ be six non
zero and non diagonal coefficients of the $i_8$-th row of $S$.
Necessarily, there exists $s\in\{2,\dots,7\}$ and
$t\in\{1,\dots,6\}$ such that $i_s=j_t$. It comes:
$$\x(X_{kl})(\e_{i_s})= s_{i_8j_t}\e_{i_8}$$
where $X_{kl}$ is in front of the $i_8$-th row. It follows that
$\e_{i_8}$ belongs to $H$ and the same method applies to show that
$\e_{i_9}$ and $\e_{i_{10}}$ are in $H$. We conclude that
$\x(X_{ij})(H)=0$, since otherwise $H$ wouldn't be proper. This
equality holds for all the $X_{ij}$'s, hence considering their sum
$\mathcal{S}$, we get $\x(\mathcal{S})(H)=0$. In other words, we
have $\forall z\in H,\,Sz=0$. We conclude that if $\x$ is not
irreducible, we must have $det(S)=0$. Using Maple yields the
equivalence $det(S)=0\Longleftrightarrow
l\in\{r,\,-r^3,\,-\unsur{r^2},\,\unsur{r^2},\,\unsur{r^7}\}$. Let's
show conversely that for any of these values of $l$, the
representation $\x$ is not irreducible. First, we show that
$$\x(g_k)\big(\underset{1\leq i<j\leq 5}{\bigcap}
Ker\,\x(X_{ij})\big)\subset\underset{1\leq i<j\leq 5}{\bigcap}
Ker\,\x(X_{ij})$$ for all $k$ in $\{1,\,\dots,\,4\}$. Let $h$ be a
vector of $F^{10}$ such that $\x(X_{ij})(h)=0$ for all $1\leq
i<j\leq 5$. Given $i$ and $j$ such that $1\leq i<j\leq 5$, using the
conjugation formulas for $X_{ij}$ and the algebra relations like in
the proof of lemma $2$, but this time on the matrices, we get
$Mat_{\mathcal{E}}(\x(X_{ij}))G_{k}(5)h=0$ for all
$k\in\{1,\,\dots,\,4\}$. So we have $(\x(X_{ij})\circ \x(g_k))(h)=0$
\emph{i.e.\/} $\x(g_k)(h)\in Ker\,\x(X_{ij})$ for all
$k\in\{1,\,\dots,\,4\}$. Hence the inclusion above is satisfied and
the $F$-vector subspace $\cap_{1\leq i<j\leq 5} Ker\,\x(X_{ij})$ of
$F^{10}$ is a $B(A_4)$-module. The next step is to show that for
each of the values of $l$ above, this space is nonzero. A vector $v$
of $F^{10}$ is in $\cap_{1\leq i<j\leq 5} Ker\,\x(X_{ij})$ if and
only if it is in the kernel of the matrix $S$. We check that:
\begin{align}
\text{If}&\; l = r, & \!\!\!\!\!\!\!\!\!\text{then} &\;\;
(r^2,\,0,\,-r,\,1,\,-r,\,0,\,0,\,0,\,0,\,0)\in Ker\,S\\
\text{If}&\; l = -r^3, & \!\!\!\!\!\!\!\!\!\!\text{then} &\;\;
(0,\,-r,\,0,\,-\unsurr,\,1,\,0,\,0,\,0,\,0,\,0)\in Ker\,S\\
\text{If}& \;l  =  -\unsur{r^2}, & \!\!\!\!\!\!\text{then} &\;\;
(-r^2,\,r^2+\unsurr,\,r,\,-\unsurr,\,1,\,0,\,0,\,-1,\,r,\,0)\in Ker\,S\\
\text{If}& \;l = \unsur{r^2}, & \!\!\!\!\!\!\!\!\!\text{then} &\;\;
(r^2,\,-r^2+\unsurr,\,-r,\,-\unsurr,\,1,\,0,\,0,\,-1,\,r,\,0)\in Ker\,S\\
\text{If}& \;l = \unsur{r^7}, &\!\!\!\!\!\!\!\!\! \text{then} &\;\;
(\unsur{r^3},\,\unsurr,\,\unsur{r^2},\,r,\,1,\,\unsurr,\,r^3,\,r^2,\,r,\,1)\in
Ker\,S
\end{align}
We conclude that if $l$ is one of these values, then $\x$ is not
irreducible. Hence we have the equivalence: $\x\; \text{is
irreducible}\;
\emph{iff\/}\;l\not\in\{r,\,-r^3,\,-\unsur{r^2},\,\unsur{r^2},\,\unsur{r^7}\}$

%\begin{eqnarray*}
%g_1(g_4e_3e_2e_1) &=& g_4e_3g_1e_2e_1\;\;\;\;\;\;\;\;\;\;\;\;\;\;\;\;\;\;\;\text{by $(1)$}\\
%&=& g_4e_3g_2e_1-m\,g_4e_3e_2e_1\,\,\,\text{by $(12)$ and the fact
%that $e_3e_1=0$}\\ &=&
%g_4e_3e_2e_1^{(12)}\;\;\;\;\;\;\;\;\;\;\;\;\;\;\;\;\;\;\;\!\,\text{by
%$(31)$}
%\end{eqnarray*}

\section{Generalization} In this section, we extend the previous
constructions to $B(A_{n-1})$, for any integer $n$. In $A_{n-1}$,
there are $\binom{n}{2}$ positive roots. As already described in
section $4$, we order the roots in the following way: incrementing
$i$, we start with $\al_i$ and list all the positive roots that have
nodes less or equal to $i$ in their support (there are $i$ of them),
in an increasing height order. This ordering allows us to "build
$A_1$ up to $A_{n-1}$". Geometrically with the tangles, we join a
node to its left neighbors starting with its adjacent neighbor and
moving to the left; we shift the right extremity of the top
horizontal line and start again. After step $i$, all the possible
pairs of nodes for the top horizontal line appear exactly once on
the first $i+1$ nodes. Hence we have listed $\binom{i+1}{2}$
tangles. The first $\binom{n-1}{2}$ positive roots correspond to the
tangles spanning $W$ in $B(A_{n-2})$ where a vertical string has
been added to the right side. Hence, proceeding inductively, we
already know the action of $g_1,\dots,\,g_{n-2}$ on these tangles,
the last vertical string on the right being left invariant with
these actions. To obtain the matrices $G_1,\,G_2,\dots,\,G_{n-2}$
inductively, it remains to compute the actions of these $g_{i}$'s on
the $n-1$ tangles
\begin{multline*} e_{n-1}e_{n-2}\dots e_1,\;\;g_{n-1}e_{n-2}\dots
e_1,\;\;g_{n-1}g_{n-2}e_{n-3}\dots
e_1\\,\dots,\;\;g_{n-1}g_{n-2}\dots g_2e_1\end{multline*} \noindent
The action by $g_1$ on the first $n-3$ tangles (in the roots
$\al_{n-1},\,\al_{n-1}+\al_{n-2},\dots,$\\$\al_{n-1}+\al_{n-2}+\dots
+\al_3$) is simply a crossing, hence a multiplication by $r$. Recall
that the top horizontal line of the tangle associated with the
positive root $\al_{n-1}+\dots +\al_{i}$ joins the nodes $i$ and
$n$. As long as the top horizontal line does not begin with either
node $i$ or node $i+1$, as in $\al_{n-1},\;\dots,\;\al_{n-1}+\dots
+\al_{i+2}$ (first $n-i-2$ tangles) or as in $\al_{n-1}+\dots
+\al_{i-1},\;\dots,\;\al_{n-1}+\dots +\al_1$ (last $i-1$ tangles),
the action by $g_i$ is a crossing and the coefficient in the matrix
is an $r$ on the diagonal. Suppose now that we are looking for the
action of a $g_i$ on the tangle whose top horizontal line starts at
node $i$ (and ends at node $n$). The left extremity of the top line
of the tangle resulting from this action is now shifted to the right
and starts at node $i+1$. In other words, the root $\al_{n-1}+\dots
+\al_i$ is sent to the root $\al_{n+1}+\dots +\al_{i+1}$, which
corresponds to a $1$ above the diagonal in the matrix of $G_i$. It
remains to look at the action of $g_i$ on the tangle $(i+1,n)$. The
node $i+1$ is sent onto the node $i$ with an undercrossing. Using
the tangle formula $(31)$, we get:
\begin{eqnarray}g_i.(i+1,n)=(i,n) + m\,e_i.(i+1,n) -m\,(i+1,n)\end{eqnarray}
where the permutation $(k,l)$ denotes the tangle in $W$ whose top
line starts with node $k$ and ends with node $l$. Acting by $e_i$
transforms the top horizontal line of the tangle into a vertical
line that still over crosses the same number of vertical strings,
that is a total of $n-i-2$ vertical strings. Its top horizontal line
now joins the nodes $i$ and $i+1$ as in $e_ie_{i-1}\dots e_1$. We
read on the equation $(69)$ that there are three coefficients in the
column of the matrix. The first one is a $1$ just below the
diagonal, the second one is $mr^{n-i-2}$ on the
$\{1+2+\dots+(i-1)+1\}$th row of the matrix, where the factor on the
right corresponds to the $n-i-2$ over crossings. The third
coefficient is a $-m$ on the diagonal. We get the following matrices
for $G_1$, $G_2,\,\dots,\,G_{n-2}$ (that we will call "matrices of
the first kind"):
\begin{center}
\begin{tabular}{|p{1in}|p{0.6in}|l|l|l}
\multicolumn{5}{l}{$\!\!\!\!\underleftrightarrow{\text{\hspace{0.42in}}\binom{n-1}{2}\text{\hspace{0.42in}}}$}\\
\cline{1-4} \bigskip \hspace{0.1in} $G_1(n-1)$ \bigskip\bigskip & \multicolumn{3}{l|}{\hspace{0.7in}$m\,r^{n-3}\;\;\;$} & $\leftarrow e_1$ \\
\hfill & \multicolumn{3}{l|}{$\!\!\!\!\,\ulra{\text{\hspace{0.21in}}n-3\text{\hspace{0.21in}}}$} & \\ \cline{1-2} \hfill & $\begin{matrix}r & & \\
& \ddots &
\\ & & r
\end{matrix}$ & \multicolumn{2}{l|}{} & \\ \cline{2-2}
\multicolumn{4}{|l|}{\hspace{1.8in}
$\begin{pmatrix} -m & 1 \\ 1 & 0 \end{pmatrix}$} & \\
\cline{1-4}
\end{tabular}
\end{center}

\begin{center}
\hspace{0.14in}\begin{tabular}{|p{1in}|p{0.4in}|l|l|l|l}
\cline{1-5} \hspace{1in} & \multicolumn{4}{l|}{}\\
\bigskip \hspace{0.1in} $G_2(n-1)$ \bigskip &
\multicolumn{4}{c|}{\hspace{0.16in}$m\,r^{n-4}$} & $\leftarrow e_2e_1$\\
\hfill & \multicolumn{4}{l|}{$\!\!\!\!\,\ulra{\text{\hspace{0.11in}}n-4\text{\hspace{0.11in}}}$} & \\
\cline{1-2} \hfill & \hspace{-0.02in} $\begin{smallmatrix} r & & \\ & \ddots & \\
& & r\end{smallmatrix}\;\;\;$ \vspace{0.005in} & \multicolumn{3}{l|}{} & \\
\cline{2-2} \multicolumn{5}{|l|}{\hspace{1.65in}$\begin{pmatrix} -m & 1 \\ 1 & 0 \end{pmatrix}$\hspace{0.7in}} & \\
\multicolumn{5}{|l|}{\hspace{2.45in}$r$} & \\ \cline{1-5}
\end{tabular}\\
\smallskip $\vdots$\\\vspace{-0.06in}$\vdots$\\\vspace{-0.06in}$\vdots$
\end{center}

\begin{center}
\hspace{0.675in}\begin{tabular}{|p{1in}|l|l|l|l} \cline{1-4}
\hspace{1in} & \multicolumn{3}{l|}{}\\
\hspace{1in} & \multicolumn{3}{l|}{}\\
%\hspace{1in} & \multicolumn{3}{l|}{}\\
\hspace{0.1in} $G_{n-2}(n-1)$ & \multicolumn{3}{l|}{}\\
\hspace{1in} & \multicolumn{3}{l|}{}\\
\hspace{1in} & \multicolumn{3}{l|}{$m$} & $\leftarrow e_{n-2}\dots e_2e_1$\\
\hspace{1in} & \multicolumn{3}{l|}{}\\
\cline{1-1} \multicolumn{4}{|l|}{\hspace{1.11in}$\!\!\begin{pmatrix} \!-m & 1\\
1 & 0
\end{pmatrix}\hspace{0.07in}\ulra{\text{\hspace{0.21in}}n-3\text{\hspace{0.14in}}}$}\\
\cline{4-4} \multicolumn{3}{|l}{\hfill}\hspace{1.66in} \vline & $\begin{matrix} r & & \\
& \ddots
& \\ & & r \end{matrix}$\\
\cline{1-4}
\end{tabular}
\end{center}

\newcounter{RM}
\setcounter{RM}{2}

\newpage\noindent We now compute the action of $g_{n-1}$ on the spanning tangles
of $W$. When $n-2$ and $n-1$ are not in the support of the
corresponding positive roots, the tangle's top line does not end
with node $n-1$ or node $n$. Thus the result of the action by
$g_{n-1}$ is a crossing of the last two vertical strings.
Consequently, the upper left square block of size $\binom{n-2}{2}$
of the matrix $G_{n-1}$ is a scalar matrix with $r$'s on the
diagonal. Now if $n-2$ is in the support of the positive root, the
top line of the associated tangle ends with node $n-1$ and a left
action by $g_{n-1}$ shifts this end to node $n$. In other words, the
positive root $\al_{n-2}+\dots +\al_i$ is sent to the positive root
$\al_{n-1}+\al_{n-2}+\dots +\al_i$. Next, if $n-1$ is in the support
of the positive root, the top line of the associated tangle ends
with node $n$. We start with the easy case of the action of
$g_{n-1}$ on the positive root $\al_{n-1}$ associated with the
tangle $e_{n-1}\dots e_1$. The result of this action is simply a
multiplication by $1/l$ to remove the loop, hence a $1/l$ on the
diagonal. If the top horizontal line crosses now $1$ to $n-2$
vertical strings, as in $\al_{n-1}+\al_{n-2}$, \dots,
$\al_{n-1}+\dots+\al_{1}$, the action by $g_{n-1}$ shifts the right
extremity of the top horizontal line to the left. The top line ends
now in node $n-1$ instead of node $n$. But to be able to use the
regular isotopy and more specifically Reidemeister's move \Roman{RM}
(see for instance \cite{MW}) that allows to separate two strings
that intersect in two overcrossings or two undercrossings, we first
need to transform the overcrossing of $g_{n-1}$ into an
undercrossing, using the tangle formula $(31)$. It yields:
$$g_{n-1}.(i,n)=(i,n-1)+\,m\,e_{n-1}.(i,n)-m(i,n)$$ The action by
$e_{n-1}$ transforms the top horizontal line into a vertical string
that overcrosses one less vertical lines as the last vertical line
is now part of a loop. In the general picture, the $n-i-2$ crossings
are under and to build a representation as before, we divide the
corresponding coefficient by the power $r^{n-i-2}$. The top line now
joins the two last nodes as in $e_{n-1}\dots e_1$. Thus, for $2\leq
j\leq n-1$, the $(1+2+\dots +(n-2)+1, 1+2+\dots +(n-2)+j)$
coefficient of the matrix $G_{n-1}$ is $\frac{m}{lr^{j-2}}$. The
matrix $G_{n-1}$ has the following shape, where the blanks must be
filled with zeros. We will call it a "matrix of the second kind":
\begin{center}
\begin{tabular}{|l|l|l|l|l}
\multicolumn{5}{l}{$\!\!\!\!\ulra{\text{\hspace{0.31in}}\binom{n-2}{2}\text{\hspace{0.31in}}}
\ulra{\text{\hspace{0.57in}}n-2\text{\hspace{0.57in}}}
\ulra{\text{\hspace{0.77in}}n-1\text{\hspace{0.77in}}}$} \\
\cline{1-3}
$\begin{matrix} r & & & \\
& r & & \\
& & \ddots & \\
& & & r
\end{matrix}$ & \hspace{1.15in} & \hspace{1in} \\
\cline{1-1}\cline{3-3} \multicolumn{2}{|l|}{} & $\begin{matrix}
\resizebox{0.1in}{1in}{o}\end{matrix}$ \hspace{0.07in} \vline \; $\begin{matrix}  1 & & & & \\
& 1 & & & \\ & & 1 & & \\
& & & \ddots & & \\ & & & & \ddots & \\ & & & & & 1
\end{matrix}$ \\
\cline{1-3} \hspace{0.55in} & \resizebox{1.2in}{0.1in}{o} & $l^{-1}$
\vline $\;\;m\,l^{-1}\; \frac{m\,l^{-1}}{r}\;\dots\;
\frac{m\,l^{-1}}{r^{n-3}}$ \\
\cline{2-3} \hspace{0.55in} & $\begin{matrix}  1 & & & & \\
& 1 & & & \\ & & 1 & & \\
& & & \ddots & & \\ & & & & \ddots & \\ & & & & & 1
\end{matrix}$ & $\begin{matrix}\resizebox{0.1in}{1in}{o}\end{matrix}$ \hspace{0.07in} \vline $\begin{matrix}  -m & & & & \\
& -m & & & \\ & & \!\!\!\!\ddots & & \\ & & & \ddots & \\ & & & &
\!\!\!\!-m
 \end{matrix}$ \\
\cline{1-3}
\end{tabular}
\end{center}

\noindent \\ We defined matrices $G_1,\,\dots,\,G_{n-1}$ that
correspond to the actions of $g_1,\,\dots,\,g_{n-1}$ on the spanning
elements of $W$. We now define matrices $E_1,\,\dots,\,E_{n-1}$ by\\
setting
$$\forall 1\leq i\leq n-1,\, E_i:=\frac{l}{m}\Big(G_i^2+m\,G_i-I_{\binom{n}{2}}\Big)$$

\noindent Our conjecture is that:
\newtheorem{Conj}{Conjecture}
\begin{Conj}\hfill
\begin{description}
\item{\underline{Case n=3}}
%\hspace{0.15in} The map on the generators of $B(A_2)$}
\begin{center}$\begin{array}{ccc}
B(A_2)& \lra & \mathcal{M}(3,F)\\
g_1,\,g_2 & \longmapsto & G_1,\,G_2\\
e_1,\,e_2 & \longmapsto & E_1,\,E_2
\end{array}\;\;\;$ is a representation. \end{center}
It is irreducible iff $l\not\in\lbrace
 -r^3,-1,1,\frac{1}{r^3}\}$.
\item{\underline{General case}}
%\hspace{0.15in} The map on the generators of $B(A_{n-1})$}
\begin{center} $\begin{array}{ccc}
B(A_{n-1})& \lra & \mathcal{M}(\binom{n}{2},F)\\
g_i & \longmapsto & G_i\\
e_i & \longmapsto & E_i
\end{array}\;\;\;$ is a representation. \end{center}
It is irreducible iff $l\not\in\lbrace
 r,-r^3,-\frac{1}{r^{n-3}},\frac{1}{r^{n-3}},\frac{1}{r^{2n-3}}\}$.
\end{description}\end{Conj}\noindent The conjecture is verified in the cases
$n=3,4,5$ (cf respective sections $3$, $4$, $5$). Computations with
Maple were done in the case $n=6$. %joined to the existing
%results for the smaller values $n=3,4,5$ led us to the conjecture
%$1$ for a general $n$.
We wrote a program in Maple that defines the matrices
$G_1,\,G_2,\,\dots,\,G_{n-1}$ for a given $n$,
%computes the corresponding $E_i$'s and the inverses of the $G_i$'s, then computes
computes the $X_{ij}$'s,
%by using the conjugation formulas,
forms their sum matrix $S$, then solves in the variable $l$ the
equation $det(S)=0$. We ran the program for $n=6$, obtaining a
similar results (but without any proof) as in the cases $n=4$ and
$n=5$ and leading us to the conjecture $1$ for a general $n$. We
will prove the conjecture in the general case, but first we need to
give a more visible expression for the announced representation.
This is the object of the next section.

\section{The Representation itself}
\newcounter{m}
\setcounter{m}{1} Let $\mathcal{V}$ denote the vector space over $F$
with spanning vectors $x_{\be}$ indexed by $\be\in\phi^{+} $, where
$\phi^{+}$ is the set of positive roots. In this section we give a
formal definition of a representation

$$\begin{array}{ccc}
B(A_{n-1}) & \lra & End_{F}(\mathcal{V})\\
g_i & \longmapsto & \n_i
\end{array}$$
%An element $\xb$ will be viewed as an appropriate linear combination
%over $F$ of the tangles associated with $\be$, as described in
%previous section.
In the case $n=4$, an element $\xb$ can be viewed as an appropriate
linear combination over $F$ of the tangles associated with $\be$ (cf
section $4$). In what follows, $Supp(\be)$ denotes the support of
the positive root $\be$ \emph{i.e.\/} the set of $k\in\{1\dots n\}$
such that the coefficient of $\al_k$ in $\be$ is nonzero ; $ht(\be)$
denotes the height of the positive root $\be$: if $\be=\sum_i
n_i\al_i$, then $ht(\be):=\sum_i n_i$. We read on the matrices of
the previous section an expression for the representation in two
special cases: $i=n-1$ on one hand (point 1. below) and $i\in\lbrace
1,...,n-2\}$ and $n-1\in Supp(\be)$ on the other hand (point 2.
below). We then deduce inductively a formula for the representation
in the case $i\in\{1,\,\dots,\,n-2\}$ (point 3. below). We have:

\begin{enumerate}
\item
If $i=n-1$, then we read on the matrix $G_{n-1}$ of previous
section:
$$\n_i(\xb)= \begin{array}{cc} \begin{matrix} (\text{\Roman{m}})\\\stepcounter{m}(\text{\Roman{m}})\\\stepcounter{m}(\text{\Roman{m}})\\\stepcounter{m}(\text{\Roman{m}})\end{matrix} & \!\!\!\!\left| \begin{array}{ccc} r\;\;\;\;\!\xb \hfill\hfill\text{if} \;(\be|\al_i)&=&0\\
l^{-1}\xb \hfill\hfill\hfill\text{if}\; (\be|\al_i)&=&1\\
x_{\be-\al_i}+\frac{m}{l\,r^{ht(\be)-2}}\;x_{\al_i}-m\xb\;\;
\hfill\text{if}\;
(\be|\al_i)&=&\frac{1}{2}\\
x_{\be+\al_i}\;\;\;\;\;\;\;\;\;\;\;\;\;\;\;\;\;\;\;\;\;\;\;\;\;\hfill\text{if}\;
(\be|\al_i)&=&\!\!\!-\frac{1}{2}
\end{array}
\right.\end{array}$$

\setcounter{m}{1}
\item If $i\in\{1,\,\dots,\,n-2\}$ and if $n-1\in Supp(\be)$, then
we read on the last $n-1$ columns of the matrices
$G_1,\,\dots,\,G_{n-2}$ of previous section:

$$\n_i(\xb)=\begin{array}{cc}\begin{matrix}
(\text{\Roman{m}})\\\addtocounter{m}{2}(\text{\Roman{m}}^{'})\\\stepcounter{m}(\text{\Roman{m}}^{'})\end{matrix}
& \!\!\!\!\left|
\begin{array}{ccc}
r\,\xb \hfill\hfill\text{if} \;(\be|\al_i)&=& 0\\
x_{\be-\al_i}\hfill\hfill\;\;\text{if}\;(\be|\al_i)&=&\frac{1}{2}\\
x_{\be+\al_i}+\,m\,r^{n-i-2}x_{\al_i}-m\,\xb\;\;
\hfill\text{if}\;(\be|\al_i)&=& \!\!\!-\frac{1}{2}
\end{array}
\right.\end{array}$$

\item $i\in\{1,\,2,\,\dots,\,n-2\}$ and no restriction
on $\be$: we use induction with the expressions for the
representation in points 1. and 2. above. More specifically, we will
always use the induction in the following way: suppose that we want
to evaluate $\n_i(n)(x_{\al_l+\dots+\al_k})$ where $l\geq k$ and
$l\leq n-2$. By induction this value has already been computed in
lower dimension (in $B(A_{n-2})$) and proceeding successively by
induction, we may in fact decrement the integer $n$ till we reach
the integer $M:=Max(i+1,l+1)$. On the matrices the picture is to
keep moving backward in the upper left corner till you can either
use the last columns of a matrix of the first kind (in the case
$M=l+1$) or the matrix of the second kind (in the case $M=i+1$) in
the suitable $B(A_{M-1})$. We note that if $i+1<l+1$ then by
definition $M=l+1$ and we have $i\leq l-1=(l+1)-2$, so that the
inductive steps are all justified and the expression for
$\n_i(l+1)(x_{\al_l+\dots+\al_k})$ may indeed be read on the last
$l$ columns of the matrix of the first kind $G_i(l+1)$. If on the
contrary $i+1\geq l+1$, then by definition $M=i+1$. The inductive
steps are again justified by the fact that $i\leq (i+2)-2$ and
$l\leq i$. The expression for $\n_i(i+1)(x_{\al_l+\dots+\al_k})$ is
now obtained by using the matrix of the second kind $G_i(i+1)$ in
$B(A_i)$. We will now apply those preliminary remarks, while dealing
with the different values for the inner product $(\be|\al_i)$.
\begin{itemize}\item If $\ps=-\frac{1}{2}$, there are two cases:
\begin{list}{\texttt{*}}{}
\item $\be=\al_{i-1}+\dots+\al_k$ with $k\leq i-1$. We have:
\begin{eqnarray*}\n_i(n)(x_{\al_{i-1}+\dots+\al_k})=\dots\!\!\!\!&=& \n_i(i+1)(x_{\al_{i-1}+\dots+\al_k})\\
&=& x_{\be+\al_i}\;\;\;\text{by using}\;\;(\text{\Roman{m}})
\end{eqnarray*}
\item $\be=\al_l+\dots+\al_{i+1}$ with $l\geq i+1$. Then we have:
\begin{eqnarray*}
\n_i(n)(x_{\al_l+\dots+\al_{i+1}})=\dots \!\!\!\!&=&\n_i(l+1)(x_{\al_l+\dots+\al_{i+1}})\\
&=& x_{\be+\al_i}+m\,r^{l-i-1}x_{\al_i}-m\,\xb
\end{eqnarray*}
The last equality is obtained by using $(\text{\Roman{m}}^{'})$ of
point 2.
\\For the last equality to hold, we must have $i\leq (l+1)-2$,
\emph{i.e.\/} $i\leq l-1$ which is true by our assumption in
$\star$.
\end{list}
\item If $\ps=\frac{1}{2}$, there are again two cases:
\begin{list}{\texttt{*}}{}
\item $\be=\al_i+\dots+\al_k$ with $k\leq i-1$. Then we have:
\begin{eqnarray*} \n_i(n)(x_{\al_i+\dots+\al_k})&=&
\n_i(i+1)(x_{\al_i+\dots+\al_k})\\
&=& x_{\be-\al_i}+\frac{m}{l\,r^{ht(\be)-2}}\;x_{\al_i}-m\,\xb
\end{eqnarray*}
\addtocounter{m}{-1}The last equality is obtained by using the
expression (\Roman{m}) of point 1.
\item $\be=\al_l+\dots+\al_i$ with $l\geq i+1$. Then we have:
\begin{eqnarray*}
\n_i(n)(x_{\al_l+\dots+\al_i})&=& \n_i(l+1)(x_{\al_l+\dots+\al_i})\\
&=& x_{\be-\al_i}\;\;\;\text{by}\;\;(\text{\Roman{m}}^{'})
\end{eqnarray*}
\end{list}
\item If $\ps=0$, then there are again two possibilities:
\begin{list}{\texttt{*}}{}
\item $i-1,\,i,\,i+1\not\in Supp(\be)$\\
Then $\n_i(\xb)=r\,\xb$
\item $\{i-1,\,i,\,i+1\}\subseteq Supp(\be)$\\
Then we can write
$\be=\al_l+\dots+\al_{i+1}+\al_i+\al_{i-1}+\dots+\al_k$, some $l\geq
i+1$ and $k\leq i-1$. It comes:
\begin{equation*}
\begin{split}
\n_i(n)(& x_{\al_l+\dots+\al_{i+1}+\al_i+\al_{i-1}+\dots+\al_k})\\
&=\dots=\n_i(l+1)(x_{\al_l+\dots+\al_{i-1}+\al_i+\al_{i+1}+\dots+\al_k})
\end{split}
\end{equation*}
We have by assumption $i\leq l-1=(l+1)-2$, so that we may apply
point 2. Further, $\be$ has height at least $l-i+2$ and there are
$(l+1)-i-2=l-i-1$ $r$'s on the diagonal of the first $l-i-1$ columns
of the last $l$ columns of the matrix $G_i(l+1)$. By previous
section, skip two more columns and get again a pattern of $r$'s on
the diagonal, which now yields:
$$\n_i(n)(x_{\al_l+\dots+\al_{i+1}+\al_i+\al_{i-1}+\dots+\al_k})=r\,\xb$$
\end{list}
Thus, we see that in both cases we have $\n_i(\xb)=r\,\xb$.
\item If $(\be|\al_i)=1$, the only possibility is $\be=\al_i$ and we
have: \addtocounter{m}{-1}
\begin{eqnarray*}
\n_i(n)(x_{\al_i})&=&\n_i(i+1)(x_{\al_i})\\
&=& \frac{1}{l}\,x_{\al_i}\;\;\;\text{as in}\; (\text{\Roman{m}})
\end{eqnarray*}
We deduce from that last point and from the equation (\Roman{m}) of
point 1. that if $(\be|\al_i)=1$, then
$\;\n_i(\xb)=\frac{1}{l}\,\xb\;$ for all $i$.
\end{itemize}
\end{enumerate}

\noindent We now gather all these results to obtain the final
definition of the endomorphisms $\n_i$'s. In what follows, let $\p$
be the total order on the positive roots, as used and described
several times before:
\begin{multline*}\al_1\p\al_2\p\al_2+\al_1\p\al_3\p\al_3+\al_2\p\al_3+\al_2+\al_1\\
\p\dots\p\al_{n-1}\p\al_{n-1}+\al_{n-2}\p\dots\p\al_{n-1}+\al_{n-2}+\dots+\al_1
\end{multline*}
It will also be convenient to read $\be\q\gamma$ for some positive
roots $\be$ and $\gamma$. By $\be\q\gamma$ we will understand
$\gamma\p\be$. We see with points 1., 2. and 3. that we always have
$\n_i(\xb)=r\,\xb$ when $\ps=0$. Similarly, by points 1. and 3. we
always have $\n_i(\xb)=l^{-1}\,\xb$ when $\ps=1$. Suppose now that
$\ps=\frac{1}{2}$. Then $\be-\al_i$ is a positive root and there are
two cases: either $\be-\al_i\p\al_i$ or $\be-\al_i\q\al_i$. The
first situation occurs in point 1. since $\be-\al_{n-1}$ starts with
$\al_{n-2}$ and is thus ranked before $\al_{n-1}$ and in point 3.
when $\be=\al_i+\dots+\al_k$ with $k\leq i-1$. In the latter case,
$\be-\al_i=\al_{i-1}+\dots+\al_k\p\al_i$. In both cases we have
$\n_i(\xb)=x_{\be-\al_i}+\frac{m}{l\,r^{ht(\be)-2}}\,x_{\al_i}-m\,\xb$.
The second situation occurs in point 2. since then
$\be=\al_{n-1}+\dots+\al_{i+1}+\al_i$ and so
$\be-\al_i=\al_{n-1}+\dots+\al_{i+1}\q\al_i$. It also occurs in
point 3. when $\be=\al_l+\dots+\al_i$ with $l\geq i+1$, since then
$\be-\al_i=\al_l+\dots+\al_{i+1}\q\al_i$. Again, in both cases the
result is the same and we have $\n_i(\xb)=x_{\be-\al_i}$. Finally
the last case $\ps=-\frac{1}{2}$ splits in turn into two cases:
$\be\q\al_i$ on one hand and $\be\p\al_i$ on the other hand. Indeed,
in point 1. and in point 3. when $\be=\al_{i-1}+\dots+\al_k$ with
$k\leq i-1$, we have $\be\p\al_i$ and $\n_i(\xb)=x_{\be+\al_i}$.
Whereas in point 2. and in point 3. when $\be=\al_l+\dots+\al_{i+1}$
with $l\geq i+1$, we have $\be\q\al_i$ and
$\n_i(\xb)=x_{\be+\al_i}+m\,r^{ht(\be)-1}\,x_{\al_i}-m\,\xb$. \\Note
that when $\ps=\frac{1}{2}$, a direct comparison between $\be$ and
$\al_i$ is not the criterion as in that case $i$ lies in the support
of $\be$, which always yields $\be$ greater than $\al_i$. In fact,
if $\ps=\unsur{2}$, $i$ is contained in the support of $\be$ and
moreover $\be$ ends with $\al_i$ or begins with $\al_i$. Comparing
$\be-\al_i$ with $\al_i$ tells the relative position of $\al_i$ in
the sum. If $\ps=-\unsur{2}$, then $\be$ "stops just before $\al_i$"
(in which case $\be\q\al_i$
) or "starts right after $\al_i$" (in which case $\be\p\al_i$).\\
We are now able to give a complete expression for the endomorphism
$\n_i$:
\begin{equation*}\n_i(\xb)=\begin{cases} r\;\;\;\;\xb & \text{if $\ps=\;\;\;0\qquad\qquad\qquad\;\;\;\;\;\;\,\;(a)$} \\
l^{-1}\,\xb & \text{if $\ps=\;\;\;1\qquad\qquad\qquad\;\;\;\;\;\;\,\;(b)$} \\
x_{\be-\al_i} & \text{if $\ps=\;\;\;\frac{1}{2}\;\;$}\&
\text{$\;\;\be-\al_i\q\;\al_i$}\;\,\;(c)\\
x_{\be-\al_i}+\;\frac{m}{l\,r^{ht(\be)-2}}\;\;\,x_{\al_i}\;\!-m\,\xb
& \text{if
$\ps=\;\;\;\frac{1}{2}\;\;$} \& \text{$\;\;\be-\al_i\p\;\al_i$}\;\,\;(d)\\
x_{\be+\al_i}+m\,r^{ht(\be)-1}\,x_{\al_i}-m\,\xb & \text{if
$\ps=-\frac{1}{2}\;\;$} \& \text{$\;\;\;\;\;\;\be\;\;\;\;\q\;\al_i$}\;\,\;(e)\\
x_{\be+\al_i} & \text{if $\ps=-\frac{1}{2}\;\;$} \&
\text{$\;\;\;\;\;\;\be\;\;\;\;\p\;\al_i$}\;\,\;(f)\end{cases}
\end{equation*}

\noindent We now define
$\n(e_i):=\frac{l}{m}(\n_i^2+m\,\n_i-id_{\mathcal{V}})$ and compute
the explicit expression of the endomorphism $\n(e_i)$. For $\ps=0$,
we use the defining relation $r^2+m\,r-1=0$. When $\ps=1$, we get:
\begin{eqnarray*}
\n(e_i)(\xb)&=& \frac{l}{m}\Big(\unsur{l^2}+\frac{m}{l}-1\Big)\,\xb\\
            &=& \Big(\unsur{ml}+1-\frac{l}{m}\Big)\,\xb\\
            &=& \Big(1-\frac{l-\unsur{l}}{m}\Big)\,\xb\\
            &=& \Big(1-\frac{l-\unsur{l}}{m}\Big)\,x_{\al_i}\;\;\text{  as in that
            case $\be=\al_i$}
\end{eqnarray*}
\noindent In the other cases, we notice that:

$$\begin{array}{ccccccccccccc} \text{if} & \nts\nts\ps &\nts =&
\unsur{2}&\nts\nts\nts,&\nts \text{then}& \be-\al_i & \text{is  a
root and}&
\nts(\be-\al_i|\al_i)& \nts = & \unsur{2}-1 & \nts = & \nts\nts\nts -\unsur{2}\\
\text{if}& \nts\nts\ps & \nts = & \nts\nts -\unsur{2}&\nts\nts,&\nts
\text{then}& \be+\al_i & \text{is  a  root  and}&
\nts(\be+\al_i|\al_i)& \nts = & \nts\nts -\unsur{2}+1 &\nts =
&\unsur{2}
\end{array}$$

\noindent The computation of the square $\n_i^2$ uses the third
relation together with the fifth relation and the fourth relation
together with the sixth relation in a very pretty way. The idea is
that a condition of application of $(c)$ (resp $(d)$) with $\be$
yields a condition of application of $(e)$ (resp $(f)$) with
$\be-\al_i$ and a condition of application of $(e)$ (resp $(f)$)
with $\be$ yields a condition of application of $(c)$ (resp $(d)$)
with $\be+\al_i$. And indeed,\\
if $\ps=\unsur{2}$ and $\be-\al_i\q\al_i$, then
$(\be-\al_i|\al_i)=-\unsur{2}$ and $\be-\al_i\q\al_i$. \\It follows
that:
\begin{eqnarray*}
\n_i^2(\xb)&=&  \n_i(x_{\be-\al_i})\qquad\qquad\qquad\qquad\qquad\qquad\qquad\,\text{by application of $(c)$ with $\be$}\\
           & = & x_{(\be-\al_i)+\al_i}+m\,r^{ht(\be-\al_i)-1}\,x_{\al_i}-m\,x_{\be-\al_i}\;\,\text{by
           application of $(e)$ with $\be-\al_i$}\\
           & = & x_{\be}+m\,r^{ht(\be)-2}\,x_{\al_i}-m\,x_{\be-\al_i}
\end{eqnarray*}
Computing $\n(e_i)(\xb)$, the first term on the right hand side of
the last equality above cancels with $-Id_{\mathcal{V}}(\xb)=-\xb$
and the third term cancels with $m\,\n_i(\xb)=m\,x_{\be-\al_i}$. It
remains the term proportional to $\al_i$ which multiplied by the
coefficient $\frac{l}{m}$ yields $l\,r^{ht(\be)-2}\,x_{\al_i}$.\\
If $\ps=\unsur{2}$ and $\be-\al_i\p\al_i$,
$(\be-\al_i|\al_i)=-\unsur{2}$ and $\be-\al_i\p\al_i$, so that:
\begin{eqnarray*}
\n_i^2(\xb) &=&
\n_i(x_{\be-\al_i}+\frac{m}{l\,r^{ht(\be)-2}}\,x_{\al_i}-m\,\xb)
\!\!\!\!\!\!\qquad\text{by application of $(d)$ with $\be$}\\
&=& \xb +
\frac{m}{l^2\,r^{ht(\be)-2}}\,x_{\al_i}-m\,\n_i(\xb)\qquad\text{by
application of $(f)$ with $\be-\al_i$ and by $(b)$}
\end{eqnarray*}
When computing $\n(e_i)$, the two terms on the extremities of the
last equality are canceled and it remains the only term
$\unsur{l\,r^{ht(\be)-2}}\,x_{\al_i}$ after multiplication by the
factor
$\frac{l}{m}$.\\
If $\ps=-\unsur{2}$ and $\be\q\al_i$, then
$(\be+\al_i|\al_i)=\unsur{2}$ and $(\be+\al_i)-\al_i\q\al_i$. Hence,
\begin{eqnarray*}
\n_i^2(\xb)&=&\n_i(x_{\be+\al_i}+m\,r^{ht(\be)-1}\,x_{\al_i}-m\,\xb\;\;\;\;\text{by
application of $(e)$ with $\be$}\\
&=&
x_{\be}+\frac{m\,r^{ht(\be)-1}}{l}\,x_{\al_i}-m\,\n_i(\xb)\qquad\text{by
application of $(c)$ with $\be+\al_i$ and by $(b)$}\\
\end{eqnarray*}
As in the previous cases, the only remaining term in $\n(e_i)(\xb)$
is the multiple of $x_{\al_i}$ and this time, the
coefficient is just $r^{ht(\be)-1}$. \\
Finally if $\ps=-\unsur{2}$ and $\be\p\al_i$, then
$(\be+\al_i|\al_i)=\unsur{2}$ and $(\be+\al_i)-\al_i\p\al_i$. It
follows that:
\begin{eqnarray*}
\n_i^2(\xb)&=&
\n_i(x_{\be+\al_i})\qquad\qquad\qquad\qquad\qquad\qquad\;\;\;\;\;\;\;\text{by
application of $(f)$ with $\be$}\\
&=&
x_{(\be+\al_i)-\al_i}+\frac{m}{l\,r^{ht(\be+\al_i)-2}}\,x_{\al_i}-m\,x_{\be+\al_i}\;\;\text{by
application of $(d)$ with $\be+\al_i$}
\end{eqnarray*}
After canceling the two terms in $\xb$ and in $x_{\be+\al_i}$, the
final result is $\unsur{r^{ht(\be)-1}}\,x_{\al_i}$. Thus we have:
$$\n(e_i)(\xb)= \left\lbrace\begin{array}{ccc} 0 & & \!\!\!\!\!\!\!\!\!\!\!\!\!\!\text{if $\ps=\;\;\;0$}
\\\\
\Big(1-\frac{l-\unsur{l}}{m}\Big)& x_{\al_i} &
\!\!\!\!\!\!\!\!\!\!\!\!\!\!\text{if $\ps=\;\;\;1$}
\\\\
l\,r^{ht(\be)-2}& x_{\al_i} & \qquad\qquad\;\;\;\;\text{if
$\ps=\;\;\;\frac{1}{2}\;\;$}\hfill \&
\text{$\;\;\be-\al_i\q\;\al_i$}\\\\
\unsur{l\,r^{ht(\be)-2}} & x_{\al_i} & \qquad\qquad\;\;\;\;\text{if
$\ps=\;\;\;\frac{1}{2}\;\;$} \hfill\&
\text{$\;\;\be-\al_i\p\;\al_i$}\\\\
r^{ht(\be)-1} & x_{\al_i} & \qquad\qquad\;\;\;\;\text{if
$\ps=-\frac{1}{2}\;\;$} \&
\text{$\;\;\;\;\;\;\be\;\;\;\;\q\;\al_i$}\\\\
\unsur{r^{ht(\be)-1}} & x_{\al_i} & \qquad\qquad\;\;\;\;\text{if
$\ps=-\frac{1}{2}\;\;$} \&
\text{$\;\;\;\;\;\;\be\;\;\;\;\p\;\al_i$}\end{array}\right.$$

\noindent By construction of $\n_i$, $G_i$ is the matrix of the
endomorphism $\n_i$ in the basis
$\mathcal{B}_{\mathcal{V}}:=(x_{\al_1},x_{\al_2},x_{\al_2+\al_1},\dots,x_{\al_{n-1}},x_{\al_{n-1}+\al_{n-2}},x_{\al_{n-1}+\al_{n-2}+\dots+\al_1})$
and by definition, $E_i$ is the matrix of the endomorphism $\n(e_i)$
in the basis $\mathcal{B}_{\mathcal{V}}$.

\begin{thm}
$$\n^{(n)}:\left.\begin{array}{ccc}
B(A_{n-1})& \lra & End_{F}(\mathcal{V})\\
g_i & \longmapsto & \n_i\\
e_i & \longmapsto & \n(e_i)
\end{array}\right.$$
is a representation of the algebra $B(A_{n-1})$ in the $F$-vector
space $\mathcal{V}$.
\end{thm}
\noindent If we can show the theorem, then by the remark above, the
maps on the generators in the conjecture 1. also define a
representation of $B(A_2)$ (resp $B(A_{n-1})$). \\\\\textsc{Proof of
the Theorem:}\\ First we must show that for any two nodes $i$ and
$j$, if $i\nsim j$, then $\n_i\n_j=\n_j\n_i$ and if $i\sim j$ then
$\n_i\n_j\n_i=\n_j\n_i\n_j$. Suppose first that $i\nsim j$. Then we
have $(\al_i|\al_j)=0$. We want to show that:
$$\begin{array}{cccc}
\n_i\n_j(x_{\al_j})& = & \n_j\n_i(x_{\al_j})&\\
\n_i\n_j(\xb)&=& \n_j\n_i(\xb) & \text{for
$\be\not\in\lbrace\al_i,\al_j\}$}
\end{array}$$
It is a direct consequence of $(a)$ and $(b)$ that the two members
of the first equality are equal to $\frac{r}{l}\xalj$. As for the
second equality, we check it by computing the common value depending
on the inner products $\ps$ and $(\be|\al_j)$. We summarize the
results in the following table:

\hspace{-0.68in}
\begin{tabular}{|c|c|c|c|l|}\hline
case & $\ps$ & $(\be|\al_j)$ & $\n_i\n_j(\xb)=\n_j\n_i(\xb)$ & rule\\
\hline
1 & $0$ & $0$ & $r^2\,\xb$ & $(a)$ \\
%\hline
2 & $0$ & $\unsur{2}$ \& $\be-\al_j\q\al_j$ & $r\,x_{\be-\al_j}$ & $(a)$ \& $(c)$\\
%\hline
3 & $0$ & $\unsur{2}$ \& $\be-\al_j\p\al_j$ & $r\,x_{\be-\al_j}+\frac{m}{l\,r^{ht(\be)-1}}\,\xalj-mr\,\xb$& $(a)$ \& $(d)$\\
%\hline
4 & $0$ & $-\unsur{2}$ \& $\be\q\al_j$ & $r\,x_{\be+\al_j}+m\,r^{ht(\be)}\xalj-mr\,\xb$ & $(a)$ \& $(e)$\\
%\hline
5 & $0$ & $-\unsur{2}$ \& $\be\p\al_j$& $r\,x_{\be+\al_j}$ & $(a)$ \& $(f)$ \\
\hline 6 & $\unsur{2}$ \& $\be-\al_i\p\al_i$ & $\unsur{2}$ \&
$\be-\al_j\q\al_j$& $x_{\be-\al_j-\al_i}+\frac{m}{l\,r^{ht(\be)-3}}\,\xali-m\,\baljm$ & $(a)$ \& $(c)$ \& $(d)$\\
%\hline
7 & $-\unsur{2}$ \& $\be\p\al_i$ & $-\unsur{2}$ \& $\be\q\al_j$& $\baljip +m\,r^{ht(\be)}\,\xalj-m\,\balip$ & $(a)$ \& $(e)$ \& $(f)$\\
\hline 8 & $\unsur{2}$ \& $\be-\al_i\q\al_i$ & $-\unsur{2}$ \& $\be\p\al_j$ & $x_{\be+\al_j-\al_i}$ & $(c)$ \& $(f)$\\
%\hline
9 & $\unsur{2}$ \& $\be-\al_i\p\al_i$ & $-\unsur{2}$ \&
$\be\q\al_j$&
\begin{minipage}[t]{2.3in} $x_{\be+\al_j-\al_i}-m\,\baljp-m\,\balim$\\
$+mr^{ht\be}\xalj+\frac{m}{lr^{ht(\be)-3}}\,\xali+m^2\xb$\end{minipage} & $(a)$ \& $(d)$ \& $(e)$\\
\hline
\end{tabular}

$ $ \\\\ The cases $1$ to $9$ are the only possible cases when
excluding $\ps=1$ or $(\be|\al_j)=1$. In the table above, $(a)$ is
often used with $(\al_i|\al_j)=0$. In the first five cases, $\ps=0$.
Since we assumed that $(\al_i|\al_j)=0$, it follows that
$(\be-\al_j|\al_i)=0$ (used in cases $2$ and $3$ to compute
$\n_i\n_j(\xb)$) and $(\be+\al_j|\al_i)=0$ (used in cases $4$ and
$5$ to compute $\n_i\n_j(\xb)$). Next, if both inner products are
equal to $\unsur{2}$, then without loss of generality, $\be$
"starts" with $\al_i$ and "ends with" $\al_j$, \emph{i.e.\/}
$\be=\al_i+\dots+\al_j$ with $i>j$. Then $\be-\al_i\p\al_i$ and
$\be-\al_j\q\al_j$. Moreover,
$(\be-\al_i|\al_j)=(\be|\al_j)=\unsur{2}$ and we read on the
expression giving $\be$ that $(\be-\al_i)-\al_j\q\al_j$. Then we may
use $(c)$ with $\be-\al_i$ and $\al_j$ for computing
$\n_j\n_i(\xb)$. Similarly, $(\be-\al_j|\al_i)=\unsur{2}$ and
$(\be-\al_j)-\al_i\p\al_i$, so we may use $(d)$ with $\be-\al_j$ and
$\al_i$ for computing $\n_i\n_j(\xb)$. If both inner products are
equal to $-\unsur{2}$, then without loss of generality, we may write
$\be=\al_{i-1}+\dots+\al_{j+1}$. We have $\al_j\p\be\p\al_i$ and
also $\be+\al_j\p\al_i$ and $\be+\al_i\q\al_j$. This allows us to
use $(f)$ with $\be+\al_j$ and $\al_i$ when computing
$\n_i\n_j(\xb)$ on one hand and $(e)$ with $\be+\al_i$ and $\al_j$
when computing $\n_j\n_i(\xb)$ on the other hand. The same methods
apply for the cases $8$ and $9$, with $\be=\al_{j-1}+\dots+\al_i$ in
$8$ and $\be=\al_i+\dots+\al_{j+1}$ in $9$.\\

Suppose now that $i\sim j$. Then we have $(\al_i|\al_j)=-\unsur{2}$.
Again by symmetry of the roles played by $i$ and $j$, we need to
show that:
$$\begin{array}{cccc}
\n_i\n_j\n_i(x_{\al_j})& = & \n_j\n_i\n_j(x_{\al_j})&\\
\n_i\n_j\n_i(\xb)&=& \n_j\n_i\n_j(\xb) & \text{for
$\be\not\in\lbrace\al_i,\al_j\}$}
\end{array}$$
The right hand side of the first equation is
$l^{-1}\n_j\n_i(\xalj)$. Next, we need to distinguish between two
cases: $\al_j\p\al_i$ and $\al_j\q\al_i$. The first case is the
short one as it uses the short expressions $(c)$ and $(f)$, which
immediately yield
$\n_i\n_j\n_i(x_{\al_j})=\n_j\n_i\n_j(x_{\al_j})=l^{-1}x_{\al_i}$.
In the second case, we have,
$\n_i(\xalj)=x_{\al_i+\al_j}+m\,\xali-m\xalj$. Since
$(\al_i+\al_j|\al_j)=1-\unsur{2}=\unsur{2}$ and
$(\al_i+\al_j)-\al_j=\al_i\p\al_j$, we get $\n_j\n_i(\xalj)=\xali
+\frac{m}{l}x_{\al_j}-m\,x_{\al_i+\al_j}+m\,x_{\al_i+\al_j}-\frac{m}{l}x_{\al_j}$,
\emph{i.e.\/} $\n_j\n_i(\xalj)=\xali$. In that case the final result
is again $l^{-1}x_{\al_i}$ for both members of the first equality.
The proof and results for the second equality are gathered in the
following table:\\

\hspace{-1.2in}
\begin{tabular}{|c|c|c|c|l|}\hline
case & $\ps$ & $(\be|\al_j)$ & $\n_i\n_j\n_i(\xb)=\n_j\n_i\n_j(\xb)$ & rule\\
\hline
1 & $0$ & $0$ & $r^3\,\xb$ & $(a)$ \\
%\hline
2 & $0$ & $\unsur{2}$ \& $\be-\al_j\q\al_j$ & $r\,x_{\be-\al_j-\al_i}$ & $(a)$ \& $(c)$\\
%\hline
3 & $0$ & $\unsur{2}$ \& $\be-\al_j\p\al_j$ &
\begin{minipage}[t]{2.7in}
$r\,x_{\be-\al_i-\al_j}-mr\,x_{\be-\al_j}-mr^2x_{\be}$\\
$+\frac{m}{lr^{ht\be-3}}\,x_{\al_i+\al_j}+\frac{m}{lr^{ht\be-2}}\,x_{\al_i}-\frac{m^2}{lr^{ht\be-3}}\,x_{\al_j}$
\end{minipage}& all of them\\
%\hline
4 & $0$ & $-\unsur{2}$ \& $\be\q\al_j$ & \begin{minipage}[t]{2.7in}$r\,x_{\be+\al_i+\al_j}-mr\,x_{\be+\al_j}-mr^2\xb$\\
$+mr^{ht\be}x_{\al_i+\al_j}+mr^{ht\be-1}x_{\al_i}-m^2r^{ht\be}x_{\al_j}$\end{minipage} & all of them\\
%\hline
5 & $0$ & $-\unsur{2}$ \& $\be\p\al_j$& $r\,x_{\be+\al_i+\al_j}$ & $(a)$ \& $(f)$ \\
\hline 6 & $\unsur{2}$ \& $\be-\al_i\p\al_i$ & $\unsur{2}$ \&
$\be-\al_j\q\al_j$& $\unsur{l}\,\xb+\frac{m}{l}\,\xalj-\frac{m}{l}\,x_{\be-\al_j}$ & all of them except $(a)$\\
\hline
7 & $\unsur{2}$ \& $\be-\al_i\q\al_i$ & $-\unsur{2}$ \& $\be\q\al_j$ & $r\,\xb+mr^{ht\be-1}\,\xali-mr\,x_{\be-\al_i}$ & $(a)$ \& $(c)$ \& $(e)$ \& $(f)$\\
%\hline
8 & $\unsur{2}$ \& $\be-\al_i\p\al_i$ & $-\unsur{2}$ \&
$\be\p\al_j$&
$r\,\xb+\frac{m}{lr^{ht\be-2}}\,\xalj-mr\,x_{\be+\al_j}$ & $(a)$ \& $(c)$ \& $(d)$ \& $(f)$\\
\hline
\end{tabular}

$ $ \\\\ The values for the inner products and the fact that $i$ and
$j$ are adjacent nodes determine in each case the only possible
expression for $\be$. In case $2$, we have
$\be=\al_l+\dots+\al_{i+1}+\al_i+\al_j$ with $l\geq i+1$. Then
$(\be-\al_j|\al_i)=\unsur{2}$ and $\be-\al_j-\al_i\q\al_i$ and
$(\be-\al_j-\al_i|\al_j)=0$. In case $3$, we have
$\be=\al_j+\al_i+\al_{i-1}+\dots+\al_k$ with $k\leq
i-1$.\\
In case $4$, we have $\be=\al_l+\dots+\al_{j+1}$ with $l\geq j+1$
and $\al_j\q\al_i$.\\
In case $5$, we have $\be=\al_{j-1}+\dots+\al_k$ with $k\leq j-1$
and $\al_i\q\al_j$.\\
If $(\be|\al_i)=(\be|\al_j)=\unsur{2}$ (case $6$), the only
possibility is $\be=\al_i+\al_j$ with without loss of generality
$\al_i\q\al_j$. We get: $$\begin{array}{ccc} \n_i\n_j\n_i(\xb) &=&
\unsur{l}\,\xb+\frac{m}{lr^{ht\be-2}}\,\xalj-\frac{m}{l}\,x_{\be-\al_j}\\
\n_j\n_i\n_j(\xb) &=&
\unsur{l}\,\xb+\frac{mr^{ht\be-2}}{l}\,\xalj-\frac{m}{l}\,x_{\be-\al_j}
\end{array}$$
Since $ht\be=2$, the two results are the same. \\The case
$(\be|\al_i)=(\be|\al_j)=-\unsur{2}$ is impossible.\\
In case $7$, we have $\be=\al_l+\dots+\al_i$, with $l>i$ and $i\sim
j$, $i>j$. Finally in case $8$, we have $\be=\al_i+\dots+\al_k$ with
$k<i$ and $j\sim i$, $j>i$.

To get a representation, there are two more things to check:\\
$\n_i\n(e_i)=l^{-1}\n(e_i)$ and $\n(e_i)\n_j\n(e_i)=l\n(e_i)$ for
any two adjacent nodes $i$ and $j$. The first equality holds since
$\n_i(\xali)=l^{-1}\xali$ and $\n(e_i)$ is always a multiple of
$\xali$. Next, given two adjacent nodes $i$ and $j$, their inner
product is equal to $-\unsur{2}$. We want to show that
$\n(e_i)\n_j(\xali)=l\xali$. Suppose first that $\al_i\p\al_j$. Then
$\n_j(\xali)=x_{\al_i+\al_j}$. Since
$(\al_i+\al_j|\al_i)=1-\unsur{2}=\unsur{2}$ and since
$\al_j\q\al_i$, we get $\n(e_i)(x_{\al_i+\al_j})=l\,\xali$. Next
suppose that $\al_i\q\al_j$. We have
$\n(e_i)(x_{\al_i+\al_j}+m\,\xalj-m\,\xali)=
(\unsur{l}+m(1-x))\xali=l\,\xali$. This ends the proof of the
theorem $3$.

Finally, it is direct to observe that $\n(e_ie_j)=0$ when $i\nsim
j$. Thus, the so built representation factors through the quotient
$B(A_{n-1})/I_2$.

\section{Proof of the Main Theorem}
\subsection{properties of the representation and the case $n=6$}
Suppose that the representation $\n$ is not irreducible. Then there
exists $\mathcal{U}$ an $F$-vector subspace of $\mathcal{V}$ such
that $\mathcal{U}\neq 0$, $\mathcal{U}\neq \mathcal{V}$ and
$\n(w)(\mathcal{U})\subseteq\mathcal{U}$ for any element $w$ of the
BMW algebra $B(A_{n-1})$. We will show that
$\n(X_{ij})(\mathcal{U})=0$ for any $1\leq i<j\leq n$. And indeed,
suppose that there exists $u\in\mathcal{U}$ such that
$\n(e_i)(u)\neq 0$. Then by the action of $\n(e_i)$ described above,
we get that $\xali$ is in $\mathcal{U}$. From there we deduce that
all the $\xb$'s are in $\mathcal{U}$ and there are three steps.
First we have $\n_{i+1}(\xali)\in\mathcal{U}$, \emph{i.e.\/}
$x_{\al_{i+1}+\al_i}\in\mathcal{U}$. By successive applications of
$(f)$, it follows that all the $x_{\al_l+\dots+\al_i}$ with $l\geq
i$ are in $\mathcal{U}$. Next, starting again from
$x_{\al_{i+1}+\al_i}$, we have
$\n_i(x_{\al_{i+1}+\al_i})=x_{\al_{i+1}}$, so that
$x_{\al_{i+1}}\in\mathcal{U}$ and proceeding inductively with step
$1$, all the $x_{\al_l}$ with $l\geq i$ are also in $\mathcal{U}$.
Third,
$\n_{i-1}(\xali)=x_{\al_i+\al_{i-1}}+m\,x_{\al_{i-1}}\,modulo\,
F\xali$ and another application of $\n_{i-1}$ yields
$\n_{i-1}(x_{\al_i+\al_{i-1}}+m\,x_{\al_{i-1}})=\xali+\frac{m}{l}\,x_{\al_{i-1}}$,
from which we derive that $x_{\al_{i-1}}$ is in $\mathcal{U}$. Again
by induction we also have
$x_{\al_1},\dots,x_{\al_{i-2}}\in\mathcal{U}$. We gather all these
informations as follows: by the steps $2$ and $3$, all the
$x_{\al_i}$'s are in $\mathcal{U}$. Furthermore, by the first step
applied to the $\xali$'s, all the $\xb$'s are in $\mathcal{U}$. This
contradicts $\mathcal{U}\neq\mathcal{V}$. We conclude that
$\n(e_i)(\mathcal{U})=0$ for all $i$. It remains to show that
$\forall j\geq i+2,\,\n(X_{ij})(\mathcal{U})=0$. Using the
conjugation formulas, an expression for $\n(X_{ij})$ is:
$$\n(X_{ij})=\n_{j-1}\dots\n_{i+1}\n(e_i)\n_{i+1}^{-1}\dots\n_{j-1}^{-1}$$
For any positive root $\be$, $\n(X_{ij})(\xb)$ is a multiple of
$\n_{j-1}\dots\n_{i+1}(\xali)$, hence a multiple of
$x_{\al_i+\dots+\al_{j-1}}$. We will remember this fact:
\begin{Prop}
\hfill\\ $\n(X_{ij})(\xb)$ is always a multiple of
$x_{\al_i+\dots+\al_{j-1}}$. \\In other words, the
$\{(1+2+\dots+j-2)+(j-i)\}$th row is the only non zero row in the
matrix $Mat_{\mathcal{B}_{\mathcal{V}}}\n(X_{ij})$.
\end{Prop}
\noindent Thus, if $\n(X_{ij})(u)\neq 0$, some $u\in\mathcal{U}$, it
comes $x_{\al_i+\dots+\al_{j-1}}\in\mathcal{U}$. As expected,
$\al_i+\dots+\al_{j-1}$ is the positive root having nodes $i$ and
$j$ joined on its top horizontal line like in $X_{ij}$. Then, by
successive applications of $(c)$, we get
$x_{\al_{j-1}}\in\mathcal{U}$. By the same arguments as above, we
deduce from that fact that all the $\xb$'s are in fact in
$\mathcal{U}$, which is a contradiction. Again we conclude that
$\n(X_{ij})(\mathcal{U})=0$. Thus, we have: $$\n\Bigg(\sum_{1\leq
i<j\leq n}X_{ij}\Bigg)(\mathcal{U})=0$$ Since $S$ is the matrix of
the endomorphism $\n(\sum_{1\leq i<j\leq n}X_{ij})$ in the basis
$\mathcal{B}_{\mathcal{V}}$, $\n$ not irreducible implies $det(S)=0$
(otherwise $\mathcal{U}$ would be trivial). A direct consequence of
that fact is that the sufficient condition on $l$ found with Maple
so that the representation is irreducible holds in the case $n=6$.
Furthermore, we obtained with Maple that for each of the values of
$l$, there exists a nonzero vector in the kernel of the matrix $S$.
Such vectors belong to the intersection $\cap Ker\;\n(X_{ij})$ by
proposition $2$. Since the $F$-vector subspace $\cap
Ker\;\n(X_{ij})$ is stable under the $\n(w)'s,\,w\in B(A_5)$ by the
same arguments as those exposed at the end of part $5$, the Main
Theorem is true in the case $n=6$.
%\subsection{The converse}
\subsection{The case $l=\unsur{r^{2n-3}}$} In this section, we show
that for $l=\unsur{r^{2n-3}}$ the representation is reducible and we
give a necessary and sufficient condition on $l$ so that there
exists a one dimensional invariant subspace of $
\V$. %the
%$B(A_{n-1})$-module $\bigcap_{1\leq i<j\leq n}\,Ker\,\n(X_{ij})$ is
%one-dimensional as a vector space over $F$.
We introduce new notations. Recall that each of the positive roots
is associated with a top line joining two nodes in a unique way. We
will denote by $w_{ij}$ ($i<j$) the root whose associated tangle has
nodes $i$ and $j$ joined on its top horizontal line. Note that
$w_{ij}$ is nothing else but $\al_i+\dots+\al_{j-1}$. We will again
denote by $w_{ij}$ the element $x_{w_{ij}}$ of $\V$. Our result is
the following:
\begin{thm} Assume $(r^2)^2\neq 1$.\\\\
\textbf{Suppose} $\mathbf{n=3}$. There exists a one dimensional
invariant subspace of $\mathcal{V}$ if and only if $l=\unsur{r^3}$
or $l=-r^3$.\\\\ If those values are distinct then, if such a space
exists, it is unique and
\\\\ \indent If $l=\unsur{r^3}$, it is spanned by
$w_{12}+r\,w_{13}+r^2\,w_{23}$
\\ \indent If $l=-r^3$,
it is spanned by $w_{12}-\unsurr\,w_{13}+\unsur{r^2}\,w_{23}$\\\\
Moreover, if $-r^3=\unsur{r^3}$, \emph{i.e\/} $r^6=-1$, then there
exists exactly two one-dimensional invariant subspaces:
$$Span_F(w_{12}+r\,w_{13}+r^2\,w_{23})\;\;\text{and}\;\;
Span_F(w_{12}-\unsurr\,w_{13}+\unsur{r^2}\,w_{23})$$ In other words,
the sum
$$Span_F(w_{12}+r\,w_{13}+r^2\,w_{23})+Span_F(w_{12}-\unsurr\,w_{13}+\unsur{r^2}\,w_{23})$$
is direct.\\\\
\textbf{Suppose} $\mathbf{n\geq 4}$. There exists a one dimensional
invariant subspace of $\mathcal{V}$ if and only if
$l=\unsur{r^{2n-3}}$. If so, it is spanned by $\sum_{1\leq s<t\leq
n}r^{s+t}w_{st}$\\
\end{thm}
\noindent\textsc{Proof of the theorem:} suppose that there exists a
one-dimensional invariant subspace $\mathcal{U}$ of $\mathcal{V}$,
say $\U$ is spanned by $u$. We have seen that
$$\mathcal{U}\subseteq\bigcap_{1\leq i<j\leq n}\,Ker\,\n(X_{ij})$$
In particular, we must have $\n(e_i)(\mathcal{U})=0$ for all $i$. By
definition of $\n(e_i)$, this implies
$(\n_i^2+m\,\n_i-id_{\mathcal{V}})(u)=0$. For each $i$, let $\la_i$
be the scalar so that $\n_i(u)=\la_iu$. It comes
$\la_i^2+m\la_i-1=0$, hence $\la_i\in\{r,-\unsur{r}\}$. Furthermore,
by using the braid relation $\n_i\n_j\n_i=\n_j\n_i\n_j$, we see that
all the $\la_i$'s for $1\leq i\leq n-1$ must take the same value $r$
or $-\unsur{r}$. Let's denote this common value by $\la$. Thus all
the $\la_i$'s are determined by $\la_1$. We want to determine wether
$\la_1=r$ or $\la_1=-\unsur{r}$.  Let's write
$$u=\sum_{1\leq i<j\leq n}\mu_{ij}w_{ij}$$
We claim that the coefficients in this sum are related in a certain
way. First we notice that an action of $\n_i$ on
$\mu_{i,k}\,w_{i,k}$ makes a term in $w_{i+1,k}$ appear with
coefficient $\mu_{i,k}$ and an action of $\n_i$ on
$\mu_{i+1,k}\,w_{i+1,k}$ makes a term in $w_{i+1,k}$ appear with
coefficient $-m \mu_{i+1,k}$. Moreover, $w_{i+1,k}$ cannot be
obtained from other $w_{st}$'s by acting with $\n_i$. Thus, if
$\n_i(u)=\la u$, the coefficients $\mu_{i,k}$ and $\mu_{i+1,k}$ must
be related by the relation:
$$
\mu_{i,k}-m\,\mu_{i+1,k}=\la\,\mu_{i+1,k}
$$
We note that if one of the coefficients $\mu_{i,k}$, $\mu_{i+1,k}$
is zero, then the other one is zero. In other words the two
coefficients are either both zero or both non zero. They are related
by \begin{equation} \mu_{i+1,k}=\la\,\mu_{i,k}\end{equation}
Similarly, by operating on $w_{l,i}$ and $w_{l,i+1}$ with $\n_i$, we
get the equation:
\begin{equation}\mu_{l,i+1}=\la\,\mu_{l,i}\end{equation} Now the
relations $(70)$ and $(71)$ show that all the coefficients
$\mu_{ij}$'s are nonzero as $u$ itself is nonzero. Moreover, up to a
multiplication by a scalar, $u$ is of the form
$$u=\sum_{1\leq i<j\leq n}\la^{i+j}w_{ij}$$
Suppose first $n=3$. Then $u$ can be written in the form
$$u=\la^3\,w_{12}+\la^4\,w_{13}+\la^5\,w_{23}$$
Let's compute $\n_1(u)$ and $\n_2(u)$:
$$\begin{array}{ccccccccc}
\n_1(u) & = & (\frac{\la^3}{l}+\,m\,\la^5) & w_{12} & + &
\la^5\,w_{13} & + & \la^6 & w_{23}\\
\n_2(u) & = & \la^4 & w_{12} & + & \la^5\,w_{13} & + &
(\frac{\la^5}{l}+\frac{m\,\la^4}{l}) & w_{23}
\end{array}$$
Since $\n_1(u)=\la\,u$ and $\n_2(u)=\la\,u$, we must have:
$$\left\lbrace
\begin{array}{ccccccc} \frac{\la^3}{l}+m\,\la^5=\la^4&
\textit{i.e.}& \unsur{l}=\la(1-m\,\la)&\textit{i.e.}&
l=\unsur{\la^3} &\text{as}& 1-m\,\la=\la^2\\
\frac{\la^5}{l}+\frac{m\,\la^4}{l}=\la^6 &\textit{i.e.}& l=
\unsur{\la}(1+\frac{m}{\la}) &\textit{i.e.}& l=\unsur{\la^3}&
\text{as}& 1+\frac{m}{\la}=\unsur{\la^2}\end{array}\right.$$
\noindent Thus, in the case $n=3$, if there exists a one dimensional
invariant subspace of $\mathcal{V}$ then $l$ must take the values
$\unsur{r^3}$ or $-r^3$. Conversely, let's consider the two vectors:
$$\begin{array}{cccccccc}
&u_r &=& w_{12} &+& r\,w_{13} &+& r^2\,w_{23}\\
&u_{-\unsurr} &=& w_{12} &-& \unsurr\,w_{13}&+&\unsur{r^2}\,w_{23}
\end{array}$$
We read on the equations giving the expressions for $\n_1(u)$ and
$\n_2(u)$ that
$$\begin{array}{cccccccccc} \text{If} &
l&=&\unsur{r^3}&\text{then} & \n_1(u_r)&=&\n_2(u_r)&=&r\;u_r\\
\text{If} & l &=& -r^3 & \text{then} &
\n_1(u_{-\unsurr})&=&\n_2(u_{-\unsurr})&=&-\unsurr\;u_{-\unsurr}
\end{array}$$
The theorem is thus proved in the case $n=3$. Suppose now $n\geq 4$.
By the above, the coefficient $\mu_{34}$ of $u$ is nonzero. An
action of $\n_1$ on $w_{34}$ is a multiplication by $r$ and an
action of $\n_1$ on the other $w_{ij}$'s does not affect the
coefficient of $w_{34}$. This last point forces $\la=r$ and $u$ can
be rewritten:
$$u=\sum_{1\leq i<j\leq n}r^{i+j}w_{ij}$$
We will now see how this expression of $u$ forces the value of $l$.
It suffices for instance to look at the action of $\n_1$ on $u$ and
the resulting coefficient in $w_{12}$. A term in $w_{12}$ appears
when $\n_1$ acts on $w_{12}$ and on the $w_{2j}$'s for $3\leq j\leq
n$. In the first case the resulting coefficient is $\frac{r^3}{l}$
and in the second case the resulting coefficient is
$m\,r^{j-3}\times r^{j+2}$. Hence we get the equation:
\begin{equation} \frac{m}{r}\sum_{j=3}^{n}(r^2)^{j} +
\frac{r^3}{l}=r^4
\end{equation}
from which we derive that $l=\unsur{r^{2n-3}}$. Conversely, for this
value of $l$ and letting $u=\sum_{1\leq i<j\leq n}r^{i+j}w_{ij}$, we
check that $\n_i(u)=r\,u$ for all integers $1\leq i\leq n-1$. We
start the proof by showing that the coefficient $r^{2i+1}$ of
$w_{i,i+1}$ is multiplied by $r$ when we act by $\n_i$. We notice
that the root $\al_i$ appears when we act by $\n_i$ on:
\begin{enumerate}
\item $\al_i$
\item $\be$ with $(\be|\al_i)=\unsur{2}$ and $\be-\al_i\p\al_i$
\item $\be$ with $(\be|\al_i)=-\unsur{2}$ and $\be\q\al_i$
\end{enumerate}
In the first case, it yields the coefficient $r^{2n-3}\times
r^{2i+1}=r^{2n+2i-2}$ \\In the second case, $\be=w_{k,i+1}$ with
$k=1,\dots, i-1$. It yields the coefficient:
$$\sum_{k=1}^{i-1}\frac{m\,r^{2n-3}}{r^{i-k-1}}\times
r^{k+i+1}=m\,r^{2n-1}\sum_{k=1}^{i-1}(r^2)^k=r^{2n}-r^{2i+2n-2}$$ In
the third case, $\be=w_{i+1,s}$ with $s=i+2,\dots, n$. It yields the
coefficient:
$$\sum_{s=i+2}^{n}m\,r^{s-i-2}\times
r^{s+i+1}=\frac{m}{r}\,\sum_{s=i+2}^{n}(r^2)^s=r^{2i+2}-r^{2n}$$
When summing the three coefficients, it only remains
$r^{2i+2}=r\times r^{2i+1}$. Thus the coefficient $r^{2i+1}$ of
$w_{i,i+1}$ is multiplied by $r$.\\ Next, given a positive root
$\be$, if none of the nodes $i-1,\,i,\,i+1$ is in the support of
$\be$ or if all three nodes $i-1,\,i,\,i+1$ are in the support of
$\be$, then it comes:
$$\n_i(x_{\be})=r\,x_{\be}$$
Thus, we only need to study the combined effect of $\n_i$ on
$w_{k,i}$, $w_{k,i+1}$, with $k\leq i-1$ on one hand and $w_{i,l}$,
$w_{i+1,l}$, with $l\geq i+2$ on the other hand. \\We have:
\begin{eqnarray}
r^{k+i}\;\n_i(w_{k,i})&=&r^{k+i}\;w_{k,i+1}\\
r^{k+i+1}\;\n_i(w_{k,i+1})&=&r^{k+i+1}\;w_{k,i}-m\,r^{k+i+1}\;w_{k,i+1}\;
modulo\, F\,\xali
\end{eqnarray}
So we get:
$$\n_i\Big(r^{k+i}\,w_{k,i}+r^{k+i+1}\,w_{k,i+1}\Big)=r^{k+i+1}\,w_{k,i}+r^{k+i+2}\,w_{k,i+1}\;modulo\,F\,\xali$$
Similarly, we have:
\begin{eqnarray}
r^{l+i}\;\n_i(w_{i,l})&=&r^{l+i}\;w_{i+1,l}\\
r^{l+i+1}\;\n_i(w_{i+1,l})&=&r^{l+i+1}\;w_{i,l}-m\,r^{l+i+1}\;w_{i+1,l}\;modulo\,
F\,\xali
\end{eqnarray}

\noindent so that:
$$\n_i\Big(r^{l+i}\,w_{i,l}+r^{l+i+1}\,w_{i+1,l}\Big)=r^{l+i+1}\,w_{i,l}+r^{l+i+2}\,w_{i+1,l}\;modulo\,F\,\xali$$

\noindent This ends the proof of Theorem $4$.

\subsection{The cases $l=\unsur{r^{n-3}}$ and $l=-\unsur{r^{n-3}}$}
In this section we show a necessary and sufficient condition on $l$
and $r$ so that there exists an irreducible $(n-1)$-dimensional
invariant subspace of the $F$-vector space $\mathcal{V}$, under the
condition that the Iwahori-Hecke algebra of the symmetric group
$Sym(n)$ with parameter $r^2$ over the field $F$ is semisimple. This
last condition is equivalent to $(r^2)^k\neq 1$ for all
$k=1,\,\dots,\,n$. Explicitly we have the two theorems:
\begin{thm}\hfill\\\\
Let $n$ be a positive integer with $n\geq 3$ and $n\neq 4$. Let's
assume that the Iwahori-Hecke agebra of the symmetric group $Sym(n)$
with parameter $r^2$ over the field $F$ is semisimple. Then, there
exists an irreducible $(n-1)$-dimensional invariant subspace of
$\mathcal{V}$ if and only if $l=\unsur{r^{n-3}}$ or
$l=-\unsur{r^{n-3}}$. \\\\
If so, it is spanned by the $v_i$'s, $1\leq i\leq n-1$, where $v_i$
is defined by the formula:
$$v_i:=\Big(\unsur{r}\,-\,\unsur{l}\Big)w_{i,i+1}\,+\sum_{k=i+2}^{n}r^{k-i-2}(w_{i,k}\,-\,\unsurr\;w_{i+1,k})\,+\;\e_l\,\sum_{s=1}^{i-1}r^{n-i-2+s}(w_{s,i}\,-\,\unsurr\;w_{s,i+1})$$
\begin{center} with \hspace{.1in}
$\begin{cases}\e_{\unsur{r^{n-3}}}\;\;=\,1\\\e_{-\unsur{r^{n-3}}}=-1\end{cases}$
\end{center}
\end{thm}
\begin{thm} (Case $n=4$) \hfill\\ Let's assume that the Iwahori-Hecke algebra of the symmetric group $sym(4)$ with parameter $r^2$ over the field $F$
is semisimple. Then:\\\\ there exists an irreducible $3$-dimensional
invariant subspace of $\mathcal{V}$ if and only if
$l\in\lb\unsurr,-\unsurr,-r^3\rb$.\\\\
If $l\in\lb -\unsurr,\unsurr\rb$, it is spanned by the vectors:
$$\begin{array}{ccccccc} v_1&=&(\unsurr-\unsur{l})w_{12} &+&
(w_{13}-\unsurr\,w_{23}) &+&
r\,(w_{14}-\unsurr\,w_{24})\\
v_2&=&(\unsurr-\unsur{l})w_{23} &+& (w_{24}-\unsurr\,w_{34}) &+&
\e_l\,r\,(w_{12}-\unsurr\,w_{13})\\
v_3 &=& (\unsurr-\unsur{l})w_{34} &+&
\e_l\,(w_{13}-\unsurr\,w_{14})&+& \e_l(w_{23}-\unsurr\,w_{24})
\end{array}$$
\begin{center} where \hspace{.1in}
$\begin{cases}\e_{\unsurr}\;\;=\,1\\\e_{-\unsurr}=-1\end{cases}$\end{center}
\noin If $l=-r^3$, it is spanned by the vectors:
$$\hspace{-.8in}\begin{array}{cccccccccccccc} v_1\eg & & r & w_{23}\;\;+&
&w_{13}\;\;+& (\unsurr+\unsur{r^3})&w_{34}\;\;-&
& w_{24}\;\;-& \unsurr & w_{14}\\
v_2\eg -r&w_{12}\;\;-& & & r^2 & w_{13}\;\;-&
\unsurr & w_{34}\;\;-&\unsur{r^2}& w_{24}\;\;+& (r+\unsurr)&w_{14}\\
v_3\eg (r+\unsur{r^3})&w_{12}\;\;+&\unsurr & w_{23}\;\;-&
&w_{13}\;\;+& & & &w_{24}\;\;-& r&w_{14}
\end{array}$$
\end{thm}
\noindent\textsc{Joint proof of the Theorems:}\\ We first recall
that if $\mathcal{U}$ is a proper invariant subspace of $\V$, it
must be annihilated by all the algebra elements $X_{ij}$'s. In
particular it is annihilated by all the $e_i$'s. Thus, the action of
the BMW algebra $B(A_{n-1})$ on the $F$-vector space $\U$ is an
Iwahori-Hecke algebra action. Further, for $n\geq 4$ and $n\neq 6$,
there are exactly two inequivalent irreducible (n-1)-dimensional
representations of the Iwahori-Hecke algebra of the symmetric group
$Sym(n)$ with parameter $r^2$ over the field $F$ and they are
respectively given by the matrices:
\begin{center} $M_1=\begin{bmatrix}
-1/r & 1/r & & & \\
 &r & & &\\
 & & r & & \\
 & & & \ddots & \\
 & & & &r\end{bmatrix},\;
M_2=\begin{bmatrix} r & & & & & \\
r & -1/r & 1/r & & & \\
  &  & r & & & \\
  &   &   &r& & \\
  &   &   & &\ddots & \\
  &   &   & &       &r
 \end{bmatrix},$\end{center}

$$M_3=\begin{bmatrix}
r& & & & & & \\
 &r& & & & & \\
 &r&-1/r&1/r& & & \\
 & & &r& & & \\
 & & & &r& &\\
 & & & & &\ddots& \\
 & & & & &      &r
\end{bmatrix},\;\dots,\; M_{n-2}=\begin{bmatrix}
r&       & & &    &   \\
 &\ddots & & &    &   \\
 &       &r& &    &   \\
 &       & &r&    &   \\
 &       & &r&-1/r&1/r\\
 &       & & &    & r
\end{bmatrix},$$
$M_{n-1}=\begin{bmatrix}
r & & & & \\
  &\ddots& &  &\\
  &      &r & &\\
  &      &  &r & \\
  &      &  &r&-1/r
\end{bmatrix}$ \\\\

\noindent and for the conjugate representation $$N_1=\begin{bmatrix}
r & -r & & & \\
 &-1/r & & &\\
 & & -1/r & & \\
 & & & \ddots & \\
 & & & &-1/r\end{bmatrix},\;
N_2=\begin{bmatrix} -1/r & & & & & \\
-1/r & r & -r & & & \\
  &  & -1/r & & & \\
  &   &   &-1/r& & \\
  &   &   & &\ddots & \\
  &   &   & &       &-1/r
 \end{bmatrix},$$

$$N_3=\begin{bmatrix}
-1/r& & & & & & \\
 &-1/r& & & & & \\
 &-1/r&r&-r& & & \\
 & & &-1/r& & & \\
 & & & &-1/r& &\\
 & & & & &\ddots& \\
 & & & & &      &-1/r
\end{bmatrix},\;\dots,$$ $$ N_{n-2}=\begin{bmatrix}
-1/r&       & & &    &   \\
 &\ddots & & &    &   \\
 &       &-1/r& &    &   \\
 &       & &-1/r&    &   \\
 &       & &-1/r&r&-r\\
 &       & & &    & -1/r
\end{bmatrix},\;
N_{n-1}=\begin{bmatrix}
-1/r & & & & \\
  &\ddots& &  &\\
  &      &-1/r& &\\
  &      &  &-1/r & \\
  &      &  &-1/r&r
\end{bmatrix}$$ \\

\noindent where all the matrices are square matrices of size $n-1$
and where the blanks are filled with zeros. Note that those two
representations are equivalent for $n=3$. We now show that the
latter representation cannot occur inside $\V$ when $n\geq 5$.
Indeed, suppose that there exists $\U$ an irreducible $(n-1)$-
dimensional invariant subspace of $\V$. Let
$(v_1,\,\dots,\,v_{n-1})$ be a basis of $\U$ in which the matrices
of the left actions of the $g_i$'s are the $N_i$'s. We have the
relations (when the indices make sense):
$$(\bigtriangledown)\left|\begin{array}{cccc} \n_t(v_i)&=&-1/r\;v_i,&\text{for}\;
t\not\in\{
i-1,\,i,\,i+1\}\\
\n_i(v_i)&=&r\;v_i&\\
\n_{i+1}(v_i)&=&-1/r(v_i+v_{i+1})&\\
\n_{i-1}(v_i)&=&-1/r\;v_i-r\;v_{i-1}&
\end{array}\right.$$

\begin{lemma}
It is impossible to have such a set of relations.
\end{lemma}
\noindent Throughout the proof of the lemma and the theorems, we
will make an extensive use of the following equalities. We have for
any node $q$:

%$\forall k\geq q+2,\;$\begin{equation}\begin{aligned}\n_q(w_{q,k}) &=  w_{q+1,k}\\
%\n_q(w_{q+1,k})&= w_{q,k}-m\,
%w_{q+1,k}\;\;\text{modulo}\;F\,x_{\al_q}\end{aligned}\end{equation}

%$\forall s\leq q-1,$\begin{equation}\begin{aligned}\n_q(w_{s,q}) &=  w_{s,q+1}\\
%\n_q(w_{s,q+1}) &= w_{s,q}-m\,
%w_{s,q+1}\;\;\text{modulo}\;F\,x_{\al_q}\end{aligned}\end{equation}\\

$$\begin{array}{ccccc} \forall k\geq q+2, & \n_q(w_{q+1,k}) &=&
w_{q,k}-m\,w_{q+1,k}\;\;\text{modulo}\;F\,x_{\al_q}& (\star)\\
\forall s\leq q-1, & \n_q(w_{s,q+1}) &=& w_{s,q}-m\,
w_{s,q+1}\;\;\text{modulo}\;F\,x_{\al_q}& (\star\star)
\end{array}$$

% We claim that $(77)$ (resp $(78)$) is the only way to &obtain
%terms in $w_{q+1,k}$, $w_{q,k}$ for a given $k\geq q+2$ (resp
%$w_{s,q}$, $w_{s,q+1}$ for a given $s\leq q-1$) when acting &by
%$\n_q$. Indeed, an action by $\n_q$ results in a linear &combination
%of the following operations:

\noindent We claim that $(\star)$ (resp $(\star\star)$) is the only
way to get a term in $w_{q,k}$ for a given $k\geq q+2$ (resp
$w_{s,q}$ for a given $s\leq q-1$) when acting by $\n_q$. Indeed, an
action by $\n_q$ results in a linear combination of the following
operations:
\begin{list}{\texttt{$(i)$}}{}
\item create a root $\al_q$
\end{list}
\begin{list}{\texttt{$(ii)$}}{}
\item add or subtract a root $\al_q$ to an existing positive root
\end{list}
\begin{list}{\texttt{$(iii)$}}{}
\item leave an element of $\V$ unchanged
\end{list}
in the way described in the defining representation. %Then,
%obviously, $w_{q,k}$ (resp $w_{s,q}$) can only be obtained from
%$w_{q+1,k}$ (resp $w_{s,q+1}$) and $w_{q+1,k}$ (resp $w_{s+1,q}$)
%can only be obtained from $w_{q,k}$ and $w_{q+1,k}$ (resp $w_{s,q}$
%and $w_{s,q+1}$).\\
With this remark, we see that $w_{q,k}$ (resp $w_{s,q}$) can only be
obtained from $w_{q+1,k}$ or $w_{q,k}$ itself (resp from $w_{s,q+1}$
or $w_{s,q}$ itself). Since $\n_q(w_{q,k})=w_{q+1,k}$ (resp
$\n_q(w_{s,q})=w_{s,q+1}$), $w_{q,k}$ (resp $w_{s,q}$) can only be
obtained from $w_{q+1,k}$ (resp from $w_{s,q+1}$). Hence, if
$\n_q(v_i)=\la\,v_i$ with $\la\in\{r,-1/r\}$ and the $\mu_{s,t}$'s
are the respective coefficients of the $w_{s,t}$'s in $v_i$, the
relations $(\star)$ and $(\star\star)$ respectively imply that:
\begin{eqnarray}
\forall k\geq q+2,\;\mu_{q+1,k}&=&\la\mu_{q,k}\\
\forall s\leq q-1,\;\mu_{s,q+1}&=&\la\mu_{s,q}
\end{eqnarray}
We will make an extensive use of these equalities. It will also be
useful to note that for any node $q\in\{1,\,\dots,\,n-1\}$, we have:
\begin{eqnarray*}
\forall k\geq q+2,\;\n_q(w_{q,k})&=&w_{q+1,k}\\
\forall s\leq q-1,\;\n_q(w_{s,q})&=&w_{s,q+1}
\end{eqnarray*}
as it has already been mentioned above. Finally, we will let the
endomorphisms of $\V$ over $F$ act on the right on $\V$ by
$$\xb.\n_i=\n_i(\xb)$$
Consider the first vector of the basis $v_1$ and recall that $n\geq
4$. The relation $v_1.\,\n_{n-1}=-\unsurr\, v_1$ implies that in
$v_1$ there are no terms in $w_{s,t}$ for integers
$s,t\in\{1,\,\dots,\,n-2\}$ such that $s<t$. Hence we may write:
\begin{equation}
v_1=\sum_{j=1}^{n-2}\m_{j,n-1}w_{j,n-1}+\sum_{j=1}^{n-1}\m_{j,n}w_{j,n}
\end{equation}
We now use the relation $v1.\,\n_3=-\unsurr\,v_1$ to get rid of more
terms in this sum. Indeed, this relation implies that there are no
terms in $w_{j,k}$ for $j\geq 5$ in $v_1$. Furthermore, the
relations: $$\left\lbrace\begin{array}{ccc}
w_{2,n-1}.\,\n_1 &=& w_{1,n-1}+mr^{n-4}\,w_{12}-m\,w_{2,n-1}\\
w_{2,n}.\,\n_1 &=& w_{1,n}+mr^{n-3}\,w_{12}-m\,w_{2,n}
\end{array}\right.$$
\noindent imply in turn that:
$$mr^{n-4}\m_{2,n-1}+mr^{n-3}\m_{2,n}=0$$
\emph{i.e}
\begin{equation}
\m_{2,n}=-\unsurr\,\m_{2,n-1}
\end{equation}
\noindent as there is no term in $w_{12}$ in $v_1$ by the first
point. Finally, an application of $(77)$ with $q=1$ and $\lambda=r$
also yields the relations between the coefficients:
\begin{eqnarray}
\m_{2,n-1} &=& r\,\m_{1,n-1}\\
\m_{2,n} &=& r\,\m_{1,n}
\end{eqnarray}

\noindent We gather all these results to get, up to a multiplication
by a scalar:
\begin{multline}
v_1=w_{1,n-1}+r\,w_{2,n-1}-\unsurr\,w_{1,n}-w_{2,n}\\
+\m_{3,n-1}\,w_{3,n-1}+\m_{4,n-1}\,w_{4,n-1}+\m_{3,n}\,w_{3,n}+\m_{4,n}\,w_{4,n}
\end{multline}
or
\begin{equation}
v_1=\m_{3,n-1}\,w_{3,n-1}+\m_{4,n-1}\,w_{4,n-1}+\m_{3,n}\,w_{3,n}+\m_{4,n}\,w_{4,n}
\end{equation}

\begin{enumerate}
\item\textbf{Suppose first that} $\mathbf{n\geq 5}$ \\\\
\noindent Then $v1.\,\n_3=-\unsurr\,v_1$ implies that $\m_{2,n}=0$.
Then the expression for $v_1$ is given by $(84)$. We will deal with
the case $n=5$ separately. Hence, assume first $n>5$.
\begin{enumerate}
\item\underline{
$n>5$}\\\\
We keep looking at the action of $\n_3$. By $(77)$ applied with
$\lambda=-\unsurr$ and $q=3$, we have:
\begin{eqnarray}
\m_{4,n} &=& -\unsurr\;\m_{3,n}\\
\m_{4,n-1} &=& -\unsurr\;\m_{3,n-1}
\end{eqnarray}

\noindent Furthermore, by the set of equalities:

$$\left\lbrace\begin{array}{ccc}
w_{4,n-1}.\,\n_3 &=& w_{3,n-1}+mr^{n-6}\,w_{34}-m\,w_{4,n-1}\\
w_{4,n}.\,\n_3 &=& w_{3,n} + mr^{n-5}\,w_{34}-m\,w_{4,n}
\end{array}\right.$$

we get by the same computation as in $(80)$ that:
\begin{equation}\m_{4,n}=\;-\unsurr\,\m_{4,n-1}\end{equation}

as there is no term in $w_{34}$ in $v_1$ (recall that $n>5$). Next,
we have $w_{3,n}.\,\n_4=r\,w_{3,n}$ and an action of $\n_4$ on the
other terms of $v_1$ never makes any term in $w_{3,n}$ appear.
Hence, $v1.\,\n_4=-\unsurr\,v_1$ forces $\m_{3,n}=0$. Then by
$(85)$, $(87)$ and $(86)$, all the coefficients of $v_1$ are zero.
In other words $v_1=0$, which is impossible as $v_1$ is a basis
vector.
\item \underline{The case $n=5$}\\\\
In this case, there are only three terms in $v_1$:
$$v_1=\m_{34}\,w_{34}+\m_{35}w_{35}+\m_{45}w_{45}$$
Let's act by $\n_3$. We have:
\begin{eqnarray}
w_{34}.\,\n_3&=&\unsur{l}\,w_{34}\\
w_{45}.\,\n_3&=&w_{35}+m\,w_{34}-m\,w_{45}
\end{eqnarray}
from which we derive by looking at the term in $w_{34}$:
$$\Big(\unsur{l}+\unsurr\Big)\,\m_{34}=-m\,\m_{45}$$
Moreover, by an equality of type $(77)$ with $q=3$ and
$\lambda=-\unsurr$, we get:
$$\m_{45}=-\unsurr\,\m_{35}$$
At this point, we need to distinguish between two cases:
\begin{enumerate}
\item If $l=-r$, we get by the two previous equalities:
$\m_{45}=\m_{35}=0$. Then, $v_1=\m_{34}\,w_{34}$ and $\m_{34}\neq
0$. Then we must have
$v_1.\,\n_4=\m_{34}\,w_{35}=-\unsurr\,\m_{34}\,w_{34}$ \emph{i.e}
$w_{35}=-\unsurr\,w_{34}$, which is a contradiction.
\item If $l\neq -r$, we get
$\m_{34}=-\frac{m}{\unsurr+\unsur{l}}\;\m_{45}$, so that $v_1$ is
proportional to
$$ w_{35}-\unsurr\,w_{45}+\frac{ml}{l+r}\,w_{34}$$
We now check that this expression for $v_1$ is compatible with
$v_1.\,\n_4=-\unsurr\,v_1$ for a certain value of $l$. We have:
$$\left\lbrace\begin{array}{ccc}
w_{35}.\,\n_4 &=& w_{34}+\frac{m}{l}\,w_{45}-m\,w_{35}\\
w_{45}.\,\n_4 &=& \!\!\!\unsur{l}\,w_{45}\\
w_{34}.\,\n_4 &=& \qquad\qquad\qquad\qquad w_{35}
\end{array}\right.$$
By looking at the coefficient of the term in $w_{45}$ in
$v_1.\,\n_4=-\unsurr\,v_1$, we must have:\\
$\unsur{r^2}=\frac{m}{l}-\unsur{lr}$ \emph{i.e} $l=-r^3$. We note
that $-r^3\neq -r$ as we have assumed that the parameter $m$ is
nonzero. Looking at the coefficients of the terms in $w_{34}$ and
$w_{35}$ in $v_1.\,\n_4=-\unsurr\,v_1$ yields in turn:
$$ -r=\frac{ml}{l+r}$$ With $l=-r^3$, this equality holds.
We replace $l$ by $-r^3$ in the expression giving $v_1$ to get:
\begin{equation}
v_1=w_{35}-\unsurr\,w_{45}-r\,w_{34}
\end{equation}
Next, we must have $v_1.\,\n_2=-\unsurr\,(v_1+v_2)$ and we have:
\begin{eqnarray}
w_{35}.\,\n_2 &=& w_{25}+mr\,w_{23}-m\,w_{35}\\
w_{45}.\,\n_2 &=& r\,w_{45}\\
w_{34}.\,\n_2 &=& w_{24}+m\,w_{23}-m\,w_{34}
\end{eqnarray}

We need to investigate the coefficients for $v_2$. The equality
$v_2.\,\n_4=-\unsurr\,v_2$ forces the coefficients of the $w_{ij}$'s
to be zero for $1\leq i<j\leq 3$. Hence we may write:
\begin{multline} v_2=\la_{14}\,w_{14}+\la_{24}\,w_{24}+\la_{34}\,w_{34}+\\
\la_{15}\,w_{15}+\la_{25}\,w_{25}+\la_{35}\,w_{35}+\la_{45}\,w_{45}\end{multline}

Moreover, since $v_2.\,\n_2=r\,v_2$, we have by $(77)$ applied with
$\la=r$ and $q=2$:
\begin{eqnarray*}
\la_{34}&=&r\,\la_{24}\\
\la_{35}&=&r\,\la_{25}
\end{eqnarray*}
and with the equalities $(91)$ and $(93)$
%\begin{eqnarray*}
%w_{34}.\,\n_2 &=& w_{24}+m\,w_{23}-m\,w_{34}\\
%w_{35}.\,\n_2 &=& w_{25}+mr\,w_{23}-m\,w_{35}
%\end{eqnarray*}
and the fact that there are no terms in $w_{23}$ in $v_1$ and $v_2$,
we also deduce that:
$$\la_{35}=-\unsurr\,\la_{34}$$

Furthermore, since there must be a term in $w_{25}$ in $v_2$, none
of these coefficients is zero. Finally, by the equalities $(91)$,
$(92)$, $(93)$, the expression for $v_1$ in $(90)$ and the fact that
$v_1.\,\n_2=-\unsurr (v_1+v_2)$, we have: $\la_{14}=\la_{15}=0$.
Thus, $v_2$ must be proportional to the vector
$$w_{24}+r\,w_{34}-w_{35}-\unsurr\,w_{25}+\la_{45}\,w_{45}$$
Next, if $v_1$ is given by $(90)$, by looking at the coefficient of
$w_{25}$ in $v_1.\,\n_2=-\unsurr\,(v_1+v_2)$, the vector $v_2$ must
be given by:
$$v_2=r^2\,w_{24}+r^3\,w_{34}-r^2w_{35}-r\,w_{25}+\la_{45}\,w_{45}$$
%Moreover, by looking in turn at the coefficient of $w_{45}$, we get:
%$$-1=\unsur{r^2}-\frac{\la_{45}}{r}\qquad\emph{i.e}\qquad\la_{45}=\unsurr +r$$
%We now get a complete expression for the vector $v_2$:
%$$v_2=r^2\,w_{24}+r^3\,w_{34}-r^2w_{35}-r\,w_{25}+\Big(\unsurr+r\Big)\,w_{45}$$
From there, the contradiction comes from the relation
$v_2.\,\n_1=-\unsurr\,v_2-r\,v_1$, as for instance the coefficient
of $w_{34}$ in $v_2.\,\n_1$ is $r^4$ while the coefficient of
$w_{34}$ in $-\unsurr\,v_2-r\,v_1$ is zero. We conclude that when
$n=5$ it is also impossible to have an irreducible $4$- dimensional
invariant subspace of $\V$ given by the matrices $N_i$'s. It remains
to study the case $n=4$.
\end{enumerate}
\end{enumerate}
\item \textbf{The case} $\mathbf{n=4}$.\\\\
By the expressions given in $(83)$ and $(84)$, we know that $v_1$ is
proportional to
$w_{13}+r\,w_{23}-\unsurr\,w_{14}-w_{24}+\m_{34}w_{34}$ or to
$w_{34}$. Assume first that $v_1$ is a multiple of $w_{34}$. The
relation $v_1.\,\n_3=-\unsurr\,v_1$ forces $l=-r$. Suppose $l=-r$
and without loss of generality $v_1=w_{34}$. We must have
$v_1.\,\n_2=-\unsurr\,(v_1+v_2)$, \emph{i.e}
\begin{equation} w_{24}+m\,w_{23}-m\,w_{34}=-\unsurr\,w_{34}-\unsurr\,v_2\end{equation} A
general form for $v_2$ is:
\begin{equation}v_2=\la_{12}\,w_{12}+\la_{13}\,w_{13}+\la_{14}\,w_{14}
+\la_{23}\,w_{23}+\la_{24}\,w_{24} +\la_{34}\,w_{34}\end{equation}
By $(95)$, we have $\la_{12}=\la_{13}=\la_{14}=0$. Next, by $(77)$
applied with $q=2$ and $\la=r$, we get $\la_{34}=r\,\la_{24}$. Also,
from
\begin{eqnarray*}
w_{34}.\,\n_2&=&w_{24}+m\,w_{23}-m\,w_{34}\\
w_{23}.\,\n_2&=&-\unsurr\,w_{23}
\end{eqnarray*}
we deduce that $$\la_{23}=\frac{m}{r+\unsurr}\,\la_{34}$$ Thus, the
last three coefficients in $v_2$ are nonzero and $v_2$ is
proportional to the vector
$$w_{24}+r\,w_{34}+\frac{mr}{r+\unsurr}\,w_{23},$$
and in fact $(95)$ forces
$$v_2=-r\Bigg(w_{24}+r\,w_{34}+\frac{mr}{\unsurr+r}\,w_{23}\Bigg)$$
to match the terms in $w_{24}$ and $w_{34}$. Now the contradiction
comes from the term in $w_{23}$: we must have:
$$m=\frac{mr}{r+\unsurr}\qquad \emph{i.e} \qquad \unsurr =0,$$ which is
impossible.\\
Thus, we must have
$v_1=w_{13}+r\,w_{23}-\unsurr\,w_{14}-w_{24}+\m_{34}w_{34}$.\\
% We
%have:
%\begin{eqnarray}
%w_{34}.\,\n_3 &=& \unsur{l}\,w_{34}\\
%w_{14}.\,\n_3 &=& w_{13}+\frac{m}{lr}\,w_{34}-m\,w_{14}\\
%w_{24}.\,\n_3 &=& w_{23}+\frac{m}{l}\,w_{34}-m\,w_{24}
%\end{eqnarray}
%The relation $v_1.\,\n_3=-\unsurr\,v_1$ forces by looking at the
%coefficient of the term in $w_{34}$:
%$$\Big(\unsur{l}+\unsurr\Big)\,\m_{34}=\frac{m}{l}\Big(1+\unsur{r^2}\Big)$$
%If $l=-r$, we get $\m_{34}=0$, so that
%$v_1=w_{13}+r\,w_{23}-\unsurr\,w_{14}-w_{24}$. Then by $(98)$,
%$(99)$ and the expression for $v_1$, we must have
%$-\frac{m}{lr^2}-\frac{m}{l}=0$, \emph{i.e} $r^2=1$, which is
%impossible as $m\neq0$.\\
%Hence $l\neq -r$ and $$\m_{34}=\frac{m\Big(\unsurr+r\Big)}{l+r}$$ so
%that
%$$v_1=w_{13}+r\,w_{23}-\unsurr\,w_{14}-w_{24}+\frac{m\Big(\unsurr+r\Big)}{l+r}\,w_{34}$$
% We use the relation $v_1.\,\n_2=-\unsurr\,(v_1+v_2)$ to get a
%contradiction. But first we need to study the vector $v_2$.
We have the set of equations:
\begin{eqnarray}
w_{13}.\,\n_2&=&w_{12}+\frac{m}{l}\,w_{23}-m\,w_{13}\\
w_{23}.\,\n_2&=&\unsur{l}\,w_{23}\\
w_{14}.\,\n_2&=&r\,w_{14}\\
w_{24}.\,\n_2&=&w_{34}\\
w_{34}.\,\n_2&=&w_{24}+m\,w_{23}-m\,w_{34}
\end{eqnarray}
\noin By looking at the coefficient of $w_{12}$ in the relation
$v1.\,\n_2=-\unsurr\,v_1-\unsurr\,v_2$, we get by using the equation
$(97)$ and the expression for $v_1$ above:
$$1=-\frac{\la_{12}}{r}\qquad\emph{i.e.\/}\qquad \la_{12}=-r$$
Consequently, $$\la_{13}=r\la_{12}=-r^2$$ by the relation
$v_2.\,\n_2=r\,v_2$ and $(78)$. By looking at the coefficient of
$w_{14}$ in the relation $v1.\,\n_2=-\unsurr\,v_1-\unsurr\,v_2$ and
by using $(99)$ and the expression for $v_1$, we get:
$$-1=-\frac{\la_{14}}{r}+\unsur{r^2}\qquad\emph{i.e.\/}\qquad
\la_{14}=r+\unsurr$$ Similarly, by studying the coefficient of
$w_{24}$ in the same relation, we get:
$$\mu_{34}=\unsurr-\frac{\la_{24}}{r}$$
Next, we use the relation $v_2.\,\n_1=-\unsurr\,v_2-r\,v_1$ together
with the set of equations:
\begin{eqnarray*}
w_{12}.\,\n_1&=&\unsur{l}\,w_{12}\\
w_{23}.\,\n_1&=&w_{13}+m\,w_{12}-m\,w_{23}\\
w_{13}.\,\n_1&=&w_{23}\\
w_{34}.\,\n_1&=&r\,w_{34}\\
w_{24}.\,\n_1&=&w_{14}+mr\,w_{12}-m\,w_{24}\\
w_{14}.\,\n_1&=&w_{24}
\end{eqnarray*}
By looking at the coefficient in $w_{13}$, we get with
$\la_{13}=-r^2$ (in $v_2$) and $\mu_{13}=1$ (in $v_1$):
$$\la_{23}=r-r=0$$
Also, by looking at the coefficient in $w_{14}$, we get
$$\la_{24}=1-\frac{\la_{14}}{r}$$
And since $\la_{14}=r+\unsurr$ from above, we actually get:
$$\la_{24}=-\unsur{r^2}$$ Then we also have
$$\la_{34}=r\,\la_{24}=-\unsurr$$ by $v_2.\,\n_2=r\,v_2$ and $(77)$.
From the expression for $\la_{24}$, it also follows that
$$\mu_{34}=\unsurr+\unsur{r^3}$$ To finish, we look at the
coefficient of $w_{12}$ in $v_2.\,\n_1=-\unsurr\,v_2-r\,v_1$. By
using the second set of equations with the adequate coefficients
$\la_{12}$, $\la_{23}$ and $\la_{24}$ for $v_2$, we immediately get:
$$\frac{-r}{l}-\frac{m}{r}=1$$
$$\emph{i.e\/}\qquad\boxed{l=-r^3}$$
We conclude that $n=4$ is the only case for which it is possible to
have an irreducible representation of degree $n-1$ inside
$\mathcal{V}$ that is equivalent to the matrix representation
defined by the matrices $N_i$'s. Moreover, we showed in that case
that $l$ must take the value $-r^3$. Furthermore, we have seen along
the proof that such a $3$-dimensional invariant subspace must be
spanned over $F$ by the vectors:
\begin{eqnarray*}
v_1&=& w_{13}+r\,w_{23}-\unsurr\,w_{14}-w_{24}+\big(\unsurr+\unsur{r^3}\big)w_{34}\\
v_2&=&
-r\,w_{12}-r^2\,w_{13}-\unsurr\,w_{34}-\unsur{r^2}\,w_{24}+\big(r+\unsurr\big)w_{14}
\end{eqnarray*}
\noin and a third linearly independent vector $v_3$ such that:
$$\left|\begin{array}{ccc}
v_3.\,\n_3\eg r\,v_3\\
v_3.\,\n_2\eg -\unsurr\,v_3-r\,v_2\\
v_2.\,\n_3\eg -\unsurr(v_2+v_3)
\end{array}\right.$$
Only two of the above relations will be useful to force the value
for the coefficients of $v_3$ in the basis of $\mathcal{V}$. More
precisely, we use the action by $g_3$ as in the first and third
relations, together with the following defining equations (where we
took care to replace $l$ by its specialization):

\begin{eqnarray*}
w_{12}.\,\n_3\eg r\,w_{12}\\
w_{23}.\,\n_3\eg w_{24}\\
w_{13}.\,\n_3\eg w_{14}\\
w_{34}.\,\n_3\eg-\unsur{r^3}\,w_{34}\\
w_{24}.\,\n_3\eg w_{23}-\frac{m}{r^3}\,w_{34}-m\,w_{24}\\
w_{14}.\,\n_3\eg w_{13}-\frac{m}{r^4}\,w_{34}-m\,w_{14}
\end{eqnarray*}

Let's write a general form for $v_3$:
$$v_3:=\g_{12}\,w_{12}+\g_{23}\,w_{23}+\g_{13}\,w_{13}+\g_{34}\,w_{34}+\g_{24}\,w_{24}+\g_{14}\,w_{14}$$

By looking at the coefficient of $w_{23}$ in the third relation, we
get\\ $\la_{24}=-\unsurr\,\g_{23}$, from which we derive
$\g_{23}=\unsurr$ by replacing $\la_{24}$ by its value. Similarly by
looking at the coefficient of $w_{14}$, still using the third
relation, we see that: $$\la_{14}=r-\unsurr\,\g_{13}\qquad
\emph{i.e\/}\qquad r+\unsurr =r-\unsurr\,\g_{13}\qquad\emph{i.e\/}
\qquad\g_{13}=-1$$ Since by the first relation and $(78)$ we have
$$\left\lb\begin{array}{ccc}
\g_{14}&=&r\,\g_{13}\\
\g_{24}&=&r\,\g_{23}
\end{array}\right. ,$$
we derive that $\g_{14}=-r$ on one hand and $\g_{24}=1$ on the other
hand. It remains to find the coefficients $\g_{12}$ and $\g_{34}$.
Still by the last relation and by looking this time at the
coefficient in $w_{12}$, we get:
$$r\la_{12}=1-\frac{\g_{12}}{r}$$
Replacing $\la_{12}$ by the value $-r$ yields: $$\g_{12}=r+r^3$$
Finally, using the third relation for the last time and looking at
the coefficient in $w_{34}$ yields:
$$\big(-\unsur{r^3}\big)\big(-\unsurr\big)-\frac{m}{r^3}\big(-\unsur{r^2}\big)-\frac{m}{r^4}\big(r+\unsurr\big)=\unsur{r^2}-\unsurr\,\g_{34},$$
which leads to $\g_{34}=0$.

When gathering all the coefficients for $v_3$, we obtain:
$$v_3=(r+r^3)\,w_{12}+\unsurr\,w_{23}-w_{13}+w_{24}-r\,w_{14}$$
Thus, if there exists an irreducible $3$-dimensional invariant
subspace of $\mathcal{V}$ whose matrix representation is equivalent
to the one defined by the matrices $N_i$'s, then it is spanned by
the vectors $v_1,v_2,v_3$ as we defined them above and $l$ must take
the value $-r^3$.

Conversely, we show that for the value $-r^3$ of $l$, the vectors
$$\left\lb\begin{array}{ccc}
v_1\eg
r\,w_{23}+w_{13}+\big(\unsurr+\unsur{r^3}\big)\,w_{34}-w_{24}-\unsurr\,w_{14}\\
v_2\eg
-r\,w_{12}-r^2\,w_{13}-\unsurr\,w_{34}-\unsur{r^2}\,w_{24}+\big(r+\unsurr\big)\,w_{14}\\
v_3\eg (r+r^3)\,w_{12}+\unsurr\,w_{23}-w_{13}+w_{24}-r\,w_{14}
\end{array}\right.$$
form a free family of vectors that satisfy to the relations
$(\bigtriangledown)$. This will prove that
$\text{Span}_F(v_1,v_2,v_3)$ is an irreducible $3$-dimensional
invariant subspace of $\mathcal{V}$. For the freedom of the family
of vectors, we note that in $v_1$ (resp $v_2$, resp $v_3$) there is
no term in $w_{12}$ (resp $w_{23}$, resp $w_{34}$). Then, if
$\la_1$, $\la_2$ and $\la_3$ are scalars such that:
$$\la_1v_1+\la_2v_2+\la_3v_3=0,$$
they must be related by: $$\left\lb\begin{array}{ccc}
\la_2&=&(1+r^2)\,\la_3\\
\la_2&=&(1+\unsur{r^2})\,\la_1\\
\la_3&=& \nts\nts\nts\nts\nts\nts -r^2\,\la_1
\end{array}\right.,$$
\noin where we used the freedom of the $w_{ij}$'s. These equations
imply that: $$(1+r^2+r^4+r^6)\la_1=0$$ Since $r^2\neq 1$ and
$(r^2)^4\neq 1$, we get $\la_1=0$. It follows that
$\la_1=\la_2=\la_3=0$. Hence the family $(v_1,v_2,v_3)$ is free. It
remains to show that the relations $(\bigtriangledown)$ are
satisfied on the vectors $v_1$, $v_2$, $v_3$. This is an easy
verification that is left to the reader.

\end{enumerate}
Let's go back to the proof of the theorems and suppose that there is
an irreducible $(n-1)$-dimensional invariant subspace $\U$ of $\V$.
When $n=4$, assume that $l\neq -r^3$. By the preceding work, there
exists a basis $(v_1,\,\dots,\,v_{n-1})$ of $\U$ in which the
matrices of the left actions of the $g_i$'s are the matrices $M_i$'s
and we have the relations, when the indices make sense:
$$(\tr)\left|\begin{array}{cccc} \n_t(v_i)&=&r\;v_i,&\text{for}\;
t\not\in\{
i-1,\,i,\,i+1\}\\
\n_i(v_i)&=&-\unsurr\;v_i&\\
\n_{i+1}(v_i)&=&r(v_i+v_{i+1})&\\
\n_{i-1}(v_i)&=&r\;v_i+\unsurr\;v_{i-1}&
\end{array}\right.$$
We will show that these relations force the values $\unsur{r^{n-3}}$
or $-\unsur{r^{n-3}}$ for $l$.\\ First of all, the relation
$v_i.\,\n_i=-\unsurr\,v_i$ implies that in $v_i$ there is no term in
$w_{s,t}$ for
$$\begin{array}{ccc}
& t\leq i-1 &(a)\\
or & s\geq i+2 & (b)\\
or & \left\lbrace\begin{array}{ll} s\leq i-1 \\ t\geq i+2
&\end{array}\right. & \begin{array}{ll} (c)\\(d) \end{array}
\end{array}$$

\noindent since those terms are all multiplied by $r$ when acting by
$\n_i$ and they cannot be obtained from other $w_{k,q}$'s with an
action by $\n_i$ by points $(i)$, $(ii)$ and $(iii)$ above. By
$(b)$, either $s\leq i-1$ or $s=i$ or $s=i+1$. In the first case, by
$(d)$, the integer $t$ must be less or equal to $i+1$. By $(a)$ the
only possibilities for $t$ are $t=i$ or $t=i+1$. When $s=i$, we may
have $t\geq i+1$ and when $s=i+1$, we may have $t\geq i+2$. Thus,
the possible values for $s$ and $t$ for $w_{s,t}$ in $v_i$ are:
$$\begin{array}{lllll}
s\leq i-1\;\text{and}\;t\in\lbrace i,\,i+1\rbrace\\
\qquad\qquad\qquad \text{or}\\
s=i\;\text{and no restriction on}\; t\; (t\geq i+1)\\
\qquad\qquad\qquad \text{or}\\
s=i+1\;\text{and no restriction on}\; t\;(t\geq i+2)
\end{array}$$

\noindent Thus, a general form for $v_i$ must be:

\begin{equation}
\begin{split}
v_i=\m_{i,i+1}\,w_{i,i+1} & +\sum_{k=i+2}^n\m_{i,k}\,w_{i,k}+\sum_{k=i+2}^n\m_{i+1,k}\,w_{i+1,k}\\
&\qquad\qquad+\sum_{s=1}^{i-1}\m_{s,i}\,w_{s,i}+\sum_{s=1}^{i-1}\m_{s,i+1}\,w_{s,i+1}
\end{split}
\end{equation}
\noindent Next, from $v_i.\n_i=-\unsurr\,v_i$, we have by $(78)$
applied with $q=i$ and $\la=-\unsurr$ on one hand and $(77)$ applied
with $q=i$ and $\la=-\unsurr$ on the other hand:
%$$\forall s\leq i-1,\;\;\left\lbrace\begin{array}{ll}
%\;\; \,\,w_{s,i}.\,\n_i  =  w_{s,i+1}\\
%w_{s,i+1}.\,\n_i  =  w_{s,i}-m\,w_{s,i+1}\;\;\text{modulo}\;
%F\,\xali
%\end{array}\right.$$
%
%\noindent and we also have:
%$$\forall k\geq i+2,\;\;\left\lbrace\begin{array}{ll}
%\;\; \,\,w_{i,k}.\,\n_i  =  w_{i+1,k}\\
%w_{i+1,k}.\,\n_i  =  w_{i,k}-m\,w_{i+1,k}\;\;\text{modulo}\;
%F\,\xali
%\end{array}\right.$$
%
%\noindent Moreover, by arguments already discussed above, these
%equations are respectively the only way to obtain terms in
%$w_{s,i}$, $w_{s,i+1}$ (resp $w_{i,k}$, $w_{i+1,k}$). Hence the
%relation $\vi.\,\nui=-\unsurr\,\vi$ implies that:
%$$\left\lbrace\begin{array}{ll}
%\m_{s,i}-m\,\m_{s,i+1}=-\unsurr\,\m_{s,i+1},\;\;\;\forall s\leq i-1\\
%\m_{i,k}-m\,\m_{i+1,k}=-\unsurr\,\m_{i+1,k},\;\;\;\forall k\geq i+2
%\end{array}\right.$$
%
%\noindent or equivalently:
$$\left\lbrace\begin{array}{ll}
\m_{s,i+1}=-\unsurr\,\m_{s,i},\;\;\;\forall s\leq i-1\\
\m_{i+1,k}=-\unsurr\,\m_{i,k},\;\;\;\forall k\geq i+2
\end{array}\right.$$

\noindent Further, to get more relations between the coefficients,
we use the relations
$$v_i.\,\n_q=r\,v_i\;\;\text{for}\;q\not\in\lbrace
i-1,i,i+1\rbrace$$ \noindent By $(77)$ applied with each
$q\not\in\lbrace i-1,\,i,\,i+1\rbrace$ and $\la=r$, we get:
\begin{equation}\forall k\geq q+2, \;\;\m_{q+1,k}=r\,\m_{q,k}\end{equation}
\noindent and by $(78)$ also applied with each $q\not\in\lbrace
i-1,\,i,\,i+1\rbrace$ and $\la=r$, we get: \begin{equation}\forall
s\leq q-1, \;\;\m_{s,q+1}=r\,\m_{s,q}\end{equation} \noindent Let's
apply $(103)$ with $q=s\leq i-2$ and $k\in\lb i,\,i+1\rb$ ($i\geq
s+2$) on one hand and $(104)$ with $q=k\geq i+2$ and $s\in\lb
i,\,i+1\rb$ ($i+1\leq k-1$) on the other hand (where we used the
same notations as in the sums of $(102)$) to get:

$$\!\!\!\forall\; 1\leq s\leq i-2,\,\left\lb\begin{array}{ccc}
\!\!\!\!\!\!\m_{s+1,i}&=&\!\!\!\!\!r\,\m_{s,i}\\
\m_{s+1,i+1}&=&r\,\m_{s,i+1}
\end{array}\right.\qquad\text{on one hand}$$

$$\forall\; i+2\leq k\leq n-1,\,\left\lb\begin{array}{ccc}
\!\!\!\!\!\m_{i,k+1}&=&\!\!\!\!\!r\,\m_{i,k}\\
\m_{i+1,k+1}&=&r\,\m_{i+1,k}
\end{array}\right.\qquad\text{on the other hand}$$

\noindent The formula $(102)$ can now be rewritten as:

\begin{equation*}
\begin{split}
v_i=\m^i\,w_{i,i+1} &+ \sum_{k=i+2}^n
\del^i\,r^{k-i-2}(w_{i,k}-\unsurr\,w_{i+1,k}) \\ &+
\sum_{s=1}^{i-1}\,\la^i\,r^{s-1}(w_{s,i}-\unsurr\,w_{s,i+1})
\end{split}
\end{equation*}

\noindent We will show that the $\del^i\;$'s for $i=1,\,\dots,\,n-2$
can be set to the value one. First, we notice that if
$v_1,\,\dots,\,v_{n-1}$ are vectors satisfying to the relations
$(\tr)$, then the vectors $\del\,v_1,\,\dots,\,\del\,v_{n-1}$ also
satisfy to the relations $(\tr)$ for any nonzero scalar $\del$.
Thus, without loss of generality, we may write:
$$v_1=\mu^1\,w_{12}+\sum_{k=3}^n
r^{k-3}(w_{1,k}-\unsurr\,w_{2,k})$$ \noindent where we set
$\del^1=1$. Next, we notice on the expression above that an action
of $\n_2$ on $v_1$ never makes a term in $w_{24}$ appear. From
there, it suffices to look at the coefficient of $w_{24}$ in the
relation
$$v_1.\n_2=r\,v_1+r\,v_2$$ to get:
$$0=-r+r\,\del^2$$
\emph{i.e.\/}
$$\del^2=1$$
\noin Let's proceed by induction and suppose that $\del^i=1$ for a
given $i$ with $2\leq i\leq n-3$. We notice that $\del^{i+1}$ is the
coefficient of $w_{i+1,i+3}$ in $v_{i+1}$. Since an action of
$\n_{i+1}$ on $v_i$ never makes a term in $w_{i+1,i+3}$ appear, by
looking at the coefficient of $w_{i+1,i+3}$ in the relation
$\n_{i+1}\,v_i=r\,v_i+r\,v_{i+1}$, we get:
$$0=-r\,\del^i+r\,\del^{i+1}$$
\emph{i.e.\/}
$$\del^{i+1}=1$$
\noin Thus, we have shown that all the $\del^i\;$'s, for $1\leq
i\leq n-2$, can be set to the value one. Hence the formula giving
the $v_i$'s can be rewritten as follows:
\begin{equation}
\begin{split}
v_i=\m^i\,w_{i,i+1} &+ \sum_{k=i+2}^n
\,r^{k-i-2}(w_{i,k}-\unsurr\,w_{i+1,k}) \\ &+
\sum_{s=1}^{i-1}\,\la^i\,r^{s-1}(w_{s,i}-\unsurr\,w_{s,i+1})
\end{split}
\end{equation}
\noin where the only two unknown coefficients are $\la^i$ and
$\mu^i$. It remains to find the coefficients $\la^i$ and $\m^i$. By
looking at the coefficient of the term $w_{i,i+1}$ in the relation
$v_i.\,\n_{i+1}=r(v_i+v_{i+1})$, we get the equations:

\begin{equation}
\forall\; 1\leq i\leq n-2,\;\;\boxed{r\,\m^i+r^i\,\la^{i+1}=1}
\end{equation}

\noindent Next, by looking at the coefficient of the same term
$w_{i,i+1}$ in the relation $v_i.\,\n_{i-1}=r\,v_i+\unsurr\,v_{i-1}$
for $i\geq 2$, we get the equations:
\begin{equation*}
\forall\; 2\leq i\leq
n-1,\;\;-m\,\m^i-\la^i\,r^{i-3}=r\,\m^i-\unsur{r^2}
\end{equation*}
which can be rewritten after multiplication by a factor $r^2$:
\begin{equation}
\forall\; 2\leq i\leq n-1,\;\;\boxed{r\,\m^i+r^{i-1}\,\la^i=1}
\end{equation}
\noindent Subtracting equalities $(106)$ and $(107)$ yields:
\begin{equation}
\forall\; 2\leq i\leq n-2,\,\la^{i+1}=\unsurr\,\la^i
\end{equation}
\noin Next, we write equality $(106)$ with $i=1$ and equality
$(107)$ with $i=2$:
\begin{eqnarray}
r\,\m^1+r\,\la^2&=&1\\
r\,\m^2+r\,\la^2&=&1 \notag
\end{eqnarray}
\noin from which we derive:
\begin{equation}
\m^1=\m^2
\end{equation}
Further, by $(107)$, we have:
\begin{equation*}
\forall\;2\leq i\leq n-2,\;\;r\,\m^{i+1}+r^i\,\la^{i+1}=1
\end{equation*}
And by using $(108)$ we get:
\begin{equation}
\forall\;2\leq i\leq n-2,\;\;r\,\m^{i+1}+r^{i-1}\,\la^i=1
\end{equation}
\noin Now by $(107)$ and $(111)$, it comes:
\begin{equation}
\forall\;2\leq i\leq n-2,\;\;\m^{i+1}=\m^i
\end{equation}
\noin Gathering $(110)$ and $(112)$, we get:
\begin{equation}
\m_1=\m_2=\!\!\!\!\!\!\!\!\hspace{0.2in}\dots=\m_{n-1}
\end{equation}
\noin %In other words, the sequence $(\m^i)_{1\leq i\leq n-1}$ is
%stationary.
Thus, the coefficient $\m^i$ of the term $w_{i,i+1}$ in the
expression giving $v_i$ does not depend on the integer $i$. We will
denote it by $\m$.\\ By $(108)$, all the $\la^i\;$'s are determined
by $\la^2$ in the following way:
\begin{equation}
\forall\;2\leq i\leq n-1,\;\;\la^i=\la^2\,\bigg(\unsurr\bigg)^{i-2}
%\forall\;2\leq i\leq n-1,\;\;\la^i=\frac{\la^2}{r^{i-2}}
\end{equation}
\noin By $(109)$, we have $\la^2=\unsurr-\m$. Thus, by determining
the coefficient $\m$, we will get a complete expression for all the
vectors $v_i$'s. Recall that $\m$ is the coefficient of $w_{i,i+1}$
in $v_i$. We look at the coefficient of $w_{i,i+1}$ in
$v_i.\,\nui=-\unsurr\,v_i$. We have:
\begin{eqnarray*}
\forall\;k\geq i+2,\;\Big[w_{i+1,k}.\,\nui\Big]_{w_{i,i+1}} &=& m\,r^{k-i-2}\\
\forall\;s\leq i-1,\;\Big[w_{s,i+1}.\,\nui \Big]_{w_{i,i+1}} &=&
\frac{m}{l\,r^{i-s-1}}
\end{eqnarray*}
\noin Hence we get the equation:
$$\frac{\m}{l}-\sum_{k=i+2}^n r^{k-i-3}\,m\,r^{k-i-2}
-\sum_{s=1}^{i-1}
\la^i\,r^{s-2}\,\frac{m}{l\,r^{i-s-1}}=-\frac{\m}{r}$$ \emph{i.e}
$$\frac{\m}{l}-\unsur{r^2}+(r^2)^{n-i-2}-\frac{\la^i}{l}\Big(\unsur{r^i}-r^{i-2}\Big)=-\frac{\m}{r}\qquad(\star)_i$$
\noin Let's write down $(\star)_2$ and $(\star)_3$:
\begin{eqnarray*}
\m\,\Big(\unsur{l}+\unsurr\Big)&=&\frac{\la^2}{l}\Big(\unsur{r^2}-1\Big)+\unsur{r^2}-(r^2)^{n-4}\qquad (\star)_2\\
\m\,\Big(\unsur{l}+\unsurr\Big)&=&\frac{\la^2}{lr}\Big(\unsur{r^3}-r\Big)+\unsur{r^2}-(r^2)^{n-5}\qquad
(\star)_3
\end{eqnarray*}
\noin where $\la^3$ has been replaced by $\frac{\la^2}{r}$. Let's
subtract these two equalities:
\begin{eqnarray*}
\!\!\!\!\!\!\!\!\frac{\la^2}{l}\Big(\unsur{r^2}-\unsur{r^4}\Big)\;=\;(r^2)^{n-4}\Big(1-\unsur{r^2}\Big)\qquad
(\star)_2-(\star)_3
\end{eqnarray*}
\noin We multiply this last equality by $\unsur{r^2}$ and then
divide it by $\unsur{r^2}-\unsur{r^4}$ (licit as $m\neq 0$) to get:
$$\la^2=l\,(r^2)^{n-3}$$
\noin We recall that $\m=\unsurr-\la^2$. Hence we get a relation
binding $\m$ and $l$ as follows:
$$\m=\unsurr-l\,(r^2)^{n-3}$$
Next, by looking at the coefficient of $w_{12}$ in
$v_1.\,\n_1=-\unsurr\,v_1$, we get:
\begin{eqnarray}
\m\,\Big(\unsur{l}+\unsurr\Big)&=&\unsur{r^2}-(r^2)^{n-3}
\end{eqnarray}
\noin as $v_1$ is only composed of one term and one sum as below:
$$v_1=\m\,w_{12}+\sum_{k=3}^n r^{k-3}(w_{1,k}-\unsurr\,w_{2,k})$$
Replacing $\m$ by its value in equality $(115)$ yields:
$l^2=\unsur{(r^2)^{n-3}}$. Hence there are two possible values for
$l$: \\
$$\begin{array}{cccccc} \text{Either} & l=\unsur{r^{n-3}} &
\text{and} & \m=\unsurr-r^{n-3} & \text{and} &
\la^2=r^{n-3}\\
\text{Or} & l=-\unsur{r^{n-3}} & \text{and} & \m=\unsurr+r^{n-3} &
\text{and} & \;\;\la^2=-r^{n-3}
\end{array}$$
\noin In both cases, we have $\m=\unsurr-\unsur{l}$. In the first
case, we get $\la^i=r^{n-i-1}$ and in the second case, we get
$\la^i=-r^{n-i-1}$, where we used the expression given in $(114)$.\\
At this point we have proven that if there exists an irreducible
$(n-1)$-dimensional invariant subspace $\U$ of $\V$, then $l$ must
take the values $\unsur{r^{n-3}}$ or $-\unsur{r^{n-3}}$. %Moreover,
%for each of these values of $l$, there must exist a basis
%$(v_1,v_2,\dots,v_{n-1})$ of $\U$ with
%$$v_i:=\Big(\unsur{r}\,-\,\unsur{l}\Big)w_{i,i+1}\,+\sum_{k=i+2}^{n}r^{k-i-2}(w_{i,k}\,-\,\unsurr\;w_{i+1,k})\,+\;\e_l\,\sum_{s=1}^{i-1}r^{n-i-2+s}(w_{s,i}\,-\,\unsurr\;w_{s,i+1})$$
%\begin{center} where \hspace{.1in}
%$\begin{cases}\e_{\unsur{r^{n-3}}}\;\;=\,1\\\e_{-\unsur{r^{n-3}}}=-1\end{cases}$
%\end{center}
\noin Conversely, given
$l\in\lb\unsur{r^{n-3}},-\unsur{r^{n-3}}\rb$, we show that the
vectors $v_i$'s, $1\leq i\leq n-1$, defined by:
\begin{equation}
v_i:=\Big(\unsur{r}\,-\,\unsur{l}\Big)w_{i,i+1}\,+\sum_{k=i+2}^{n}r^{k-i-2}(w_{i,k}\,-\,\unsurr\;w_{i+1,k})\,+\;\e_l\,\sum_{s=1}^{i-1}r^{n-i-2+s}(w_{s,i}\,-\,\unsurr\;w_{s,i+1})\end{equation}
\begin{center} where \hspace{.1in}
$\begin{cases}\e_{\unsur{r^{n-3}}}\;\;=\,1\\\e_{-\unsur{r^{n-3}}}=-1\end{cases}$
\end{center}
form a free family of vectors and satisfy to the relations (when the
indices make sense):
$$\left|\begin{array}{cccc} v_i.\,\n_t&=&r\;v_i,&\text{for}\;
t\not\in\{
i-1,\,i,\,i+1\}\\
v_i.\nui &=&-\unsurr\;v_i&\\
v_i.\,\n_{i+1} &=&r(v_i+v_{i+1})&\\
v_i.\,\n_{i-1}&=&r\;v_i+\unsurr\;v_{i-1}&
\end{array}\right.$$

Let's compute $v_i.\,\n_t$ for $t<i-1$. All the terms in $(116)$ are
multiplied by a factor $r$ when acting by $\n_t$, except possibly:
\begin{equation}
\begin{aligned}
\e_l\;r^{n-i-2+t}&(w_{t,i}-\unsurr\,w_{t,i+1})\\
\e_l\;r^{n-i-2+t+1}&(w_{t+1,i}-\unsurr\,w_{t+1,i+1})
\end{aligned}
\end{equation}
\noin By definition of $\n_t$, we have:
\begin{eqnarray}
w_{t,i}.\,\n_t &=& w_{t+1,i}\\
w_{t,i+1}.\,\n_t &=& w_{t+1,i+1}\\
w_{t+1,i}.\,\n_t &=& w_{t,i}+m\,r^{i-t-2}\,w_{t,t+1}-m\,w_{t+1,i}\\
w_{t+1,i+1}.\,\n_t &=&
w_{t,i+1}+m\,r^{i-t-1}\,w_{t,t+1}-m\,w_{t+1,i+1}
\end{eqnarray}

\noin When we act by $\n_t$ on $v_i$, all the terms $w_{s,t}$'s
composing $v_i$ are multiplied by $r$ except the four mentioned
above. Hence we read on the expressions $(118)$, $(119)$, $(120)$
and $(121)$ that $w_{t,i}$ is obtained from and only from
$w_{t+1,i}$. And similarly, $w_{t,i+1}$ is obtained from and only
from $w_{t+1,i+1}$. Now we read on $(117)$ and $(120)$ that the
coefficient of $w_{t,i}$ is multiplied by $r$ when acting by $\n_t$
on $v_i$. Similarly, by $(117)$ and $(121)$, the coefficient of
$w_{t,i+1}$ is multiplied by $r$ when acting by $\n_t$ on $v_i$.
Similarly, we read on $(117)$, $(118)$ and $(120)$ that after acting
by $\n_t$ on $v_i$, the coefficient of $w_{t+1,i}$ is
$\e_l\,r^{n-i-2+t}(1-mr)$, \emph{i.e}
$\e_l\,r^{n-i+t}=r\times\e_l\,r^{n-i-2+t+1}$. And we read on
$(117)$, $(119)$ and $(121)$ that the coefficient of $w_{t+1,i+1}$
is $-\frac{\e_l}{r}\,r^{n-i-2+t}(1-mr)$, \emph{i.e}
$-\frac{\e_l}{r}\,r^{n-i+t}=r\times\e_l\,r^{n-i-2+t+1}\Big(-\unsurr\Big)$.
We conclude that all the terms composing $v_i$ are multiplied by a
factor $r$ when we act by $\n_t$ on $v_i$. And we note that the
terms in $w_{t,t+1}$ appearing in $(120)$ and $(121)$ cancel each
other with the adequate coefficients $1$ and $-\unsurr$ of
$w_{t+1,i}$ and $w_{t+1,i+1}$. Thus, we have shown that:
$$\forall\,t\leq i-2,\,v_i.\,\n_t=r\,v_i$$
\noin Similarly, given $t\geq i+2$, all the terms in $(116)$ are
multiplied by a factor $r$ when acting by $\n_t$, except possibly:
\begin{equation}
\begin{aligned}
r^{t-i-2}&(w_{i,t}-\unsurr\,w_{i+1,t})\\
r^{t-i-1}&(w_{i,t+1}-\unsurr\,w_{i+1,t+1})
\end{aligned}
\end{equation}
\noin We compute the action of $\n_t$ on the four terms $w_{i,t}$,
$w_{i+1,t}$, $w_{i,t+1}$ and $w_{i+1,t+1}$:
\begin{eqnarray}
w_{i,t}.\,\n_t &=& w_{i,t+1}\\
w_{i+1,t}.\,\n_t &=& w_{i+1,t+1}\\
w_{i,t+1}.\,\n_t &=& w_{i,t}+\frac{m}{l\,r^{t-i-1}}\,w_{t,t+1}-m\,w_{i,t+1}\\
w_{i+1,t+1}.\,\n_t &=&
w_{i+1,t}+\frac{m}{l\,r^{t-i-2}}\,w_{t,t+1}-m\,w_{i+1,t+1}
\end{eqnarray}
\noin We read on the equalities $(123)-(126)$ that when acting by
$\n_t$, the term $w_{i,t}$ (resp $w_{i,t+1}$) can be obtained from
and only from the term $w_{i,t+1}$ (resp $w_{i+1,t+1}$). By $(122)$
and $(125)$, the coefficient of $w_{i,t}$ is multiplied by a factor
$r$ when acting by $\n_t$ on $v_i$ and by $(122)$ and $(126)$, the
one of $w_{i+1,t}$ is also multiplied by a factor $r$. Like above,
we read on $(122)$, $(123)$ and $(125)$ that the coefficient of
$w_{i,t+1}$ is multiplied by $r$ when acting by $\n_t$ on $v_i$ and
we read on $(122)$, $(124)$ and $(126)$ that the coefficient of
$w_{i+1,t+1}$ is also multiplied by a factor $r$. Thus, all the
terms $w_{s,t}$'s in $(116)$ are in fact multiplied by $r$ and we
have:
$$\forall\,t\geq i+2,\,v_i.\,\n_t=r\,v_i$$
\noin We now show that $v_i.\,\n_i=-\unsurr\,v_i$.

We have for all $k\geq i+2$ and all $s\leq i-1$:
$$\left\lb\begin{array}{ccccc}
w_{i,i+1}&.&\nui &=& \unsur{l}\,\;\;\;\,w_{i,i+1}\\
w_{i,k}&.&\nui &=& \qquad\qquad\qquad\qquad\qquad\qquad\!\! w_{i+1,k}\\
w_{i+1,k}&.&\nui &=& w_{i,k}+m\,r^{k-i-2}\,w_{i,i+1}-m\,w_{i+1,k}\\
w_{s,i}&.&\nui &=& \qquad\qquad\qquad\qquad\qquad\qquad\!\!\! w_{s,i+1}\\
w_{s,i+1}&.&\,\nui &=&
w_{s,i}+\;\frac{m}{l\,r^{i-s-1}}\,\;\,w_{i,i+1}-m\,w_{s,i+1}
\end{array}\right.$$
\noin We read on the equations above that the coefficient of
$w_{i,i+1}$ after an action by $\nui$ on $v_i$ is:
$$\Big(\unsurr-\unsur{l}\Big)\times\unsur{l}-\frac{m}{r}\,\sum_{k=i+2}^n
(r^{k-i-2})^2-\frac{m\e_l}{lr}.\,
r^{n+1}.\,\sum_{s=1}^{i-1}(r^{s-i-1})^2$$ \emph{i.e.}
$$\Big(\unsurr-\unsur{l}\Big)\times\unsur{l}-\unsur{r^2}\Big[1-r^{2n-2i-2}+\frac{\e_l}{l}(r^{n-2i+1}-r^{n-1})\Big]$$
\noin where we used the relation $\frac{m}{1-r^2}=\unsurr$. Next, we
note that:
$$\forall\,l\in\Big\lb\unsur{r^{n-3}},-\unsur{r^{n-3}}\Big\rb,\;\frac{\e_l}{l}=r^{n-3}$$
\noin Hence, the term in the bracket is nothing else but
$$1-r^{2n-4}=(1-r^{n-2})(1+r^{n-2})=\Big(1-r\,\frac{\e_l}{l}\Big)\Big(1+r\,\frac{\e_l}{l}\Big)$$
\noin Then, "distributing" the coefficient $\unsur{r^2}$ inside
these two factors yields the product
$\Big(\unsurr-\frac{\e_l}{l}\Big)\Big(\unsurr+\frac{\e_l}{l}\Big)$.
Since $\e_l\in\lb -1,1\rb$, this product is
$\Big(\unsurr-\unsur{l}\Big)\Big(\unsurr+\unsur{l}\Big)$. Thus, the
coefficient of $w_{i,i+1}$ after an action by $\nui$ on $v_i$ is:
$$\Big(\unsurr-\unsur{l}\Big)\unsur{l}-\Big(\unsurr+\unsur{l}\Big)=\boxed{-\unsurr\bigp\unsurr-\unsurl\bigpd}$$
\noin So the coefficient of $w_{i,i+1}$ after an action by $\nui$ on
$v_i$ is multiplied by a factor $-\unsurr$. Let's look at the other
terms. For each $i+2\leq k\leq n$, the coefficient of the term in
$w_{i+1,k}$ is given by $r^{k-i-2}(1+\frac{m}{r})=r^{k-i-4}$, hence
is multiplied by $-\unsurr$. Similarly, for each $1\leq s\leq i-1$,
the coefficient of $w_{s,i+1}$ is, after an action by $\nui$,
$\e_lr^{n-i-2+s}(1+\frac{m}{r})=\e_l\,r^{n-i+s-4}$, hence is
multiplied by $-\unsurr$. Moreover, it can be read directly on the
expression $(116)$ and the third (resp the fifth) equation above
that the coefficient of $w_{i,k}$ (resp $w_{s,i}$) gets multiplied
by a factor $-\unsurr$. Thus, we have:
$$v_i.\,\nui=-\unsurr\,v_i$$
\noin Let's now compute $v_i.\,\n_{i+1}$. When acting by $\n_{i+1}$,
most of the terms $w_{s,t}$'s in $(116)$ are multiplied by $r$, with
the exception of:
\begin{itemize}
\item $w_{i,i+1}$
\item $w_{i,i+2}$, $w_{i+1,i+2}$ and $w_{i+1,k}$ for $i+3\leq k\leq n$, in the
first sum
\item $w_{s,i+1}$ for $1\leq s\leq i-1$, in the second sum
\end{itemize}
\noin We compute the action of $\n_{i+1}$ on these terms:
\begin{eqnarray}
w_{i,i+1}.\,\n_{i+1} &=& w_{i,i+2}\\
&& \notag \\
w_{i,i+2}.\,\n_{i+1} &=&
w_{i,i+1}+\frac{m}{l}\,w_{i+1,i+2}-m\,w_{i,i+2}\\
w_{i+1,i+2}.\,\n_{i+1} &=& \unsurl\,w_{i+1,i+2}\\
w_{i+1,k}.\,\n_{i+1} &=& w_{i+2,k} \qquad (i+3\leq k\leq n)\\
&& \notag \\
w_{s,i+1}.\,\n_{i+1} &=& w_{s,i+2} \qquad\;\;\;\;\;\,\, (1\leq s\leq
i-1)
\end{eqnarray}
\noin From there,
$$\begin{array}{lll}\left.\begin{array}{lll}\text{the coefficient
of}\; w_{i,i+1}\;\text{in}\;
v_i.\,\n_{i+1}\;\text{is}\;1,\\
\text{the coefficient of}\;w_{i,i+1}\;\text{in}\;v_i \;\text{is}\;
\unsurr-\unsurl,\\
\text{the coefficient of}\; w_{i,i+1}\;\text{in}\;
v_{i+1}\;\text{is}\;\e_l\,r^{n-3}=\e_l\,\frac{\e_l}{l}=\unsurl
\end{array}\right\rb 1=r\bigp\unsurr-\unsurl\bigpd+\frac{r}{l}\\\\
\left.\begin{array}{lll}\text{the coefficient of}\;
w_{i,i+2}\;\text{in}\;
v_i.\,\n_{i+1}\;\text{is}\;r-\unsurl,\\
\text{the coefficient of}\;w_{i,i+2}\;\text{in}\;v_i \;\text{is}\;
1,\\
\text{the coefficient of}\; w_{i,i+2}\;\text{in}\;
v_{i+1}\;\text{is}\;-\unsur{lr}
\end{array}\right\rb r-\unsurl=r-r\times\unsur{lr}\\\\
\left.\begin{array}{lll}\text{the coefficient of}\;
w_{i+,i+2}\;\text{in}\;
v_i.\,\n_{i+1}\;\text{is}\;\frac{m}{l}-\unsur{rl}=-\frac{r}{l},\\
\text{the coefficient of}\;w_{i+1,i+2}\;\text{in}\;v_i \;\text{is}\;
-\unsurr,\\
\text{the coefficient of}\; w_{i+1,i+2}\;\text{in}\;
v_{i+1}\;\text{is}\;\unsurr-\unsurl
\end{array}\right\rb -\frac{r}{l}=-1+r\bigp\unsurr-\unsurl\bigpd
\end{array}$$
\noin From now on we do not need to worry about these terms anymore.
We study the action of $\n_{i+1}$ on the first sum of $v_i$ for
indices $k$ such that $i+3\leq k\leq n$. By the equation $(130)$, an
action by $\n_{i+1}$ on $v_i$ makes a term in $w_{i+2,k}$ appear
with coefficient the one of $w_{i+1,k}$ in $v_i$, that is
$-\frac{r^{k-i-2}}{r}$. Hence, we get a term:
$$r\sum_{k=i+3}^n(r^{k-i-3}\,w_{i+1,k}-\unsurr\,w_{i+2,k}),$$
\noin where we need to add
$$-\sum_{k=i+3}^n r^{k-i-2}\,w_{i+1,k}$$

\noin The first expression is the first sum in $v_{i+1}$, multiplied
by a factor $r$. The latter expression is a part of the first sum in
$r\,v_i$ with indices $k$ greater than $i+3$. Since all the terms
$w_{i,k}$'s for $i+3\leq k\leq n$ are multiplied by $r$ when acting
by $\n_{i+1}$, we get a sum:
$$r\sum_{k=i+3}^n r^{k-i-2}(w_{i,k}-\unsurr\,w_{i+1,k})$$
\noin Next, we study the action of $\n_{i+1}$ on the second sum of
$v_i$. Acting by $\n_{i+1}$ makes a term in $w_{s,i+2}$ appear by
equation $(131)$, with a coefficient $-\e_l\,\frac{r^{n-i-2+s}}{r}$.
Hence we get a sum:
\begin{equation}r\e_l\,\sum_{s=1}^{i-1}
r^{n-i-3+s}(w_{s,i+1}-\unsurr\,w_{s,i+2}) \end{equation} \noin where
we need to add
$$-\e_l\sum_{s=1}^{i-1} r^{n-i-2+s}\,w_{s,i+1}$$
\noin This latter expression is a part of the second sum in
$r\,v_i$. And since all the terms $w_{s,i}$'s, for $1\leq s\leq i-1$
get multiplied by $r$ when acting by $\n_{i+1}$, we actually get the
whole second sum
$$r\,\e_l\,\sum_{s=1}^{i-1}r^{n-i-2+s}(w_{s,i}-\unsurr\,w_{s,i+1})$$
\noin of $r\,v_i$.\\
The sum $(132)$ is the second sum in $r\,v_{i+1}$, minus the term
$$\frac{r}{l}\,(w_{i,i+1}-\unsurr\,w_{i,i+2}),$$
\noin which corresponds to $s=i$. Since those terms in $w_{i,i+1}$
and $w_{i,i+2}$ have already been processed above, we may now
conclude that:
$$v_i.\,\n_{i+1}=r(v_i+v_{i+1})$$
\noin Let $2\leq i\leq n-1$. It remains to show that
$$v_i.\n_{i-1}=r\,v_i+\unsurr\,v_{i-1}$$
\noin When we act by $\n_{i-1}$ on $v_i$, all the terms $w_{s,t}$ in
$(116)$ get multiplied by $r$, except
\begin{itemize}
\item
$w_{i,i+1}$
\item
$w_{i,k}$, each $i+2\leq k\leq n$, in the first sum
\item
$w_{i-1,i}$, $w_{s,i}$ for $1\leq s\leq i-2$ and $w_{i-1,i+1}$ in
the second sum
\end{itemize}
\noin Let's compute the action of $\n_{i-1}$ on these terms.
\begin{eqnarray}
w_{i,i+1}.\,\n_{i-1} &=& w_{i-1,i+1}+m\,w_{i-1,i}-m\,w_{i,i+1}\\&& \notag \\
\forall \,k\geq i+2,\;\;\, w_{i,k}.\,\n_{i-1} &=&
\!w_{i-1,k}+m\,r^{k-i-1}\,w_{i-1,i}-m\,w_{i,k}\\&& \notag \\
w_{i-1,i}.\,\n_{i-1} &=& \unsurl\,w_{i-1,i}\\
\forall\, s\leq i-2,\;\;\, w_{s,i}.\,\n_{i-1} &=&
w_{s,i-1}+\frac{m}{l\,r^{i-s-2}}\,w_{i-1,i}-m\,w_{s,i}\\
w_{i-1,i+1}.\,\n_{i-1} &=& w_{i,i+1}
\end{eqnarray}
\noin We see with the equations $(133)$ and $(137)$ that the
coefficient of $w_{i,i+1}$ in $v_i.\,\n_{i-1}$ is:
$$-m\bigp\unsurr-\unsurl\bigpd-\e_l.\frac{r^{n-3}}{r}$$
\noin We recall from before that:
$$r^{n-3}=\frac{\e_l}{l}$$
\noin Thus, $$\e_l.r^{n-3}=\frac{\e_l^2}{l}=\unsurl$$ \noin We will
use this equality extensively. Therefore, the coefficient of
$w_{i,i+1}$ in $v_i.\,\n_{i-1}$ is in fact:
$$\boxed{1-\frac{r}{l}-\unsur{r^2}}$$
\noin We now compute the coefficient of $w_{i-1,i}$ in
$v_i.\,\n_{i-1}$. There are several contributions coming from four
different sources: \begin{list}{\texttt{*}}{}\item the term
$w_{i,i+1}$ with coefficient $m\bigp\unsurr-\unsurl\bigpd$
\item the terms $w_{i,k}$'s for $i+2\leq k\leq n$ with coefficient
\begin{equation*}\begin{split} \sum_{k=i+2}^n m\,r^{k-i-1}r^{k-i-2} &= m\,r^{-2i-3}\,\sum_{k=i+2}^n (r^2)^k\\
&= m\,r^{-2i-3}\,r^{2i+4}\,\frac{1-(r^2)^{n-i-1}}{1-r^2}\\
&= 1-r^{2n-2i-2}
\end{split}\end{equation*}
\item the term $w_{i-1,i}$ with coefficient $\unsur{l^2}$
\item the terms $w_{s,i}$'s for $1\leq s\leq i-2$ with coefficient
\begin{equation*}\begin{split}
\e_l\sum_{s=1}^{i-2}\frac{m}{lr^{i-s-2}}\,r^{n-i-2+s} &=
r^{2n-3}r^{-2i}\,m\sum_{s=1}^{i-2}(r^2)^s\\
&=r^{2n-3}r^{-2i}\,r\,(1-(r^2)^{i-2})\\
&=r^{2n-2i-2}-(r^{n-3})^2\\
&=r^{2n-2i-2}-\unsur{l^2}
\end{split}\end{equation*}
\end{list}
\noin The sum of all these contributions is:
$$1+m\bigp\unsurr-\unsurl\bigpd=\unsur{r^2}+\frac{r}{l}-\unsur{rl}=\boxed{\unsurr
\bigp\unsurr-\unsurl\bigpd+\frac{r}{l}}$$

 \noin Now we
have:
$$\begin{array}{lll}\left.\begin{array}{lll}\text{the coefficient
of}\; w_{i,i+1}\;\text{in}\;
v_i.\,\n_{i-1}\;\text{is}\;1-\frac{r}{l}-\unsur{r^2},\\
\text{the coefficient of}\;w_{i,i+1}\;\text{in}\;v_i \;\text{is}\;
\bigp\unsurr-\unsurl\bigpd,\\
\text{the coefficient of}\; w_{i,i+1}\;\text{in}\;
v_{i-1}\;\text{is}\;-\unsurr
\end{array}\right\rb 1-\frac{r}{l}-\unsur{r^2}=r\bigp\unsurr-\unsurl\bigpd+\unsurr\bigp-\unsurr\bigpd\\\\
\left.\begin{array}{lll}\text{the coefficient of}\;
w_{i-1,i}\;\text{in}\; v_i.\,\n_{i-1}\;\text{is}\;\unsurr
\bigp\unsurr-\unsurl\bigpd+\frac{r}{l},\\
\text{the coefficient of}\;w_{i-1,i}\;\text{in}\;v_i \;\text{is}\;
\unsurl,\\
\text{the coefficient of}\; w_{i-1,i}\;\text{in}\;
v_{i-1}\;\text{is}\;\unsurr-\unsurl
\end{array}\right\rb \unsurr
\bigp\unsurr-\unsurl\bigpd+\frac{r}{l}=r.\unsurl+\unsurr.\bigp\unsurr-\unsurl\bigpd\\\\
\left.\begin{array}{lll}\text{the coefficient of}\;
w_{i-1,i+1}\;\text{in}\;
v_{i-1}.\,\n_{i-1}\;\text{is}\;\unsurr-\unsurl,\\
\text{the coefficient of}\;w_{i-1,i+1}\;\text{in}\;v_i \;\text{is}\;-\unsur{lr},\\
\text{the coefficient of}\; w_{i-1,i+1}\;\text{in}\;
v_{i-1}\;\text{is}\;1
\end{array}\right\rb\unsurr-\unsurl=r.\big(-\unsur{lr}\big)+\unsurr\times 1
\end{array}$$
\noin From now on, we do not worry about the terms in $w_{i,i+1}$,
$w_{i-1,i}$ and $w_{i-1,i+1}$ anymore. We look at the equations
$(134)$ and $(136)$. By $(134)$, the term $w_{i,k}$ gets multiplied
by $r-\unsurr$. Since $w_{i+1,k}.\,\n_{i-1}=r\,w_{i+1,k}$, we get
the whole first sum for $v_i$, multiplied by a factor $r$ and a term
$$-\unsurr\,\sum_{k=i+2}^n r^{k-i-2}\,w_{i,k},$$
\noin which is part of the global sum:
$$\unsurr\sum_{k=i+2}^n r^{k-i-1}(w_{i-1,k}-\unsurr\,w_{i,k})$$
\noin of $\unsurr\,v_{i-1}$, where the terms in $w_{i-1,i+1}$ and
$w_{i,i+1}$, corresponding to $k=i+1$, have been omitted. These
terms were already processed above. But the action of $\n_{i-1}$ on
$v_i$ also makes a term in $w_{i-1,k}$ appear (still by $(134)$),
with coefficient $r^{k-i-2}$. Thus we get the whole first sum for
$v_{i-1}$, multiplied by a factor $\unsurr$. Similarly, by $(136)$,
the term $w_{s,i}$ gets multiplied by $r-\unsurr$ for each $1\leq
s\leq i-2$. And since $w_{s,i+1}$ gets multiplied by $r$ when acting
by $\n_{i-1}$, we get the whole second sum
$$\e_l\sum_{s=1}^{i-2}r^{n-i-2+s}(w_{s,i}-\unsurr\,w_{s,i+1})$$
\noin of $r\,v_i$, minus the terms in $w_{i-1,i}$ and $w_{i-1,i+1}$
corresponding to the integer $s=i-1$. Those terms have already been
processed above. Next, the part $-\unsurr\,w_{s,i}$ yields a sum:
$$\frac{\e_l}{r}\sum_{s=1}^{i-2}r^{n-i-1+s}\big(-\unsurr\big)w_{s,i},$$
\noin which is a part of the second sum:
$$\frac{\e_l}{r}\sum_{s=1}^{i-2}r^{n-i-1+s}(w_{s,i-1}-\unsurr\,w_{s,i})$$
\noin of $\unsurr\,v_{i-1}$. Since the action of $\n_{i-1}$ on $v_i$
also makes a term in $w_{s,i-1}$ appear with coefficient the one of
$w_{s,i}$ in $v_i$ by $(136)$, that is
$\frac{\e_l}{r}\,r^{n-i-1+s}$, we get the whole second sum of
$\unsurr\,v_{i-1}$.\\
All the preceding results lead to conclude that
$v_i.\,\n_{i-1}=r\,v_i+\unsurr\,v_{i-1}$. We have now shown that the
$v_i$'s satisfy to the announced relations.

Let's show that the family $(v_1,\,v_2,\,\dots,\,v_{n-1})$ is free.
In what follows, we will denote by $\IH_{F,r^2}(n)$ the
Iwahori-Hecke algebra of the symmetric group $Sym(n)$ with parameter
$r^2$ over the field $F$. For large values of the integer $n$, we
will make use of a result of James on the dimensions of the
irreducible representations of the symmetric group $Sym(n)$. For
$n\geq 7$, James states in \cite{GDJ} that an irreducible
$K\,Sym(n)$-module, where $K$ is a field for characteristic zero, is
either one of the Specht modules $S^{(n)}$, $S^{(1^n)}$,
$S^{(n-1,1)}$, $S^{(2,1^{n-2})}$ or has dimension greater than
$n-1$. Further for $n=5$, if $d$ is the degree of an irreducible
representation of $Sym(5)$ over $K$, then $d\in\{1,4,5,6\}$. Hence
this fact also holds for $n=5$. In characteristic zero, when the
Iwahori-Hecke algebra of the symmetric group $Sym(n)$ is semisimple,
the degrees of its irreducible representations are the same as the
degrees of the irreducible representations of the symmetric group
(See for instance \cite{IHA}). Thus, for $n=5$ or $n\geq 7$, if
$\text{dim Span}_F(v_1,\,\dots,\,v_{n-1})=k<n-1$, then
$\text{Span}_F(v_1,\,\dots,\,v_{n-1})$ is either one dimensional or
must contain a one dimensional invariant subspace. Then
$l=\unsur{r^{2n-3}}$ by Theorem $4$. But we have assumed that
$l\in\{\unsur{r^{n-3}},-\unsur{r^{n-3}}\}$. Hence we get a
contradiction since $r^{2n-3}=\e\,r^{n-3}$ with $\e\in\{1,-1\}$
would imply $r^{2n}=1$, which is forbidden when $\iha$ is
semisimple. Similarly for $n=3$, if $v_1$ and $v_2$ are linearly
dependent, then the invariant vector space spanned by $v_1$ and
$v_2$ over the field $F$ has dimension one. This forces $l=-r^3$ or
$l=\unsur{r^3}$ by Theorem $4$. But we have assumed that $l\in\lb
1,-1\}$ in this case and the condition of semisimplicity
$(r^2)^3\neq 1$ prevents $r^3=1$ or $r^3=-1$ from happening. Thus,
in the case $n=3$, the vectors $v_1$ and $v_2$ are linearly
independent. Let's now deal with the case $n=4$. We want to show
that $(v_1,v_2,v_3)$ is free. By definition, we have:
$$\left\lb\begin{array}{ccccccc}
v_1 &=&
(\unsurr-\unsur{l})\,w_{12} & + & (w_{13}-\unsurr\,w_{23}) & + & r(w_{14}-\unsurr\,w_{24})\\
v_2 &=& (\unsurr-\unsur{l})\,w_{23} & + &
(w_{24}-\unsurr\,w_{34}) & + & r\,(w_{12}-\unsurr\,w_{13})\,\e_l\\
v_3 &=& (\unsurr-\unsur{l})w_{34} & + &
(w_{13}-\unsurr\,w_{14})\,\e_l & + & r(w_{23}-\unsurr\,w_{24})\,\e_l
\end{array}\right.$$
Suppose $\la_1$, $\la_2$ and $\la_3$ are scalars such that:
\begin{equation} \la_1\,v_1+\la_2\,v_2+\la_3\,v_3=0 \end{equation}
\noindent Since the $w_{ij}$'s are linearly independent over $F$, we
derive the set of equations:
\begin{eqnarray}
\la_1-\la_2\,\e_l+\la_3\,\e_l &=& 0\\
r\,\la_1-\frac{\e_l}{r}\,\la_3 &=& 0\\
-\la_1+\la_2-\la_3\,\e_l &=& 0\\
-\frac{\la_2}{r}+\la_3\big(\unsurr-\unsur{l}\big) &=& 0
\end{eqnarray}

\noindent $(139)$, $(140)$, $(141)$ and $(142)$ are respectively
obtained by equaling to zero the respective coefficients of
$w_{13}$, $w_{14}$, $w_{24}$ and $w_{34}$ in the relation $(138)$.
Let's first deal with $l=-\unsurr$ and $\e_l=-1$. By $(139)$ and
$(141)$, we have: $$\la_3-\la_2=\la_2+\la_3,$$ from which we derive
$\la_2=0$ since our base field has characteristic zero. Then by
$(129)$, we get $\la_1=\la_3$ and by $(140)$, we get
$(r+\unsurr)\,\la_3=0$. Since we have assumed that $(r^2)^2\neq 1$
by semisimplicity of the Iwahori-Hecke algebra $\IH_{F,r^2}(4)$, we
cannot have $r^2=-1$. Hence we get $\la_3=0$ and so
$\la_1=\la_2=\la_3=0$.\\
Suppose now $l=\unsurr$ and $\e_l=1$. By $(140)$, we have
$\la_1=\frac{\la_3}{r^2}$ and by $(139)$ or $(141)$, we get:
$\la_2=(1+\unsur{r^2})\,\la_3$. Then by $(142)$, we get:
$$\Big(\unsur{r^3}+r\Big)\la_3=0$$
\noindent Since we have assumed that $(r^2)^4\neq 1$, it is
impossible to have $(r^2)^2=-1$. It follows that $\la_3=0$, and by
the above relations binding $\la_1$ and $\la_3$ on one hand and
$\la_2$ and $\la_3$ on the other hand, we also get $\la_1=\la_2=0$.
This achieves the proof of the fact that the family of vectors
$(v_1,v_2,v_3)$ is free.

%Suppose the $v_i$'s are
%not linearly independent. Then there exists a family
%$(\la_1,\,\dots,\,\la_{n-1})$ of scalars in
%$F\setminus\lb(0,\,\dots,\,0)\rb$ such that:
%$$\la_1v_1+\la_2v_2+\dots+\la_{n-1}v_{n-1}=0$$
%\noin By applying successively $\n_1$, $\n_2$, $\dots,$ $\n_{n-1}$
%to both members of this equation, we get $n-1$ equations that we
%gather in the matrix equation:
%$$\hspace{-1.4in}\begin{array}{cccc}
%\left(\begin{array}{cccccccc} -\frac{\la_1}{r}+\frac{\la_2}{r} &
%r\,\la_2
%&\dots & \dots & r\,\la_j & \dots & r\,\la_{n-2} & r\,\la_{n-1}\\
%r\,\la_1 & r\,\la_1-\frac{\la_2}{r}+\frac{\la_3}{r}& \dots & \dots &
%r\,\la_j & \dots & r\,\la_{n-2} & r\,\la_{n-1}\\
%r\,\la_1 & r\,\la_2 & \ddots & &\vdots & & \vdots &\vdots \\
%\vdots&\vdots& & \ddots & \vdots & & \vdots & \vdots\\
%& & & & r\,\la_{j-1}-\frac{\la_j}{r}+\frac{\la_{j+1}}{r} & & \vdots
%&
%\vdots\\
%\vdots &\vdots & & & \vdots & &
%r\,\la_{n-3}-\frac{\la_{n-2}}{r}+\frac{\la_{n-1}}{r} & r\,\la_{n-1}\\
%r\,\la_1 & r\,\la_2 & & & r\,\la_j & & r\,\la_{n-2} &
%r\,\la_{n-2}-\frac{\la_{n-1}}{r}
%\end{array}\right)
%&\!\!\!\!\!\!\!\left(\begin{array}{llllll}v_1\\v_2\\
%\,\,\vdots\\ \,\,\vdots\\ \,\,\vdots\\v_{n-2}\\v_{n-1}
%\end{array}\!\!\right) &\!\!\!\!\!\!=&\!\!\!\!\!0
%\end{array}$$
%\begin{center} [proof to be inserted, and $\IH_{F,r^2}(Sym(n))$ to be defined]\end{center}
 We
are now able to conclude: the vector space
$\text{Span}_F(v_1,\dots,v_{n-1})\subset\V$ is $(n-1)$-dimensional,
is invariant under the action of the $g_i$'s and is a
$\IH_{F,r^2}(n)$-module since it is a proper nontrivial invariant
subspace of $\V$. Then by the relations described above, it is also
irreducible.

We note that for $n=6$ the sufficient condition on $l$ and $r$ so
that there exists an irreducible $5$-dimensional invariant subspace
of $\V$ holds in Theorem $5$, since simple computations of the
dimensions of the irreducible $F\,Sym(6)$-modules show that there is
no irreducible representation of $\IH_{F,r^2}(6)$ of degree $d$ with
$1<d<5$. This forces the family of vectors $(v_1,v_2,v_3,v_4,v_5)$
to be free by the same argument as in the cases $n=5$ and $n\geq 7$.
However, for the whole theorem to be true in the case $n=6$, we need
to investigate further about the irreducible representations of $
\IH_{F,r^2}(6)$ of degree $5$ corresponding to the partition $(3,3)$
and its conjugate partition $(2,2,2)$. An idea is to reduce somehow
to the case $n=5$. We will use the branching rule for the
restriction of the Specht modules of the Iwahori-Hecke algebras that
are semisimple, as it is described in Corollary $6.2$ of \cite{IHA}.
Assuming that $\ih(6)$ is semisimple, we have:
$$S^{(3,3)}\da_{\ih(5)}\,\simeq\, S^{(3,2)}\qquad\qquad
S^{(2,2,2)}\da_{\ih(5)}\,\simeq\, S^{(2,2,1)}$$ Moreover, we have
the following fact:
\newtheorem{Fact}{Fact}
\begin{Fact} Suppose $\ih(5)$ is semisimple. Then, up to equivalence, the two irreducible matrix representations
of degree $5$ of $\ih(5)$ are respectively defined by the matrices
$P_1$, $P_2$, $P_3$, $P_4$ and $Q_1$, $Q_2$, $Q_3$, $Q_4$ given by:
$$P_1:=\begin{bmatrix}
r & & & &\\
  &r& & &\\
  & &r& &\\
1 & &-r^2&-\unsurr&\\
  &1&    &        &-\unsurr\\
\end{bmatrix},\; P_2:=\begin{bmatrix}
-\unsurr& & & 1 & \\
        &-\unsurr&1& & 1\\
        & & r& & \\
        & &  &r&\\
        & &  & &r
\end{bmatrix}$$
$$P_3:=\begin{bmatrix}
r & & & &\\
  &r& & &\\
  &1 &-\unsurr & &\\
1 & & &-\unsurr&-r^2\\
  & & &        &r\\
\end{bmatrix},\; P_4:=\begin{bmatrix}
&1 &-r &  & \\
1 &r-\unsurr&1& & \\
        & & r& & \\
        & & -r^2 & &1\\
        & &  r   & 1 &r-\unsurr
\end{bmatrix}$$
and for the conjugate representation:
$$Q_1:=\begin{bmatrix}
-\unsurr & & & &\\
  &-\unsurr& & &\\
  & &-\unsurr& &\\
1 & &-\unsur{r^2}&r&\\
  &1&    &        &r\\
\end{bmatrix},\; Q_2:=\begin{bmatrix}
r& & & 1 & \\
        &r&1& & 1\\
        & &-\unsurr & & \\
        & &  &-\unsurr&\\
        & &  & &-\unsurr
\end{bmatrix}$$
$$Q_3:=\begin{bmatrix}
-\unsurr & & & &\\
  &-\unsurr& & &\\
  &1 & r & &\\
1 & & &r&-\unsur{r^2}\\
  & & &        &-\unsurr\\
\end{bmatrix},\; Q_4:=\begin{bmatrix}
&1 &\unsurr &  & \\
1 &r-\unsurr&1& & \\
        & &-\unsurr& & \\
        & & -\unsur{r^2} & &1\\
        & &  -\unsurr   & 1 &r-\unsurr
\end{bmatrix}$$

\noin where the blanks are to be filled with zeros.
\end{Fact}
\textsc{Proof of the Fact:} since the second matrix representation
is the conjugate of the first one where we replaced $r$ by
$-\unsurr$, it suffices to show that the first set of matrices
defines a representation of $\ih(5)$ that is irreducible. It is a
direct verification that that the so-defined matrices $P_i$'s,
$i=1\dots 4$, satisfy to the braid relations and to the Hecke
algebra relations $P_i^2+m\,P_i=I_5$ where $I_5$ is the identity
matrix of size $5$. Let us now show that this Iwahori-Hecke algebra
matrix representation is irreducible. Assuming $\ih(5)$ is
semisimple, a non irreducibility of the representation means that
there exists a nonzero vector $w$ in $F^5$ and scalars
$\la_1,\la_2,\la_3,\la_4$ in $F$ such that for all $i=1,\dots,4$ we
have:
\begin{equation*}
P_i\,w=\la_i\,w\qquad\qquad(\star_i)
\end{equation*}
Let $w_1,w_2,w_3,w_4,w_5$ be the coordinates of $w$ in the canonical
basis of $F^5$. We write all the equations provided by $(\star_1)$,
$(\star_2)$, $(\star_3)$ and $(\star_4)$.
$$
(\star_1)\left\lbrace\begin{array}{cccc} r\,w_1&=&\la_1\,w_1 & (E1)\\
r\,w_2&=&\la_1\,w_2&(E2)\\
r\,w_3&=&\la_1\,w_3& (E3)\\
w_1-r^2\,w_3-\unsurr\,w_4&=&\la_1\,w_4& (E4)\\
w_2-\unsurr\,w_5&=&\la_1\,w_5& (E5)\\
\end{array}\right.$$ $$ (\star_2)\left\lbrace\begin{array}{cccc}
-\unsurr\,w_1+w_4&=&\la_2\,w_1&(E6)\\
-\unsurr\,w_2+w_3+w_5&=&\la_2\,w_2&(E7)\\
r\,w_3&=&\la_2\,w_3&(E8)\\
r\,w_4&=&\la_2\,w_4&(E9)\\
r\,w_5&=&\la_2\,w_5&(E10)
\end{array}\right.$$
$$
\,(\star_3)\left\lbrace\begin{array}{cccc} r\,w_1&=&\la_3\,w_1 & (E11)\\
r\,w_2&=&\la_3\,w_2&(E12)\\
w_2-\unsurr\,w_3&=&\la_3\,w_3& (E13)\\
w_1-\unsurr\,w_4-r^2\,w_5&=&\la_3\,w_4& (E14)\\
r\,w_5&=&\la_3\,w_5& (E15)\\
\end{array}\right.$$ $$\;\;\;\;\;\,\,(\star_4)\left\lbrace\begin{array}{cccc}
w_2-r\,w_3&=&\la_4\,w_1&(E16)\\
w_1+(r-\unsurr)\,w_2+w_3&=&\la_4\,w_2&(E17)\\
r\,w_3&=&\la_4\,w_3&(E18)\\
-r^2\,w_3+w_5&=&\la_4\,w_4&(E19)\\
r\,w_3+w_4+(r-\unsurr)w_5&=&\la_4\,w_5&(E20)
\end{array}\right.$$
\\
We will show that it is impossible to have such a set of relations.
Our first step is to show that these relations force all the
$\la_i$'s to take the same value $r$. From there, we show that all
the $w_i$'s must then be zero, which is a contradiction. First,
since all the $\la_i$'s must be equal to $r$ or all the $\la_i$'s
must be equal to $-\unsurr$, it suffices to show that $\la_1=r$.
\\\\Now, if $\la_1\neq r$, we get:$\;\;$
\begin{tabular}[c]{|c|c|}
\hline Equation used & Result\\\hline $(E_1)$&$w_1=0$\\\hline
$(E_2)$&$w_2=0$\\\hline $(E_3)$&$w_3=0$\\\hline
$(E_6)$&$w_4=0$\\\hline $(E_7)$&$w_5=0$\\\hline
\end{tabular}\\\\
\noin Since we get that all the coordinates of the vector $w$ must
be equal to zero, we see that $\la_1$ must equal $r$, so that all
the $\la_i$'s must in fact be equal to $r$. Now, from $(E_4)$ and
$(E_{14})$, we have:
\begin{eqnarray*}
\bigp r+\unsurr\bigpd\,w_4&=&w_1-r^2\,w_3\\
\bigp r+\unsurr\bigpd\,w_4&=&w_3-r^2\,w_5
\end{eqnarray*}
from which we derive that: \begin{equation*} w_3=w_5 \end{equation*}
Using this equality, we get from $(E_{13})$ that:
\begin{equation*}
w_2=\bigp r+\unsurr\bigpd\,w_5
\end{equation*}
and from $(E_7)$ that:
\begin{equation*}
\bigp r+\unsurr\bigpd\,w_2=2\,w_5
\end{equation*}
Gathering these two equalities, it yields:
$$\bigp r+\unsurr\bigpd^2\,w_5=2\,w_5$$
Since we have assumed that $\ih(5)$ is semisimple, we have
$(r^2)^4\neq 1$. Thus, it is impossible to have $\bigp
r+\unsurr\bigpd^2=2$. Hence we get $w_5=w_3=0$. Then it comes
$w_2=0$ by $(E_{13})$ and $w_4=0$ by $(E_{20})$ and $w_1=0$ by
$(E_{17})$. It is a contradiction to have $w=0$, hence our
so-defined matrix representations are irreducible. We use these two
irreducible $5$-dimensional matrix representations to show the
following theorem:
\newtheorem{t-l}{Result}
\begin{t-l}\hfill\\
Suppose $n=5$. Then there exists an irreducible $5$-dimensional
invariant subspace of $\V$ if and only if $l=r$.
\end{t-l}
\noin\textsc{Proof:} suppose that there exists an irreducible
$5$-dimensional invariant subspace of $\V$, say $\W$. Then $\W$ is
an irreducible $\ih(5)$-module. Then there must exist a basis
$(w_1,w_2,w_3,w_4,w_5)$ of $\W$ in which the matrices of the left
action of the $g_i$'s, $i=1,\dots,4$, are either the $P_i$'s or the
$Q_i$'s. We show that the first possibility leads to force the value
$r$ for $l$. Suppose such vectors exist. We read on the matrices
$P_1$ and $P_4$ that:
\begin{eqnarray*}
g_1.w_4&=&-\unsurr\,w_4\\
g_3.w_4&=&-\unsurr\,w_4
\end{eqnarray*}
These two relations imply that in $w_4$ there are no terms in
$w_{34}$, $w_{45}$, $w_{35}$ and in $w_{12}$, $w_{15}$, $w_{25}$. In
what follows, we will denote by $\mu_{ij}^{(k)}$ the coefficient of
$w_{ij}$ in $w_k$. Thus, we have:
$$w_4=\mu_{23}^{(4)}\,w_{23}+\mu_{13}^{(4)}\,w_{13}+\mu_{24}^{(4)}\,w_{24}+\mu_{14}^{(4)}\,w_{14}$$
Furthermore, we read on the matrix $P_2$ that:
$$g_2.w_4=r\,w_4+w_1$$
By looking at the coefficient of $w_{23}$ in this equality and by
using the expression for $w_4$, we get:
$$\unsur{l}\,\mu_{23}^{(4)}+\frac{m}{l}\,\mu_{13}^{(4)}=r\,\mu_{23}^{(4)}+\mu_{23}^{(1)}\;(\star)$$
We will show that $\mu_{23}^{(1)}=0$. To this aim, we look at the
coefficient of $w_{13}$ in
$$g_1.w_1=r\,w_1+w_4$$
to obtain:
$$\mu_{23}^{(1)}=r\,\mu_{13}^{(1)}+\mu_{13}^{(4)}$$
Further we look at the coefficient of $w_{13}$ in
$g_3.w_1=r\,w_1+w_4$ to get the equation
$\mu_{14}^{(1)}=r\,\mu_{13}^{(1)}+\mu_{13}^{(4)}$. Since
$g_2.w_1=-\unsurr\,w_1$, in $w_1$ there is no term in $w_{14}$.
Hence we get:
$$\mu_{13}^{(4)}=-r\,\mu_{13}^{(1)}$$
Mixing the two equalities now yields
$$\mu_{23}^{(1)}=0,$$
as desired. Now $(\star)$ can be rewritten as:
$$\unsur{l}\,\mu_{23}^{(4)}+\frac{m}{l}\,\mu_{13}^{(4)}=r\,\mu_{23}^{(4)}\;(\star\star)$$
Since $\mu_{14}^{(1)}=\mu_{23}^{(1)}=0$, we derive from
$$g_1.w_1=r\,w_1+w_4$$
that
$$\mu_{24}^{(1)}=\mu_{14}^{(4)}$$
and we derive from
$$g_3.w_1=r\,w_1+w_4$$
that
$$\mu_{24}^{(1)}=\mu_{23}^{(4)}$$
Now it comes: $$\mu_{23}^{(4)}=\mu_{14}^{(4)}$$ Moreover, since
$g_3.w_4=-\unsurr\,w_4$, we have by $(78)$ applied with $s=1$, $q=3$
and $\la=-\unsurr$ that:
$$\mu_{14}^{(4)}=-\unsurr\,\mu_{13}^{(4)}$$
It follows that
$$\mu_{23}^{(4)}=-\unsurr\,\mu_{13}^{(4)}$$
Plugging this value for $\mu_{23}^{(4)}$ into $(\star\star)$ yields:
$$\bigp\frac{m}{l}-\unsur{lr}\bigpd\,\mu_{13}^{(4)}=-\mu_{13}^{(4)}\;\;(\star\star\star)$$
Furthermore, the coefficients in $w_4$ are related in a certain way
that prevents $\mu_{13}^{(4)}$ from being zero. Indeed, as seen
along the way, we have:
$$\left|\begin{array}{ll}
\mu_{23}^{(4)}=\;\;\;\mu_{14}^{(4)}\\
\mu_{14}^{(4)}=-\unsurr\,\mu_{13}^{(4)}
\end{array}\right.$$
Further, with $g_3.w_4=-\unsurr\,w_4$, we add a new equation by
$(78)$ applied with $s=2$, $q=3$ and $\la=-\unsurr$:
$$\mu_{24}^{(4)}=-\unsurr\,\mu_{23}^{(4)}$$
It now appears clearly from these relations that if
$\mu_{13}^{(4)}=0$, then the other three coefficients in $w_4$ are
also zero. This is naturally impossible as $w_4$ is a basis vector.
Thus, $(\star\star\star)$ reduces to:
$$\frac{m}{l}-\unsur{lr}=-1,$$
\emph{i.e\/}
$$\boxed{l=r}$$
Conversely, suppose $l=r$.\\
Let $w4:=w_{23}-\unsurr\,w_{24}+w_{14}-r\,w_{13}$\\
and set
\begin{eqnarray*}
w5:=g_4.w4=r\,w_{23}-r^2\,w_{13}-\unsurr\,w_{25}+w_{15}\\
w1:=g_2.w4-r\,w4=-\unsurr\,w_{34}+w_{24}-r\,w_{12}+w_{13}\\
w2:=g_4.w1=-\unsurr\,w_{35}+w_{25}-r^2\,w_{12}+r\,w_{13}\\
w3:=g_3.w2-r\,w2=-r^2\,w_{13}+r\,w_{14}-\unsurr\,w_{45}+w_{35}
\end{eqnarray*}
\newtheorem{Claim}{Claim}
\begin{Claim}
If $l=r$, then the vectors $w1,w2,w3,w4,w5$ span an irreducible
$5$-dimensional invariant subspace of $\V$ over $\Q(r)$.
\end{Claim}
\noin\textsc{Proof:} it is a direct verification that the family of
vectors $w1,w2,w3,w4,w5$ are linearly independent over $\Q(r)$ and
that they satisfy to the relations:

$$\begin{array}{l}\left|\begin{array}{ccc}
g_1.w1\eg r\,w1+w4\\
g_2.w1\eg -\unsurr\,w1\\
g_3.w1\eg r\,w1+w4\end{array}\right.\\\\
\left|\begin{array}{ccc}
g_1.w2\eg r\,w2+w5\\
g_2.w2\eg -\unsurr\,w2\\
g_4.w2\eg (r-\unsurr)\,w2+w1
\end{array}\right.\\\\
\left|\begin{array}{l}
g_1.w3\;\; = \;\;r\,w3-r^2\,w4\\
g_2.w3\;\; = \;\;r\,w3+w2\\
g_3.w3\;\; = \;\;-\unsurr\,w3\\
g_4.w3\;\; = \;\;-r\,w1+w2+r\,w3-r^2\,w4+r\,w5\end{array}\right.\\\\
\left|\begin{array}{ccc} g_1.w4\eg-\unsurr\,w4\\
g_3.w4\eg -\unsurr\,w4
\end{array}\right.\\\\
\left|\begin{array}{l} g_1.w5\;\;= \;\;-\unsurr\,w5\\
g_2.w5\;\;=\;\; r\,w5+w2\\
g_3.w5\;\;=\;\; r\,w5-r^2\,w4\\
g_4.w5\;\;=\;\; (r-\unsurr)\,w5+w4
\end{array}\right.
\end{array}$$
\\
\noin By Fact $1$ the claim is hence proven. To finish the proof of
Result $1$, it suffices to show that the second irreducible
representation of degree $5$ described by the matrices $Q_i$'s
cannot occur. Suppose that there exists a basis
$(v_1,v_2,v_3,v_4,v_5)$ of $\W$ in which the matrices of the left
action of the $g_i$'s, $i=1,\dots,4$, are the $Q_i$'s. In what
follows, we will denote by $\la_{ij}^{(k)}$ the coefficient of
$w_{ij}$ in $v_k$. By the relations
\begin{eqnarray*}
g_3.v_4&=&r\,v_4\\
g_1.v_4&=&r\,v_4
\end{eqnarray*}
we get the relations between the coefficients in $v_4$:
\begin{eqnarray*}
\la_{45}^{(4)}&=&r\,\la_{35}^{(4)}\\
\la_{14}^{(4)}&=&r\,\la_{13}^{(4)}\\
\la_{24}^{(4)}&=&r\,\la_{23}^{(4)}\\
\la_{25}^{(4)}&=&r\,\la_{15}^{(4)}\\
\la_{24}^{(4)}&=&r\,\la_{14}^{(4)}\\
\la_{23}^{(4)}&=&r\,\la_{13}^{(4)}
\end{eqnarray*}
Next, since $g_4.v_4=v_5$ and
$g_3.v_5=-\unsurr\,v_5-\unsur{r^2}\,v_4$, we get
$$g_3\,g_4\,v_4=-\unsurr\,g_4\,v_4-\unsur{r^2}\,v_4$$ Looking at the
term in $w_{12}$ in this equation yields:
$$r^2\,\la_{12}^{(4)}=-\la_{12}^{(4)}-\unsur{r^2}\,\la_{12}^{(4)}$$
Since $(r^2)^3\neq 0$, we have $r^2+\unsur{r^2}+1\neq 0$. Thus we
get $\la_{12}^{(4)}=0$. By looking at the coefficient of $w_{12}$ in
$g_3.v_1=-\unsurr\,v_1+v_4$, we now get:
$r\,\la_{12}^{(1)}=-\unsurr\,\la_{12}^{(1)}$. Thus we also have
$\la_{12}^{(1)}=0$. Further since $\la_{13}^{(1)}=r\,\la_{12}^{(1)}$
by the relation $g_2.v_1=r\,v_1$, it follows that
$\la_{13}^{(1)}=0$. We then look at the coefficient of $w_{13}$ in
$g_2.v_4=-\unsurr\,v_4+v_1$ to get:
$-m\la_{13}^{(4)}=-\unsurr\,\la_{13}^{(4)}$ where we used that
$\la_{12}^{(4)}=0$. It comes $\la_{13}^{(4)}=0$. From there we
derive from the set of relations above that:
$$\la_{13}^{(4)}=\la_{23}^{(4)}=\la_{24}^{(4)}=\la_{14}^{(4)}=0$$
Thus we have:
$$v_4=\la_{34}^{(4)}\,w_{34}+r\,\la_{35}^{(4)}\,w_{45}+\la_{35}^{(4)}\,w_{35}+r\,\la_{15}^{(4)}\,w_{25}+\la_{15}^{(4)}\,w_{15}$$
Let's look at the term in $w_{25}$ in
$g_3\,g_4\,v_4=-\unsurr\,g_4\,v_4-\unsur{r^2}\,v_4$. It comes:
$$-mr^2\,\la_{15}^{(4)}=-\unsurr\,\bigp-mr\,\la_{15}^{(4)}\bigpd-\unsur{r^2}\,r\,\la_{15}^{(4)}$$
\emph{i.e\/} $$-r^3\,\la_{15}^{(4)}=0$$ Hence $\la_{15}^{(4)}=0$ and
there are only three terms left in $v_4$:
$$v_4=\la_{34}^{(4)}\,w_{34}+r\,\la_{35}^{(4)}\,w_{45}+\la_{35}^{(4)}\,w_{35}$$
Let's now look at the term in $w_{34}$ in the relation
$g_1.v_1=-\unsurr\,v_1+v_4$. We get:
$$r\la_{34}^{(1)}=-\unsurr\,\la_{34}^{(1)}+\la_{34}^{(4)}$$
\emph{i.e\/}
$$\la_{34}^{(4)}=\bigp r+\unsurr\bigpd\,\la_{34}^{(1)}$$
On the other hand, by looking at the coefficient of $w_{34}$ in
$g_2.\,v_4=-\unsurr\,v_4+v_1$, we get:
$$-m\,\la_{34}^{(4)}=-\unsurr\,\la_{34}^{(4)}+\la_{34}^{(1)}$$
\emph{i.e\/}
$$\la_{34}^{(4)}=\unsurr\,\la_{34}^{(1)}$$
The two relations binding $\la_{34}^{(4)}$ and $\la_{34}^{(1)}$ now
yield $\la_{34}^{(4)}=\la_{34}^{(1)}=0$. Thus, there are only two
terms left in $v_4$. Explicitly we have:
$$v_4=r\,\la_{35}^{(4)}\,w_{45}+\la_{35}^{(4)}\,w_{35}$$
But by looking at the coefficient of $w_{34}$ in $g_3.v_4=r\,v_4$ we
have:
$$m\,\la_{45}^{(4)}=0$$
Then $v_4$ would be zero, a contradiction.

Let's go back to the proof of Theorem $5$. When $n=6$, suppose that
there exists a $5$-dimensional invariant subspace $\W$ of $\V$ with
$$\W\simeq\,S^{(3,3)}\;\;\text{or}\;\;\W\simeq\,S^{(2,2,2)}$$
Then
$$\W\da_{\ih(5)}\simeq\,S^{(3,2)}\;\;\text{or}\;\;\W\da_{\ih(5)}\simeq\,S^{(2,2,1)}$$
Then there must exist a basis of $\W$ in which the matrices of the
left action by the $g_i$'s are the $P_i$'s or the $Q_i$'s. We first
place ourself in the first situation and adapt the proof of Result
$1$. Let $(w_1,w_2,w_3,w_4,w_5)$ be a basis of $\W$ in which the
matrices of the left action by the $g_i$'s are the $P_i$'s. Since
the relation $g_1.w_4=-\unsurr\,w_4$ (resp $g_3.w_4=-\unsurr\,w_4$)
forces that in $w_4$ there is no term in $w_{36},w_{46},w_{56}$
(resp in $w_{26},w_{16}$), we see that the shape of $w_4$ is still
the same. Then the presence of a sixth node does not modify the
arguments of Result $1$ and we get by the exact same arguments as in
the proof of Result $1$ that $l=r$. Furthermore $\W$ must be spanned
by the linearly independent vectors $w_1,w_2,w_3,w_4,w_5$ provided
in the proof of Result $1$. But $\W$ is a submodule of $\V$. Hence
$g_5.\W\subseteq\W$. This is not compatible with the spanning
vectors above.\\
Second, still following the proof in Result $1$, we show that it is
impossible to have a basis $(v_1,v_2,v_3,v_4,v_5)$ of $\W$ in which
the matrices of the left action by the $g_i$'s are the $Q_i$'s.
Suppose such vectors can be found. Taking the same notations as in
Result $1$ and using the same arguments as in the proof of Result
$1$ with the addition of a sixth node, we must have:
\begin{multline*}v_4=r\la_{35}^{(4)}\,w_{45}+\la_{35}^{(4)}\,w_{35}+\\
\la_{56}^{(4)}\,w_{56}+ r\la_{36}^{(4)}\,w_{46}+
\la_{36}^{(4)}\,w_{36}+ r\la_{16}^{(4)}\,w_{26}+
\la_{16}^{(4)}\,w_{16}\end{multline*} We simplify further the shape
of $v_4$ by looking at the terms in $w_{12}$ in $g_1.v_4=r\,v_4$.
Such terms appear only when $g_1$ acts on $w_{26}$, with coefficient
$m\,r^3$. Since there is no term in $w_{34}$ in $v_4$, we simply
get:
$$m\,r^3\la_{26}^{(4)}=0$$
Hence, $\la_{26}^{(4)}=\la_{16}^{(4)}=0$.\\
Further, by looking at the coefficient of $w_{34}$ in
$g_3.v_4=r\,v_4$, we obtain:
$$m\,r\la_{35}^{(4)}+m\,r^2\,\la_{36}^{(4)}=0$$
\emph{i.e\/}
$$\la_{36}^{(4)}=-\unsurr\,\la_{35}^{(4)}$$
Furthermore, we claim that $v_4$ cannot be a multiple of $w_{56}$.
Indeed, if so, in $-\unsurr\,g_4.v_4-\unsur{r^2}\,v_4$ there is no
term in $w_{36}$, but in $g_3g_4.v_4$, there is one. This would
contradict
$$g_3g_4.v_4=-\unsurr\,g_4.v_4-\unsur{r^2}\,v_4$$
Thus, without loss of generality, we may set $\la_{35}^{(4)}=1$ and
we rewrite $v_4$ as follows:
$$v_4=r\,w_{45}+w_{35}+\la_{56}^{(4)}\,w_{56}-w_{46}-\unsurr\,w_{36}$$
We now derive from $v_1=g_2.v_4+\unsurr\,v_4$ that:
\begin{multline*}v_1=
(1+r^2)\,w_{45}+r\,w_{35}+w_{25}\\+\bigp
r+\unsurr\bigpd\,\la_{56}^{(4)}\,w_{56}-\bigp
r+\unsurr\bigpd\,w_{46}-w_{36}-\unsurr\,w_{26}
\end{multline*}
If we look at the coefficient of $w_{26}$ in
$g_3.v_1=-\unsurr\,v_1+v_4$ we get
$$-1=\unsur{r^2},$$
which is a contradiction as $(r^2)^2\neq 1$. \\
Thus, we conclude that when $n=6$ it is impossible to have an
irreducible $5$-dimensional submodule of $\V$ isomorphic to
$S^{(3,3)}$ or isomorphic to $S^{(2,2,2)}$. Rather, it must be
isomorphic to $S^{(5,1)}$ or $S^{(2,1^4)}$. Then, as seen before,
one of these representations leads to $l\in\lbrace
\unsur{r^3},-\unsur{r^3}\rbrace$ while its conjugate representation
cannot occur. Hence, the proof of Theorem $5$ is now complete.

\subsection{Proof of the Main Theorem}
In this section we prove the theorem:
\begin{thm}
Assume $\iha$ is semisimple. Then,\\\\
$\n^{(3)}$ is irreducible if and only if $l\not\in\lb
-r^3,1,-1,\unsur{r^3}\}$\\\\
For $n\geq 4$, $\n^{(n)}$ is irreducible if and only if $l\not\in\lb
r,-r^3,\unsur{r^{n-3}},-\unsur{r^{n-3}},\unsur{r^{2n-3}}\}$\\
\end{thm}
\noindent\textsc{Proof:} we already proved the theorem for the small
values $n=3,4,5,6$ in the previous sections. For larger values of
$n$, we will use the following two statements of James which can be
found in \cite{GDJ}:
\begin{Prop}
Let $K$ be a field of characteristic zero.\\
For all $n\geq 7$, every irreducible $K\,Sym(n)$-module is either
isomorphic to one of the Specht modules $S^{(n)}$, $S^{(1^n)}$,
$S^{(n-1,1)}$, $S^{(2,1^{n-2})}$ or has dimension greater than n-1.
\end{Prop}
\begin{Prop} Let $K$ be a field of characteristic zero.\\
For all $n\geq 9$, every irreducible $K\,Sym(n)$-module is either
isomorphic to one of the Specht modules $S^{(n)}$, $S^{(n-1,1)}$,
$S^{(n-2,2)}$, $S^{(n-2,1,1)}$ or their conjugates, or has dimension
greater than $\frac{(n-1)(n-2)}{2}$.
\end{Prop}
\noindent For a proof of these two facts, we refer the reader to
\cite{GDJ}. We will use an immediate corollary:
\begin{Cor}
Let $\mathcal{D}$ be an irreducible $F\,Sym(n)$-module with $n=7$ or
$n\geq 9$, where $F$ is a field of characteristic zero. Then, there
are two possibilities:
$$\begin{array}{lll}
\text{either} & dim\,\D\in\lb 1,n-1,\frac{n(n-3)}{2},\frac{(n-1)(n-2)}{2}\rb\\
&\\ \text{or} & dim\,\D\,>\,\frac{(n-1)(n-2)}{2}
\end{array}$$
\end{Cor}
\noindent\textsc{Proof of the Corollary:} the dimensions of the
irreducible $F\,Sym(n)$-modules $S^{(n-2,1,1)}$, $S^{(n-2,2)}$ (and
their respective conjugates) are given by the number of standard
Young tableau of shapes $(n-2,1,1)$ and $(n-2,2)$. In a standard
tableau, the numbers increase down each row and down each column.
For the shape $(n-2,2)$, the two possible configurations are the
following:
%\begin{center}

$$\begin{array}{ll}
%\longleftarrow\hline\longrightarrow\\
%\vspace{.01in}\overline{
\begin{tabular}[c]{|c|c|p{.35in}|p{.1in}|}
\hline 1 & 2 & \dots\dots & \\
\hline
\end{tabular}
%}
\\
\begin{tabular}[c]{|p{.075in}|p{.075in}|}
 & \\\hline\end{tabular}
\end{array}$$
\noindent where you need to pick two integers out of the $(n-2)$
remaining ones, to fill in the two boxes on the second row. The
other possible configuration is:
$$\begin{array}{ll}
\begin{tabular}[c]{|c|c|p{.35in}|p{.1in}|}
\hline 1 & 3 & \dots\dots & \\
\hline
\end{tabular}\\
\begin{tabular}[c]{|c|p{.075in}|}
2 & \\\hline\end{tabular}
\end{array}$$
\noindent And there are $(n-3)$ possible choices for the node
$(2,2)$. Thus, the total number of standard Young tableau of shape
$(n-2,2)$ is
$$\binom{n-2}{2}+(n-3)=\frac{n(n-3)}{2}$$
\noindent Hence we have $\text{dim}\,S^{(n-2,2)}=\frac{n(n-3)}{2}$.
And there are $\binom{n-1}{2}$ standard Young tableaux of shape
$(n-2,1,1)$, so $\text{dim}\,S^{(n-2,1,1)}=\frac{(n-1)(n-2)}{2}$.
%\end{center}

Now for $n\geq 9$, the corollary is nothing else but a reformulation
of Proposition $4$ with $\text{dim}\, S^{(n-2,2)}=\frac{n(n-3)}{2}$
and $\text{dim}\,S^{(n-2,1,1)}=\frac{(n-1)(n-2)}{2}$. By Proposition
$3$, the irreducible $F\,Sym(7)$-modules have dimension $1,6$ or
dimension greater than $6$. In the table below, we computed the
dimensions of the irreducible $F\,Sym(7)$-modules that have
dimensions greater than $6$. We used the Hook formula. Here, the
conjugates of the Specht modules mentioned in the table have been
omitted since they have the same dimensions by the Hook formula.
% As can be read in the table, We read in the table. As illustrated by the table
As illustrated by the table, the irreducible $sym(7)$-modules have
dimensions either $1$ or $6$ or $\frac{7\ti4}{2}=14$ or
$\frac{6\ti5}{2}=15$ or dimensions greater than $15$. This achieves
the proof of the corollary. Note that Proposition $4$ fails in the
case $n=7$, as for instance $\text{dim}\,S^{(4,3)}=14$, while
Corollary $3$ holds in this case.

\newpage
\begin{tabular}[c]{|p{1.5in}|c|c|}
\hline class of irreducible $Sym(7)$-module &
hooklengths & dimension \\
\hline $\qquad\qquad S^{(5,2)}$ & $\begin{array}{llll}\\
\begin{tabular}[c]{|c|c|c|c|c|}
\hline 6 & 5 & 3 & 2 & 1 \\
\hline
\end{tabular}\\
\begin{tabular}[c]{|c|c|}
2 & 1 \\ \hline\end{tabular}\\ \begin{array}{l}\end{array}
\end{array}$
& $f^{(5,2)}=\frac{7!}{6\ti5\ti4\ti3}=14$\\
\hline $\qquad\qquad S^{(5,1,1)}$ & $\begin{array}{lllll}\\
\begin{tabular}[c]{|c|c|c|c|c|}
\hline 7 & 4 & 3 & 2 & 1\\
\hline\end{tabular}\\
\begin{tabular}[c]{|c|c|}
2\\ \hline\end{tabular}\\
\begin{tabular}[c]{|c|c|}
1\\ \hline\end{tabular}\\ \begin{array}{l}\end{array}
\end{array}$ & $f^{(5,1,1)}=\frac{7!}{7\ti4\ti3\ti2\ti2}=15$\\
\hline $\qquad\qquad S^{(4,3)}$ & $\begin{array}{llll}\\
\begin{tabular}{|c|c|c|c|}
\hline 5&4&3&1\\\hline\end{tabular}\\
\begin{tabular}{|c|c|c|}
3&2&1\\\hline\end{tabular}\\\begin{array}{l}\end{array}\end{array}$
& $f^{(4,3)}=\frac{7!}{5\ti4\ti3\ti2\ti3}=14$\\\hline $\qquad\qquad
S^{(4,2,1)}$ & $\begin{array}{lllll}\\
\begin{tabular}{|c|c|c|c|}
\hline 6&4&2&1\\\hline\end{tabular}\\
\begin{tabular}{|c|c|} 3&1\\\hline\end{tabular}\\
\begin{tabular}{|c|}
1\\\hline\end{tabular}\\\begin{array}{l}\end{array}\end{array}$ &
$f^{(4,2,1)}=\frac{7!}{6\ti4\ti3\ti2}=35$\\\hline $\qquad\qquad
S^{(4,1,1,1)}$ & $\begin{array}{llllll} \\
\begin{tabular}{|c|c|c|c|}
\hline 7&3&2&1\\\hline\end{tabular}\\
\begin{tabular}{|c|} 3\\\hline\end{tabular}\\
\begin{tabular}{|c|}
2\\\hline\end{tabular}\\
\begin{tabular}{|c|}
1\\\hline\end{tabular}\\\begin{array}{l}\end{array}\end{array}$ &
$f^{(4,1,1,1)}=\frac{7!}{7\ti6\ti3\ti2}=20$\\\hline $\qquad\qquad
S^{(3,3,1)}$ & $\begin{array}{lllll}\\ \begin{tabular}{|c|c|c|}
\hline 5 & 3 & 2\\ \hline\end{tabular}\\
\begin{tabular}{|c|c|c|} 4 & 2 & 1\\\hline\end{tabular}\\
\begin{tabular}{|c|} 1\\ \hline\end{tabular}\\
\begin{array}{l}\end{array}\end{array}$ & $f^{(3,3,1)}=\frac{7!}{5\ti4\ti3\ti2\ti2}=21$\\ \hline
\end{tabular}
\\\\\\
\noindent Finally, we note that Corollary $3$ is not true for $n=8$,
since $\frac{8\ti5}{2}=20$ and $\text{dim}\,S^{(4,4)}=14$.

Let's go back to the proof of the theorem. There are two parts in
the proof. First of all we need to show that for each of the values
$r$, $-r^3$, $\unsur{r^{n-3}}$, $-\unsur{r^{n-3}}$,
$\unsur{r^{2n-3}}$ of $l$, the representation $\n^{(n)}$ is
reducible. Second we need to show that if $\n^{(n)}$ is reducible,
then it forces one of the values $r$, $-r^3$, $\unsur{r^{n-3}}$,
$-\unsur{r^{n-3}}$, $\unsur{r^{2n-3}}$ for the parameter $l$. The
first point is achieved by exhibiting a non trivial proper invariant
subspace inside $\V$ as it was already done in the case
$l=\unsur{r^{2n-3}}$ (cf Theorem $4$) and in the cases
$l\in\lb\unsur{r^{n-3}},-\unsur{r^{n-3}}\rb$ (cf Theorems $5$ and
$6$). We leave the cases $l=r$ and $l=-r^3$ for later. For the
second part of the proof, we will use our preliminary remarks and
results of G.D. James. We proceed by induction on $n$ to show the
following property:
$$(\P_n):\nn\;\text{is reducible}\Rightarrow l\in\Big\lb r,-r^3,\unsur{r^{n-3}}
,-\unsur{r^{n-3}},\unsur{r^{2n-3}}\Big\rb$$ First we show that
$(\P_7)$ and $(\P_8)$ hold. Suppose $\n^{(7)}$ is reducible. Let
$\W$ be a proper non trivial invariant subspace of $\V$, with $\W$
irreducible. Then $\W$ is an irreducible $\IH_{F,r^2}(7)$-module. By
Corollary $3$, we have $\text{dim}\,\W\in\lb 1,6,14,15\}$ or
$\text{dim}\,\W>15$. If $\di\,\W=1$, Theorem $4$ implies that
$l=\unsur{r^{11}}$. Also, if $\di\,\W=6$, Theorem $5$ implies that
$l\in\lb\unsur{r^4},-\unsur{r^4}\rb$. Let $\V_0$ be the vector
subspace of $\V$ over $F$, spanned by the $w_{ik}$'s for $1\leq
i<k\leq 6$. The vector space $\V_0$ is a $B(A_5)$-module with the
action provided by the restriction of $\n^{(7)}$ to $B(A_5)$, but is
not a $B(A_6)$-module since for instance
$g_6.w_{5,6}=w_{5,7}\not\in\V_0$. In particular, $\V_0\neq\W$.
Further, if $\V_0\subseteq\W$, then $\xalone\in\W$. By
$\S\;\mathbf{8.1}$ this would imply that $\W=\V$, which is
impossible. Thus, considering the intersection of vector spaces
$\W\cap\V_0$, we get:
$$0\subseteq\W\cap\V_0\subset\V_0$$
\indent $\clubsuit$ Suppose $0=\W\cap\V_0$. Since
$\W\oplus\V_0\subseteq\V$, we must have
$\di_F\W+\di_F\V_0\leq\di_F\V$,
$$\emph{\/i.e}\qquad\qquad\di_F\W\leq\binom{7}{2}-\binom{6}{2}=6$$
\noindent In this situation, we get, as seen above
$l\in\big\lb\unsur{r^{11}},-\unsur{r^4},\unsur{r^4}\big\rb$.

$\clubsuit\clubsuit$ Otherwise, we have $0\subset\W\cap\V_0\subset\V_0$.\\
Since $\W$ is also a $B(A_5)$-module with an action given by the
restriction $\n^{(7)}\downarrow_{B(A_5)}$, the vector space
$\W\cap\V_0$ can be viewed as a $B(A_5)$-module. By forgetting the
last vertical line in the tangles representing the spanning elements
$w_{i,k}$'s, $1\leq i<k\leq 6$ of $\V_0$, the inclusions above mean
exactly that $\n^{(6)}$ is reducible. Then, by the case $n=6$, we
get $l\in\lb r,-r^3,\unsur{r^3},-\unsur{r^3},\unsur{r^9}\rb$. \\

We now use the same technique to still get more informations, but
this time with the $F$-vector space $\V_1:=\langle w_{i,k}|1\leq
i<k\leq 5\rangle_{_{F}}$. With the action given by
$\n^{(7)}\da_{_{B(A_4)}}$ the vector subspace $\V_1$ of $\V$ is a
$B(A_4)$-module. For the same reasons as above, it is impossible to
have $\V_1\subseteq\W$, hence we have:
$$0\subseteq\W\cap\V_1\subset\V_1$$
\noindent And there are again two cases:\\

$\spadesuit$ If $\W\cap\V_1=\lb0\rb$, then $\di\W\leq
\binom{7}{2}-\binom{5}{2}=11$. Then by the above,
$\di\W\in\lb1,6\rb$ and
$l\in\lb\unsur{r^4},-\unsur{r^4},\unsur{r^{11}}\}$

$\spadesuit\spadesuit$ Otherwise, $0\subset\W\cap\V_1\subset\V_1$.
By forgetting the last two vertical lines in the tangles
representing the spanning elements $w_{i,k}'s\;1\leq i<k\leq 5$ of
$\V_1$, we read on the inclusions above that $\n^{(5)}$ is
reducible. By the case $n=5$, it yields $l\in\lb
r,-r^3,\unsur{r^2},-\unsur{r^2},\unsur{r^7}\rb$.\\

Suppose now that
$l\not\in\lb\unsur{r^4},-\unsur{r^4},\unsur{r^7}\rb$. Then points
$\clubsuit$ and $\spadesuit$ cannot happen. By points
$\clubsuit\clubsuit$ and $\spadesuit\spadesuit$ we get:
$$\left\lb\begin{array}{lll}
l\in\{r,-r^3,\unsur{r^3},-\unsur{r^3},\unsur{r^9}\rb\\
\qquad\qquad\&\\
l\in\{r,-r^3,\unsur{r^2},-\unsur{r^2},\unsur{r^7}\rb
\end{array}\right.$$
\noindent As usual let $\e$ take the value $+1$ or $-1$. We have the
equivalences:
\begin{eqnarray}
r^2=\e\,r^3 &\Lra& \;r\,=\e\\
r^2=\e\,r^9 &\Lra& r^7=\e\\
r^7=\e\,r^3 &\Lra& r^4=\e\\
r^7=\text{ }\,\,r^9 &\Lra& r^2=1
\end{eqnarray}
\noindent Since $r^2\neq 1$ as $m\neq 0$ and since $(r^2)^7\neq 1$
and $(r^2)^4\neq 1$ by semisimplicity of $\IH_{F,r^2}(7)$, none of
$(143)-(146)$ can occur. Thus, $l$ must take the values $r$ or
$-r^3$.\\

We conclude that $l\in\lb
r,-r^3,\unsur{r^4},-\unsur{r^4},\unsur{r^{11}}\rb$, as expected.\\

\noindent Let's now prove $(\P_8)$. Taking the same notations as
before, suppose $\V$ is reducible and let $\W$ be an irreducible
submodule of $\V$. By $\S\;\mathbf{8.1}$ and the relations defining
the BMW algebra, $\W$ is also an irreducible
$\IH_{F,r^2}(8)$-module. If $\di\W=1$, then $l=\unsur{r^{13}}$ by
Theorem $4$ and if $\di\W=7$, then
$l\in\lb\unsur{r^5},-\unsur{r^5}\rb$ by Theorem $5$. Otherwise,
$\di\W\geq 8$ by Proposition $3$, a piece of the work of James in
\cite{GDJ}. Still taking the same notations as before, we define the
vector spaces $\V_0$ and $\V_1$:
$$\begin{array}{l}
\V_0=\big\langle w_{i,k}|1\leq i<k\leq 7\big\rangle_{_F}\\
\V_1=\big\langle w_{j,s}|1\leq j<s\leq 6\big\rangle_{_F}
\end{array}$$
If $\W\cap\V_0=\{0\}$, then it forces
$\di\W\leq\binom{8}{2}-\binom{7}{2}=7$, so that
$l\in\{\unsur{r^{13}},-\unsur{r^5},\unsur{r^5}\}$ by the above.
Otherwise, we have $0\subset\W\cap\V_0\subset\V_0$. If we forget
about the last node (number eight), we see with these inclusions
that $\n^{(7)}$ is reducible. Then by $(\P_7)$, it yields $l\in\lb
r,-r^3,\unsur{r^4},-\unsur{r^4},\unsur{r^{11}}\rb$.\\
If $\W\cap\V_1=\{0\}$, then it forces
$\di\W\leq\binom{8}{2}-\binom{6}{2}=13$. At this point, we need to
study wether or not there could be an irreducible
$\IH_{F,r^2}(8)$-module of dimension $d$ with $8\leq d\leq 13$. The
answer to that question is no as shown by the table below. We
referred to the Appendix table in \cite{GDJ}. James gives a
polynomial lower bound in $n$ for the dimension of an irreducible
$K\,Sym(n)$-module depending on the shape
$\bar{\mu}=(\mu_2,\mu_3,\dots)$ of the partition
$\mu=(\mu_1,\mu_2,\dots)\vdash n$ of $n$ and on the characteristic
of the field. In characteristic zero, these bounds are reached for
any $n$. We summarize the results in the case $n=8$ and $char(K)=0$
in the table below. In this case, a complete list of Specht modules
of dimension greater or equal to $8$ is:
$$S^{(6,2)},\,S^{(6,1,1)},\,S^{(5,3)},\,S^{(5,2,1)},\,
S^{(5,1,1,1)},\,S^{(4,4)},\,S^{(4,3,1)},\,S^{(4,2,2)},\,
S^{(4,2,1,1)},\,S^{(3,3,2)},$$ and their respective conjugates (when
these are not self-conjugate), as illustrated in the table of our
appendix. Hence we see that the different possibilities for
$\bar{\mu}$ are:
$$(2),\,(1^2),\,(3),\,(2,1),\,
(1^3),\,(4),\,(3,1),\,(2,1,1),\,(2,2),\,(3,2)$$ \noindent For each
of these shapes corresponding to the rows of height $\geq 2$ in the
Ferrers diagram, the dimension of the corresponding class of
irreducible $K\,Sym(8)$-module is:
\begin{center}
\begin{tabular}[c]{|c|c|c|}
\hline Specht module & $\bar{\mu}$ & dimension\\
\hline $S^{(6,2)}$ & $(2)$ & $20$\\
\hline $S^{(6,1,1)}$ & $(1^2)$ & $21$\\
\hline $S^{(5,3)}$ & $(3)$ & $28$\\
\hline $S^{(5,2,1)}$ & $(2,1)$ & $64$\\
\hline $S^{(5,1,1,1)}$ & $(1^3)$ & $35$\\
\hline $S^{(4,4)}$ & $(4)$ & $14$\\
\hline $S^{(4,3,1)}$ & $(3,1)$ & $70$\\
\hline $S^{(4,2,1,1)}$ & $(2,1^2)$ & $90$\\
\hline $S^{(4,2,2)}$ & $(2^2)$ & $56$\\
\hline $S^{(3,3,2)}$ & $(3,2)$ & $42$\\
\hline
\end{tabular}
\end{center}

The last row of the array was obtained by the Hook formula as the
appendix table of James gives results only for partitions
$\bar{\mu}\vdash m$ with $m\leq 4$. From the table, we gather the
information that the next smallest degree of an irreducible
representation of the Iwahori-Hecke algebra $\IH_{F,r^2}(8)$ is
$14$. Thus, $\di\W\leq 13$ implies in fact that $\di\W\in\lb
1,7\rb$, and this forces
$l\in\lb\unsur{r^{13}},\unsur{r^5},-\unsur{r^5}\rb$ by Theorem $4$
and Theorem $5$. On the other hand, if $\W\cap\V_1\neq\lb0\rb$, then
it comes $0\subset\W\cap\V_1\subset\V_1$, the last inclusion holding
for the same reasons as before. The $F$-vector space $\W\cap\V_1$ is
a $B(A_5)$-module with the action given by
$\n^{(8)}\da_{_{B(A_5)}}$. When the BMW algebra $B(A_5)$ acts on the
tangles representing the elements $w_{i,k}$'s, $1\leq i\leq6$ of
$\V_1$, it leaves their last two vertical lines invariant. Hence, we
may as well forget about them and we read on the inclusions above
that $\n^{(6)}$ is reducible. By the case $n=6$, it follows that
$l\in\lb r,-r^3,\unsur{r^3},-\unsur{r^3},\unsur{r^9}\rb$.\\
Suppose now that
$l\not\in\lb\unsur{r^13},\unsur{r^5},-\unsur{r^5}\rb$. Then, by the
previous considerations, we get:
$$\left\lb\begin{array}{lll}
l\in\{r,-r^3,\unsur{r^3},-\unsur{r^3},\unsur{r^9}\rb\\
\qquad\qquad\&\\
l\in\{r,-r^3,\unsur{r^4},-\unsur{r^4},\unsur{r^{11}}\rb
\end{array}\right.$$
\noindent Again, we cannot have:
$$\left\lb\begin{array}{ccc}
r^3&=&\e\,r^4\\
r^9&=&r^{11}
\end{array}\right.$$
\noindent because $r\not\in\lb 1,-1\rb$ as $m\neq 0$. Also, we
cannot have:
$$\left\lb\begin{array}{ccc}
r^3&=&\e\,r^{11}\\
r^9&=&\e\,r^4
\end{array}\right.$$
\noindent since $(r^2)^4\neq 1$ and $(r^2)^5\neq 1$ by
semisimplicity of $\IH_{F,r^2}(8)$. Then, the only possibilities
are: $l=r$ or $l=-r^3$. In summary, $l\in\lb
r,-r^3,\unsur{r^5},-\unsur{r^5},\unsur{r^{13}}\rb$. This achieves
the proof of $(\P_8)$.

We are now in a position to deal with the general case. Let $n$ be
any integer greater or equal to $9$ and suppose that $(\P_{n-2})$
and $(\P_{n-1})$ hold. Suppose that $\n^{(n)}$ is reducible and let
$\W$ be a proper nontrivial invariant subspace of $\V$. The
$F$-vector subspace $\W$ must satisfy to:
$$\left\lb\begin{array}{ll}
0\subset\W\subset\V\\ \forall 1\leq i\leq
n-1,\;\n_i(\W)\subseteq\W\qquad\qquad(\bigstar)
\end{array}
\right.$$ \noindent Without loss of generality, let's assume that
$\W$ is irreducible. We have the lemma:
\begin{lemma}
Let $\G$ be a $F$-vector subspace of $\V$. Then,\\
$\G$ is a $B(A_{n-1})$-module if and only if $\G$ is a
$\IH_{F,r^2}(n)$-module.
\end{lemma}
\noindent\textsc{Proof of the lemma:} Suppose first that $\G$ is a
$B(A_{n-1})$-module. By $\S\,\mathbf{8.1}$ and the defining relation
$(3)$ of the BMW algebra, the $B(A_{n-1})$-action on $\G$ is a
Iwahori-Hecke algebra action by $\IH_{F,r^2}(n)$ since $e_i.\G=0$
for all $1\leq i\leq n-1$. Conversely, if $\G$ is an
$\IH_{F,r^2}(n)$-module, it is obviously a $B(A_{n-1})$-module since
the expression defining $e_i$ is polynomial in $g_i$.\\

\noindent By the lemma, $\W$ is an irreducible
$\IH_{F,r^2}(n)$-module. Then by Corollary $3$, we have
$\text{dim}\,\W\in\lb
1,n-1,\frac{n(n-3)}{2},\frac{(n-1)(n-2)}{2}\rb$ or
$\text{dim}\,\W>\frac{(n-1)(n-2)}{2}$. $$\begin{array}{l} \text{If}
\;\text{dim}\,\W=1,\;\text{then}\;l=\unsur{r^{2n-3}}\;\text{by Theorem}\;4.\\
\text{If}\;\text{dim}\,\W=n-1,\;\text{then}\;
l\in\lb\unsur{r^{n-3}},-\unsur{r^{n-3}}\rb\;\text{by Theorem}
5.\end{array}$$ \noindent Otherwise,
$\text{dim}\,\W=\frac{n(n-3)}{2}$ or
$\text{dim}\,\W=\frac{(n-1)(n-2)}{2}$ or $\dw>\frac{(n-1)(n-2)}{2}$.
\unitlength=1in
%$$\begin{array}{cccccccccc}
%& \stackrel{l=\unsur{r^{2n-3}}}{\dna}& &\stackrel{l=\pm\unsur{r^{n-3}}}{\dna} & & & & & \stackrel{>}{\lra}&\\
%\hrf&|&\hrf&|&\hrf&|&\hrf&|&&\\
%&1&&n-1&&\frac{n(n-3)}{2}&&\frac{(n-1)(n-2)}{2}&&\dw
%\end{array}$$

$$\begin{array}{l}
\begin{picture}(.6,.01)\put(.6,.01){\line(-1,0){.6}}\end{picture}\!\!\!\stackrel{\stackrel{l=\unsur{r^{2n-3}}}{\dna}}{\begin{picture}(.25,0.01)\put(.25,.01){\line(-1,0){.25}}\end{picture}\!\!\!|\!\!\!\begin{picture}(.25,.01)\put(.25,.01){\line(-1,0){.25}}\end{picture}}
\!\!\!\begin{picture}(.5,.01)\put(.5,.01){\line(-1,0){.5}}\end{picture}\!\!\!\!\!\stackrel{\stackrel{l=\pm\unsur{r^{n-3}}}{\dna}}{\begin{picture}(.25,.01)\put(.25,.01){\line(-1,0){.25}}\end{picture}\!\!\!|\!\!\!\begin{picture}(.25,.01)\put(.25,.01){\line(-1,0){.25}}\end{picture}}
\!\!\!\!\!\begin{picture}(.85,.01)\put(.85,.01){\line(-1,0){.85}}\end{picture}\!\!\!|\!\!\!\begin{picture}(1,.01)\put(1,.01){\line(-1,0){1}}\end{picture}\!\!\!|
\!\!\!\stackrel{\stackrel{>}{\lra}}{\begin{picture}(.9,.01)\put(.9,.01){\line(-1,0){.9}}\end{picture}}\\
\hspace{.5in}\;\;\;\;\;\;1\hspace{.4in}\qquad\!\!
n-1\hspace{.65in}\nts\nts\!\frac{n(n-3)}{2}\hspace{.6in}\nts\nts\nts\nts\nts\!\frac{(n-1)(n-2)}{2}\hspace{.4in}\!\!\!\dw
\end{array}$$

\noindent As for $n=7$ and $n=8$, let's define the two $F$-vector
spaces:
$$\begin{array}{l}
\V_0=\big\langle w_{i,k}|1\leq i<k\leq n-1\big\rangle_{_F}\\
\V_1=\big\langle w_{j,s}|1\leq j<s\leq n-2\big\rangle_{_F}
\end{array}$$
\noindent of respective dimensions $\binom{n-1}{2}$ and
$\binom{n-2}{2}$. If $\W$ contains $\V_0$ or $\V_1$, then $\W$
contains $\xalone$ and by $\S\mathbf{8.1}$, all the $\xb$'s are in
fact in $\W$. Then $\W$ is the whole space $\V$, which is
impossible. Thus we have the inclusions of vector spaces:
$$\begin{array}{l}
0\subseteq\W\cap\V_0\subset\V_0\\
0\subseteq\W\cap\V_0\subset\V_1
\end{array}$$
\noindent If $\W\cap\V_0=\{0\}$, then $\W\bigoplus\V_0$ and we have:
\begin{eqnarray}
\notag \dw+\text{dim}\,\V_0&\leq&\text{dim}\,\V\\ \notag
\emph{\/i.e}\qquad\qquad\qquad \dw&\leq&\binom{n}{2}-\binom{n-1}{2}\\
\emph{\/i.e}\qquad\qquad\qquad \dw&\leq& n-1
\end{eqnarray}
\noindent Similarly, if $\W\cap\V_1=\{0\}$, then $\W\bigoplus\V_1$
and we get:
\begin{eqnarray}
\notag \dw+\text{dim}\,V_1&\leq&\text{dim}\,\V\\\notag
\emph{\/i.e}\qquad\qquad\qquad \dw&\leq&\binom{n}{2}-\binom{n-2}{2}\\
\emph{\/i.e}\qquad\qquad\qquad \dw&\leq& 2n-3
\end{eqnarray}
\noindent Since the inequality $2n-3<\frac{n(n-3)}{2}$ holds as soon
as $n\geq 7$, so in particular for any $n\geq 9$, we see in each
case and as illustrated by the figure, that there are only two
possibilities for the dimension of $\W$: either $\dw=1$ or
$\dw=n-1$. Thus, each of the conditions $(147)$ and $(148)$ implies
that $l\in\lb\unsur{r^{n-3}},-\unsur{r^{n-3}},\unsur{r^{2n-3}}\rb$.
Suppose now that $l$ does not take any of these values. Then we have
$\W\cap\V_0\neq\{0\}$ and $\W\cap\V_1\neq\{0\}$, so that:

$$\begin{array}{l}
0\subset\W\cap\V_0\subset\V_0\\
0\subset\W\cap\V_0\subset\V_1
\end{array}$$

\noindent By definition of the representation, $\n^{(n-1)}$ is a
representation of $B(A_{n-2})$ in $\V_0$. Still by definition of the
representation, for each $i=1,\dots,n-2$ and for all basis vector
$\xb$ in $\V_0$, the vectors $\n^{(n-1)}(g_i)(\xb)$ and
$\n^{(n)}(g_i)(\xb)$ only depend on the inner product $\ps$ and we
have:
\begin{equation*}
\n^{(n-1)}(g_i)(\xb)=\n^{(n)}(g_i)(\xb)=\n_i(\xb)
\end{equation*}
\noindent Let $z\in\W\cap\V_0$. Then $\forall \,1\leq i\leq
n-2,\;\n^{(n-1)}(g_i)(z)=\n^{(n)}(g_i)(z)\in\W$ since $\W$ is a
$B(A_{n-1})$-submodule of $\V$. Thus,
$$\forall\;i=1,\dots,n-2,\,\n^{(n-1)}(g_i)(\W\cap\V_0)\subseteq\W\cap\V_0$$
\noindent Then $\n^{(n-1)}$ is reducible and $(\P_{n-1})$ implies
that $l\in\lb
r,-r^3,\unsur{r^{n-4}},-\unsur{r^{n-4}},\unsur{r^{2n-5}}\rb$. \\
Similarly, $\n^{(n-2)}$ is a representation of $B(A_{n-3})$ in
$\V_1$ and we have:
$$\forall 1\leq i\leq n-3,\,\forall
z\in\W\cap\V_1,\,\n^{(n-2)}(g_i)(z)=\n_i(z)\in\W\;\text{by}\;(\bigstar)$$
\noindent Thus, we have: $$\forall
i=1,\dots,n-3,\,\n^{(n-2)}(g_i)(\W\cap\V_1)\subseteq\W\cap\V_1,$$
\noindent so that the representation $\n^{(n-2)}$ is reducible. Then
$(\P_{n-2})$ implies that $$l\in\bigg\lb
r,-r^3,\unsur{r^{n-5}},-\unsur{r^{n-5}},\unsur{r^{2n-7}}\bigg\rb$$
\noin In summary, if
$l\not\in\lb\unsur{r^{n-3}},-\unsur{r^{n-3}},\unsur{r^{2n-3}}\rb$,
then: $$l\in\bigg\lb
r,-r^3,\unsur{r^{n-4}},-\unsur{r^{n-4}},\unsur{r^{2n-5}}\bigg\rb\;\text{
and}\;l\in\bigg\lb
r,-r^3,\unsur{r^{n-5}},-\unsur{r^{n-5}},\unsur{r^{2n-7}}\bigg\rb$$
\noin At this stage, we recall that we made the assumption that
$\IH_{F,r^2}(n)$ is semisimple. By Corollary $3.44$ in \cite{IHA},
the smallest integer $e$ such that $$1+r^2+\dots+(r^2)^{e-1}=0$$
($e=\infty$ if no such integer exists) must be greater than $n$.
Since $m\neq 0$ implies $r^2\neq 1$, this condition is in fact
equivalent to: "the smallest integer $e$, if it exists, such that
$(r^2)^e=1$ is greater than $n$". Thus, we must have: $$(r^2)^k\neq
1,\,\forall k=1,\dots,n$$ We use this constraint on $r$ to show that
the two conditions on $l$ above force $l\in\lb r,-r^3\rb$. Indeed,
the condition
$$\begin{array}{lll}
l\in\{r,-r^3,\unsur{r^{n-4}},-\unsur{r^{n-4}},\unsur{r^{2n-5}}\rb\\
\qquad\qquad\&\\
l\in\{r,-r^3,\unsur{r^{n-5}},-\unsur{r^{n-5}},\unsur{r^{2n-7}}\rb
\end{array}$$
and the impossibility to have any of the following equalities below
$$\begin{array}{cccccccccc}
\frac{\e}{r^{n-4}}&=&\frac{\e^{'}}{r^{n-5}}&\Lra& r\in\lb
-1,1\rb &\text{impossible as}& m&\neq& 0\\
\frac{\e}{r^{n-4}}&=&\unsur{r^{2n-7}}&\Lra&
r^{n-3}=\e &\text{impossible as}& (r^2)^{n-3}&\neq& 1\\
\unsur{r^{2n-5}}&=&\frac{\e}{r^{n-5}}&\Lra& r^n=\e&\text{impossible
as}& (r^2)^n&\neq& 1\\
\unsur{r^{2n-5}}&=&\unsur{r^{2n-7}}&\Lra& r^2=1&\text{impossible as}
& m&\neq& 0
\end{array}$$

\noin imply that $l\in\lb r,-r^3\rb$. We conclude that if $\n^{(n)}$
is reducible, then $l\in\lb
r,-r^3,\unsur{r^{n-3}},-\unsur{r^{n-3}},\unsur{r^{2n-3}}\rb$. Hence,
$$(\P_{n-2})\&(\P_{n-1})\Longrightarrow (\P_n)$$
\noin Since $(\P_7)\&(\P_8)$ hold, $(\P_n)$ holds for all $n\geq 9$.

To complete the proof of the Theorem, it remains to show that for
$l=r$ and $l=r^3$, the representation $\n^{(n)}$ is reducible. We
know from before that $\cap_{1\leq i<j\leq n}Ker\,\n(X_{ij})$ is a
proper invariant subspace of $\V$. If we can show that for $l=r$ and
for $l=r^3$, this subspace is non zero, then $\n^{(n)}$ is reducible
for each of these values of $l$. For $n=4,5$ and for each value of
$l=r,-r^3$ we exhibited a non zero vector belonging to $\cap_{1\leq
i<j\leq n}Ker\,\n(X_{ij})$, which showed the reducibility of
$\n^{(n)}$ in these cases. For bigger values of $n$, we have the
result:
\begin{Prop} Let $n$ be an integer with $n\geq 5$.
$$\begin{array}{ccccccc}
\text{If}& l&=&r,&
\X:=r^2\,w_{12}-r\,w_{13}+w_{34}-r\,w_{24}&\in&\underset{1\leq i<j\leq n}{\bigcap}Ker\,\n(X_{ij})\\
\text{If}& l&=&-r^3,& \qquad\;\;\,
\Y:=-r\,w_{23}-\unsurr\,w_{34}+w_{24}&\in&\underset{1\leq i<j\leq
n}{\bigcap}Ker\,\n(X_{ij})
\end{array}$$
\end{Prop}
\noin \textsc{Proof of the Proposition:} for $n=5$, the vectors of
the proposition are the ones annihilating the matrix $S$ of
$\S\,\mathbf{5}$. By construction, $S$ is the matrix of the
endomorphism $\sum_{1\leq i<j\leq n}\n(X_{ij})$ in the basis $\B$.
By Proposition $2$ of $\S\,\mathbf{8.1}$, the matrix
$Mat_{\B}\n(X_{ij})$ is the matrix whose
$\{(1+2+\dots+j-2)+(j-i)\}$-th row is the one of
$Mat_{\B}\sum_{1\leq i<j\leq n}\n(X_{ij})$ and with zeros elsewhere.
Hence Proposition $5$ holds for $n=5$. Moreover, to show in the
general case that $\X,\Y\in\cap_{1\leq i<j\leq n}Ker\,\n(X_{ij})$,
it suffices to show that each row of the sum matrix annihilates
$\X,\Y$. Further, since the last $\binom{n}{2}-5$ coordinates of the
vectors $\X$ and $\Y$ are zero, it suffices to check that the five
first coordinates of each row of the sum matrix annihilate the
column vectors $(r^2,0,-r,1,-r)$ and $(0,-r,0,-\unsurr,1)$.
Furthermore, by the remark above, the first $\binom{5}{2}$ rows of
the sum matrix correspond to the only nonzero row of the matrices
$Mat_{\B}\n(X_{i,j})$ for $1\leq i<j\leq 5$. Since an expression for
$\n(X_{ij})$ is:
$$\n(X_{ij})=\n_{j-1}\dots\n_{i+1}\n(e_i)\n_{i+1}^{-1}\dots\n_{j-1}^{-1},$$
and since $\n_i^{-1}(\xb)$ only depends on the inner product
$(\al_i|\be)$ as in the expression of $\n_i^{-1}$ below:
$$\nts\nts\nts\nts\nts\nts\nts\nts\nts\nts
\n_i^{-1}(\xb)=\left\lbrace\begin{array}{ccccccccc}
\unsurr\,\xb&\text{if}&\ps&=&0&&&&\\
l\,\xb&\text{if}&\ps&=&1&&&&\\
x_{\be-\al_i}-ml\,r^{ht(\be)-2}\,\xali+m\,\xb&\text{if}&\ps&=&\unsur{2}&\&
&\be-\al_i&\succ &\al_i\\
x_{\be-\al_i}&\text{if}&\ps&=&\unsur{2}&\& &\be-\ali&\p&\ali\\
x_{\be+\ali}&\text{if}&\ps&=&-\unsur{2}&\&&\be&\succ &\ali\\
x_{\be+\ali}-\frac{m}{r^{ht(\be)-1}}\,\xali+m\,\xb&\text{if}&\ps&=&-\unsur{2}&\&&\be&\p&\ali
\end{array}\right.,$$ the action of $\n^{(n)}(X_{ij})$, $j\leq 5$ on
$\xalone,x_{\al_2},x_{\al_1+\al2},x_{\al_3},x_{\al_3+\al_2}$ is the
same as the action of $\n^{(5)}(X_{ij})$, $j\leq 5$ on the same
vectors. Thus, we only need to check that the five first coordinates
of the last $\binom{n}{2}-\binom{5}{2}$ rows of the sum matrix
annihilates the column vectors $(r^2,0,-r,1,-r)$ and
$(0,-r,0,-\unsurr,1)$. In other words, we need to check that for all
the BMW algebra elements $X_{ij}$'s with $1\leq i<j\leq n$ and
$j\geq
6$:\\\\
If $l=r$,
\begin{equation*}
\nts\nts\nts\nts\nts\nts\nts\nts\nts
r^2\,[\n(X_{ij})(\xalone)]_{_{w_{ij}}}-r\,[\n(X_{ij})(x_{\al_1+\al_2})]_{_{w_{ij}}}+[\n(X_{ij})(x_{\al_3})]_{_{w_{ij}}}-r\,[\n(X_{ij})(x_{\al_3+\al_2})]_{_{w_{ij}}}=0\;\!\;\;(\divideontimes)
\end{equation*}
If $l=-r^3$,
\begin{equation*}
-r\,[\n(X_{ij})(x_{\al_2})]_{_{w_{ij}}}-\unsurr\,[\n(X_{ij})(x_{\al_3})]_{_{w_{ij}}}+[\n(X_{ij})(x_{\al_2+\al_3})]_{_{w_{ij}}}=0\;\qquad\qquad\qquad\;\;\;\;(\divideontimes\divideontimes)
\end{equation*}
For a proof of these two equalities, we shall use the table of the
appendix that gives a complete description of how the $X_{ij}$'s act
on the $w_{sk}$'s or the straightforward computations below that use
the conjugation formulas for the $X_{ij}$'s. First we have the
lemma:
\begin{lemma}\hfill\\
For all $1\leq i<j-1\leq n$ with $j\geq 6$ and $i\geq 5$, we have:
$$\forall\be\in\lb\al_1,\al_2,\al_2+\al_1,\al_3,\al_3+\al_2\rb,\,\n(X_{ij})(\xb)=\n_{j-1}\dots\n_{i+1}\n(e_i)\n_{i+1}^{-1}\dots\n_{j-1}^{-1}(\xb)=0$$
For all $i\geq 5$, we have:
$$\nts\nts\nts\nts\nts\nts\forall\be\in\lb\al_1,\al_2,\al_2+\al_1,\al_3,\al_3+\al_2\rb,\,\n(X_{i,i+1})(\xb)=\n(e_i)(\xb)=0$$
\end{lemma}
\noin\textsc{Proof of the lemma:} for $1\leq i<j-1\leq n$ and $j\geq
6$ and $i\geq 5$, it suffices to notice that:
$$\text{Supp}(\be)\cap\lb i-1,i,i+1,\dots,j-1,j\rb=\emptyset$$
\noin Hence, the vector
$\n(e_i)\n_{i+1}^{-1}\dots\n_{j-1}^{-1}(\xb)$ is zero. As for the
second equality, it simply comes from the fact that $\ps=0$. \hfill
$\square$\\\\ In what follows, we shall denote by $\phi_1^{+}$ the
subset of $\phi^{+}$ composed of the positive roots
$\al_1,\al_2,\al_2+\al_1,\al_3,\al_3+\al_2$.\\\\ It remains to show
$(\divideontimes)$ and $(\divideontimes\divideontimes)$ for $i\in\lb
1,2,3,4\rb$ and $j\geq 6$. Let's start with $i=4$. We have for all
$\be\in\phi_1^{+}$:
\begin{eqnarray*}
\n(X_{4j})(\xb)&=&\n_{j-1}\dots\n_5\n(e_4)\n_5^{-1}\dots\n_{j-1}^{-1}(\xb)\\
               &=&\unsur{r^{j-5}}\,\n_{j-1}\dots\n_5\n(e_4)(\xb)\\
               &=&\left\lb\begin{array}{ccccc}
               0&\text{if}&\be&\in&\lb\al_1,\al_2,\al_1+\al_2\rb\\
               \!\!\!\!\!\!\unsur{r^{j-5}}\n_{j-1}\dots\n_5(x_{\al_4})&\text{if}&\be&=&\nts\nts\nts\nts\nts\nts\nts\nts\nts\nts\nts\al_3\\
               \unsur{r^{j-5}}\n_{j-1}\dots\n_5(x_{\al_4}).\unsurr&\text{if}&\be&=&\nts\nts\nts\nts\nts\nts\nts\nts\nts\nts\nts\al_3+\al_2\end{array}\right.
\end{eqnarray*}
\noin By the game of the coefficients in $(\divideontimes)$ and
$(\divideontimes\divideontimes)$, we conclude that
$(\divideontimes)$ and $(\divideontimes\divideontimes)$ hold when
$(i,j)=(4,j)$ with $j\geq 6$.\\
Next, we have for all $\be\in\phi_1^{+}$ and all $j\geq 6$:
\begin{eqnarray*}
\n(X_{3j})(\xb)&=&\n_{j-1}\dots\n_4\n(e_3)\n_4^{-1}\dots\n_{j-1}^{-1}(\xb)\\
               &=&\unsur{r^{j-5}}\,\n_{j-1}\dots\n_4\n(e_3)\n_4^{-1}(\xb)\\
               &=&\left\lb\begin{array}{cccccc}
               0&&\text{if}&\be&=&\al_1\\
               \unsur{r^{j-4}}&\n_{j-1}\dots\n_4(x_{\al_3})&\text{if}&\be&=&\al_2\\
               \unsur{r^{j-3}}&\n_{j-1}\dots\n_4(x_{\al_3})&\text{if}&\be&=&\al_2+\al_1\\
               \nts\nts\nts\nts\nts\nts\nts\nts\nts\nts\nts\nts\nts\nts\,\unsur{r^{j-5}}\unsur{l}&\n_{j-1}\dots\n_4(x_{\al_3})&\text{if}&\be&=&\al_3\\
               \unsur{r^{j-5}}\,m\,(\unsur{l}-\unsur{r})&\n_{j-1}\dots\n_4(x_{\al_3})&\text{if}&\be&=&\al_3+\al_2\\
               \end{array}\right.
\end{eqnarray*}
\noin When $l=r$, the last term above is zero. Thus, the only
nonzero terms in $(\divon)$ are
$-r\,\unsur{r^{j-3}}\n_{j-1}\dots\n_4(x_{\al_3})$ and
$\unsur{r^{j-5}}\unsurr\n_{j-1}\dots\n_4(x_{\al_3})$ and they cancel
each other. As for $(\divon\divon)$ when $l=-r^3$, all the powers in
$\unsurr$ cancel each other, as in the equation:
$$-r\,\unsur{r^{j-4}}+\unsurr\,\unsur{r^{j-5}}\,\unsur{r^3}+\unsur{r^{j-5}}\Big(\unsurr-r\Big)\Big(-\unsur{r^3}-\unsurr\Big)=0$$
\noin More computations show that for all $\be\in\phi_1^{+}$ and for
all $j\geq 6$,
\begin{eqnarray*}
\n(X_{2j})(\xb)&=&\n_{j-1}\dots\n_3\n(e_2)\n_3^{-1}\dots\n_{j-1}^{-1}(\xb)\\
               &=&\unsur{r^{j-5}}\,\n_{j-1}\dots\n_3\n(e_2)\n_3^{-1}\n_4^{-1}(\xb)\\
               &=&\left\lb\begin{array}{cccccc}
               \unsur{r^{j-3}}&\n_{j-1}\dots\n_3(x_{\al_2})&\text{if}&\be&=&\al_1\\
               \unsur{r^{j-4}}.\unsur{l}&\n_{j-1}\dots\n_3(x_{\al_2})&\text{if}&\be&=&\al_2\\
               \unsur{r^{j-4}}\,m\big(\unsur{l}-\unsurr\big)&\n_{j-1}\dots\n_3(x_{\al_2})&\text{if}&\be&=&\al_2+\al_1\\
               0&&\text{if}&\be&=&\al_3\\
               \unsur{r^{j-5}}.\unsur{l}&\n_{j-1}\dots\n_3(x_{\al_2})&\text{if}&\be&=&\al_3+\al_2\\
               \end{array}\right.
\end{eqnarray*}
\noin The zero comes from the fact that:
\begin{equation}\n_3^{-1}\n_4^{-1}(x_{\al_3})=x_{\al_4}\end{equation} Moreover,
if $l=r$ the third row is also zero. Since
$$r^2\,\unsur{r^{j-3}}-r\,\unsur{r^{j-5}}\unsurr=0,$$ we see that the
first and the last term on the left hand side of the equality
$(\divon)$ cancel each other. Also, the fact that:
$$r\,\unsur{r^{j-4}}\unsur{r^3}-\unsur{r^{j-5}}\unsur{r^3}=0$$
and the computations above imply that $(\divon\divon)$ holds.\\
Let's now study the action of $\n(e_1)\n_2^{-1}\n_3^{-1}\n_4^{-1}$
on an $\xb$ with $\be\in\phi_1^{+}$.\\ We have:
\begin{eqnarray*}
\n(e_1)\n_2^{-1}(\xalone)&=&\n(e_1)(x_{\al_2+\al_1}-m\,x_{\al_2}+m\,\xalone)\\
                         &=&\unsur{l}\,\xalone\\
\n(e_1)\n_2^{-1}\n_3^{-1}(x_{\al_2})&=&\n(e_1)\n_2^{-1}(x_{\al_2+\al_3}-m\,x_{\al_3}+m\,x_{\al_2})\\
                                    &=&\n(e_1)(x_{\al_3})\\
                                    &=&0\\
\n(e_1)\n_2^{-1}\n_3^{-1}(x_{\al_2+\al_1})&=&\n(e_1)\n_2^{-1}(x_{\al_3+\al_2+\al_1}-\frac{m}{r}\,x_{\al_3}+m\,x_{\al_2+\al_1})\\
&=&\n(e_1)(\unsurr\,x_{\al_3+\al_2+\al_1}-\frac{m}{r}\,x_{\al_2+\al_3}+m\,\xalone)\\
&=&\unsur{l}\,\xalone\\
\n(e_1)\n_2^{-1}\n_3^{-1}\n_4^{-1}(x_{\al_3})&=&\n(e_1)\n_2^{-1}(x_{\al_4})\;\;\text{by $(149)$}\\
&=&\unsurr\,\n(e_1)(x_{\al_4})\\
&=&0\\
\n(e_1)\n_2^{-1}\n_3^{-1}\n_4^{-1}(x_{\al_3+\al_2})&=&\n(e_1)\n_2^{-1}\n_3^{-1}(x_{\al_4+\al_3+\al_2}-\frac{m}{r}\,x_{\al_4}+m\,x_{\al_3+\al_2})\\
&=&\n(e_1)\n_2^{-1}(\unsurr\,x_{\al_4+\al_3+\al_2}-\frac{m}{r}x_{\al_3+\al_4}+m\,x_{\al_2})\\
&=&\unsurr\,\n(e_1)(x_{\al_3+\al_4})\\
&=&0
\end{eqnarray*}
\noin Therefore, we have for all $\be\in\phi_1^{+}$ and all $j\geq
6$:
\begin{eqnarray*}
\n(X_{1j})(\xb)&=&\n_{j-1}\dots\n_2\n(e_1)\n_2^{-1}\dots\n_{j-1}^{-1}(\xb)\\
               &=&\unsur{r^{j-5}}\,\n_{j-1}\dots\n_2\n(e_1)\n_2^{-1}\n_3^{-1}\n_4^{-1}(\xb)\\
               &=&\left\lb\begin{array}{cccccc}
               \unsur{r^{j-5}}\,\unsur{r^2}\,\unsur{l}&\n_{j-1}\dots\n_2(x_{\al_1})&\text{if}&\be&=&\al_1\\
               0&&\text{if}&\be&=&\al_2\\
               \unsur{r^{j-5}}\,\unsurr\,\unsur{l}&\n_{j-1}\dots\n_2(x_{\al_1})&\text{if}&\be&=&\al_2+\al_1\\
               0&&\text{if}&\be&=&\al_3\\
               0&&\text{if}&\be&=&\al_3+\al_2\\
               \end{array}\right.
\end{eqnarray*}
\noin Then, we see that all the terms of the left hand side of
equality $(\divon\divon)$ are zero. And we read on the first
expression and on the third expression above that
$$r^2\,[\n(X_{1,j})(x_{\al_1})]_{_{w_{1j}}}-r\,[\n(X_{1,j})(x_{\al_1+\al_2})]_{_{w_{1j}}}=0$$
Then $(\divon)$ also holds. This achieves the proof of Proposition
$5$ and the conjecture is thus proven when the Iwahori-Hecke algebra
$\IH_{F,r^2}(n)$ is semisimple.

\section{More properties of the representation when $l\in\lb
r,-r^3\rb$}

In this part we still assume that the Iwahori-Hecke algebra
$\IH_{F,r^2}(n)$  for the adequate integer $n$ is semisimple and we
show more properties of the representation $\n^{(n)}$. We focus on
the cases when the representation is reducible with $l\in\lb
r,-r^3\rb$. We introduce new notations: we will denote by $K(n)$ the
intersection module
$$K(n):=\bigcap_{1\leq i<j\leq n}\,Ker\,\n^{(n)}(X_{ij}),$$ and by $k(n)$ its
dimension as a vector space over $F$. We first study some properties
of reducibility or irreducibility of $K(n)$ depending on the values
of $l$ and $r$.
\subsection{Proprieties of reducibility/irreducibility of the
$K(n)$'s} We have the Proposition:
\begin{Prop}\hfill\\\\
Let $n$ be an integer with $n\geq 4$. \\
When $l=r$, $K(n)$ is always irreducible.\\\\
Let $n$ be an integer with $n\geq 3$.\\
When $l=-r^3$ and $n\neq 8$, there are two cases:
$$\begin{cases}
r^{2n}\neq -1 & \text{and}\;\, K(n) \;\,\text{is irreducible}\\
r^{2n}=-1 & K(n)\;\, \text{is reducible and}\;\, k(n)\geq\dbw
\end{cases}$$
(Case n=8) When $l=-r^3$, there are two cases:\\\\
$\;\;\;\text{                }\begin{cases}
r^{16}\neq -1 & \text{and}\;\, K(8) \;\,\text{is irreducible}\\
r^{16}=-1 & K(8)\;\, \text{is reducible and}\;\, k(8)\in\lb
15,21,22\rb
\end{cases}$\\
\end{Prop}

\textsc{Proof:}
\begin{list}{\texttt{$\star$}}{}
\item Let's first assume that $n\geq 5$ and $n\neq 8$. Suppose
$l=r$. The case $l=-r^3$ will be treated separately later on. We
know that for $n\geq 5$ and $n\neq 8$, the irreducible
$\IH_{F,r^2}(n)$-modules have dimension
$1,\,n-1,\,\frac{n(n-3)}{2},\frac{(n-1)(n-2)}{2}$ or dimension
greater than $\frac{(n-1)(n-2)}{2}$. Suppose $K(n)$ is reducible. By
semisimplicity of $\IH_{F,r^2}(n)$, the $\IH_{F,r^2}(n)-$ module
$K(n)$ would decompose as a direct sum $K_1(n)\ds K_2(n)$ with
$K_1(n)$ irreducible. With the same notations as before, let's write
$k_1(n)$ for the dimension of $K_1(n)$ over $F$ and $k_2(n)$ for the
dimension of $K_2(n)$ over $F$. Suppose now $k_1(n)\geq
\frac{(n-1)(n-2)}{2}$. Then it comes:
$$k_2(n)<\frac{n(n-1)}{2}-\frac{(n-1)(n-2)}{2}=n-1$$ By Theorem $4$,
this would force $l=\unsur{r^{2n-3}}$, which contradicts our
assumption $l=r$. Furthermore by theorem $5$, it is impossible to
have $k_1(n)=n-1$ and by theorem $4$, it is also impossible to have
$k_1(n)=1$. Then, the only remaining possibility for $k_1(n)$ is
$k_1(n)=\frac{n(n-3)}{2}$. It follows that:
$$k_2(n)+\frac{n(n-3)}{2}<\frac{n(n-1)}{2},$$
$$\emph{i.e\/}\qquad\qquad\qquad k_2(n)<n\leq\cil$$ This last set of inequalities now
forces $k_2(n)\in\lb 1,n-1\rb$, which is again a contradiction.
Hence $K(n)$ is irreducible in this case.
\item Let's now deal with the case $n=3$. Here we assume that $l=-r^3$. We have $k(3)\in\lb
1,2\rb$. \begin{itemize} \item[*] If $r^6\neq -1$. \\\\Suppose
$k(3)=2$. Since $l=-r^3$, we know that there exists a
one-dimensional invariant subspace of $\V$. Theorem $4$ also states
that this one-dimensional invariant subspace is unique with our
assumption on $r$. Moreover, it is contained in $K(3)$. Then it must
have a one-dimensional summand, which is impossible by uniqueness of
the one-dimensional invariant subspace in that case. Hence we must
have $k(3)=1$. Then $K(3)$ is irreducible and as a matter of fact,
$K(3)$ is the one-dimensional invariant subspace of $\V$.\\
\item[*] If $r^6=-1$. \\\\\
Then, by Theorem $4$, there exists two one-dimensional invariant
subspaces of $\V$ whose sum is direct. Since these spaces are
distinct, $k(3)$ cannot equal $1$. Hence $k(3)=2$, $K(3)$ is
reducible and $K(3)$ is a direct sum of these two one-dimensional
invariant subspaces.
\end{itemize}
\item The case $n=4$.\\\\
First we prove the following Theorem:\\
\begin{t-l}\hfill\\\\
Suppose $\n^{(4)}$ is reducible and let $\W$ be an irreducible
submodule of $\V$. \\If $\text{dim}\,\W=2$, then $l=r$ and $\W$ is
spanned over $F$ by the vectors:
$$\begin{array}{ll}
(w_{13}-\unsurr\,w_{23})-\unsurr\,(w_{14}-\unsurr\,w_{24})\\
(w_{12}-\unsurr\,w_{13})-\unsurr\,(w_{24}-\unsurr\,w_{34})
\end{array}$$
Conversely, if $l=r$, the two linearly independent vectors above are
stable under the action of $g_1,g_2,g_3$.
\end{t-l}
\textsc{Proof:} By lemma $6$, $\W$ is an irreducible $2$-dimensional
$\IH_{F,r^2}(4)$-module. Let's set $H_1$, $H_2$ and $H_3$ to be the
matrices:
$$ H_1:=\begin{bmatrix}
-\unsurr & 1\\
0 & r\\
\end{bmatrix},\; H_2:=\begin{bmatrix}
r &0\\
1&-\unsurr\\
\end{bmatrix},\; H_3:=\begin{bmatrix}
-\unsurr & 1\\
0 & r\\ \end{bmatrix}$$

By the symmetry of the roles played by the scalars $r$ and
$-\unsurr$ and the choice $H_3=H_1$, it suffices to check that:
$$\begin{array}{ll}
H_1^2+m\,H_1=I\\
H_1H_2H_1=H_2H_1H_2
\end{array},$$
(where $I$ is the identity matrix of size $2$) to see that these
matrices define a matrix representation of degree $2$ of the
Iwahori-Hecke algebra $\IH_{F,r^2}(4)$. Furthermore, this
representation is irreducible. And indeed, if the representation
were reducible, there would exist a nonzero vector $v=(v_1,v_2)$ of
$F^2$ and scalars $\la_1$ and $\la_2$ such that:
$$\left\lb\begin{array}{ccc}
H_1.v &=& \la_1\,v\\
H_2.v &=& \la_2\,v \end{array}\right.$$ By using the definitions for
$H_1$, we get from the first equation:
$$\left\lb\begin{array}{ccc}
-\unsurr\,v_1+v_2&=&\la_1\,v_1\\
r\,v_2&=&\la_1\,v_2\\
\end{array}\right. ,$$
from which we derive the new set of equations:
$$\left\lb\begin{array}{ccc}
(r+\unsurr)\,v_1&=&v_2\\
\la_1&=&r\\
\end{array}\right.\qquad\text{or}\qquad\left\lb\begin{array}{ccc}
v_2&=&0\\
\la_1&=&-\unsurr
\end{array}\right.$$

Since it is visible on the matrices that in $H_2.v=\la_2\,v$, the
scalar $v_1$ has been replaced by the scalar $v_2$, we also get the
set of equations:

$$\left\lb\begin{array}{ccc}
(r+\unsurr)\,v_2&=&v_1\\
\la_2&=&r\\
\end{array}\right.\qquad\text{or}\qquad\left\lb\begin{array}{ccc}
v_1&=&0\\
\la_2&=&-\unsurr
\end{array}\right.$$ where $v_1$ has been replaced by $v_2$ and $\la_1$ by
$\la_2$.

From there, we see that it is impossible to have $v_2=0$ or $v_1=0$
because since $(r^2)^2\neq 1$, this would yield $v_1=v_2=0$. Hence
the braces on the right hand side are to be excluded. Now we get:
$$v_1=\bigp r+\unsurr\bigpd^2\,v_1\;\;\text{and}\;\;
v_1\neq 0,$$ so that $\bigp r+\unsurr\bigpd^2=1$,
$$\emph{i.e\/}\qquad r+\unsurr=1\;\;\text{or}\;\;r+\unsurr=-1$$
Let's solve the quadratics $r^2-r+1=0$ and $r^2+r+1=0$. Both have
discriminant $(i\sqrt{3})^2$. Thus we get the solutions:
$$r\in\bigg\lb\unsur{2}+i\frac{\sqrt{3}}{2},\unsur{2}-i\frac{\sqrt{3}}{2},-\unsur{2}+i\frac{\sqrt{3}}{2},-\unsur{2}-i\frac{\sqrt{3}}{2}\bigg\rb$$
Equivalently, $$r\in\lb
e^{i\frac{\pi}{3}},e^{-i\frac{\pi}{3}},e^{2i\frac{\pi}{3}},e^{4i\frac{\pi}{3}}\rb$$
Then it comes:
$$r^2\in\lb
e^{2i\frac{\pi}{3}},e^{-2i\frac{\pi}{3}},e^{4i\frac{\pi}{3}},e^{8i\frac{\pi}{3}}\rb$$
Now it is visible that $(r^2)^3$=1, which is forbidden. Thus, the so
given representation of $\ih(4)$ of degree $2$ is irreducible. Then
there exists a basis $(v_1,v_2)$ in $\W$ such that:
$$\begin{array}{ccccc}
\left\lb\begin{array}{l} \n_1\,v_1=-\unsurr\,v_1\;\;\;\;\;\,(\divon)_1\\
\n_1\,v_2=v_1+r\,v_2\;(\divon)_2
\end{array}\right. &, &\left\lb\begin{array}{l} \n_2\,v_1=r\,v_1+v_2\;(\divon)_3\\
\n_2\,v_2=-\unsurr\,v_2\;\;\;\;\;\,(\divon)_4
\end{array}\right. &, &\left\lb\begin{array}{l} \n_3\,v_1=-\unsurr\,v_1\;\;\;\;\;\,(\divon)_5\\
\n_3\,v_2=v_1+r\,v_2\;(\divon)_6
\end{array}\right.
\end{array}$$
By $\eq1$ (resp $\eq5$), there is no term in $w_{34}$ (resp
$w_{12}$) in $v_1$. Similarly, by $\eq4$, there is no term in
$w_{14}$ in $v_2$. Next, let the $\la_{ij}$'s be the coefficients of
the $w_{ij}$'s in $v_1$ and let the $\mu_{ij}$'s be the coefficients
of the $w_{ij}$'s in $v_2$. By $\eq1$ and $(77)$ applied with $q=1$
and $\la=-\unsurr$, we get:
$$\la_{23}=-\unsurr\,\la_{13}\;\;\&\;\;\la_{24}=-\unsurr\,\la_{14}$$
By $\eq5$ and $(78)$ applied with $q=3$ and $\la=-\unsurr$, we get:
$$\la_{24}=-\unsurr\,\la_{23}\;\;\&\;\;\la_{14}=-\unsurr\,\la_{13}$$
By $\eq4$ and respectively $(77)$ and $(78)$ applied with $q=2$ and
$\la=-\unsurr$, we get:
$$\mu_{34}=-\unsurr\,\mu_{24}\;\;\&\;\;\mu_{13}=-\unsurr\,\mu_{12}$$
Gathering these relations between the coefficients, we obtain:
\begin{eqnarray*}
v_1&=& w_{13}-\unsurr\,w_{23}+\unsur{r^2}\,w_{24}-\unsurr\,w_{14}\\
v_2&\al&w_{12}-\unsurr\,w_{13}+\mu\,(w_{24}-\unsurr\,w_{34})+\mu^{'}\,w_{23}
\end{eqnarray*}
where $\mu$ and $\mu^{'}$ are scalars to determine. For that we use
the mixed relations $\eq2$, $\eq3$, $\eq6$. First and foremost,
$\eq3$ sets the coefficient of proportionality of $v_2$ to be one.
Also, by looking at the coefficient of $w_{24}$ in $\eq3$, we get
$\mu=-\unsurr$. Further, by looking at the coefficient of $w_{23}$
in $\eq6$ and replacing $\mu$ by its value, we get:
$$-\unsurr=-\unsurr+r\,\mu^{'}$$
It follows that $\mu^{'}=0$. Thus, if $\W$ is an irreducible
two-dimensional invariant subspace of $\V$, then it is spanned by
the vectors
\begin{eqnarray}
v_1&=&w_{13}-\unsurr\,w_{23}+\unsur{r^2}\,w_{24}-\unsurr\,w_{14}\\
v_2&=&w_{12}-\unsurr\,w_{13}-\unsurr\,w_{24}+\unsur{r^2}\,w_{34}
\end{eqnarray}
Furthermore, looking at the coefficient of $w_{12}$ in $\eq2$
yields:
$$\unsurl-m=r,$$
from which we derive $\qquad\unsurl=\unsurr$ $\qquad$ \emph{i.e\/}
$\qquad\boxed{l=r}$\\\\
Conversely, if $l=r$, it is a direct verification that the vectors
$v_1$ and $v_2$ defined by $(150)$ and $(151)$ satisfy to all the
relations $\eq i$ with $i=1,\,\dots,6$. Thus, those linearly
independent vectors span an irreducible two-dimensional invariant
subspace of $\V$.\\
We have the immediate corollary:
\begin{Cor} Let $n=4$. Then, there exists an irreducible $2$-dimensional invariant subspace of
$\V$ if and only if $l=r$.
\end{Cor}
\textsc{Proof:} contained in the above.

$\qquad$ \textit{Aparte:} we note at this stage that for
$l\in\lb\unsurr,-\unsurr\rb$, the freedom of the family of vectors
$(v_1,v_2,v_3)$ of Theorem $6$ is a direct consequence of Theorem
$4$ and Corollary $4$. Indeed, if the family $(v_1,v_2,v_3)$ is not
free, then $\di Span_F(v_1,v_2,v_3)=1$ or $\di
Span_F(v_1,v_2,v_3)=2$. In the first case, we get $l=\unsur{r^5}$ by
Theorem $4$, which is impossible. In the second case, the existence
of an irreducible $2$-dimensional invariant subspace would force
$l=r$, also impossible. We note that the same argument does not hold
for the spanning vectors of an irreducible $3$-dimensional invariant
subspace when $l=-r^3$ since we could have $-r^3=\unsur{r^5}$.
Freedom needs to be shown by hands in that case.

We have a corollary of Corollary $4$:
\begin{Cor}
When $l=r$, the intersection module $K(4)$ is the irreducible
$2$-dimensional invariant subspace of $\V$.
\end{Cor}
\textsc{Proof:} by Result $2$, an irreducible $2$-dimensional
invariant subspace of $\V$ when it exists must be unique. Also, as a
matter of fact, it must be contained in $K(4)$. Suppose $l=r$ and
$k(4)>2$. Then $k(4)\in\lb 3,4,5\rb$. Assume first $k(4)=3$. If
$K(4)$ were not irreducible, it would contain a $2$-dimensional
submodule that has a one-dimensional summand. Then
$l=\unsur{r^5}=r$, i.e $(r^2)^3=1$: absurd. So $K(4)$ is
irreducible, $3$-dimensional and $l\in\lb\unsurr,-\unsurr,-r^3\rb$
by Theorem $6$. None of these possibilities is compatible with
$l=r$. Hence $k(4)\neq 3$. If $k(4)=4$, then $K(4)$ is not
irreducible as its dimension is not $1$, $2$ or $3$. Then by
uniqueness in Result $2$ it must be a direct sum of a
$1$-dimensional submodule and an irreducible $3$-dimensional one.
This is again impossible. The other remaining possibility is
$k(4)=5$. Then the only possibility for $K(4)$ is to decompose as a
direct sum of an irreducible $3$-dimensional invariant subspace and
an irreducible two-dimensional one, which for the same reasons as
above has to be excluded. We must conclude that $k(4)=2$, which
achieves the proof of the corollary.

Let's go back to the proof of Proposition $6$ for $n=4$. First, by
Corollary $5$, half of the work is done in the case $n=4$. Indeed,
the corollary says that when $l=r$, the submodule $K(4)$ of $\V$ is
irreducible. When $l=-r^3$, we know from Theorem $6$ of the
existence of a (unique) irreducible $3$-dimensional invariant
subspace of $\V$. Then $k(4)\geq 3$. We attempt to exclude turn by
turn the two possibilities $k(4)=4$ and $k(4)=5$ and summarize the
usual arguments in the following table:
\begin{center}
\begin{tabular}{|c|c|c|c|}
\hline k(4)& K(4)& value for $l$ & Contradiction\\\hline 4&
$3\bigoplus 1$ & $l=\unsur{r^5}=-r^3$ & NONE, apparently\\\hline 5&
$3\bigoplus 2$ & $l=r=-r^3$ & $(r^2)^2\neq 1$\\\hline
\end{tabular}
\end{center}
We need to investigate further in the case $k(4)=4$. If $r^8\neq
-1$, then it is impossible to have $k(4)=4$. Since it is also
impossible to have $k(4)=5$, it comes $k(4)=3$ and $K(4)$ is
irreducible. If on the contrary $r^8=-1$, then we have
$l=-r^3=\unsur{r^5}$. The fact that $l=\unsur{r^5}$ forces the
existence of a one-dimensional invariant subspace of $\V$. Since we
have seen that $k(4)\geq 3$, this one-dimensional invariant subspace
is not $K(4)$, but is contained in $K(4)$. Then by semisimplicity of
$\ih(4)$, it has a summand in $K(4)$. This summand cannot be
two-dimensional (otherwise $l=r$, impossible with $l=-r^3$). Since
by the table above, $k(4)<5$, this summand must in fact be
$3$-dimensional. It forces $k(4)=4$ and $K(4)$ is a direct sum of
the unique one-dimensional invariant subspace of $\V$ and of the
unique irreducible $3$-dimensional invariant subspace of $\V$. This
terminates the case $n=4$.
\item Let's go back to the general case with $n\geq 5$ and $n\neq 8$
when $l=-r^3$. First if $r^{2n}\neq -1$ then it is impossible to
have $l=-r^3=\unsur{r^{2n-3}}$. Thus, by the same arguments as in
the case $l=r$, $K(n)$ is irreducible. If now $r^{2n}=-1$, we have
$l=-r^3=\unsur{r^{2n-3}}$. By Theorem $4$, there exists a
one-dimensional invariant subspace in $\V$. Since by Proposition
$5$, the intersection $K(n)\cap\V_0$ is non-trivial as the vector
$\Y$ is in $K(n)$ for every $n\geq 5$ and since by choice of $l$ and
$r$ and the previous study, the module $K(n-1)$ is irreducible, we
have the inclusions:
$$\begin{array}{l}
0\subset K(n)\cap\V_0\subseteq K(n)\\
0\subset K(n)\cap\V_0 = K(n-1)
\end{array}$$
If $k(n)=1$, then we must have $K(n)\cap\V_0=K(n)=K(n-1)$. But still
by choice of $l$ and $r$, $K(n-1)$ cannot be one-dimensional, hence
a contradiction. Thus it is impossible to have $k(n)=1$ and the
one-dimensional invariant subspace of $\V$, say
$\V_{1,\unsur{r^{2n-3}}}$ must be strictly contained in $K(n)$. This
proves that $K(n)$ is reducible. Moreover, by semisimplicity of
$\ih(n)$, the module $\pal$ has a summand in $K(n)$. This summand
cannot be or contain an irreducible $(n-1)$-dimensional invariant
subspace (otherwise $l\in\lb\unsur{r^{n-3}},-\unsur{r^{n-3}}\rb$ by
Theorem $5$, a contradiction with $l=-r^3$). Since
$$\chl-\dbw=n-1,$$ this summand must be irreducible, of dimension $\cil$, $\dbw$ or dimension greater than $\dbw$.
Hence $k(n)\geq \dbw$.

\item The case $n=8$.\\\\
If $\W$ is an irreducible submodule of $\V$, then
$$\dw\in\lb 1,7,14,20,21\rb$$
\begin{itemize}
\item[*] $\underline{\text{Assume}\;\; l=r}$\\\\
The argument is the same as in the first general case, except that
there is one more possibility for a dimension between $7$ and $21$.
Explicitly, suppose $K(8)$ is reducible. Since we assumed that
$\ih(8)$ is semisimple, there exists $K_1(8)$ and $K_2(8)$
submodules of $K(8)$ such that $K(8)=K_1(8)\bigoplus K_2(8)$.
Without loss of generality, $K_1(8)$ is irreducible. Like in the
first case, if $k_1(8)=21$, then $k_2(8)<28-21=7$. Then $K_2(8)$ is
irreducible and one-dimensional. It follows that $l=\unsur{r^{13}}$
by Theorem $4$. Since $(r^2)^7\neq 1$, we get a contradiction. Also,
since $(r^2)^6\neq 1$ and $(r^2)^7\neq 1$, we cannot have $k_1(8)\in
\lb 1,7\rb$. Hence $k_1(8)\in\lb 14,20\rb$. Suppose first
$k_1(8)=14$. Then $k_2(8)<28-14=14$. This is impossible by the same
arguments as before. On the other hand, if $k_1(8)=20$, then
$k_2(8)<28-20=8$. Then $K_2(8)$ is irreducible and has dimension $1$
or $7$, which is impossible. We conclude that $K(8)$ is irreducible
when $l=r$.\\
\item[*] $\underline{\text{Assume}\;\;l=-r^3}$\\
\begin{itemize}
\item[a)] If $r^{16}\neq -1$ then $\unsur{r^{13}}\neq -r^3$. Since
we also have $(r^2)^8\neq 1$, the same arguments as for $l=r$ yield
the irreducibility of $K(8)$ in that case.
\item[b)] If $r^{16}=-1$, then $l=-r^3=\unsur{r^{13}}$. By Theorem
$4$, there exists a unique one-dimensional invariant subspace of
$\V$. If it was $K(8)$, we would have $K(8)=K(8)\cap\V_0=K(7)$, the
last equality holding by irreducibility of $K(7)$ for these values
of $l$ and $r$. Then $K(7)$ is also one-dimensional which forces
$l=\unsur{r^{11}}$: impossible. We conclude that $K(8)$ is reducible
and contains a one-dimensional invariant subspace that has a summand
in $K(8)$. This summand cannot be one-dimensional or contain any
one-dimensional submodule by uniqueness in Theorem $4$. This summand
may also not be $7$-dimensional since it would then be irreducible
and we would have $l\in\lb\unsur{r^5},-\unsur{r^5}\rb$. This is not
compatible with $l=-r^3$ and $(r^2)^8\neq 1$. Neither can it contain
an irreducible $7$-dimensional invariant subspace. Since it has
dimension less than $28$, it must be $14$, $20$ or $21$-dimensional.
This ends the proof of Proposition $6$.
\end{itemize}
\end{itemize}
\end{list}

\noin We deduce from the Proposition some properties of inclusions
of the $K(n)$'s. Let us first introduce a few more notations. We
extend the definitions of $\V_0$ and $\V_1$ and define $\V_{n-s}$ to
be the $F$-vector space spanned by the $w_{ij}$'s with $1\leq
i<j\leq s-1$. Precisely,
\begin{definition}
\begin{equation}
\V_{n-s}:=\text{Span}_F(w_{ij}|1\leq i<j\leq s-1)\qquad\forall
s=4,\dots,n
\end{equation}
And for better clarity in the notations, we set:
$$\left\lb\begin{array}{cccc}
\V^{(s)}&:=&\V_{n-s-1}&\forall s=3,\dots,n-1\\
\V^{(n)}&:=&\!\!\!\!\!\!\!\!\!\!\!\!\!\!\V&
\end{array}\right.$$
\end{definition}
\noin We will use either of these definitions depending on the
context. The second definition is somehow nicer as it has a meaning
independently of the choice of integer $n$. We have the
Propositions:

\begin{Prop}
Suppose $l=r$. Assume $\ih(n)$ is semisimple with $n\geq 4$. Then
$$K(n)\supset K(n-1)$$
\end{Prop}
\noin\textsc{Proof of Proposition $7$}\\\\
The inclusion uses the irreducibility of the $K(n)$'s for $n\geq 4$,
when $l=r$ and the fact that they are nontrivial. Once this is done,
we recover the result of Proposition $5$. To show that the
inclusion is strict uses this same result.\\

\noin \textit{The proof itself:} suppose $l=r$ and let $n$ be an
integer with $n\geq 5$. We consider the intersection
$K(n)\cap\V^{(n-1)}$. If this intersection was trivial, we would
have $k(n)\leq n-1$, which is impossible. Since
$K(n)\cap\V^{(n-1)}\subset\V^{(n-1)}$ (otherwise $K(n)=\V$ by $\S
8$), $K(n)\cap\V^{(n-1)}$ is a proper non-trivial invariant subspace
of $\V^{(n-1)}$ and we must have:
\begin{equation}
0\subset K(n)\cap\V^{(n-1)}\subseteq K(n-1)
\end{equation}
By irreducibility of $K(n)$ for every $n\geq 4$ when $l=r$, we now
get:
\begin{equation}
K(n)\cap\V^{(n-1)}=K(n-1)
\end{equation}
Thus, we read: $$K(n-1)\subseteq K(n)$$ Further, if the equality
holds between the two sets, then,
\begin{equation}
K(n)\subset\V^{(n-1)}
\end{equation}
We show that $(155)$ is a contradiction. We recall from Proposition
$5$ that the vector $\X$ defined by:
$$\X=r^2\,w_{12}-r\,w_{13}+w_{34}-r\,w_{24}$$
belongs to $K(n)$. Since $K(n)$ is a $B(A_{n-1})$-module, the vector
$$(g_{n-1}\dots g_5g_4)\,.\,\X=
r^{n-2}\,w_{12}-r^{n-3}\,w_{13}+w_{3n}-r\,w_{2n}$$ also belongs to
$K(n)$. As a matter of fact, this vector is in $K(n)$, but is not in
$\V^{(n-1)}$, hence the contradiction in $(155)$. We conclude that
$K(n)\supset K(n-1)$. It remains to show that $K(3)\subset K(4)$.
But by Theorem $7$, the representation $\n^{(3)}$ is irreducible
when $l=r$. Moreover, $K(3)$ is a submodule of $\V$ and is not $\V$
itself. Thus $K(3)$ must be trivial. Finally, $K(4)\neq 0$ as it
contains for instance the vector
$x(4)=r^2\,w_{13}+w_{24}-r\,w_{14}-r\,w_{23}$. Thus, when $l=r$, we
have:
$$0=K(3)\subset K(4)$$ This ends the proof of Proposition $7$.

\begin{Prop}\hfill\\
In Proposition $7$ it suffices to assume that $\ih(n-1)$ is
semisimple.
\end{Prop}

\noin Indeed, in the proof above, we considered the intersection
module $K(n)\cap\V^{(n-1)}$. Consider instead the intersection
$K(n)\cap K(n-1)$ and use Proposition $5$ to claim that this
intersection is non trivial for $n\geq 6$. When $n=5$, use for
instance that $x(4)\in K(4)\cap K(5)$. Then by irreducibility of
$K(n-1)$, we have $K(n-1)\cap K(n)=K(n-1)$, from which we derive
that $K(n)\supseteq K(n-1)$. Since $(g_{n-1}\dots g_4).\X\in
K(n)\setminus K(n-1)$, we have $K(n)\supset K(n-1)$. \\\\
We now study some properties of inclusions of the $K(n)$'s when
$l=-r^3$.
\begin{Prop} Let $n$ be any integer with $n\geq 4$.
\textbf{Suppose $l=-r^3$ and suppose $\ih(n-1)$ is semisimple.}\\\\
$1)$ If $r^{2(n-1)}\neq -1$, then $K(n-1)\subset K(n)$.\\
$2)$ $(n=4)$ If $r^6=-1$, then we have:
$$K(3)\not\subseteq K(4)\subset K(5)\subset K(6)$$
$3)$ If $r^{2(n-1)}=-1$ for some $n\geq 5$ then we have the towers
of inclusion:
$$K(3)\subset\dots\subset K(n-1)\not\subseteq K(n)\subset\dots\subset K(2(n-1))$$
\end{Prop}

\noin\textsc{Proof:} let's first prove $1)$. Following the same
scheme as in the proof of Proposition $8$, we show the inclusion
$K(n-1)\subset K(n)$. For $n\geq 6$, the intersection $K(n)\cap
K(n-1)$ is nontrivial as $\Y$ belongs to all the $K(s),\;s\geq 5$ by
Proposition $5$. When $n=5$, $y(4)\in K(5)\cap K(4)$. Also, when
$n=4$, $-r\,w_{12}-\unsurr\,w_{23}+w_{13}\in K(4)\cap K(3)$. Thus,
$0\subset K(n)\cap K(n-1)\subseteq K(n-1)$. By irreducibility of
$K(n-1)$ when $l=-r^3$, $\ih(n-1)$ is semisimple and $r^{2(n-1)}\neq
-1$, it follows that $K(n-1)\subseteq K(n)$. Since
$$(g_{n-1}\dots g_4).\Y=-r\,w_{23}-\unsurr\,w_{3n}+w_{2n}\in K(n)\setminus K(n-1),$$
we see that $K(n-1)\subset K(n)$. This ends the proof of point $1)$.
Let's prove points $2)$ and $3)$.
\begin{Claim}
Suppose $r^{2n}=-1$ for some $n\geq 3$ and assume that $\ih(n)$ is
semisimple. Then $r^{2k}\neq -1,\;\forall
k\not\in(2\mathbb{N}+1)\,n$.
\end{Claim}
\textsc{Proof:} Suppose $r^{2k}=-1$, some $k\neq n$. The proof is in
two steps.
\begin{itemize}
\item First $k\leq 2n$: we have
$$-1=r^{2n}=r^{2k}r^{2(n-k)}$$
If $r^{2k}=-1$, it comes $r^{2(n-k)}=1$, a contradiction with
$\ih(n)$ semisimple as $k-n\leq n$.
\item If $k\geq 2n$, let's divide $k$ by $2n$:
$$k=t\times 2n+s\;\text{with}\;0\leq s\leq 2n-1$$
It comes $-1=r^{2k}=r^{2s}(r^{2n})^{2t}$. Then $r^{2s}=-1$ with
$0\leq s\leq 2n-1$. It forces $s=n$ by the first point. Then
$k=(1+2t)n$. Hence $k\in(2\mathbb{N}+1)\,n$.
\end{itemize}
We deduce from the claim a lemma:
\begin{lemma}
If $r^{2n}=-1$ and $\ih(n)$ is semisimple, then $\ih(2n-1)$ is
semisimple (but $\ih(2n)$ is not).
\end{lemma}
\textsc{Proof of the lemma:} first if $r^{2n}=-1$, then
$(r^2)^{2n}=1$, so $\ih(2n)$ is not semisimple. What the lemma says
is that $2n$ is the first positive integer $k$ for which $r^{2k}=1$.
In other words, $r^2\neq 1,\dots,(r^2)^{2n-1}\neq 1$, so that
$\ih(2n-1)$ is semisimple. Suppose that there exists an integer
$k_0$ such that $0<k_0<2n$ and $r^{2k_0}=1$. Then it comes:
$$-1=r^{2k_0}r^{2n-2k_0}=r^{2(n-k_0)}$$ so that $$r^{2(n-k_0)}=-1$$
Now $0<|n-k_0|<n$ implies that $|n-k_0|\not\in(2\mathbb{N}+1)n$.
Then by the claim we have $r^{2|n-k_0|}\neq -1$, which yields a
contradiction. Thus, $\forall 0<k\leq 2n-1, \;r^{2k}\neq 1$. This
ends the proof of the lemma.\\

\noin Suppose $n\geq 5$ and let's go back to the proof of
Proposition $9$: by the lemma $\ih(2n-3)$ is semisimple. Also, by
the claim, $r^{2k}\neq -1,\;\forall k\not\in(2\mathbb{N}+1)(n-1)$.
In particular, for all positive integer $k$ such that $k\leq 2n-3$
and $k\neq n-1$, we have $r^{2k}\neq -1$. By applying point $1)$ of
Proposition $9$ to all these admissible $k$'s, we get:
\begin{center}
$K(3)\subset\dots\subset K(n-1)$\\
$K(n)\subset\dots\subset K(2n-2)$
\end{center}
It remains to show that $K(n-1)\not\subseteq K(n)$. When $l=-r^3$
and $r^{2(n-1)}=-1$, there exists a one-dimensional invariant
subspace of $\V^{(n-1)}$, say $\V_{1,n-1}$. Moreover, by Theorem
$4$, $$\V_{1,n-1}=Span_F(\sum_{1\leq s<t\leq n-1}r^{s+t}\,w_{s,t})$$
We have:
\begin{eqnarray*}e_{n-1}.\,\sum_{1\leq s<t\leq
n-1}r^{s+t}\,w_{s,t}&=&\sum_{i=1}^{n-2}r^{i+n-1}\,e_{n-1}.\,w_{i,n-1}\\
&=&\sum_{i=1}^{n-2}\frac{r^{i+n-1}}{r^{n-i-2}}\,w_{n-1,n}\\
&=&r\,\sum_{i=1}^{n-2}(r^2)^i\,w_{n-1,n}\\
&=&r\,\frac{1-(r^2)^{n-2}}{1-r^2}\,w_{n-1,n}\\
\end{eqnarray*}

\noin Since $(r^2)^{n-2}\neq 1$, we see that $e_{n-1}.\sum_{1\leq
s<t\leq n-1}r^{s+t}\,w_{s,t}\neq 0$. Hence $\V_{1,n-1}\not\subseteq
K(n)$ and a fortiori, $K(n-1)\not\subseteq K(n)$. This ends the
proof of point $3)$. Let's prove $2)$. If $l=-r^3$ and $r^6=-1$, we
know from before that
$$K(3)=Span_F(w_{12}+r\,w_{13}+r^2\,w_{23})\bigoplus
Span_F(w_{12}-\unsurr\,w_{13}+\unsur{r^2}\,w_{23})$$ Since
$e_3.(w_{12}+r\,w_{13}+r^2\,w_{23})=(1+r^2)\,w_{34}\neq 0$, we see
that $$K(3)\not\subseteq K(4)$$ Let's show that $K(4)\subset K(5)$.
When $l=-r^3$, $y(4)\in K(4)\cap K(5)$, hence $K(4)\cap K(5)\neq 0$.
Moreover $K(4)$ is irreducible as $r^8\neq -1$ (otherwise $r^2=1$)
and $r^8\neq 1$ (otherwise $r^2=-1$). Thus $K(4)\cap K(5)= K(4)$,
which implies $K(4)\subseteq K(5)$. From $y(4)\in K(4)$ we derive
$g_4.y(4)\in K(5)$ id est
$-r^3\,w_{12}-\unsur{r^2}w_{35}+w_{15}-r\,w_{23}\in K(5)$. This
element is not in $K(4)$. Hence we have $K(4)\subset K(5)$. Finally,
by Proposition $5$, $K(5)\cap K(6)\neq 0$ and since $r^{10}\neq -1$
(as otherwise $(r^2)^2=1$) and $r^{10}\neq 1$ (as otherwise
$(r^2)^2=1$), we know from Proposition $6$ that $K(5)$ is
irreducible. Thus, $K(5)\cap K(6)=K(5)$ and $K(5)\subseteq K(6)$.
Further, by the usual argument,
$$g_5.\,\Y\in K(6)\setminus K(5)$$ Thus $K(5)\subset K(6)$.
%Gathering
%these results, we have shown that:
In summary,
$$\begin{array}{l}\text{When}\;\; l=-r^3,\;(r^2)^2\neq 1\;\;\text{and}\;\; r^6=-1,\;\;\text{we have:}\\
K(3)\not\subseteq K(4)\subset K(5)\subset K(6)
\end{array}$$

\noin This achieves the proof of Proposition $9$. In the next
section, we prove two more theorems that finish characterizing the
dimensions of the invariant subspaces of $\V$ when the
representation is reducible.
\subsection{Properties of dimension of the $K(n)$'s}
Our first main result is the following:
\begin{Prop}
Let $n\geq 4$. When $l=r$, $k(n)=\cil$
\end{Prop}
\textsc{Proof:} first we show that when $l=r$, we have
$k(n)\geq\cil$. We will deal with the case $n=8$ separately. By
Corollary $3$, for $n=7$ or $n\geq 9$, the irreducible
$\ih(n)$-modules have dimensions $1,n-1,\cil,\dbw$ or dimension
greater than $\dbw$. Hence in the case $n=7$ or $n\geq 9$, if
$k(n)<\cil$, then there exists an irreducible $(n-1)$-dimensional
invariant subspace of $\V$ or there exists a one-dimensional
invariant subspace of $\V$, which forces $l\in\lb
\mfl,\jca,-\jca\rb$. None of this value is compatible with $l=r$.
Thus, $k(n)\geq\cil$ in these cases. When $n\in\lb 5,6\rb$, the
story is similar as Corollary $3$ still holds in these cases. As for
$n=4$, $k(4)=2$ by Corollary $5$. Suppose now $n=8$. The irreducible
$\ih(8)$-modules have dimensions $1,7,14,20,21$ or dimension greater
or equal to $28=\text{dim}\,\V^{(8)}$. By Proposition $7$, we know
that $K(8)\supset K(7)$. By the above, we have $k(7)\geq 14$. Hence
it comes $k(8)>14$. Moreover, by Proposition $6$, $K(8)$ is
irreducible when $l=r$. Then $k(8)\geq 20$ and we are done with all
the cases. We will prove the other inequality geometrically, by
using the tangles. Let's introduce a few notations. Let $T(n)$
denote the matrix of the sum endomorphism $\n(n):=\sum_{1\leq
i<j\leq n} \n^{(n)}(X_{ij})$ as in the definition below:
\begin{definition}\hfill
$$\begin{array}{l}\n(n):=\sum_{1\leq i<j\leq n}
\n^{(n)}(X_{ij})\\T(n):=Mat_{\mathcal{B}_{\V^{(n)}}}\n(n)\end{array}$$
\end{definition}
\newtheorem{Prop-Lemma}{Property}
\begin{Prop-Lemma}
$K(n)=Ker\,\n(n)$
\end{Prop-Lemma}
\noin\textsc{Proof:} by Proposition $2$,
$Mat_{\B_{\V^{(n)}}}\n^{(n)}(X_{ij})$ is the matrix whose
$[\binom{j-1}{2}+(j-i)]$-th row is $T(n)$'s and with zeros
elsewhere. It follows that
$$\bigcap_{1\leq i<j\leq
n}\,\text{Ker}\,\n^{(n)}(X_{ij})=\text{Ker}\,\bigp\sum_{1\leq
i<j\leq n}\n^{(n)}(X_{ij})\bigpd$$ With our notations,
$$K(n)=\text{Ker}\,(\n(n))$$
Since $rk(\n(n))=rk(T(n))$, we get the equality on the dimensions:
$$k(n)+rk(T(n))=\chl$$
Thus, to show that $k(n)\leq\cil$, it suffices to show that
$rk(T(n))\geq n$.
\begin{Prop-Lemma}
When $l=r$ and $n\geq 4$, $rk(T(n))\geq n$
\end{Prop-Lemma}
\textsc{Proof:} when $n=4$, since $k(4)=2$, $rk(T(4))=6-2=4$. When
$n=5$ we recall from Result $1$ that there exists an irreducible
$5$-dimensional invariant subspace. Since $K(5)$ is irreducible by
Proposition $6$, we get $k(5)=5$, hence $rk(T(5))=10-5=5$. Since the
square submatrix of $T(6)$ composed of the first ten rows and first
ten columns is $T(5)$, to show that $rk(T(6))\geq 6$, an idea
consists of extracting from $T(5)$ a square submatrix of size $5$
(such a matrix exists as $rk(T(5))=5$) and try and build from it an
invertible square submatrix of $T(6)$ of size $6$ by adding a sixth
subrow and a sixth subcolumn. We notice that in each last five rows
of $T(6)$, there are six zeros amongst the first ten coefficients.
These are indicated in bold below:

\hspace{-1.7in}
$$\nts\nts\nts\nts\nts\nts\nts T(6)=\begin{pmatrix}
\begin{tabular}{c|c}
$T(5)$ & $\begin{array}{ccccc} &&&&\\&&&&\\&&&&\\&&&&\\&&&&\\&&$\Huge{$\ast$}$&&\\&&&&\\&&&&\\&&&&\\&&&&\end{array}$\\
\hline $\begin{array}{cccccccccc}
\ze\ze\ze\ze\ze\ze 1&\unsurr&\unsur{r^2}&\unsur{r^3}\\
\ze\ze\ze\unsurr&\unsur{r^2}&\unsur{r^3}&\unsurr&\ze\ze\mathbf{0}\\
\ze\unsur{r^2}&\unsur{r^3}&\unsur{r^2}&\ze\ze\ze\unsurr&\ze\mathbf{0}\\
\unsur{r^3}&\unsur{r^3}&\ze\ze\unsur{r^2}&\ze\ze\ze\unsurr&\mathbf{0}\\
\unsur{r^4}&\ze\unsur{r^3}&\ze\ze\unsur{r^2}&\ze\ze\ze\unsurr
\end{array}$ & $\begin{array}{ccccc} 2&\unsurr&\unsur{r^2}&\unsur{r^3}&\unsur{r^4}\\r&2&\unsurr&\unsur{r^2}&\unsur{r^3}\\ r^2&r&2&\unsurr&\unsur{r^2} \\r^3&r^2&r&2&\unsurr\\r^4 &r^3&r^2&r&2\end{array}$
\end{tabular}
\end{pmatrix}$$

\noin Since the coefficients in the lower right part of the matrix
are all non-zero, life would be wonderful if there existed five
columns in $T(5)$ aligned on one of these rows of zeros that made a
square submatrix of $T(5)$ invertible. Then, adding this subrow of
five zeros and any subcolumn amongst the last five columns of $T(6)$
would still make the determinant of the resulting extended submatrix
non-zero. Unfortunately things don't happen this way. We wrote a
program in Maple that computes the determinant of all the square
submatrices of $T(5)$ (See appendix $D$). The first five numbers
that are printed in the output correspond to the columns and the
five ones that follow, to the rows. We fix an admissible column and
see if for one of the rows we get a nonzero determinant. To do that
we scroll the bar down onto an admissible column (the first one
being $[1,2,3,4,5]$ and the last one being $[5,6,7,9,10]$ and there
are $30$ of these) and look for a row that would give a nonzero
coefficient at the end. The result is negative. The mathematical
explanation for that is unknown to the author. We need to ask for
slightly less and content ourselves with five columns aligned on
only four zeros. Let's pick the first four zeros of row $11$ and
notice that it could be nice to have a one as the last non-zero
coefficient. Thus, we investigate about the determinants of the
submatrices of $T(5)$ with the pattern of columns $\C:=[1,2,3,4,7]$.
The first (in the lexicographic order) admissible $5$-tuple of rows
is actually $\Rrond:=[1,2,3,4,7](=\C)$. And in fact picking $\C$ and
$\Rrond$ is also a natural choice as, by our program, it is for
these subcolumns and subrows that the first nonzero determinant
arises in lexicographic order (when ordering the columns first).\\\\
\textit{A few notations:}
\begin{itemize}
\item Given a matrix $A$,
$\text{submatrix}(A,[i_1,\,\dots,\,i_s],[j_1,\,\dots,\,j_s])$
denotes the square submatrix of size $s$ of $A$ with subrows
$i_1,\dots,i_s$ and subcolumns $j_1,\dots,j_s$, where we followed
Maple's notations.
\item Given $\Rrond$ a subset of rows,
$c_{i_k,\Rrond}$ denotes the extracted $i_k$-th column with rows
$\Rrond$.\\
\end{itemize}
Let $M:=\sm(T(5),[1,2,3,4,7],[1,2,3,4,7])$.\\
We have $det(M)=\frac{(r^2+1)^2}{r^2}\neq0$, as $(r^2)^2\neq 1$.\\\\
Hence the family of columns
$(c_{1,\Rrond},c_{2,\Rrond},c_{3,\Rrond},c_{4,\Rrond},c_{7,\Rrond})$
is free. A fortiori, the subfamily
$(c_{1,\Rrond},c_{2,\Rrond},c_{3,\Rrond},c_{4,\Rrond})$ is free.

$$\begin{array}{l}
\text{Let}\;\; M^{'}:=\sm(T(6),[1,2,3,4,7],[1,2,3,4,12])\\
det(M^{'})=0\;\;\text{(see Appendix $D$)}
\end{array}$$

\noin Then
$(c_{1,\Rrond},c_{2,\Rrond},c_{3,\Rrond},c_{4,\Rrond},c_{12,\Rrond})$
is not free, but
$(c_{1,\Rrond},c_{2,\Rrond},c_{3,\Rrond},c_{4,\Rrond})$ is. Hence,
the column $c_{12,\Rrond}$ is a linear combination of columns
$c_{1,\Rrond},c_{2,\Rrond},c_{3,\Rrond},c_{4,\Rrond}$. Thus, we
don't modify the determinant of $M$ by adding to $c_{7,\Rrond}$ a
multiple of $c_{12,\Rrond}$. The $12$th column of $T(6)$ is the
vector $\n^{(6)}(w_{46})$. The coefficient on the $11$th row is
$[\n^{(6)}(w_{46})]_{x_{\al_5}}=[\n^{(6)}(X_{56})(w_{46})]_{x_{\al_5}}=\unsurr$
by $(TL)_1$ of Appendix $C$ with $l=r$. We recall that the
coefficient of the $7$th column and the $11$th row is a one, while
the other coefficients of the extracted matrix on the same row are
$0$'s. We now do the operation on the columns
$$c_{7,\Rrond@[\binom{5}{2}+1]}\leftarrow\;c_{7,\Rrond@[\binom{5}{2}+1]}-r\,c_{12,\Rrond@[\binom{5}{2}+1]}$$
in order to make a fifth zero appear, where we used the symbol $@$
for concatenation of lists. By doing so the determinant of $M$ is
unchanged. Let's consider the matrix
$$M^{''}:=\sm(T(6),[1,2,3,4,7,11],[1,2,3,4,7,12])$$
$M^{''}$ is a square submatrix of $T(6)$ of size $6$ and
$$\text{det}(M^{''})=\underset{\neq 0}{\underbrace{\text{det}(M)}}\times\unsurr\neq 0$$
Thus, we have exhibited an invertible submatrix of $T(6)$ of size
$6$. This shows that $rk(T(6))\geq 6$. Furthermore, we have seen
that when $l=r$, we have for all $n\geq 4$ that $k(n)\geq\cil$. This
is equivalent to $rk(T(n))\leq n$. Thus, we have $rk(T(6))\leq 6$.
Hence $rk(T(6))=6$. Next, for $n\geq 7$ we inductively build
invertible submatrices $S(n)$ of $T(n)$ of size $n$ in the following
way: $$\begin{array}{l}
S(5):=\sm(T(5),\Rrond_5,\C_5)\;\text{with}\;\C_5=\Rrond_5:=[1,2,3,4,7](=\Rrond)\\\\
S(n):=\sm(T(n),\Rrond_{n-1}@[\binom{n-1}{2}+(n-5)],\C_{n-1}@[\binom{n-1}{2}+(n-4)])
\end{array}$$
$S(n)$ is built from $S(n-1)$ by adding the extracted
$[\binom{n-1}{2}+(n-5)]$-th subrow and the extracted
$[\binom{n-1}{2}+(n-4)]$-th subcolumn. In other words, $S(n)$ is
built from $S(n-1)$ by adding the subrow corresponding to $X_{5n}$
and the subcolumn corresponding to $w_{4n}$ as shown on Figure $1$.
In the table below, we gathered the results of the actions of
$X_{12}$, $X_{23}$, $X_{13}$, $X_{34}$ and $X_{45}$ on the vectors
$w_{4,n}$'s for $n\geq 6$. To calculate the coefficient of the
action of $X_{34}$ on $w_{4k}$, we used $(SR)_{k-4}$ with $i=3$ and
$j=4$. It yields the coefficient $r^{k-5}$; to calculate the
coefficient of the action of $X_{45}$ on $w_{4k}$, we used
$(TR)_{k-5}$ with $l=r$.
It also yields the coefficient $r^{k-5}$.\\
\begin{center}
\begin{tabular}[c]{c|c|c|c|c|c}
& $w_{46}$ & $w_{47}$ & $\dots$ &
$w_{4n}$&\dots\\\cline{1-3}\cline{5-5} $X_{12}$& $0$
&$0$&&$0$&\\\cline{1-3}\cline{5-5}
$X_{23}$&$0$&$0$&&$0$&\\\cline{1-3}\cline{5-5} $X_{13}$
&$0$&$0$&\dots&$0$&\dots\\\cline{1-3}\cline{5-5}
$X_{34}$&$r$&$r^2$&&$r^{n-5}$&\\\cline{1-3}\cline{5-5}
$X_{45}$&$r$&$r^2$&&$r^{n-5}$&\\\cline{1-3}\cline{5-5}
\end{tabular}
\end{center}
It appears clearly on the table that
$c_{\binom{n-1}{2}+(n-4),\Rrond}$ is a multiple of $c_{12,\Rrond}$.
Since, as seen above, $c_{12,\Rrond}$ is a linear combination of the
extracted columns
$c_{1,\Rrond},c_{2,\Rrond},c_{3,\Rrond},c_{4,\Rrond}$ of $T(5)$,
each extracted column $c_{\binom{n-1}{2}+(n-4),\Rrond}$ of $T(n)$,
$n\geq 7$, is also a linear combination of columns
$c_{1,\Rrond},c_{2,\Rrond},c_{3,\Rrond},c_{4,\Rrond}$. Thus, for any
$n\geq 7$, we won't modify the determinant of $S(5)(=M)$ by adding
to its $7$-th column $c_{7,\Rrond}$ a multiple of
$c_{\binom{n-1}{2}+(n-4),\Rrond}$. As it appears on the Figure below,\\

\hspace{-1.2in} $\begin{array}{cc}
\begin{array}{l}
\\
\\
\\
\\
\\
\\
\\
\\
X_{56}\\
\\
\\
\\
\\
\\
X_{57}\\
\\
\\
\\
\\
\\
\\
\\
\\
\\
X_{5n}\\
\\
\\
\end{array}
& \begin{array}{l}\begin{array}{ccccccccccccccccccccccccc} &&&&&&&&&&\;\;\;\,w_{46}&&&&&w_{47}&&&&&&\nts\nts\nts\nts w_{4n}&&\\
&\!1&\,2&\,\,3&\;4&&\;\;\,\;7&&&10&\;\;\;12&&15&&&\,\!18
 &&\!\!\dots&\binom{n-1}{2}&&&\binom{n-1}{2}+(n-4)&&\!\!\!\nts\binom{n}{2}
\end{array}\\
$\begin{tabular}[c]{c|c|c|c|} $\begin{array}{cccccccccc}
\sq&\sq&\sq&\sq&&&\sq&&&\\
\sq&\sq&\sq&\sq&&&\sq&&&\\
\sq&\sq&\sq&\sq&&&\sq&&&\\
\sq&\sq&\sq&\sq&&&\sq&&&\\
&&&&&&&&&\\
&&&&&&&&&\\
\sq&\sq&\sq&\sq&&&\sq&&&\\
&&&&&&&&&\\
&&&&&&&&&\\
&&&&&&&&&\\
\end{array}$
& $\begin{array}{ccccc} &0&&&\\
&0&&& \\
&0&&& \\
&r&&& \\
&&&& \\
&&&& \\
&r&&& \\
&&&& \\&&&& \\&&&&
\end{array}$
& $\begin{array}{ccccccccc} &&0&&&&&&\\
&&0&&&\\&&0&&&\\&&r^2&&&&&&\\
&&&&&&&&\\
&&&&&&&&\\
&&r^2&&&&&&\\
&&&&&&&&\\&&&&&&&&\\&&&&&&&&
\end{array}$
& $\begin{array}{cccccccccc}
&&&&0&&&&&\\
&&&&0&&&&&\\&&&&0&&&&&\\&&&&r^{n-5}&&&&&\\
&&&&&&&&&\\&&&&&&&&&\\
&&&&r^{n-5}&&&&&\\
&&&&&&&&&\\&&&&&&&&&\\&&&&&&&&&\\
\end{array}$\\
\cline{1-1} \end{tabular}$\\
$\begin{tabular}[c]{c|c|c|}$\begin{array}{cccccccccc}
\,0&\,0&\,\,0&\,\,0&&&1&&&\\
&&&&&&&&&\\
&&&&&&&&&\\
&&&&&&&&&\\
&&&&&&&&&\\
\end{array}$
$\begin{array}{cccccc} &&\,\unsurr&&&\\
&&&&& \\
&&&&& \\
&&&&& \\
&&&&&
\\
\end{array}$\hspace{0.0025in}
&
$\begin{array}{ccccccccc} &&0&&&&&&\,\;\\
&&&&&&&&\,\;\\&&&&&&&&\,\;\\&&&&&&&&\,\;\\
&&&&&&&&\,\;\\
\end{array}$
&\hspace{-0.02in}$\begin{array}{ccccccccccc}
\!&&&&\,\,0&&&&&&\,\,\,\\
\!&&&&\,\,&&&&&&\,\,\,\\\!&&&&\,\,&&&&&&\,\,\,\\\!&&&&\,\,&&&&&&\,\,\,\\
\!&&&&\,\,&&&&&&\,\,\,
\end{array}$\\
\cline{1-1}
\end{tabular}$
\\
$\begin{tabular}[c]{c|c|} $\begin{array}{cccccccccccccccccccccccc}
&&&&&&&&&&&&&&&&&&&&&&&\!\!\!\hspace{-0.02in}\\
\,0&\,0&\,\,0&\,\,0&&&\unsurr&&&&&\;\;\;0&&&&&&&\;\;\unsurr&&&&&\!\!\!\hspace{-0.02in}\\
&&&&&&&&&&&&&&&&&&&\!\!\ddots&&&&\!\!\!\hspace{-0.02in}\\
&&&&&&&&&&&&&&&&&&&&\,\ddots&&&\!\!\!\hspace{-0.02in}\\
&&&&&&&&&&&&&&&&&&&&&&&\!\!\!\hspace{-0.02in}\\
&&&&&&&&&&&&&&&&&&&&&&&\!\!\!\hspace{-0.02in}\\
\end{array}$\hspace{-0.2in}&\hspace{-0.08in}
$\begin{array}{ccccccccccc}\!\!&&&&\,\,&&&&&&\;\hspace{0.02in}\\
\!\!&&&&\,\,\,0&&&&&&\;\hspace{0.02in}\\
\!\!&&&&\,\,&&&&&\;\\\!\!&&&&\,\,&&&&&\\\!\!&&&&\,\,&&&&&&\;\hspace{0.02in}\\
\!\!&&&&\,\,&&&&&\;\\\!\!&&&&\,\,&&&&&&\;\hspace{0.02in}\\
\end{array}$\\\cline{1-1}
\end{tabular}$\\
$\begin{tabular}[c]{c|}
$\begin{array}{cccccccccccccccccccccccccccccccccccccc}
&&&&&&&&&&&&&&&&&&&&&&&&&&&&&&&&&&&&&\;\,\hspace{0.01in}\\
&&&&&&&&&&&&&&&&&&&&&&&&&&&&&&&&&&&&&\;\,\hspace{0.01in}\\
&&&&&&&&&&&&&&&&&&&&&&&&&&&&&&&&&&&&&\;\,\hspace{0.01in}\\
&&&&&&&&&&&&&&&&&&&&&&&&&&&&&&&&&&&&&\;\,\hspace{0.01in}\\
\,0&\,0&\,\,0&\,\,0&&\;\;\unsur{r^{n-6}}&&&&&&0&&&&&&&\!\!\!0&&&&&&&&&&&&&\!\!\,\unsurr&&&&&&\;\,\hspace{0.01in}\\
&&&&&&&&&&&&&&&&&&&&&&&&&&&&&&&&&&&&&\;\,\hspace{0.01in}\\
&&&&&&&&&&&&&&&&&&&&&&&&&&&&&&&&&&&&&\;\,\hspace{0.01in}\\
&&&&&&&&&&&&&&&&&&&&&&&&&&&&&&&&&&&&&\;\,\hspace{0.01in}\\
&&&&&&&&&&&&&&&&&&&&&&&&&&&&&&&&&&&&&\;\,\hspace{0.01in}\end{array}$\\
\cline{1-1}
\end{tabular}$\\\\
Figure\;\;1
\end{array}
\end{array}$

\noin we have:
$$\begin{array}{l}
\forall n\geq 7,\,X_{5n}.w_{4n}=\unsurr\,w_{5n}\qquad\qquad\text{(Rule $(TL)_1$ with $l=r$)}  \\
\forall n\geq
7,\,X_{5n}.w_{45}=\unsur{r^{n-6}}\,w_{5n}\qquad\qquad\;\;\text{(Rule
$(SL)_1$ with $i=5$ $\&$ $j=n$)}\\
\forall n\geq 7,\,X_{5n}.w_{4j}=0\;\text{for all $6\leq j\leq
n-1$}\;\text{BECAUSE $l=r$ (cf Rule $(CL)$)}\\
\forall n\geq 7,\,X_{5n}.w_{s,t}=0\;\text{for
all}\;(s,t)\in\lb(1,2),(2,3),(1,3),(3,4)\rb\\
\forall n\geq 7,\,X_{5k}.w_{4j}=0,\,\forall 6\leq k\leq
n-1,\,\forall k+1\leq j\leq n
\end{array}$$
Hence, doing the following operation $\mathcal{O}_n$ on $S(n)$ with:
$$\mathcal{O}_s:\C_{7,\Rrond_s}\leftarrow
\C_{7,\Rrond_s}-\sum_{k=0}^{n-6}\unsur{r^{k-1}}\,\C_{\binom{k+5}{2}+(k+2),\Rrond_s}$$
transforms $S(n)$ into a matrix by blocks with $0$'s in the lower
left quadrant, a square matrix
$M_{\mathcal{O}_5}=S(5)_{\mathcal{O}_5}$ of size $5$ obtained from
$M$ by doing the operation $\mathcal{O}_5$ and whose determinant
equals the one of $M$ in the upper left quadrant and an upper
triangular square matrix of size $n-5$ with $\unsurr$'s over the
diagonal in the lower right quadrant:
$$S(n)_{\mathcal{O}_n}=\begin{array}{l}\begin{bmatrix}
$\begin{tabular}[c]{c|cccc} &\\
 $S(5)_{\mathcal{O}_5}$&\Huge{$\ast$}\\&\\\hline
\Huge{$0$}&$\begin{array}{cccc}
\unsurr&&&\\&\ddots&\qquad\text{\huge{$\ast$}}&\\\text{\huge{$0$}}&&\ddots&\\&&&\unsurr
\end{array}$
\end{tabular}$
\end{bmatrix}\\
\begin{array}{cc}
&\\
\!\overlra{\qquad 5\qquad} &\!\!\!\!\overlra{\qquad\qquad
n-5\qquad\qquad}
\end{array}
\end{array}$$
Then it comes:
\begin{eqnarray*}\text{det}(S(n))\;=\;\dete(S(n)_{\mathcal{O}_n})&=&\dete(S(5)_{\mathcal{O}_5}).\left(\unsurr\right)^{n-5}\\&=&\dete(M_{\mathcal{O}_5}).\left(\unsurr\right)^{n-5}\\
&=&\dete(M).\left(\unsurr\right)^{n-5}\\&\neq& 0\end{eqnarray*}
Hence, $S(n)$ is a square submatrix of size $n$ of $T(n)$ that is
invertible. This shows that $rk(T(n))\geq n$ for all $n\geq 7$. And
since we know from before that for all $n\geq 4$, $rk(T(n))\leq n$
when $l=r$, we have in fact:
$\text{When}\;l=r,\;rk(T(n))=n\;\text{for all}\;n\geq 7$. This is
the same as saying that $\text{When}\;l=r,\;k(n)=\cil\;\text{for
all}\;n\geq 7$. Since as we have seen along the proof, we also have
$rk(T(4))=4$ (\emph{i.e\/} $k(4)=2=\frac{4.1}{2}$), $rk(T(5))=5$
(\emph{i.e\/} $k(5)=5=\frac{5.2}{2}$) and $rk(T(6))=6$ (\emph{i.e\/}
$k(6)=9=\frac{6.3}{2}$) , this ends the proof of Proposition $10$.
Moreover we will remember the following fact, equivalent to
Proposition $10$:
\begin{Fact}
When $l=r$, the rank of $T(n)$ is $n$ for all $n\geq 4$.
\end{Fact}
\noin We have two direct consequences of Propostion $10$. Unless
otherwise mentioned, $\ih(n)$ is assumed to be semisimple.
%In what follows, we will denote by $c_{i,\Rrond}$ the column
%extracted from the $i$-th column with rows $\Rrond$. Following
%Maple's notations, we will denote by
%$submatrix(A,[i_1,\,\dots,\,\i_s],[j_1,\,\dots,\,j_s$ the square
%submatrix of $A$ of size $s$ with subrows $i_1,\dots,i_s$ and
%subcolumns $j_1,\dots,j_s$.

%of Using Maple, we extracted from the $10\times10$ matrix $T(5)$ a
%square submatrix of size $5$ that is invertible

\begin{Cor}
Let $n$ be an integer with $n\geq 4$. If $l=r$, then there exists an
irreducible $\cil$-dimensional invariant subspace of $\V$.
\end{Cor}

\textsc{Proof:} by Proposition $6$, $K(n)$ is irreducible. By
Proposition $10$, $K(n)$ is $\cil$-dimensional. This ends the proof.

\begin{Cor} Let $n$ be an integer with
$n\geq 5$. If there exists an irreducible $\dbw$-dimensional
invariant subspace of $\V$, then $l=-r^3$.
\end{Cor}

\textsc{Proof:} Suppose $n\geq 6$. Let $\W$ be an irreducible
$\dbw$-dimensional invariant subspace. Since for $n\geq 6$, we have
$2n-3\leq\cil<\cil+1=\dbw$, it comes $\dw>2n-3$, hence
$\W\cap\V_1\neq\lb 0\rb$. Also, $n-1<\dbw$ implies that
$\W\cap\V_0\neq\lb 0\rb$. Thus, by theorem $7$ we get: $l\in\lb
r,-r^3,\unsur{r^{2(n-2)-3}},\unsur{r^{n-5}},-\unsur{r^{n-5}}\rb$ and
$l\in\lb
r,-r^3,\unsur{r^{2(n-1)-3}},\unsur{r^{n-4}},-\unsur{r^{n-4}}\rb$.
This only leaves the possibility $l\in\lb r,-r^3\rb$ by
semisimplicity of $\ih(n)$. Now $\dw=\dbw$ implies that
$k(n)\geq\cil+1$ and when $l=r$ we have $k(n)=\cil$ by Proposition
$10$. Hence the only possibility for $l$ is in fact $l=-r^3$. \\\\
\textit{The case $n=5$}\\\\
Suppose that there exists an irreducible $6$-dimensional invariant
subspace of $\V$. Then we claim that there cannot exist an
irreducible $5$-dimensional invariant subspace of $\V$. Indeed let's
name them $\V_5$ and $\V_6$. We have $\V_5\subseteq K(5)$ and
$\V_6\subseteq K(5)$. The existence of an irreducible
$5$-dimensional invariant subspace of $\V$ implies $l=r$ by Result
$1$; by Proposition $6$, $K(5)$ is irreducible as $l=r$. Then
$K(5)=\V_5=\V_6$, contradiction. Hence there cannot exist an
irreducible $5$-dimensional invariant subspace. Hence $l\neq r$. But
$\n^{(5)}$ is reducible, hence $l\in\lb -r^3,
\unsur{r^2},-\unsur{r^2},\unsur{r^7}\rb$. Since $6>4$,
$\V_6\cap\V^{(4)}\neq 0$. Hence, $\n^{(4)}$ is reducible and
$l\in\lb -r^3,\unsurr,-\unsurr,\unsur{r^5}\rb$. Then $l$ must take
one of the two values $-r^3$ or $\unsur{r^7}$. Since $l\neq r$, we
have $k(5)\in\lb 6,7\rb$. If $l=\unsur{r^7}$ then there exists a
one-dimensional invariant subspace $\V_1$ of $\V^{(5)}$ which is in
direct sum with $\V_6$. Hence $k(5)=7$ and $K(5)=\V_1\oplus\V_6$. If
$l\neq \unsur{r^7}$, then $l=-r^3$ and $k(5)=6$. Then $K(5)$ is
irreducible and $K(5)=\V_6$. Let's push the study a little bit more.
We have $\di(\V_6\cap\V^{(4)})=6+6-\di(\V_6+\V^{(4)})\geq 2$ and
$\V_6\cap\V^{(4)}\subseteq K(4)$. Hence, $k(4)\geq 2$. Moreover,
since $l\neq r$, $k(4)\in\lb 3, 4\rb$. Suppose $k(4)=4$. Then $K(4)$
is a direct sum of an irreducible $3$-dimensional invariant subspace
of $\V^{(4)}$ and of a one-dimensional one. Then
$l=\unsur{r^5}=-r^3$. So $r^8=-1$. Then $r^{10}\neq -1$ and
$l=-r^3$. Hence $K(5)$ is irreducible. Moreover $l\neq \unsur{r^7}$.
Thus, if $l=\unsur{r^7}$, then $k(4)=3$ and $K(4)$ is irreducible as
$l\neq r$. Also it follows that $l\in\lb -\unsurr,-r^3\rb$, hence
either $r^6=-1$ or $r^{10}=-1$. We conclude that Corollary $7$ holds
for $n=5$ in the case $r^6\neq -1$. Also, we deduce from this
discussion that if $K(4)$ is reducible, then $l\neq \unsur{r^7}$.
Hence $l=-r^3$ and $K(5)$ is irreducible. Thus, $K(4)$ and $K(5)$
cannot be simultaneously reducible. Suppose now that $r^6=-1$. We
will show that $l=-r^3$, and this will prove that Corollary $7$
holds in fact in any case. We have $r^{10}\neq -1$. Thus,
$\unsur{r^7}\neq -r^3$. Suppose $l\neq -r^3$. Then $l=\unsur{r^7}$.
Then as seen above, $K(5)=\V_1\oplus\V_6$. Then $k(5)=7>4$ and
$K(4)$ is irreducible as for instance $K(5)$ is reducible. It
follows that $K(5)\cap\V^{(4)}=K(4)$. In particular, we have
$K(4)\subseteq K(5)$. Also, by the above, when $l=\unsur{r^7}$, we
have $k(4)=3$. Thus, $K(4)$ is irreducible, $3$-dimensional. Since
$l=\unsur{r^7}=-\unsurr$, we know from Theorem $6$ that $K(4)$ is
spanned over $F$ by the vectors:
\begin{eqnarray*}
v_1&=&\bigp\unsurr+r\bigpd\,w_{12}+(w_{13}-\unsurr\,w_{23})+r\,(w_{14}-\unsurr\,w_{24})\\
v_2&=&\bigp\unsurr+r\bigpd\,w_{23}+(w_{24}-\unsurr\,w_{34})-r(w_{12}-\unsurr\,w_{13})\\
v_3&=&\bigp\unsurr+r\bigpd\,w_{34}-(w_{13}-\unsurr\,w_{14})-(w_{23}-\unsurr\,w_{24})
\end{eqnarray*}

\noin Let's compute the action of $X_{35}$ on $v_2$. By using the
table in Appendix $C$, where we replaced $l$ by $-\unsurr$, we have
the equalities:

$$\begin{array}{ccccc}
X_{35}.w_{34}&=&-r&x_{\al_3+\al_4} &\text{by (ML)}_1\\
X_{35}.w_{13}&=&\unsur{r^2}&x_{\al_3+\al_4} &\!\!\text{by (SL)}_2\\
X_{35}.w_{23}&=&\unsurr&x_{\al_3+\al_4} &\!\!\text{by (SL)}_1\\
X_{35}.w_{24}&=&(\unsurr-r)(-\unsurr-r)&x_{\al_3+\al_4} &\;\;\;\text{by (CL)}_{(1,1)}\\
\end{array}$$

\noin Then it comes:
\begin{equation}
X_{35}.v_2=\bigp r+\unsurr\bigpd^2\,x_{\al_3+\al_4}
\end{equation}

\noin Since $K(4)\subseteq K(5)$, we must have $\bigp
r+\unsurr\bigpd^2=0$, a contradiction. Hence it is impossible to
have $l=\unsur{r^7}$ when $r^6=-1$. Then $l$ must take the value
$-r^3$, which ends all the cases. Conversely, if $l=-r^3$, we claim
that there exists an irreducible $6$-dimensional invariant subspace
inside $\V^{(5)}$. Indeed, by Theorem $7$, $\n^{(5)}$ is reducible,
hence there exists an irreducible $B(A_4)$-submodule of $\V$. It
cannot be $5$-dimensional by Result $1$ ($l\neq r$); nor can it be
$4$-dimensional ($l\not\in\lb \unsur{r^2},-\unsur{r^2}\rb$). If
$r^{10}\neq -1$, it cannot be one-dimensional either. Then it must
be $6$-dimensional. If $r^{10}=-1$, then we recall from Proposition
$6$ that $k(5)\geq 6$. Then there must exist an irreducible
$6$-dimensional submodule as well. Hence the Theorem:
\begin{thm}\hfill\\
Suppose $n=5$. There exists an irreducible $6$-dimensional invariant
subspace of $\V$ if and only if $l=-r^3$.
\end{thm}

Finally we note that Corollary $7$ does not hold for $n=4$.
Indeed, since $(r^2)^4\neq 1$, we have $\unsurr\neq -r^3$ and $-\unsurr\neq -r^3$. And when $l\in\lb\unsurr,-\unsurr\rb$, there exists an irreducible $\frac{3.2}{2}$-dimensional invariant subspace of $\V^{(4)}$ by Theorem $6$. \\

We have a corollary of Corollary $7$:
\begin{Cor}\hfill\\
Let $n$ be an integer with $n\geq 4$. \\
$i)$ If there exists an irreducible $\cil$-dimensional invariant subspace of $\V$, then there does not exist any $\dbw$-dimensional invariant subspace of $\V$.\\
$ii)$ If there exists an irreducible $\dbw$-dimensional invariant subspace of $\V$, then there does not exist any $\cil$-dimensional invariant subspace of $\V$.\\
\end{Cor}

\textsc{Proof:} we prove both points at the same time. Let's first assume $n\geq 5$. Suppose that there exists both an irreducible $\cil$-dimensional invariant subspace $\V_{\cil}$ of $\V$ and an irreducible $\dbw$ invariant subspace $\V_{\dbw}$ of $\V$. By Corollary $7$, we have $l=-r^3$. If $r^{2n}\neq -1$, $K(n)$ is irreducible by Proposition $6$, hence a contradiction. If $r^{2n}=-1$, then on one hand there exists a one-dimensional invariant subspace $\V_1$ of $\V$, which is in direct sum with $\V_{\dbw}$. Hence we have $\V_1\oplus\V_{\dbw}\subseteq K(n)$. Then $k(n)\geq \dbw+1$. On the other hand, $\V_{\cil}$ has a summand, say $S$ in $K(n)$. Let's denote by $s$ its dimension. Since $k(n)\geq \cil+2$, we must have $s\geq 2$. Since $l=-r^3$, there does not exist any $(n-1)$-dimensional invariant subspace inside $\V$. Hence $s\geq 2$ implies in fact $s>n-1$. Now if $s\geq n$, we have $\di(\V_{\cil}\oplus S)\geq n+\cil=\chl$: a contradiction since $k(n)<\chl$. This ends the proof in the case $n\geq 5$. It remains to do the case $n=4$. In this case, if there exists an irreducible $2$-dimensional invariant subspace of $\V^{(4)}$, then $l=r$ by Result $2$; if there exists an irreducible $3$-dimensional invariant subspace of $\V^{(4)}$, then $l\in\lb -r^3,\unsurr,-\unsurr\rb$ by Theorem $6$. Then it is impossible to have an irreducible $3$-dimensional invariant subspace and an irreducible $2$-dimensional one. \\

When $l=-r^3$, it appears that it is not as easy to show properties
on the rank of the matrix $T(n)$. However, we have the nice result:

\begin{Prop}
Suppose $\ih(n-1)$ is semisimple, $l=-r^3$ and $r^{2(n-1)}\neq -1$.
Then,
$$k(n)\geq k(n-1)+(n-2)$$
\end{Prop}

\noin\textsc{Proof of the proposition:} by proposition $9$, point
$1)$, we know under these assumptions that $K(n-1)\subset K(n)$.
Thus, any vector that annihilates the matrix $T(n-1)$ annihilates
the matrix $T(n)$. Let $v_1,v_2,\dots,v_{k(n-1)}$ be a basis of
$K(n-1)$. Define $(n-2)$ linearly independent vectors by:
$$V_k=w_{k+1,n}-r\,w_{k,n}+r^{n-k}\,w_{k,k+1},\;\;k=1\dots n-2$$

\begin{Claim}\hfill\\
$v_1,v_2,\dots,v_{k(n-1)},V_1,\dots,V_{n-2}\,$ are
$\,k(n-1)+(n-2)\,$ linearly independent vectors of $K(n)$
\end{Claim}

\textsc{Proof of the claim:} we want to show that the $X_{ij}$'s
annihilate the $V_k$'s for $k=1,\dots,n-2$. First, let's show that
the last $(n-1)$ rows of the matrix $T(n)$ annihilate these vectors.
To that aim, we compute:

\begin{equation*}
\forall\,1\leq j,k\leq n-1,\;\;[X_{j,n}.w_{k,n}]_{_{w_{j,n}}}=\begin{cases} -r^{k-j-2}&\text{if $k<j$}\\
-r^{k-j+2}&\text{if $k>j$}\\
-\unsur{r^2}-r^2&\text{if $k=j$}
\end{cases}
\end{equation*}

\noin Let's fix a row of the matrix $T(n)$ that corresponds to the
action of a $X_{jn}$ with $1\leq j\leq n-1$. We want to multiply
this row by the vector $V_k$. First, we let $X_{j,n}$ act on the
vectors $w_{k+1,n}$ and $w_{k,n}$ and multiply the first resulting
coefficient by $1$ and the second resulting coefficient by $-r$,
then add the two coefficients.
\begin{itemize}
\item[*] If $k>j$, we get $-r^{k-j+3}+r^{k-j+3}=0$
\item[*] If $k<j-1$, we get $-r^{k-j-1}+r^{k-j-1}=0$
\item[*] If $k=j-1$, we get $-\unsur{r^2}-r^2+\unsur{r^2}=-r^2$
\item[*] If $k=j$, we get $-r^3+\unsurr+r^3=\unsurr$
\end{itemize}

\noin Now fix a $k\in\lb 1,\dots,n-2\rb$. We have:
$$X_{jn}.(w_{k+1,n}-r\,w_{k,n})=0\;\;\text{except for}\;\;j\in\lb k,k+1\rb$$
\noin Also, we have:
$$\begin{array}{cccc}
[X_{k+1,n}.w_{k,k+1}]_{_{w_{k+1,n}}}&=&\unsur{r^{n-k-2}}&\\
\begin{bmatrix}X_{k,n}.w_{k,k+1}\end{bmatrix}_{_{w_{k,n}}} &=& -\unsur{r^{n-k+1}}&\\
\begin{bmatrix} X_{j,n}.w_{k,k+1}\end{bmatrix}_{_{w_{j,n}}}&=&0&\text{if $j\not\in\lb k,k+1\rb$}
\end{array}$$

\noin Thus, we get:
$$\forall\; 1\leq j\leq n-1,\;\forall\; 1\leq k\leq n-2,\; X_{jn}.(w_{k+1,n}-r\,w_{k,n}+r^{n-k}\,w_{k,k+1})=0$$
\noin This shows that the last $(n-1)$ rows of the matrix $T(n)$
annihilate the vectors $V_k$'s, $k=1,\dots,n-2$. We will now show
that the whole matrix $T(n)$ annihilates in fact these vectors.
Given two positive integers $s$ and $t$ such that $1\leq s<t\leq
n-1$, we have:

$$\begin{bmatrix}X_{s,t}.\,w_{k,n}\end{bmatrix}_{_{w_{s,t}}}\;\;\;=\begin{cases}
-r^{n-t+2}&\text{if $s=k$}\\
r^{n-s-2}&\text{if $t=k$}\\
(-r^3-r)(r^{k-s+n-t-1}-r^{k-s+n-t-3})&\text{if $s<k$ and $t>k$}\\
0&\text{otherwise}
\end{cases}$$

$$\begin{bmatrix}X_{s,t}.\,w_{k+1,n}\end{bmatrix}_{_{w_{s,t}}}=\begin{cases}
-r^{n-t+2}&\;\;\text{if $s=k+1$}\\
r^{n-s-2}&\;\;\text{if $t=k+1$}\\
(-r^3-r)(r^{k-s+n-t}-r^{k-s+n-t-2})&\;\;\text{if $s<k+1$ and $t>k+1$}\\
-r^2-\unsur{r^2}&\;\;\text{if $(s,t)=(k,k+1)$}\\
0&\;\;\text{otherwise}
\end{cases}$$

\hspace{-0.275in}
$\begin{bmatrix}X_{s,t}.\,w_{k,k+1}\end{bmatrix}_{_{w_{s,t}}}=\begin{cases}
-\unsur{r^{t-k+1}}&\;\;\,\qquad\qquad\qquad\qquad\qquad\;\;\;\;\text{if $s=k$ and $t>k+1$}\\
\unsur{r^{t-k-2}}&\;\;\,\qquad\qquad\qquad\qquad\qquad\;\;\;\;\text{if $s=k+1$ and $t>k+1$}\\
r^{k-s-1}&\;\;\,\qquad\qquad\qquad\qquad\qquad\;\;\;\;\text{if $t=k$ and $s<k$}\\
-r^{k-s+2}&\;\;\,\qquad\qquad\qquad\qquad\qquad\;\;\;\;\text{if $t=k+1$ and $s<k$}\\
0&\;\;\,\qquad\qquad\qquad\qquad\qquad\;\;\;\;\text{otherwise}
\end{cases}\\\\$

\noin From these equalities, we derive:
\begin{eqnarray}
\forall\,t>k+1,\,X_{k,t}.V_k&=&
(-r^3-r)(r^{n-t}-r^{n-t-2})+r\,r^{n-t+2}-r^{n-t-1}\notag\\
&=&0\\
\forall\,s<k,\,X_{s,k}.V_k&=&-r\,r^{n-s-2}+r^{n-k}\,r^{k-s-1}\notag\\
&=&0\\
\forall\,s<k,\,X_{s,k+1}.V_k&=&
r^{n-s-2}+(r^4+r^2)(r^{n-s-2}-r^{n-s-4})-r^{n-k}\,r^{k-s+2}\notag\\
&=&0\\
\forall\,t>k+1,\,X_{k+1,t}.V_k&=&-r^{n-t-2}+r^{n-k}\,r^{k-t+2}\notag\\
&=&0
\end{eqnarray}

\noin Equalities $(157)$, $(158)$, $(159)$ and $(160)$ show that all
the $X_{ij}$'s with $1\leq i<j\leq n-1$ annihilate the vectors
$V_k$'s, $1\leq k\leq n-2$. Thus, we have shown that the first
$\binom{n-1}{2}$ rows of the matrix $T(n)$ annihilate the vectors
$V_1,\dots,V_{n-2}$. And with the work from before, all the rows of
$T(n)$ annihilate in fact the vectors $V_1,\dots,V_{n-2}$. Thus,
these vectors belong to $K(n)$. Since we picked the vectors
$v_1,v_2,\dots,v_{k(n-1)}$ to form a basis of $K(n-1)$ and since
$K(n-1)$ is contained in $K(n)$, the latter vectors also belong to
$K(n)$. Finally, it is visible that all the participating vectors in
Claim $3$ are linearly independent. thus, we conclude that the claim
holds and Proposition $11$ as well.

\begin{Cor}
Let $n\geq 4$. Assume that $r^{n+1}\neq -1,\,r^{n+2}\neq
-1,\,\dots,\,r^{2n}\neq -1$. \\If there exists an irreducible
$\cil$-dimensional invariant subspace of $\V$, then $l=r$.
\end{Cor}

\textsc{Proof:} we proceed by induction on $n$. When $n=4$, if there
exists an irreducible $2$-dimensional invariant subspace of $\V$,
then $l=r$ without any further assumption on $r$ by Result $2$. When
$n=5$, if there exists an irreducible $5$-dimensional invariant
subspace of $\V$ then $l=r$ still without any further assumption on
$r$ by Result $1$. Let $n$ be an integer with $n\geq 6$. Let's name
$W$ the irreducible $\cil$-dimensional invariant subspace of $\V$.
Consider the intersection $W\cap\V^{(n-1)}$. Since
$\text{dim}\,W=\cil>n-1$ for any $n$ greater or equal to $6$, this
intersection is nontrivial. Hence we get $l\in\lb
r,-r^3,\unsur{r^{n-4}},-\unsur{r^{n-4}},\unsur{r^{2n-5}}\rb$ as
$\n^{(n-1)}$ is then reducible. Also, from $W\cap\V^{(n-1)}\subseteq
K(n-1)$, we derive the inequality on the dimensions over $F$:
$$k(n-1)\geq \cil-(n-1)=\frac{n^2-5n}{2}+1=\frac{(n-1)(n-4)}{2}-1\qquad\qquad (\mathcal{I})_n$$
We note that $\frac{n^2-5n}{2}+1> n-2$ is equivalent to $(n-1)(n-6)>
0$, which is itself equivalent to $\cil>2n-3$. Let's first deal with
the case $n\geq 7$, assuming that Corollary $9$ holds for $n=6$.
Then $l\in\lb
r,-r^3,-\unsur{r^{n-5}},\unsur{r^{n-5}},\unsur{r^{2n-7}}\rb$, as
$\n^{(n-2)}$ is reducible ($\dw>2n-3$). Then $l\in\lb r,-r^3\rb$. We
want to show that $l=r$. Suppose $l=-r^3$. Since $r^{2n}\neq -1$, we
have seen in Proposition $6$ that $K(n)$ is irreducible. Then
$K(n)=W$ and $k(n)=\cil$. If $r^{2(n-1)}\neq -1$, then an
application of Proposition $11$ with $k(n)=\cil$ yields:
$$k(n-1)\leq\cil-(n-2)=\frac{n^2-5n}{2}+2=\frac{(n-1)(n-4)}{2}$$
Also, if $r^{2(n-1)}\neq -1$ and $l=-r^3$, then by Proposition $6$,
$K(n-1)$ is irreducible. Suppose first $n\neq 9$. Then $n-1\neq 8$,
hence, $$n-2<\frac{(n-1)(n-4)}{2}-1\leq k(n-1)\leq
\frac{(n-1)(n-4)}{2}$$ forces in fact $k(n-1)=\frac{(n-1)(n-4)}{2}$,
as there is no irreducible invariant subspace of $\V^{(n-1)}$ of
dimension strictly between $(n-2)$ and $\frac{(n-1)(n-4)}{2}$. Also,
if $n=9$ \emph{i.e\/} $n-1=8$, then the inequality above reads
$7<19\leq k(8)\leq 20$. Since the degrees of the irreducible
representations of $\ih(8)$ are $1,7,14,20,21$ and degrees higher or
equal to $28$, the case $n=9$ is not different. Thus we get
$k(n-1)=\frac{(n-1)(n-4)}{2}$. Then by induction hypothesis, the
existence of an irreducible $\frac{(n-1)(n-4)}{2}$-dimensional
invariant subspace of $\V^{(n-1)}$ implies $l=r$. This is a
contradiction with $l=-r^3$. Hence our hypothesis $l=-r^3$ was
absurd and $l=r$ as announced. It remains to study the case $n=6$.
When $n=6$, we find another argument to claim that the intersection
$W\cap\V^{(4)}$ is non trivial. Indeed, it suffices to notice that
$\text{dim}\,W+\text{dim}\,\V^{(4)}=9+6=15=\text{dim}\,\V^{(6)}$.
Then, if the sum $W+\V^{(4)}$ is direct, it comes
$W\oplus\V^{(4)}=\V^{(6)}$. Acting with $e_5$ on both sides yields
$e_5.\V^{(6)}=0$, which is a contradiction. Hence $W\cap\V^{(4)}\neq
0$ and again in the case $n=6$ we must have $l\in\lb r,-r^3\rb$. In
particular, since $r^{10}\neq -1$, this implies that $K(5)$ is
irreducible. Further, the inequality $(\mathcal{I})_6$ reads
$k(5)\geq 4$. Thus, we have $k(5)\in\lb 4,5,6\rb$. Now $k(5)=4$ is
to eliminate as the existence of an irreducible $4$-dimensional
invariant subspace of $\V^{(5)}$ would force
$l\in\lb\unsur{r^2},-\unsur{r^2}\rb$. Hence $k(5)\in\lb 5,6\rb$. We
show that in any case this forces $l=r$. And indeed if $k(5)=5$, it
is automatic with Result $1$. Suppose now $k(5)=6$. We show that it
is not possible to have $l=-r^3$. If $l=-r^3$, since $r^{10}\neq
-1$, an application of Proposition $11$ with $n=6$ yields:
$$k(6)\geq k(5)+4$$
Hence it comes $k(6)\geq 6+4=10$. But since $l=-r^3$ and $r^{12}\neq
-1$, $K(6)$ is irreducible by Proposition $6$. Hence $K(6)=W$ and
$k(6)=9$: contradiction. Thus, we must have $l=r$. This ends the
proof of Corollary $9$. The next Corollary specifies the cases
$r^{2n}=-1$ and uses the result of Corollary $9$.

%\\ Suppose now that $r^{2(n-1)}=-1$. Then, by Proposition $6$, $K(n-1)$ is reducible and moreover $k(n-1)\in\lb %\frac{(n-2)(n-3)}{2},\frac{(n-2)(n-3)}{2}+1\rb$.
%Also since $r^{2(n-1)}=-1$, certainly $r^{2(n-2)}\neq -1$ (otherwise $r^2=1$). Then by Proposition $11$, we have:
%$$k(n-1)\geq k(n-2)+(n-3)$$ Thus $$k(n-2)\leq \frac{(n-2)(n-3)}{2}+1-(n-3)=$$

\begin{Cor}
Let $n$ be an integer with $n\geq 4$. If $r^{2n}=-1$ and there
exists an irreducible $\cil$-dimensional invariant subspace of $\V$,
then
$$\begin{array}{ccccc} \text{Either}& k(n)=\cil,& k(n-1)=\frac{(n-1)(n-4)}{2}&\text{and}&l=r\\
\text{Or}&
k(n)=\dbw,&k(n-1)=\frac{(n-2)(n-3)}{2}&\text{and}&\;\;l=-r^3\end{array}$$
\end{Cor}

\textsc{Proof:} Let's first assume $n\geq 6$. Following the proof of
Corollary $9$, the existence of an irreducible $\cil$-dimensional
invariant subspace of $\V$ forces $l\in\lb
r,-r^3,\unsur{r^{n-4}},-\unsur{r^{n-4}},\unsur{r^{2n-5}}\rb$ and
$l\in\lb r,
-r^3,\unsur{r^{n-5}},-\unsur{r^{n-5}},\unsur{r^{2n-7}}\rb$, which
forces in turn $l\in\lb r,-r^3\rb$. If $l=r$, we already know from
Proposition $10$ that $k(n-1)=\frac{(n-1)(n-4)}{2}$ and $k(n)=\cil$.
If $l=-r^3$, since $r^{2n}=-1$, there exists a one-dimensional
invariant subspace of $\V$, say $\V_1$. By hypothesis, there also
exists an irreducible $\cil$-dimensional invariant subspace, say
$\V_{\cil}$, of $\V$. The two vector spaces must be in direct sum as
$\V_{\cil}$ is irreducible. Hence we have
$\V_1\oplus\V_{\cil}\subseteq K(n)$. If this inclusion is strict,
then $k(n)>\dbw$ and $\V_1\oplus\V_{\cil}$ has a summand in $K(n)$
whose dimension is strictly less than $n-1$ (recall that
$\chl-\dbw=n-1$). This is a contradiction. Hence
$K(n)=\V_1\oplus\V_{\cil}$ and $k(n)=\dbw$. To complete the proof,
it remains to show that $k(n-1)=\frac{(n-2)(n-3)}{2}$. Since
$\dbw>n-1$, we have $K(n)\cap\V^{(n-1)}\neq 0$. Hence
$K(n)\cap\V^{(n-1)}\subseteq K(n-1)$. It follows on the dimensions
over $F$ that:
$$k(n-1)\geq \dbw-(n-1)=\frac{(n-1)(n-4)}{2}$$
Here comes the use of Corollary $9$: since $r^{2n}=-1$, we have
$r^n\neq -1,\,r^{n+1}\neq -1,\,\dots,\,r^{2(n-1)}\neq -1$. In
particular $r^{2(n-1)}\neq -1$ implies that $K(n-1)$ is irreducible.
If $k(n-1)=\frac{(n-1)(n-4)}{2}$, a licit application of Corollary
$9$ yields $l=r$, a contradiction with $l=-r^3$. Thus, the
inequality above can be bettered:
$$k(n-1)\geq\frac{(n-2)(n-3)}{2}$$
As for the other way, again since $r^{2n}=-1$, we have
$r^{2(n-1)}\neq -1$. Hence we may apply Proposition $11$. It
provides us with the inequality: $\dbw\geq k(n-1)+(n-2)$,
\emph{i.e\/}
$$k(n-1)\leq \frac{(n-2)(n-3)}{2}$$
Gathering the two inequalities finally yields:
$$k(n-1)=\frac{(n-2)(n-3)}{2}$$
Finally the Corollary is true when $n=5$ and in that case we know from before that $l=r$ and $k(5)=5$ and $k(4)=2$. As for $n=4$, we also have $l=r$ and $k(4)=2$. Moreover, it is true that $k(3)=0$ since for $l=r$ the representation $\n^{(3)}$ is irreducible by Theorem $7$.\\

Joining the results of Proposition $6$ and Corollary $10$ adds a bit
of information to Proposition $6$:

\begin{Cor}
Let $n$ be an integer with $n\geq 4$. If $l=-r^3$ and $r^{2n}=-1$,
then either there exists an irreducible $\cil$-dimensional invariant
subspace of $\V$ and $k(n)=\dbw$ or $k(n)>\dbw$
\end{Cor}

\textsc{Proof:} first the corollary holds for $n=4$ and $n=5$.
Indeed, for $n=4$, if $l=-r^3$, there does not exist any irreducible
$2$-dimensional invariant subspace of $\V$ (otherwise $l=r$ by
Result $2$). On the contrary, there exists an irreducible
$3$-dimensional invariant subspace of $\V$ by Theorem $6$ and there
exists a one-dimensional invariant subspace of $\V$ since
$l=-r^3=\unsur{r^5}$ ($r^8=-1$ by hypothesis). Those two invariant
subspaces are in direct sum. Hence $k(4)>3$. When $n=5$, there does
not exist any irreducible $5$-dimensional invariant subspace of $\V$
(otherwise $l=r$ by Result $1$). By Theorem $8$ there exists an
irreducible $6$-dimensional invariant subspace of $\V$. This
subspace must be in direct sum with the existing one-dimensional
one. Hence $k(5)>6$. Suppose now $n\geq 6$. If there exists an
irreducible $\cil$-dimensional invariant subspace of $\V$, then by
Corollary $10$, case $l=-r^3$, we have $k(n)=\dbw$. Suppose that
there does not exist any irreducible $\cil$-dimensional invariant
subspace of $\V$. Since $l=-r^3$ and $r^{2n}=-1$, we know from
Proposition $6$ that $K(n)$ is reducible and $k(n)\geq\dbw$. Since
with our assumptions there must exist a one-dimensional invariant
subspace $\V_1$ of $\V$, if we had $k(n)=\dbw$, then this
one-dimensional invariant subspace of $\V$ would have an irreducible
$\cil$-dimensional summand in $K(n)$, impossible by hypothesis.
Thus, $k(n)\geq\dbw+1$.
\begin{Cor}
Let $n$ be an integer with $n\geq 4$. If $l=-r^3$ and $r^{2n}=-1$,
then either there exists an irreducible $\cil$-dimensional invariant
subspace of $\V$ and $k(n)=\dbw$ or there exists an irreducible
$\dbw$-dimensional invariant subspace of $\V$ and $k(n)=1+\dbw$.
\end{Cor}
\textsc{Proof:} by Corollary $11$, it suffices to prove that if
there does not exist an irreducible $\cil$-dimensional invariant
subspace of $\V$, then there exists an irreducible
$\dbw$-dimensional invariant subspace of $\V$ and $k(n)=1+\dbw$.

Let's first deal with the case $n=4$. In this case, there does not
exist any irreducible $2$-dimensional invariant subspace of $\V$
(otherwise $l=r$, impossible with $l=-r^3$). Further, by choice of
$l$ and $r$, there exists an irreducible $3$-dimensional invariant
subspace of $\V$ and a one-dimensional invariant subspace of $\V$
and their sum is direct. Then by uniqueness of the one-dimensional
invariant subspace of $\V$, it is impossible to have $k(4)=5$. Hence
$k(4)=4$. This ends the case $n=4$.

When $n=5$, by choice of $l$ and $r$, we know that there exists an
irreducible $6$-dimensional invariant subspace of $\V$ (cf Theorem
$8$) and there exists a one-dimensional invariant subspace of $\V$.
Their sum is direct. It forbids $k(5)=8$ or $k(5)=9$. Hence
$k(5)=7$. Also, since $l\neq r$, there does not exist any
irreducible $5$-dimensional invariant subspace of $\V$ by Result
$1$, hence we are done with the case $n=5$.

Suppose now $n\geq 6$ and suppose that there does not exist any
irreducible $\cil$-dimensional invariant subspace of $\V$. By
Corollary $11$, we have $k(n)\geq 1+\dbw$. In particular $k(n)>n-1$,
hence the intersection $K(n)\cap\V^{(n-1)}$ is non-trivial. Further,
since $r^{2n}=-1$ by hypothesis, we have $r^{2(n-1)}\neq -1$. Hence
by Proposition $6$, $K(n-1)$ is irreducible. Thus we get
$K(n)\cap\V^{(n-1)}=K(n-1)$. From there, we have:
$$k(n-1)=k(n)+\frac{(n-1)(n-2)}{2}-\text{dim}(K(n)+\V^{(n-1)})$$
Hence, we get:
$$k(n-1)\geq k(n)-(n-1)$$
that we rewrite:
$$k(n)\leq k(n-1)+(n-1)$$
Since we already know that $k(n)\geq 1+\dbw$, it suffices to show
that $k(n-1)=\frac{(n-2)(n-3)}{2}$ to get the desired result. To
that aim, we show a lemma:

%In particular $k(n)>2n-3$, hence the intersection
%$K(n)\cap\V^{(n-2)}$ is non-trivial. Further, since $r^{2n}=-1$ by
%hypothesis, we have $r^{2(n-2)}\neq -1$. Hence by Proposition $6$,
%$K(n-2)$ is irreducible. Thus we get $K(n)\cap\V^{(n-2)}=K(n-2)$.
%From there, we have:
%$$k(n-2)=k(n)+\frac{(n-2)(n-3)}{2}-\text{dim}(K(n)+\V^{(n-2)})$$
%We notice that it is impossible to have $K(n)+\V^{(n-2)}=\V^{(n)}$.
%Indeed, if such an equality hold, we would get $e_{n-1}.\V^{(n)}=0$,
%which is absurd. It follows that $\text{dim}(K(n)+\V^{(n-2)})<\chl$.
%Hence we get the inequality
%$$k(n-2)>k(n)+\frac{(n-2)(n-3)}{2}-\chl$$
%which can still be rewritten as:
%$$k(n)-k(n-2)<2n-3$$
%Next, we notice that:
%$$\dbw+2-\frac{(n-3)(n-4)}{2}=2n-3$$
%Thus, if $k(n-2)=\frac{(n-3)(n-4)}{2}$, then it forces
%$k(n)=\dbw+1$. We will show that it is true with our assumptions
%that $k(n-2)=\frac{(n-3)(n-4)}{2}$. For that, we need a lemma.

\begin{lemma}
Let $n\geq 3$. If $r^{n+1}\neq -1,r^{n+2}\neq -1,\dots,r^{2n}\neq
-1$ and $l=-r^3$, then $k(n)=\dbw$.
\end{lemma}
\textsc{Proof of the lemma:} when $n=3$ and $l=-r^3$ and $r^6\neq
-1$, there exists a unique one-dimensional invariant subspace of
$\V$. Moreover, there does not exists any irreducible
$2$-dimensional invariant subspace as $l\not\in\lb -1,1\rb$. Thus we
have $k(3)=1$. When $n=4$ and $l=-r^3$, there exists an irreducible
$3$-dimensional invariant subspace of $\V$. Moreover, since $l\neq
r$, there does not exist any irreducible $2$-dimensional invariant
subspace of $\V$. Hence $k(4)\neq 5$. Also since $r^8\neq -1$, there
does not exist any one-dimensional invariant subspace of $\V$, hence
$k(4)\neq 4$. Thus, we have $k(4)=3$. Let's also do the case $n=5$.
By Theorem $8$, there exists an irreducible $6$-dimensional
invariant subspace of $\V$. Hence $k(5)\geq 6$. Since there does not
exist any one-dimensional invariant subspace of $\V$, $k(5)$ cannot
equal $7,8$ or $9$. Thus we have $k(5)=6$. Let $n\geq 6$. Under the
assumptions on $l$ and $r$, we know that there does not exist any
one-dimensional invariant subspace of $\V$ by Theorem $4$, there
does not exist any irreducible $(n-1)$-dimensional invariant
subspace of $\V$ by Theorem $5$ and there does not exist any
irreducible $\cil$-dimensional invariant subspace of $\V$ by
Corollary $9$. Moreover, by Proposition $6$ when $l=-r^3$ and
$r^{2n}\neq -1$, we know that $K(n)$ is irreducible. We recall that
when $n=6,7$ or $n\geq 9$, the irreducible representations of
$\ih(n)$ have degree $1$, $n-1$, $\cil$, $\dbw$ or a degree greater
than $\dbw$. Hence, for these values of $n$, we must have
$k(n)\geq\dbw$. When $n=8$, we must have $k(8)\in\lb 14,20,21\rb$.
Since $k(7)\geq 15$ and since $K(7)\subset K(8)$ by Proposition $9$
(as $r^{14}\neq -1$), it is impossible to have $k(8)=14$. Hence, the
case $n=8$ is not exceptional and we must have $k(8)=21$. In any
case, we have $k(n)\geq \dbw$. We will show conversely that
$k(n)\leq \dbw$. It is equivalent to show that the rank of the
matrix $T(n)$ is greater or equal to $n-1$. We have the lemma:
\begin{lemma} Let $n$ be an integer with $n\geq 5$. Suppose $r^{2n}\neq -1$. Then,
$$l=-r^3\Longrightarrow rk(T(n))\geq n-1$$
\end{lemma}
\textsc{Proof of the lemma:} we computed with Maple the determinant
of the submatrix of $T(5)$, composed of subcolumns $1,3,4,7$ and
subrows $1,3,4,7$. We define:
$$S(5):=\text{submatrix}(T(5),[1,3,4,7],[1,3,4,7])$$ in Maple notations. The value of the determinant of $S(5)$ is
$\frac{1+r^4+r^8+r^{12}+r^{16}}{r^8}$. For future reference, define
$\mathcal{R}_5=\C_5:=[1,3,4,7]$. Next, given $n\geq 6$, we
inductively build from $S(n-1)$ a submatrix $S(n)$ of $T(n)$ by
adding the subrow correponding to the action of $X_{n-1,n}$ and the
subcolumn corresponding to the vector $w_{4,n}$. Explicitly, we
have:
$$S(n):=\text{submatrix}(T(n),\mathcal{R}_{n-1}@[\binom{n-1}{2}+1]),\mathcal{C}_{n-1}@[\binom{n-1}{2}+(n-4)])$$

And on Figure $2$ is how the matrix looks like. We will show that
when $r^{2n}\neq -1$, the determinant of this matrix is nonzero. In
fact we have:
\begin{Prop}
$$\text{det}(S(n))=(-1)^{n+1}\Bigg[\frac{1+r^4+\dots+r^{4(n-1)}}{r^{8+\frac{n(n-5)}{2}}}\Bigg]$$
\end{Prop}

\noin\textsc{Proof of the Proposition:} by construction, there are
only two nonzero coefficients on each of the $(\binom{k}{2}+1)$-th
rows of the matrix $S(n)$ for $k=5,\dots,n-1$ and they are
respectively given by:

\begin{equation}
[X_{k,k+1}.w_{4,k}]_{_{w_{k,k+1}}}=\unsur{r^{k-5}}\qquad\text{by
$(SL)_{k-4}$}
\end{equation}
\begin{equation}
\;\;[X_{k,k+1}.w_{4,k+1}]_{_{w_{k,k+1}}}=-\unsur{r^{3}.r^{k-5}}\qquad\text{by
$(TL)_{k-4}$}
\end{equation}

\noin These coefficients are the ones corresponding to the columns
$\C_{\binom{k-1}{2}+(k-4)}$ and $\C_{\binom{k}{2}+(k-3)}$. In
particular, doing the operation
$$\C_{\binom{n-2}{2}+(n-5)}\leftarrow\C_{\binom{n-2}{2}+(n-5)}+r^3\,\C_{\binom{n-1}{2}+(n-4)}$$ on the columns makes a zero
appear on the last row, making all the coefficients of the last row
zero except $-\unsur{r^3.r^{n-6}}$. Hence the determinant of $S(n)$
is:

$$\text{det}(S(n))=-\unsur{r^3}\,\unsur{r^{n-6}}\,\text{det}(\tilde{S}(n-1))$$

\noin where $\tilde{S}(n-1)$ is obtained from $S(n-1)$ by replacing
$\C_{\binom{n-2}{2}+(n-5)}(\mathcal{R}_{n-1})$ with
$\C_{\binom{n-2}{2}+(n-5)}(\mathcal{R}_{n-1})+r^3\,\C_{\binom{n-1}{2}+(n-4)}(\mathcal{R}_{n-1})$.

\noin It comes:
\begin{multline*}\text{det}(\tilde{S}(n-1))=\text{det}(S(n-1))\\
+r^3\,\text{det}(\C_1(\Rrond_{n-1}),\dots,\C_{\binom{n-3}{2}+(n-6)}(\Rrond_{n-1}),\widehat{\C_{\binom{n-2}{2}+(n-5)}}(\Rrond_{n-1}),\C_{\binom{n-1}{2}+(n-4)}(\Rrond_{n-1}))\end{multline*}

\noin The second determinant in the sum above is:
$$(-1)(-\unsurr)(-\unsur{r^2})\dots(-\unsur{r^{n-7}})\,\text{det}([1,3,4,7],[1,3,4,\binom{n-1}{2}+(n-4)])$$

\noin We computed with Maple $\text{det}([1,3,4,7],[1,3,4,12])$ and
found the value $r^9$. Then,
$$\text{det}([1,3,4,7],[1,3,4,\binom{n-1}{2}+(n-4)])=r^{n+3}$$
\noin Thus, we get:
$$\text{det}(S(n))=-\unsur{r^3}\,\unsur{r^{n-6}}\,\text{det}\,S(n-1)-\unsur{r^{n-6}}\,(-1)^n\,\frac{r^{n+3}}{r^{\frac{(n-6)(n-7)}{2}}}$$

\noin Let's proceed by induction on $n$ and assume that Proposition
$12$ holds for $S(n-1)$. Then, replacing $\text{det}(S(n-1))$ by its
value yields the new equality:
$$\text{det}(S(n))=(-1)^{n+1}\,\unsur{r^{n-6}}\Bigg[\frac{1+r^4+\dots+r^{4(n-2)}}{r^{11+\frac{(n-1)(n-6)}{2}}}+\frac{r^{n+3}}{r^{\frac{(n-6)(n-7)}{2}}}\Bigg]$$
\noin And by reducing to the same denominator, we get:
$$\text{det}(S(n))=(-1)^{n+1}\Bigg[\frac{1+r^4+\dots+r^{4(n-1)}}{r^{8+\frac{n(n-5)}{2}}}\Bigg]$$

\noin As mentioned at the beginning of the proof of Lemma $10$, the
formula in Proposition $12$ holds for $n=5$. Furthermore, we deduce
from it the value for $\text{det}(S(6))$. Indeed, by adding to the
$7$-th column $r^3$ times the $12$-th column, we see that:
\begin{eqnarray*}
\text{det}(S(6))&=&-\unsur{r^3}\times\Big(\text{det}(S(5))+r^3\,\text{det}([1,3,4,7],[1,3,4,12])\Big)\\
         &=&-\unsur{r^3}\times\Big(\frac{1+r^4+r^8+r^{12}+r^{16}}{r^8}+r^3.\,r^9\Big)\\
         &=&-\;\frac{1+r^4+r^8+r^{12}+r^{16}+r^{20}}{r^{11}}
\end{eqnarray*}
So Proposition $12$ also holds for $n=6$. Then, by induction,
Proposition $12$ holds for every $n\geq 5$.

$\begin{array}{cc}\hspace{-1.2in}
\begin{array}{l}
\\
\\
\\
\\
\\
\\
\\
\\
\\
\\
X_{56}\\
\\
\\
\\
\\
X_{67}\\
%\\
%\\
%\\
%\\
%\\
\\
\\
\\
\\
\\
X_{n-1,n}\\
\\
\\
\\
\\
\\
\\
\\
\end{array}
& \begin{array}{l}\begin{array}{ccccccccccccccccccccccccc} &&&&&&&&&&\;\;\;\,w_{46}&&&&&w_{47}&&&&&&\nts\nts\nts\nts w_{4n}&&\\
&\!1&\!2&\!3&\,\,4&&\;\;\,\;7&&&10&\;\;12&&15&&&\,\!18
 &&\!\!\dots&\binom{n-1}{2}&&&\binom{n-1}{2}+(n-4)&&\!\!\!\nts\binom{n}{2}
\end{array}\\
$\begin{tabular}[c]{c|c|c|c|} $\begin{array}{cccccccccc}
\sq&&\sq&\sq&&&\sq&&&\\
&&&&&&&&&\\
\sq&&\sq&\sq&&&\sq&&&\\
\sq&&\sq&\sq&&&\sq&&&\\
&&&&&&&&&\\
&&&&&&&&&\\
\sq&&\sq&\sq&&&\sq&&&\\
&&&&&&&&&\\
&&&&&&&&&\\
&&&&&&&&&\\
\end{array}$
& $\begin{array}{ccccc} &0&&&\\
&&&& \\
&0&&& \\
&r&&& \\
&&&& \\
&&&& \\
&r&&& \\
&&&& \\&&&& \\&&&&
\end{array}$
& $\begin{array}{ccccccccc} &&0&&&&&&\\
&&&&&\\&&0&&&\\&&r^2&&&&&&\\
&&&&&&&&\\
&&&&&&&&\\
&&r^2&&&&&&\\
&&&&&&&&\\&&&&&&&&\\&&&&&&&&
\end{array}$
& $\begin{array}{cccccccccc}
&&&&0&&&&&\\
&&&&&&&&&\\&&&&0&&&&&\\&&&&r^{n-5}&&&&&\\
&&&&&&&&&\\&&&&&&&&&\\
&&&&r^{n-5}&&&&&\\
&&&&&&&&&\\&&&&&&&&&\\&&&&&&&&&\\
\end{array}$\\
\cline{1-1} \end{tabular}$\\
$\begin{tabular}[c]{c|c|c|}$\begin{array}{cccccccccc}
\,0&&\,0&\,\,0&&&1&&&\\
&&&&&&&&&\\
&&&&&&&&&\\
&&&&&&&&&\\
&&&&&&&&&\\
\end{array}$\hspace{-0.022in}
$\begin{array}{ccccc}
&&\!\!-\unsur{r^3}&&\hspace{-.01in}\\
&&&& \hspace{-.01in}\\
&&&&\hspace{-.01in}\\
&&&&\hspace{-.01in} \\
&&&&\hspace{-.01in}
\\
\end{array}$\hspace{0.085in}
& $\begin{array}{ccccccccc}
&&0&&&&&&\,\,\\
&&&&&&&&\,\;\\&&&&&&&&\,\;\\&&&&&&&&\,\,\\
&&&&&&&&\,\,\\
\end{array}$
&\hspace{-0.2in}$\begin{array}{ccccccccccc}
%&&&&\;\;\;\;\;&&&&&&\;\,\hspace{.005in}\\
&&&&\;\;\;\;\;0&&&&&&\;\,\hspace{.005in}\\
&&&&&&&&&&\\&&&&\;\;\;\;\;&&&&&&\\&&&&\;\;\;\;\;&&&&&&\;\,\hspace{.005in}\\
&&&&\;\;\;\;\;&&&&&&\;\,\hspace{.005in}
\end{array}$\\
\cline{1-1}
\end{tabular}$
\\
$\begin{tabular}[c]{c|c|} $\begin{array}{cccccccccccccccccccccccc}
%&&&&&&&&&&&&&&&&&&&&&&&\!\!\\
\,0&&\,\,0&\,\,0&&&0&&&&&\;\;\;\,\unsurr&&&&&&&\!\!\!-\unsur{r^3}\,\unsurr&&&&&\hspace{-.12in}\\
&&&&&&&&&&&&&&&&&&&&&&&\hspace{-.12in}\\
&&&&&&&&&&&&&&&&&&&\ddots&&&&\hspace{-.12in}\\
&&&&&&&&&&&&&&&&&&&&&&&\hspace{-.12in}\\
&&&&&&&&&&&&&&&&&&&&&&&\hspace{-.12in}\\
\end{array}$\hspace{.085in}& $\begin{array}{ccccccccccc}
%&&&&\,\,&&&&&&\hspace{-.045in}\\
&&&&\,\,\,0&&&&&&\hspace{-.045in}\\
&&&&\,\,&&&&&\;\\\!\!&&&&\,\,&&&&&\\\!\!&&&&\,\,&&&&&&\hspace{-.045in}\\
&&&&\,\,&&&&&\;\\\!\!&&&&\,\,&&&&&&\hspace{-.045in}\\
\end{array}$\\\cline{1-1}
\end{tabular}$\\
$\begin{tabular}[c]{c|}
$\begin{array}{cccccccccccccccccccccccccccccccccccccc}
%&&&&&&&&&&&&&&&&&&&&&&&&&&&&&&&&&&&&&\hspace{-.26in}\\
%&&&&&&&&&&&&&&&&&&&&&&&&&&&&&&&&&&&&&\hspace{-.26in}\\
%&&&&&&&&&&&&&&&&&&&&&&&&&&&&&&&&&&&&&\hspace{-.26in}\\
%&&&&&&&&&&&&&&&&&&&&&&&&&&&&&&&&&&&&&\hspace{-.26in}\\
\,0&&\,\,0&\,\,0&&\;\;0&&&&&&0&&&&&&&\;\,0&&&&&&&&&&&&&\!\!\!\!\!-\unsur{r^3}\,\unsur{r^{n-6}}&&&&&\hspace{-.26in}\\
&&&&&&&&&&&&&&&&&&&&&&&&&&&&&&&&&&&&&\hspace{-.26in}\\
&&&&&&&&&&&&&&&&&&&&&&&&&&&&&&&&&&&&&\hspace{-.26in}\\
&&&&&&&&&&&&&&&&&&&&&&&&&&&&&&&&&&&&&\hspace{-.26in}\\
&&&&&&&&&&&&&&&&&&&&&&&&&&&&&&&&&&&&&\hspace{-.26in}\\
&&&&&&&&&&&&&&&&&&&&&&&&&&&&&&&&&&&&&\hspace{-.26in}\\
&&&&&&&&&&&&&&&&&&&&&&&&&&&&&&&&&&&&&\hspace{-.26in}\\
&&&&&&&&&&&&&&&&&&&&&&&&&&&&&&&&&&&&&\hspace{-.26in}\\\end{array}$\\
\cline{1-1}
\end{tabular}$\\\\
Figure\;\;2
\end{array}
\end{array}\\$

This achieves the proof of Lemma $9$. Let's go back to the proof of
the Corollary. Since when $r^{2n}=-1$, we have $$r^n\neq
-1,\,r^{n+1}\neq -1,\,\dots,\,r^{2(n-1)}\neq -1,$$ by Lemma $9$, we
get $k(n-1)=\frac{(n-2)(n-3)}{2}$. As already explained above, it
follows that $k(n)=1+\dbw$. Then there exists an irreducible
$\dbw$-dimensional invariant subspace of $\V$. This ends the proof
of Corollary $12$. Thus, Proposition $6$ can be slightly bettered
and rewritten in the following way:
\begin{Prop}\hfill\\\\
Let $n$ be an integer with $n\geq 4$. \\
When $l=r$, $K(n)$ is always irreducible.\\\\
Let $n$ be an integer with $n\geq 3$.\\
When $l=-r^3$, there are two cases:\\\\
$1)\;\;r^{2n}\neq -1\;\; \text{and}\;\, K(n) \;\,\text{is irreducible}$\\
$2)\;\;r^{2n}=-1 \;\;\text{and}\;\, K(n)\;\, \text{is reducible}$.
Moreover, when $n\geq 4$,
\begin{list}{\texttt{$a)$}}{}\item Either there exists an irreducible $\cil$-dimensional invariant subspace, $K(n)$ is the direct sum of an irreducible $\cil$ dimensional invariant subspace and of the unique one-dimensional invariant subspace and\\ $k(n)=\dbw$.\end{list}
\begin{list}{\texttt{$b)$}}{}\item
Or there exists an irreducible $\dbw$-dimensional invariant subspace
$K(n)$ is the direct sum of an irreducible $\dbw$-dimensional
invariant subspace and of the unique one-dimensional invariant
subspace and \\$k(n)=1+\dbw$.
\end{list}
\textit{(Case $n=3$)} There exists exactly two one-dimensional
invariant subspaces of $\V$ and $K(3)$ is the direct sum of these
two one-dimensional invariant subspaces.
\end{Prop}
\noin In the case $n=5$, we have a Corollary of this Proposition:
\begin{Cor}
If $l=-r^3$ and if $r^{10}=-1$, then:
$$K(5)+\V^{(4)}=\V^{(5)}$$
\end{Cor}
\noin\textsc{Proof:} Suppose that $r^{10}=-1$. Then $r^8\neq -1$, so
$K(4)$ is irreducible. If $K(5)+\V^{(4)}$ is strictly contained in
$\V^{(5)}$, then its dimension as a vector space over $F$ is less
than $10$. Then it comes:
$$\text{dim}(K(5)\cap\V^{(4)})=k(5)+6-\text{dim}(K(5)+\V^{(4)})\geq
k(5)+6-9=k(5)-3$$ Since by Proposition $13$, we know that
$k(5)\in\lb 6,7\rb$, the fact that $k(5)>4$ implies that the
intersection $K(5)\cap\V^{(4)}$ is nonzero. Thus by irreducibility
of $K(4)$, we have $K(5)\cap\V^{(4)}=K(4)$ and the inequality above
reads:
$$k(4)\geq k(5)-3$$
Moreover, when $l=-r^3$, there exists an irreducible $3$-dimensional
invariant subspace of $\V^{(4)}$. Since $K(4)$ is irreducible, it
must be $K(4)$. Hence $k(4)=3$. Thus, the inequality above becomes:
$$k(5)\leq 6$$
Then $k(5)=6$. Then by point $2)a)$ of Proposition $13$, there must
exist an irreducible $5$-dimensional invariant subspace of
$\V^{(5)}$. But this forces $l=r$ by Result $1$: a contradiction.
Thus, the vector spaces $K(5)+\V^{(4)}$ and $\V^{(5)}$ have the same
dimension and they are actually equal.

%Let's now deal with the case $r^{10}\neq -1$. Then, by Lemma $10$, we have $k(5)\leq 6$. Moreover, since $K(5)$ is irreducible, it has dimension $1$, %$4$, $5$ or $6$. Since $r^{10}\neq -1$, dimension $1$ is to be excluded. Because $l=-r^3$ and thus $l\not\in\lb \unsur{r^2},-\unsur{r^2},r\rb$, it is %not possible to have $k(5)\in\lb 4,5\rb$ either. This leaves the only possibility $k(5)=6$ for the dimension of $k(5)$. Next, by the same arguments as %in the first case above, we see that $k(4)=6+6-\text{dim}(K(5)+\V^{(4)})$, id est
%$$\text{dim}(K(5)+\V^{(4)})=12-k(4)$$
%By inspection, $k(4)$ cannot be five. So $k(4)\in\lb 3,4\rb$. Thus we get
\noin In fact Corollary $13$ generalizes to each $n$ by noticing
that the linearly independent set of vectors
$$\boxed{S:=\V^{(n-1)}\cup\lb V_1,\,\dots,\,V_{n-2}\rb}$$ of $K(n)+\V^{(n-1)}$ of cardinality $\chl-1$ does not span the vector space $K(n)+\V^{(n-1)}$ when $r^{2n}=-1$.
Explicitly, we will prove the following Proposition:
\begin{Prop}
Let $n$ be an integer with $n\geq 5$. Suppose $l=-r^3$ and
$r^{2n}=-1$. Then,
$$K(n)+\V^{(n-1)}=\V^{(n)}$$
\end{Prop}
\noin \textsc{Proof.} Suppose that
$K(n)+\V^{(n-1)}=\text{Span}_{F}(V_1,\dots,V_{n-2})\oplus
\V^{(n-1)}$. If $l=-r^3$ and $r^{2n}=-1$, then $l=\unsur{r^{2n-3}}$.
Hence, there exists an irreducible $1$-dimensional invariant
subspace of $\V^{(n)}$. Moreover, by Theorem $4$, it is spanned by
$$\sum_{1\leq s<t\leq n}r^{s+t}\,w_{st}$$ Thus, if the equality
above holds, the vector
$$w_{1,n}+r\,w_{2,n}+r^2\,w_{3,n}+\dots+r^{n-2}\,w_{n-1,n}$$ must be
a linear combination with coefficients in $F$ of the $(n-2)$ vectors
$w_{k+1,n}-r\,w_{k,n}$ where $k=1,\dots,n-2$. A direct verification
shows right away that this is impossible.
Then the set $S$ is a linearly independent set of $K(n)+\V^{(n-1)}$ of cardinality $\chl-1$, which does not span $K(n)+\V^{(n-1)}$. This shows that the dimension of $K(n)+\V^{(n-1)}$ must be greater than $\chl-1$, hence must in fact equal $\chl$, the dimension of $\V^{(n)}$. Thus, $K(n)+\V^{(n-1)}=\V^{(n)}$, as announced.\\\\
Proposition $14$ now allows us to give a more accurate version of
Proposition $13$, point $2)$, as it shows that point $a)$ cannot
occur.
\begin{Prop}\hfill\\\\
Let $n$ be an integer with $n\geq 4$. \\
When $l=r$, $K(n)$ is always irreducible.\\\\
Let $n$ be an integer with $n\geq 3$.\\
When $l=-r^3$, there are two cases:\\\\
$1)\;\;r^{2n}\neq -1\;\; \text{and}\;\, K(n) \;\,\text{is irreducible}$\\
$2)\;\;r^{2n}=-1 \;\;\text{and}\;\, K(n)\;\, \text{is reducible}$.
\\\begin{list}{\texttt{}}{}\item Moreover, when $n\geq 4$, $K(n)$ is
the direct sum of an irreducible $\dbw$-dimensional invariant
subspace and of the unique one-dimensional invariant subspace and
$k(n)=1+\dbw$.\item \textit{(Case $n=3$)} There exists exactly two
one-dimensional invariant subspaces of $\V$ and $K(3)$ is the direct
sum of these two one-dimensional invariant subspaces.\end{list}
\end{Prop}
\textsc{Proof:} Point $1)$ was already proven in the proof of
Proposition $13$ or its former version, and there is nothing more to
add. Hence, let's prove the new version of point $2)$. Assume first
$n\geq 5$. Suppose that $l=-r^3$ and $r^{2n}=-1$. By Proposition
$13$, point $1)$, $K(n-1)$ is irreducible, as under these
assumptions on $l$ and $r$, we have $l=-r^3$ and $r^{2(n-1)}\neq
-1$. Moreover, by Proposition $13$ point $2)$ this time, $k(n)$ is
"big enough" so that $K(n)\cap\V^{(n-1)}$ cannot be trivial.
Explicitly we have $\dbw>n-1$ for all $n\geq 5$.
$$\begin{array}{cccc}\text{Hence we have:} & K(n)\cap\V^{(n-1)}&=&K(n-1)\\
 \text{By Proposition $14$, we also have:}& K(n)+\V^{(n-1)}&=&\V^{(n)}\end{array}$$
Both equalities yields the equality on the dimensions:
\begin{eqnarray*}k(n-1)&=&k(n)+\dbw-\chl\\
                       &=&k(n)-(n-1)
\end{eqnarray*}
Thus, we get: \begin{equation}k(n)=k(n-1)+(n-1)\end{equation}
Following the result of Proposition $13$, only two cases are
possible: $2)a)$ or $2)b)$. Suppose $2)a)$ holds. Then there exists
an irreducible $\cil$-dimensional invariant subspace of $\V^{(n)}$.
By Corollary $10$ we must have
$$k(n)=\dbw\qquad\&\qquad k(n-1)=\frac{(n-2)(n-3)}{2},$$
a contradiction with $(163)$. Hence the hypothesis $2)a)$ was absurd and $2)b)$ holds: there exists an irreducible $\dbw$-dimensional invariant subspace of $\V^{(n)}$ and $k(n)=1+\dbw$. To achieve the proof, it remains to deal with the case $n=4$. In that case, there exists a one-dimensional invariant subspace and an irreducible $3$-dimensional one, and they must be in direct sum. Hence, again, situation $2)b)$ holds. \\
In turn, Corollary $10$ can be rewritten:
\begin{Cor}
Let $n$ be an integer with $n\geq 4$. Suppose $r^{2n}=-1$. If there
exists an irreducible $\cil$-dimensional invariant subspace of $\V$,
then $l=r$.
\end{Cor}
\subsection{A Proof of Theorems $C$ and $D$}
We are now in a position to give a complete characterization for the
irreducible representations. We gather our two main results in the
following theorems:
\begin{thm}
Let $n$ be an integer with $n\geq 4$. There exists an irreducible
$\cil$-dimensional invariant subspace of $\V$ if and only if $l=r$.
\end{thm}
\begin{thm}
Let $n$ be an integer with $n\geq 5$.
There exists an irreducible $\dbw$-dimensional invariant subspace of $\V$ if and only if $l=-r^3$.\\
\end{thm}

Before we start the joint proof of these two theorems, let's gather
some known facts from earlier. In chronological order, we have the
following results:
\begin{itemize}
\item For $n\geq 4$: if $l=r$, then there exists an irreducible $\cil$-dimensional invariant subspace of $\V$.\hfill \textit{(Corollary $6$)}
\item For $n\geq 5$: if there exists an irreducible $\dbw$-dimensional invariant subspace of $\V$, then $l=-r^3$.\hfill\textit{(Corollary$7$)}
\item For $n\geq 4$: $\V$ cannot contain both an irreducible $\cil$-dimensional invariant subspace and an irreducible $\dbw$-dimensional one. \\    \textit{(Corollary $8$)}
\item For $n\geq 6$ and $r^{2(n-1)}\neq -1$: if $l=-r^3$, then $k(n)\geq k(n-1)+(n-2)$. \textit{(Proposition $11$)}
\item For $n\geq 4$ and $r^{n+1}\neq -1,\,r^{n+2}\neq -1,\,\dots,\,r^{2n}\neq -1$: if there exists an irreducible $\cil$-dimensional invariant subspace of $\V$ then $l=r$.\hfill\textit{(Corollary $9$)}
    \item For $n\geq 5$ and $r^{2n}\neq -1$: if $l=-r^3$, then $k(n)\leq\dbw$.\\\textit{(Lemma $10$)}
    \item For $n\geq 3$ and $r^{2n}\neq -1$: if $l=-r^3$, then $K(n)$ is irreducible \\ \textit{(Proposition $15$)}
    \item For $n\geq 4$ and $r^{2n}=-1$: if $l=-r^3$, then there exists an irreducible $\dbw$-invariant subspace of $\V$. Moreover, $k(n)=\dbw+1$.\hfill\textit{(Proposition $15$)}
\item For $n\geq 4$ and $r^{2n}=-1$: if there exists an irreducible $\cil$-dimensional invariant subspace of $\V$, then $l=r$.\hfill\textit{(Corollary $14$)}
\end{itemize}

\noin\textsc{The proof itself:} Given $n\geq 5$, it remains to show that if $r^{2n}\neq -1$ and $r^{2s}=-1$ for some integer $\frac{n+1}{2}\leq s\leq n-1$, then the following two statements hold:\\\\
\textit{statement $1$}: if there exists an irreducible $\cil$-dimensional invariant subspace of $\V$, then $l=r$. $\,(S1)$\\
\textit{statement $2$}: if $l=-r^3$, then there exists an irreducible $\dbw$-dimensional invariant subspace of $\V$. $\,(S2)$\\

We will prove by induction that:
$$\forall n\geq 5, (\mathcal{P}_n)$$
where $$(\mathcal{P}_n): \text{if $r^{2n}\neq -1$ and $r^{2s}=-1$,
some integer $\frac{n+1}{2}\leq s\leq n-1$, then $(S1)$ and
$(S2)$}$$

First, $(\matp_5)$ holds: $(S1)$ is true by Result $1$; $(S2)$ is
true by Theorem $8$. Let $n$ be an integer with $n\geq 6$ and
suppose that $(\matp_k)$ holds for all $5\leq k\leq n-1$. Let's
first deal with the case $r^{2s}=-1$ for some integer
$\frac{n+1}{2}\leq s<n-1$. In particular, we have $r^{2(s+1)}\neq
-1$ and $r^{2s}=-1$ with $s+1<n$. By induction hypothesis,
$(\matp_{s+1})$ then holds. Since we assume $l=-r^3$ in $(S2)$,
there exists an irreducible $\frac{s(s-1)}{2}$-dimensional invariant
subspace of $\V^{(s+1)}$. Moreover, since for $l=-r^3$ and
$r^{2(s+1)}\neq -1$, $K(s+1)$ is irreducible, we get
$k(s+1)=\frac{s(s-1)}{2}$. Then we show by induction on $l$ that:
\begin{equation}\forall s+2\leq l\leq n,\,k(l)=\frac{(l-1)(l-2)}{2}\end{equation}
For $l=s+2$, since $r^{s+1}\neq -1$ and $s\geq 4$, we may apply
point number $4$ with $n=s+2$. It yields:
$$k(s+2)\geq k(s+1)+s=\frac{s(s-1)}{2}+s=\frac{s(s+1)}{2}$$
Moreover, by point number $6$ above, with $n=s+2$, we also have:
$$k(s+2)\leq \frac{s(s+1)}{2}$$ It follows that
$k(s+2)=\frac{s(s+1)}{2}$. Let $l$ be an integer with $s+3\leq l\leq
n$ and suppose equation $(164)$ holds for the integer $l-1$. Again,
since $r^{2(l-1)}\neq -1$, point number $4$ forces:
$$k(l)\geq k(l-1)+(l-2)=\frac{(l-2)(l-3)}{2}+(l-2)=\frac{(l-1)(l-2)}{2}$$
And since $r^{2l}\neq -1$, point $6$ again forces: $$k(l)\leq
\frac{(l-1)(l-2)}{2}$$ Thus, equation $(164)$ holds for each
$s+1\leq l\leq n$ and in particular holds for $n$. Then $K(n)$ is
irreducible and $\dbw$-dimensional, so that $(S2)$ holds. Now, if
there exists an irreducible $\cil$-dimensional invariant subspace of
$\V$, then $l$ cannot equal $-r^3$ by $(S2)$ and point number $3$.
Also, if $k(n)\geq\cil$, then $K(n)\cap\V^{(n-1)}\neq\lb 0\rb$,
otherwise $k(n)\leq n-1$, but $\cil\geq n$ for every $n\geq 5$.
Assume first $n\geq 7$. Then, if $k(n)\geq \cil$, then
$K(n)\cap\V^{(n-2)}\neq\lb 0\rb$, otherwise $k(n)\leq 2n-3$, but
$\cil>2n-3$ as soon as $n\geq 7$. Then it comes:
$$l\in \Big\lb r,\unsur{r^{2n-5}},\unsur{r^{n-4}},-\unsur{r^{n-4}}\Big\rb\qquad\&\qquad l\in\Big\lb r,\unsur{r^{2n-7}},\unsur{r^{n-5}},-\unsur{r^{n-5}}\Big\rb$$
This only leaves the possibility $l=r$. It remains to deal with the
case $n=6$. When $n=6$, $s$ must be $4$ and the only possibilities
for $l$ are $l=r$ or $l=\unsur{r^9}$. We will show that the second
possibility for $l$ is to be excluded. And indeed, if there exists a
one-dimensional invariant subspace of $\V^{(6)}$, and an irreducible
$9$-dimensional one, these must be in direct sum, which forces
$k(6)\geq 10$ and in fact $k(6)=10$. Since $\unsur{r^9}$ must equal
$\unsur{r^2}$ or $-\unsur{r^2}$, there must exist a unique
irreducible $4$-dimensional submodule of $\V^{(5)}$ that is the only
submodule of $\V^{(5)}$. Thus, $k(5)=4$. Now, from the inclusion
$$K(6)\cap\V^{(5)}\subseteq K(5),$$
we derive on the dimensions:
$$k(5)\geq k(6)-5,$$
so that $$k(6)\leq k(5)+5=4+5=9,$$
a contradiction with $k(6)=10$ as mentioned above. Hence the case $n=6$ is no exception and $(S1)$ also holds in that case. \\

Suppose now $r^{2n}\neq -1$ and $r^{2(n-1)}=-1$ and let's show
$(S1)$ and $(S2)$ under these assumptions. Let's first do it for
$n=6$. We assume that $r^{12}\neq -1$ and $r^{10}=-1$ and try and
show $(S1)$. Suppose there exists an irreducible $9$-dimensional
invariant subspace of $\V^{(6)}$ and suppose $l=-r^3$. Let's first
determine $k(5)$ and $k(6)$. Since $l=-r^3$ and $r^{10}=-1$, there
exists a one-dimensional invariant subspace of $\V^{(5)}$. Moreover,
since $l=-r^3$, there also exists an irreducible $6$-dimensional
invariant subspace of $\V^{(5)}$. Moreover, there cannot exist any
irreducible $4$-dimensional invariant subspace of $\V^{(5)}$ when
$l=-r^3$ and there cannot exist any irreducible $5$-dimensional
invariant subspace as well. Hence we have $k(5)=7$. Since
$r^{12}\neq -1$ and $l=-r^3$, we also know from points $6$ and $7$
that $k(6)\in\lb 9,10\rb$. Since there exists an irreducible
$9$-dimensional invariant subspace of $\V^{(6)}$ by hypothesis, if
$k(6)$ were equal to $10$, there would also exist a one-dimensional
submodule of $\V^{(6)}$, which would force $l=\unsur{r^9}$. But when
$l=-r^3$ and $r^{12}\neq -1$, this is impossible. Hence we have
$k(6)=9$. Consider now the intersection $K(6)\cap\V^{(5)}$. The
$\ih(5)$-module $K(6)\cap\V^{(5)}$ is contained in
$K(6)\da_{\ih(5)}$. By semisimplicity of $\ih(5)$, there exists an
$\ih(5)$-submodule $S$ of $K(6)\da_{\ih(5)}$ which is a summand for
$K(6)\cap\V^{(5)}$:
\begin{equation}K(6)\cap\V^{(5)}\oplus S=K(6)\da_{\ih(5)}\end{equation}
Let's study the dimension of $K(6)\cap\V^{(5)}$. First, we have the
inequality:
$$\di (K(6)\cap\V^{(5)})\geq 9+10-15=4$$
If the dimension of $K(6)\cap\V^{(5)}$ were $4$, $K(6)\cap\V^{(5)}$
would have a $3$-dimensional summand in $K(5)$. This is impossible;
if the dimension of $K(6)\cap\V^{(5)}$ were $5$, $K(6)\cap\V^{(5)}$
would have a $2$-dimensional summand in $K(5)$. This is impossible
by uniqueness of the one-dimensional invariant subpsace of
$\V^{(5)}$. If now $\di (K(5)\cap\V^{(5)})=6$, then by $(165)$ the
dimension of $S$ would be $3$. We show that this is impossible.
\begin{lemma}
Under the hypothesis $l=-r^3$ and $r^{10}=-1$, it is impossible to
have:
$$\left\lb\begin{array}{l}
S\subseteq K(6)\\
\text{$S$ is a $\ih(5)$-module}\\
\di\,S=3
\end{array}\right.$$
\end{lemma}
\noin\textsc{Proof of the lemma:} it suffices to show it when
$\di\,S=1$. We leave the proof for later in a more general setting
(cf Lemma $13$).\\\\
Assuming the lemma holds, the dimension of $K(6)\cap\V^{(5)}$ must
hence be $7$. But then
$$K(6)\cap\V^{(5)}=K(5),$$ which implies in particular:
$$K(5)\subseteq K(6)$$
But if $K(5)\subseteq K(6)$, the inequality of point $4$ above
becomes true, although $r^{10}=-1$ and we must have:
$$k(6)\geq k(5)+4$$
With $k(6)=9$ and $k(5)=7$, this inequality yields a contradiction.\\\\
\underline{\textit{Partial conclusion}}: we have proven that if there exists an irreducible $9$-dimensional invariant subspace of $\V^{(6)}$, then $l$ cannot equal $-r^3$. Then $l\in\lb\unsur{r^9},r\rb$.\\
We have the general lemma:
\begin{lemma}
Let $n$ be an integer with $n\geq 5$. If $r^{2(n-1)}=-1$ and there
exists an irreducible $\cil$-dimensional invariant subspace of $\V$,
then $l\in\lb r,-r^3\rb$.
\end{lemma}
\noin\textsc{Proof of the lemma:} Let's denote by $W$ the
irreducible $\cil$-dimensional invariant subspace of $\V$. For
$n\geq 5$, we have $\cil>n-1$, hence $W\cap\V^{(n-1)}$ is not
trivial and so $l\in\lb
r,-r^3,\unsur{r^{2n-5}},\unsur{r^{n-4}},-\unsur{r^{n-4}}\rb$. Also
since $\n^{(n)}$ is reducible we have: $l\in\lb
r,-r^3,\unsur{2n-3},\unsur{r^{n-3}},-\unsur{r^{n-3}}\rb$. Thanks to
the hypothesis $r^{2(n-1)}=-1$, it is impossible to have
$$\unsur{r^{2n-3}}=\frac{\epsilon}{r^{n-4}}$$
as otherwise $r^{2(n+1)}=1$ and $r^{2(n-1)}=-1$ would force $r^2=-1$, which is excluded. Since it is also impossible to have $\unsur{r^{n-3}}=\frac{\epsilon}{r^{2n-5}}$ and $\unsur{r^{n-3}}=\frac{\epsilon}{r^{n-4}}$, we see that $l\in\lb r,-r^3\rb$.\\
In fact, we have the more general proposition:
\begin{Prop}
Let $n$ be an integer with $n\geq 5$. Lemma $12$ holds even without
the assumption $r^{2(n-1)}=-1$: if there exists an irreducible
$\cil$-dimensional invariant subspace of $\V$, then $l\in\lb
r,-r^3\rb$.
\end{Prop}
\noin\textsc{Proof of the proposition:} by an argument repeated many
times in the past, the existence of an irreducible
$\cil$-dimensional invariant subspace of $\V$ implies that $l\in\lb
r,-r^3,\unsur{r^{2n-3}}\rb$. If $l=\unsur{r^{2n-3}}$, then there
must exist a one-dimensional invariant subspace of $\V$. Denote it
by $\V_1$ and denote the irreducible $\cil$-dimensional invariant
subspace of $\V$ by $\V_{\cil}$.
$$\text{Let}\;\tilde{\V}:=\V_1\oplus\V_{\cil}$$
$\tilde{\V}$ is $\dbw$-dimensional. The equation $$\dbw>2n-3$$ is
equivalent to $$n^2-7n+8>0,$$ thus holds for every $n$. It follows
that $\tilde{\V}\cap\V^{(n-2)}$ is not trivial, so that $l\in\lb
r,-r^3,\unsur{r^{2n-7}},\unsur{r^{n-5}},-\unsur{r^{n-5}}\rb$ and
$l\in\lb
r,-r^3,\unsur{r^{2n-5}},\unsur{r^{n-4}},-\unsur{r^{n-4}}\rb$, which
leaves the only possibilities $l\in\lb r,-r^3\rb$ for $l$.

Let's go back to our case $n=6$. Applying the proposition with $n=6$ excludes the possibility $l=\unsur{r^9}$. Thus, we have shown $(S1)$ when $n=6$ and $r^{10}=-1$ but $r^{12}\neq -1$. From there we easily deduce $(S2)$ under the same assumptions. Indeed, if $l=-r^3$, the representation $\n^{(6)}$ is reducible, hence $\V^{(6)}$ must have an irreducible submodule, which is also an $\ih(6)$-module by the old lemma $6$. The irreducible representations of $\ih(6)$ have degrees $1,5,9$ or $10$. The hypothesis $r^{12}\neq -1$ forbids to have $l=-r^3=\unsur{r^9}$. So there cannot exist any one-dimensional invariant subspace of $\V^{(6)}$. Since it is also not possible to have $-r^3\in\lb\unsur{r^3},-\unsur{r^3}\rb$, there does not exist any irreducible $5$-dimensional invariant subspace of $\V^{(6)}$ by Theorem $5$. Moreover, the existence of an irreducible $9$-dimensional submodule of $\V^{(6)}$ would force $l=r$ as we just saw in the proof of $(S1)$. Then the only remaining possibility is that there exists an irreducible $10$-dimensional invariant subspace of $\V^{(6)}$. Hence $(S2)$ holds for $n=6$, $r^{10}=-1$ and $r^{12}\neq -1$. We note that by showing the statements $(S1)$ and $(S2)$ under these conditions and together with all the previous considerations, we have actually shown the theorems $9$ and $10$ for $n=6$.\\

To finish the proof of $(\matp_n)$ by induction, let's go back to
the general case under the assumptions $r^{2n}\neq -1$ and
$r^{2(n-1)}=-1$.
%After a brief thought, we see that our new induction hypothesis is Theorem $9$ and $10$ for $n-1$.
First let's deal with the proof of $(S1)$. Suppose there exists an
irreducible $\cil$-dimensional invariant subspace of $\V^{(n)}$. By
Proposition $16$, we know that $l\in\lb r,-r^3\rb$ and we want to
show that $l=r$. We will follow the same path as in the case $n=6$.
Suppose $l=-r^3$.
%When $l=-r^3$, we have the following informations about $\V^{(n-1)}$:
%\begin{list}{\texttt{(i)}}{}
%\item There exists a one-dimensional invariant subspace.
%\end{list}
%\begin{list}{\texttt{(ii)}}{}
%\item There exists an irreducible $\frac{(n-2)(n-3)}{2}$-dimensional invariant subspace.
%\end{list}
%\begin{list}{\texttt{(iii)}}{}
%\item There cannot exist any $\frac{(n-1)(n-4)}{2}$-dimensional invariant subspace.
%\end{list}
%\begin{list}{\texttt{(i$\upsilon$)}}{}
%\item There cannot exist any $(n-2)$-dimensional invariant subspace.
%\end{list}
Since $l=-r^3$ and $r^{2(n-1)}=-1$, we have by point number $8$ that
$k(n-1)=1+\frac{(n-2)(n-3)}{2}$. With $l=-r^3$ and $r^{2n}\neq -1$,
we also know from point $6$ and $7$ that $k(n)\in\lb \cil,\dbw\rb$
(we will see later on in the proof that this is still true when
$n=8$ and admit it for now). Since there exists an irreducible
$\cil$-dimensional invariant subspace and there does not exist any
one-dimensional invariant subspace, we must have $k(n)=\cil$. Let
$S$ be a summand for $K(n)\cap\V^{(n-1)}$ in $K(n)\da_{\ih(n-1)}$:
\begin{equation}K(n)\cap\V^{(n-1)}\oplus S=K(n)\da_{\ih(n-1)}\end{equation}
We have $$\di(K(n)\cap\V^{(n-1)})\geq\frac{(n-1)(n-4)}{2}-1$$ Then
we must have $$\di(K(n)\cap\V^{(n-1)})\geq \frac{(n-2)(n-3)}{2}$$
Moreover,
$$K(n)\cap\V^{(n-1)}\subseteq K(n-1)$$
implies that $$\di(K(n)\cap\V^{(n-1)})\leq 1+\frac{(n-2)(n-3)}{2}$$
If $\di(K(n)\cap\V^{(n-1)})=1+\frac{(n-2)(n-3)}{2}$, we get
$K(n)\cap\V^{(n-1)}=K(n-1)$ which would imply that $K(n-1)\subseteq
K(n)$. Then it comes:
$$k(n)\geq k(n-1)+(n-2)=\frac{(n-2)(n-3)}{2}+(n-2)+1=\frac{(n-2)(n-1)}{2}+1,$$
a contradiction. Hence we must have
$\di(K(n)\cap\V^{(n-1)})=\frac{(n-2)(n-3)}{2}$. But then, from
equation $(166)$, the $\ih(n-1)$-module $S$ would have dimension
$n-3$. By James'result in Proposition $3$ of the thesis, $S$ must
then contain a one-dimensional $\ih(n-1)$-submodule, say $U$. When
$n=7$, we also use the fact that there is no irreducible
representation of $\ih(6)$ of degree between $1$ and $5$.
\begin{lemma}
Let $n$ be an integer with $n\geq 5$. Suppose $l=-r^3$ and
$r^{2(n-1)}=-1$. In $K(n)$, there does not exist any one-dimensional
$\ih(n-1)$-module.
\end{lemma}
\textsc{Proof of the lemma:} suppose such a module $U$ exists and
let $$u=\sum_{1\leq i<j\leq n} \mu_{ij}\,w_{ij}$$ be a spanning
vector of $U$ over $F$. By the same arguments as in the proof of
Theorem $4$, we must have: $$\forall 1\leq i\leq
n-2,\,\n_i(u)=\la\,u\;\;\text{where}\;\;\la\in\big\lb
r,-\unsurr\big\rb$$ It follows from these relations that for every
node $i$ with $1\leq i\leq n-2$, we have:
\begin{eqnarray}
\forall k\geq i+2,\,\mu_{i+1,k}&=&\la\,\mu_{i,k}\\
\forall l\leq i-1,\,\mu_{l,i+1}&=&\la\,\mu_{l,i}
\end{eqnarray}

\noin From there, we see that if one of the coefficients $\mu_{st}$,
some $1\leq s<t\leq n-1$, is zero then all of the coefficients
$\mu_{ij}$ for $1\leq i<j\leq n-1$ are zero. Suppose we are in this
situation. Then $u$ reduces to:

$$u=\sum_{i=1}^{n-1}\mu_{in}\,w_{in}$$

But we have:
\begin{eqnarray*}
e_1.\,w_{1n}&=&-r^n\,w_{12}\\
e_1.\,w_{2n}&=&r^{n-3}\,w_{12}\\
e_1.\,w_{jn}&=&0\qquad\qquad\forall 3\leq j\leq n-1
\end{eqnarray*}

Then $u$ would not be annihilated by $e_1$, which is impossible. If
none of the coefficients for the $w_{ij}$'s with $1\leq i<j\leq n-1$
are zero, in particular $\mu_{34}$ is nonzero. Then by the same
argument as in the proof of Theorem $4$, case $n\geq 4$, the
coefficient $\la$ must be $r$ and not $-\unsurr$. Thus, $u$ can be
written as:

$$u=\sum_{1\leq i<j\leq n-1}
r^{i+j}\,w_{ij}+\mu\,\sum_{i=1}^{n-1}\,r^{i+n}\,w_{in}$$

\noin Like we did in the proof of Theorem $4$ (cf equation $(72)$),
let's look at the action of $\n_1$ on $u$ and the resulting
coefficient in $w_{12}$. This time we get:
\begin{equation}
\frac{m}{r}\,\sum_{j=3}^{n-1}(r^2)^j+\mu\,\frac{m}{r}\,r^{2n}-1=r^4
\end{equation}
After simplifying this expression and replacing $r^{2(n-1)}$ by its
value $-1$, we obtain:
$$\mu\,\frac{m}{r}\,r^{2n}=0$$
Then it comes $\mu=0$, so that $U$ is in fact spanned by
$$u=\sum_{1\leq i<j\leq n-1}\,r^{i+j}\,w_{ij}$$
We note that we recover the fact that $l=\unsur{r^{2n-5}}$ (when
$r^{2(n-1)}=-1$, we can check that $\unsur{r^{2n-5}}=-r^3$). To
conclude, it suffices now to look at the action of $e_{n-1}$ on $u$.
We have for every $i$ with $1\leq i\leq n-1$:
$$e_{n-1}.\,w_{i,n-1}=\unsur{r^{n-i-2}}\,w_{n-1,n}$$
Then, it comes:
$$e_{n-1}.u=r\,\sum_{i=1}^{n-2}(r^2)^i\;w_{n-1,n}=r^3\,\frac{1-(r^2)^{n-2}}{1-r^2}\;w_{n-1,n}$$
Then, $e_{n-1}.u$ is nonzero as $r^{2(n-2)}\neq 1$. Then $U$ is not
contained in $K(n)$: a contradiction. On the way, we recovered the
result from Proposition $9$ that when $n\geq 5$, $l=-r^3$ and
$r^{2(n-1)}=-1$, $K(n-1)\nsubseteq K(n)$, where we used the same
action on the same vector.\\
Now Lemma $13$ holds and we have hence proven that it is impossible
to have $l=-r^3$. Thus, $l=r$ and $(S1)$ holds. Let's finally prove
$(S2)$. If $l=-r^3$, then $\n^{(n)}$ is reducible. Moreover, if
$r^{2n}\neq -1$, we know that $K(n)$ is irreducible by point $7$ and
$k(n)\leq \dbw$ by point $6$. Suppose first $n=7$ or $n\geq 9$. Then
an irreducible $\ih(n)$-module has dimension $1,n-1,\cil,\dbw$ or
dimension greater than $\dbw$. Since $l=-r^3$ and $r^{2n}\neq -1$,
$K(n)$ cannot have dimension $1$. Neither can it have dimension
$n-1$ (as otherwise $l\in\lb\unsur{r^{n-3}},-\unsur{r^{n-3}}\rb$).
Then it must have dimension $\cil$ or $\dbw$. If it had dimension
$\cil$, then we would have $l=r$ by $(S1)$. So we see that $K(n)$
must have dimension $\dbw$. Thus $(S2)$ is proven for $n=7$ or
$n\geq 9$. The case $n=8$ is in fact not different, but needs to be
slightly adapted. Recall that the irreducible representations of
$\ih(8)$ have dimensions $1,7,14,20,21$ or dimensions greater than
or equal to $28$. Still by points $6$ and $7$, we have $K(8)$ is
irreducible with $k(8)\leq 21$. Then $k(8)\in\lb 14,21\rb$ by the
same arguments as in the cases $n=7$ or $n\geq 9$. If $k(8)=14$,
then $\di (K(8)\cap\V^{(7)})\geq 7$. By point number $9$, the
existence of an irreducible $14$-dimensional invariant subspace of
$\V^{(7)}$ must be excluded. Also, by semisimplicity of $\ih(8)$, it
is impossible to have $-r^3\in\lb\unsur{r^4},-\unsur{r^4}\rb$. Then
$K(8)\cap\V^{(7)}$ must have dimension $1+15=16$, which is also the
dimension of $K(7)$ when $l=-r^3$ and $r^{14}=-1$. Then we get
$K(8)\cap\V^{(7)}=K(7)$ and in particular it comes $K(7)\subseteq
K(8)$. This leads to a contradiction as any spanning vector for the
one-dimensional invariant subspace of $\V^{(7)}$ does not belong to
$K(8)$. So we have shown that it is impossible to have $k(8)=14$ and
the only remaining possibility is thus to have $k(8)=21$. This shows
$(S2)$ in the case $n=8$ and $r^{16}\neq -1$ and $r^{14}=-1$.\\\\
To summarize, we have proven that $(S1)$ and $(S2)$ hold when
$r^{2n}\neq -1$ and $r^{2(n-1)}=-1$. We have also proven that $(S1)$
and $(S2)$ hold when $r^{2n}\neq -1$ and $r^{2s}=1$ for some
$\frac{n+1}{2}\leq s<n-1$ (and that is where we used induction).
Then $(\matp_n)$ holds for every $n\geq 5$ and Theorem $9$ and $10$
are thus entirely proven as soon as the Iwahori-Hecke algebra
$\ih(n)$ is assumed to be semisimple.\\

We have the immediate Corollary:
\begin{Cor}
Let $n$ be an integer with $n\geq 3$. Suppose $l=-r^3$.\\
\indent If $r^{2n}\neq -1$, then $k(n)=\dbw$\\
\indent If $r^{2n}=-1$, then $k(n)=1+\dbw$
\end{Cor}

\noin\textsc{Proof:} the second point is point $2)$ of Proposition
$15$ when $n\geq 4$. When $n=3$, we know from Theorem $4$ that there
exists exactly two one-dimensional invariant subspaces of $\V$ and
that their sum is direct. Then $k(3)\geq 2$ and since $k(3)\neq 3$,
in fact $k(3)=2$. As for the first point, Lemma $10$ for $n\geq 5$
yields the inequality $k(n)\leq \dbw$. Since when $l=-r^3$ and
$n\geq 5$ there exists an irreducible $\dbw$-invariant subspace by
Theorem $10$, we also have $k(n)\geq \dbw$, so that $k(n)=\dbw$. It
remains to deal with the cases $n=3$ and $n=4$. When $n=3$, $l=-r^3$
and $l\neq\unsur{r^3}$, there exists a unique one-dimensional
invariant subspace of $\V^{(3)}$. The uniqueness forces $k(3)=1$.
When $n=4$ and $l=-r^3$, there exists an irreducible $3$-dimensional
invariant subspace of $\V^{(4)}$ by Theorem $6$ and there does not
exist any one-dimensional invariant subspace of $\V^{(4)}$ as
$r^8\neq -1$. Then $k(4)\in\lb 3,5\rb$. Moreover, since $l\neq r$ it
is impossible to have $k(4)=5$. Hence $k(4)=3$, as announced.

\section{The Uniqueness Theorems and a Complete Description of the
Invariant Subspace of $\V^{(n)}$ when $l=r$.}

In this part, we prove theorems of uniqueness and describe the
irreducible $\cil$-dimensional invariant subspace of $\V^{(n)}$ when
$l=r$.
% and the irreducible $\dbw$-invariant subspace of $\V^{(n)}$
%when $l=-r^3$. First, we show the following theorem:

\begin{thm}\hfill\\
Let $n$ be an integer with $n\geq 3$. We assume that $\n^{(n)}$ is
reducible and exclude the case when $l=-r^3$ and $r^{2n}= -1$. Then,
there exists a unique non-trivial proper invariant subspace in
$\V^{(n)}$.
\end{thm}

\textsc{Proof:} Assume first $n\geq 5$. Since $\n^{(n)}$ is
reducible, $l$ must take one of the values $\unsur{r^{n-3}},
-\unsur{r^{n-3}}, \unsur{r^{2n-3}}, r \;\text{or}\; -r^3$ by the
Main Theorem. Moreover, with the assumptions $r^{2n}\neq -1$ and
$\ih(n)$ is semisimple, all these values are distinct. When $l=r$ or
$l=-r^3$ and $r^{2n}\neq -1$, one of the results of Proposition $15$
is that $K(n)$ is irreducible. Since any proper invariant subspace
of $\V^{(n)}$ must be contained in $K(n)$ and $K(n)$ is non-trivial
for these values of $l$ and $r$ by Proposition $5$, the vector space
$K(n)$ must in fact be the only non-trivial proper invariant
subspace of $\V^{(n)}$. Moreover, when $l=r$, $K(n)$ is
$\cil$-dimensional by Proposition $10$ and when $l=-r^3$ and
$r^{2n}\neq -1$, $K(n)$ is $\dbw$-dimensional as we just saw in
Corollary $15$. Suppose now that $l=\unsur{r^{2n-3}}$. We know by
Theorem $4$ that there exists a one-dimensional invariant subspace
of $\V^{(n)}$. Further, by Theorem $5$, (resp $9$, resp $10$), there
cannot exist any irreducible $(n-1)$ (resp $\cil$, resp $\dbw$)
dimensional invariant subspace of $\V^{(n)}$. Then in the cases
$n=5$ and $n=6$, this unique one-dimensional invariant subspace of
$\V^{(n)}$ must be the only invariant subspace of $\V^{(n)}$. Let's
now consider an integer $n$ with $n=7$ or $n\geq 9$. Then, if $K(n)$
is not one-dimensional, its dimension must be greater than or equal
to $2+\dbw$. It forces both $K(n)\cap\V^{(n-1)}\neq 0$ and
$K(n)\cap\V^{(n-2)}\neq 0$. Then, $l\in\lb r,-r^3\rb$, which is
impossible. Again, $K(n)$ is the unique invariant subspace of
$\V^{(n)}$ in that case. When $n=8$, if $k(8)\geq 15$, then again,
$k(8)>13$ and $k(8)>7$, so that $l\in\lb r,-r^3\rb$, a
contradiction. Hence we conclude again that $k(8)=1$ in that case.
Finally, suppose that $l\in\lb\unsur{r^{n-3}},-\unsur{r^{n-3}}\rb$.
In these cases, we know that there exists a unique
$(n-1)$-dimensional invariant subspace of $\V^{(n)}$ and there does
not exist any irreducible $1$ or $\cil$ or $\dbw$-dimensional
invariant subspace of $\V^{(n)}$. Again, we will deal with the case
$n=8$ apart. When $n\in\lb 5,6\rb$, as before, we may conclude
immediately. If now $n=7$ or $n\geq 9$, suppose that there exists an
irreducible $s$-dimensional invariant subspace of $\V^{(n)}$ with
$s\geq 1+\dbw$. That would force $l\in\lb r,-r^3\rb$, which is
impossible. Thus, the unique $(n-1)$-dimensional invariant subspace
of $\V^{(n)}$ is the only invariant subspace of $\V^{(n)}$. Finally,
in the case $n=8$, there cannot exist any irreducible
$14$-dimensional invariant subspace of $\V^{(8)}$ for the same
reasons as before. Hence the unique $7$-dimensional invariant
subspace of $\V^{(8)}$ is the only invariant subspace of $\V^{(8)}$.
This ends the proof of Theorem $11$. Let's now describe the unique
invariant subspace of $\V^{(n)}$ when $l=r$. We have the theorem:

\begin{thm}\hfill\\

\noin \underline{Assume $l=r$}.\begin{itemize}\item When $n=4$, the
unique invariant subspace $K(4)$ of $\V^{(4)}$ is spanned by the two
linearly independent vectors:

$$\begin{array}{ccc}
w_1^{(4)}&:=&(w_{14}-\unsurr\,w_{24})+(w_{23}-r\,w_{13})\\
w_2^{(4)}&:=&(w_{24}-\unsurr\,w_{34})+(w_{13}-r\,w_{12})
\end{array}$$

\item When $n\geq 5$, the unique invariant subspace $K(n)$ of $\V^{(n)}$ is
built inductively as a direct sum of the unique invariant subspace
$K(n-1)$ of $\V^{(n-1)}$ and of an $(n-2)$-dimensional vector space
spanned by the vectors:

$$\begin{array}{cccccc}
w_1^{(n)}&:=&w_{1,n}-\unsurr\,w_{2,n}&+&r^{n-4}\,(w_{23}-r\,w_{13})&\\
w_k^{(n)}&:=&w_{k,n}-\unsurr\,w_{k+1,n}&+&r^{n-4}\,(w_{1,k+1}-r\,w_{1,k}),&
2\leq k\leq n-2
\end{array}$$

\end{itemize}
\end{thm}

\noin\textsc{Proof of the Theorem:} the case $n=4$ is contained in
Result $2$. Let's deal with $n\geq 5$. When $n=5$, Claim $1$ of Part
$8.3$ provides us with a spanning set of vectors for the unique
invariant subspace of $\V^{(5)}$. Up to a reordering, we read that
these vectors are $w_1^{(4)}$, $w_2^{(4)}$, the spanning vectors of
the unique invariant subspace of $\V^{(4)}$ and the three vectors
$w_1^{(5)}$, $w_2^{(5)}$ and $w_3^{(5)}$. Hence the Theorem holds in
that case. We will proceed by induction. Given an integer $n$ with
$n\geq 6$, let's assume that the Theorem holds for the integer
$n-1$. So $K(n-1)$ is spanned by all the vectors $w_s^{(t)}$'s of
the Theorem with $4\leq t\leq n-1$ and $1\leq s\leq t-2$. When
$l=r$, we know that $K(n)$ is the unique invariant subspace of
$\V^{(n)}$ and that it contains $K(n-1)$, the unique invariant
subspace of $\V^{(n-1)}$. Moreover, it appears that the $(n-2)$
vectors $w_i^{(n)},\;\; i=1,\dots,n-2$, span an $(n-2)$-dimensional
subspace of $\V^{(n)}$ that is in direct sum with $K(n-1)$. Since we
notice that
$$k(n)=\cil=\frac{(n-1)(n-4)}{2}+(n-2)=k(n-1)+(n-2),$$
it will suffice to show
%that this $(n-2)$-dimensional subspace of
%$\V^{(n)}$ is in fact a subspace of $K(n)$. Hence it remains to show
that the $(n-2)$ vectors $w_i^{(n)},\,i=1,\dots,n-2$, belong to
$K(n)$. To do that, we follow the same steps as in the proof of
Claim $3$. First, let $j$ be an integer with $1\leq j\leq n-1$. We
compute with the tables of Appendix $C$ used with $l=r$:
\begin{equation*}
\forall 1\leq k\leq
n-1,\,[X_{j,n}.w_{k,n}]_{_{w_{j,n}}}=\begin{cases}
r^{k-j}&\text{if $k\neq j$}\\
2&\text{if $k=j$}
\end{cases}
\end{equation*}

\noin It follows immediately that:

$$\begin{array}{cccc} \big\lbrack
X_{j,n}.(w_{k,n}-\unsurr\,w_{k+1,n})\big\rbrack_{_{w_{j,n}}}&=&0&\text{if
$j\not\in\lb
k,k+1\rb$}\\
\big\lbrack X_{k,n}.(w_{k,n}-\unsurr\,w_{k+1,n})\big\rbrack_{_{w_{k,n}}}&=&1&\\
\big\lbrack
X_{k+1,n}.(w_{k,n}-\unsurr\,w_{k+1,n})\big\rbrack_{_{w_{k+1,n}}}&=&-\unsurr&
\end{array}$$

\noin When $l=r$, the actions $(CL)_{(s,t)}$ of Appendix $C$ are all
trivial. Then we have for any integer $k$ with $2\leq k\leq n-2$:
\begin{eqnarray*}
\big[X_{1,n}.(w_{1,k+1}-r\,w_{1,k})\big]_{_{w_{1,n}}}&=&\unsur{r.r^{n-k-2}}-\frac{r}{r.r^{n-k-1}}=0\\
\big[X_{j,n}.(w_{1,k+1}-r\,w_{1,k})\big]_{_{w_{j,n}}}&=&0,\;\begin{array}{l}\text{for
all $j$ such that
$\left\lb\begin{array}{l}2\leq j\leq n-1\\ j\not\in\lb k,k+1\rb\end{array}\right.$}\end{array}\\
\big[X_{k,n}.(w_{1,k+1}-r\,w_{1,k})\big]_{_{w_{k,n}}}&=&-\unsur{r^{n-4}}\\
\big[X_{k+1,n}.(w_{1,k+1}-r\,w_{1,k})\big]_{_{w_{k+1,n}}}&=&\unsur{r^{n-3}}
\end{eqnarray*}

Also, we have,
\begin{eqnarray*}
\big[X_{j,n}.(w_{23}-r\,w_{13})\big]_{_{w_{j,n}}}&=&0,\;\,\text{forall $j\not\in\lb 1,2,3\rb$}\\
\big[X_{1,n}.(w_{23}-r\,w_{13})\big]_{_{w_{1,n}}}&=&-\unsur{r^{n-4}}\\
\big[X_{2,n}.(w_{23}-r\,w_{13})\big]_{_{w_{2,n}}}&=&\unsur{r^{n-3}}\\
\big[X_{3,n}.(w_{23}-r\,w_{13})\big]_{_{w_{3,n}}}&=&\unsur{r^{n-4}}-\frac{r}{r^{n-3}}=0
\end{eqnarray*}

\noin Now it appears that the $(n-1)$ rows of the matrix $T(n)$
annihilate the $(n-2)$ vectors $w_i^{(n)},\,i=1,\dots,n-2$. To
complete the proof, we want to show that the $X_{s,t}$'s with $1\leq
s<t\leq n-1$ all annihilate these vectors. This verification is left
to the reader.\\

We now study in detail the case $l=-r^3$. In Theorem $11$, we
excluded the case when $l=-r^3$ and $r^{2n}=-1$. In this case, we
know that there exists both an irreducible $\dbw$-invariant subspace
and a one-dimensional invariant subspace of $\V^{(n)}$ and we have
seen that $k(n)=1+\dbw$. We will show the uniqueness of the
irreducible $\dbw$-invariant subspace when $n\geq 4$.

\begin{thm}
Let $n$ be an integer with $n\geq 4$. When $l=-r^3$ and $r^{2n}=-1$,
there exists a unique irreducible $\dbw$-invariant subspace of
$\V^{(n)}$.
\end{thm}

\noin\textsc{Proof of the Theorem.} We must show the uniqueness. For
$n=4$, the uniqueness is contained in the proof of Theorem $6$.
Suppose now $n\geq 5$ and let $\W$ be an irreducible
$\dbw$-dimensional invariant subspace of $\V^{(n)}$. When $n\geq 5$,
we have $\dbw>n-1$, so that the intersection $\W\cap\V^{(n-1)}$ is
nontrivial. Then, by irreducibility of $K(n-1)$, we get:
$$\W\cap\V^{(n-1)}=K(n-1)$$
In particular, we have $K(n-1)\subseteq W$. Then there exists an
$\ih(n-1)$-module $S$ that is a summand for $K(n-1)$ in $\W$:
$$\W=K(n-1)\oplus S$$ To prove the uniqueness part in Theorem $13$, it suffices to
show that this summand is unique. Let's analyze the situation.
First, since $k(n-1)=\frac{(n-2)(n-3)}{2}$, we see that $S$ must be
$(n-2)$-dimensional. Then $S$ is either irreducible or is a direct
sum of one-dimensional submodules. We show that the latter case may
not occur. To do so, we slightly adapt the proof of Lemma $13$,
assuming $r^{2n}=-1$ instead of $r^{2(n-1)}=-1$. Under this new
assumption, equation $(169)$ now yields $\mu=1$. Then a
one-dimensional submodule of $S$ must in fact be the unique
one-dimensional invariant subspace of $\V^{(n)}$. Thus, it is
impossible to have $S$ decompose as a direct sum of one-dimensional
submodules. So $S$ is an irreducible $(n-2)$-dimensional
$\ih(n-1)$-module that is contained in $K(n)$. To conclude, it
suffices to prove the following Proposition:
\begin{Prop}
Suppose $l=-r^3$ and $r^{2n}=-1$. In $K(n)$ there exists a unique
irreducible $\ih(n-1)$-module of dimension $(n-2)$.
\end{Prop}
\textsc{Proof of the Proposition:} The existence is provided by the
module $S$ of above. We then adapt the proofs of Theorems $5$ and
$6$ and show the uniqueness of such a module. We proceed step by
step following the same proofs. We leave the cases $n=5$ and $n=7$
for later. Those are special and require more attention. However,
the arguments that we expand below hold in these cases.

First, let $(v_1,v_2,\dots,v_{n-2})$ be a basis of $S$ and suppose
the $v_i$'s satisfy to the relations $(\triangle)$. By the relation
$\n_iv_i=-\unsurr\,v_i$ for every $i=1,\dots,n-2$, equation $(102)$
still holds. By using the key relations $(77)$ and $(78)$ we further
get:
\begin{multline} v_i=\mu^iw_{i,i+1}+\la^i\sum_{s=1}^{i-1}r^{s-1}(w_{s,i}-\unsurr\,w_{s,i+1})\\+
\delta^i\sum_{k=i+2}^{n-1}r^{k-i-2}(w_{i,k}-\unsurr\,w_{i+1,k})+\eta^i(w_{i,n}-\unsurr\,w_{i+1,n})\end{multline}
As before, all the $\delta^i$ are equal, say to $\delta$. Moreover,
similar arguments show that the $\eta^i$ 's are related by the
relations: $$\eta^{i+1}=\unsurr\,\eta^i,\;\;\forall 1\leq i\leq
n-3,$$ and thus are determined by $\eta^1$. Further, $(106)$
becomes:
\begin{equation}\forall 1\leq i\leq n-3,\; r\,\mu^i+r^i\,\la^{i+1}=\delta\end{equation}
Furthermore, $(107)$ becomes:
\begin{equation}\forall 2\leq i\leq n-2,\; r\,\mu^i+r^{i-1}\,\la^i=1\end{equation}

Next, looking at the coefficient of $w_{1,i+1}$ in
$\n_{i+1}.v_i=r(v_i+v_{i+1})$ yields:
$$\forall 2\leq i\leq n-3,\,\la^{i+1}=\unsurr\,\la^i$$
Thus, by $(171)$ and $(172)$, we get $\delta=1$. Now $(172)$ after a
change of indices can be expressed as:
\begin{equation} \forall 1\leq i\leq n-3,\; r\,\mu^{i+1}+r^i\,\la^{i+1}=1\end{equation}
Thus, confronting $(171)$ with $\delta=1$ and $(173)$, we get:
\begin{equation}\forall 1\leq i\leq n-3,\;\mu^{i+1}=\mu^i\end{equation}
Gathering all the informations above, there are only two
coefficients to determine: $\la^2$ and $\eta^1$. By looking at the
coefficient of $w_{12}$ in $\n_1.v_1=-\unsurr\,v_1$, we get:
\begin{equation}
\mu\bigg(\unsurr-\unsur{r^3}\bigg)=m\,\eta^1\,r^{n-4}+\unsur{r^2}+\unsur{r^8}
\end{equation}
Next, since $e_{n-1}.v_{n-2}=0$, we get the equation on the
coefficients:
$$\mu+\la^{n-2}\sum_{s=1}^{n-3}r^{s-1}\bigg(-\unsurr\bigg)\unsur{r^{n-s-2}}+\eta^{n-2}\bigg(-\unsur{r^3}-\unsurr\bigg(-r^2-\unsur{r^2}\bigg)\bigg)=0$$
After replacing $\la^{n-2}$ by $\unsur{r^{n-4}}\,\la^2$ and
$\eta^{n-2}$ by $\unsur{r^{n-3}}\,\eta^1$ and simplifying this
expression, we obtain:
\begin{equation}
\bigg(\frac{r^2+r^6}{r^2-1}\bigg)\mu+\frac{\eta^1}{r^{n-4}}=\frac{1+r^6}{r(r^2-1)}
\end{equation}
Computation of the determinant of the following system obtained from
equations $(175)$ and $(176)$
$$(\star)\left\lb\begin{array}{ccccc}
\bigg(\unsurr-\unsur{r^3}\bigg)\,\mu&-&m\,r^{n-4}\,\eta^1&=&\unsur{r^2}-\unsur{r^8}\\
(r^2+r^6)\,\mu&+&\frac{r^2-1}{r^{n-4}}\,\eta^1&=&\frac{1+r^6}{r}
\end{array}\right.$$
show that it is zero if and only if $r^6=-1$. For $n\geq 6$, the
semisimplicity of $\ih(n)$ prevents it from happening. For $n=5$,
since we assumed $r^{10}=-1$, having $r^6=-1$ would then force
$r^4=1$ which is impossible as $(r^2)^2\neq 1$. Thus, the system
$(\star)$ has a unique couple solution $(\mu,\eta^1)$. Then $\la^2$
is also uniquely determined by $\la^2=\unsurr-\mu$. Thus, all the
coefficients of the $v_i$'s are uniquely determined and to finish
the proof of uniqueness when $n\geq 5$ and $n\not\in\lb 5,7\rb$, it
will suffice to show that it is impossible to have a basis of
vectors $v_1,\dots,v_{n-2}$ of $S$ that satisfy to the set of
relations $(\bigtriangledown)$. We use the proof of Theorem $5$ and
adapt it. Equation $(79)$ giving an expression of $v_1$ becomes:
\begin{equation}
v_1=\sum_{j=1}^{n-3}\mu_{j,n-2}w_{j,n-2}+\sum_{j=1}^{n-2}\mu_{j,n-1}w_{j,n-1}+\mu_{n-2,n}w_{n-2,n}+\mu_{n-1,n}w_{n-1,n}
\end{equation}
Moreover, like before, all the $\mu_{j,k}$ for $j\geq 5$ are zero,
so that when $n\geq 7$, the two "added" terms on the right hand side
of $(177)$ are no longer there. Next, $(80)$ becomes:
\begin{equation}
\mu_{2,n-1}=-\unsurr\,\mu_{2,n-2},
\end{equation}
and $(80)$ and $(81)$ becomes:
\begin{eqnarray}
\mu_{2,n-1}&=&r\,\mu_{1,n-1}\\
\mu_{2,n-2}&=&r\,\mu_{1,n-2}
\end{eqnarray}
Since we assumed $n\geq 6$, by $\n_3.v_1=-\unsurr\,v_1$, we get
$$\mu_{2,n-1}=0$$
Then, when $n\geq 7$, an expression for $v_1$ is:
\begin{equation}
v_1=\mu_{3,n-2}w_{3,n-2}+\mu_{4,n-2}w_{4,n-2}+\mu_{3,n-1}w_{3,n-1}+\mu_{4,n-1}w_{4,n-1}
\end{equation}
and when $n=6$, an expression for $v_1$ is:
\begin{equation}
v_1=\mu_{3,4}w_{3,4}+\mu_{3,5}w_{3,5}+\mu_{4,5}w_{4,5}+\mu_{4,6}w_{4,6}
\end{equation}
Let's first deal with the case $n\geq 7$. Equations $(85)$, $(86)$
and $(87)$ respectively become:
\begin{eqnarray}
\mu_{4,n-1}&=&-\unsurr\,\mu_{3,n-1}\\
\mu_{4,n-2}&=&-\unsurr\,\mu_{3,n-2}\\
\mu_{4,n-1}&=&-\unsurr\,\mu_{4,n-2}
\end{eqnarray}
Like before, we show that $\mu_{3,n-1}=0$. Then $v_1$ is zero: a
contradiction. \\\\
{\underline{\textit{Important partial conclusion}}}: at this stage,
Proposition $17$ and Theorem $13$ are proven for every $n\geq 8$.
\\\\
Let's go back to the case $n=6$ and the expression given in $(182)$.
Since $\mu_{46}=-\unsurr\mu_{36}$ and $\mu_{36}=0$, there is in fact
no term in $w_{46}$. Hence we have:
$$v_1=w_{35}-\unsurr\,w_{45}-r\,w_{34}$$
Again, we have to study $v_2$. We must add a priori to $(94)$ two
more terms, namely $\la_{46}w_{46}$ and $\la_{56}w_{56}$. But with
$\n_2.v_1=-\unsurr\,(v_1+v_2)$, we see that we may as well withdraw
them. The end of the proof is then the same as before, leading to a
contradiction.\\\\
{\underline{\textit{New intermediate conclusion}}}: Proposition $17$ and Theorem $13$ are proven for every $n\geq 4$ and $n\not\in\lb 5,7\rb$.\\\\
\textsc{The case $n=5$}: two questions arise: \\$1)$ Can there exist
a basis $(w_1,w_2,w_3)$ of $S$, where the $w_i$'s satisfy to the set
of relations $(\bigtriangledown)$ ?\\
$2)$ If the answer to the first question is positive, is the linear
span of the $w_i$'s equal to the linear span of the $v_i$'s, where
the $v_i$'s are the vectors arising from the conjugate
representation $(\triangle)$ ?\\\\
The branching rule provides us with a no answer to these questions.
And indeed, when $n=5$, the unique irreducible representation of
$\ih(5)$ of degree $6$ is self-conjugate. Then $\W$ must be
isomorphic to the Specht module $S^{(3,1,1)}$, whose restriction to
$\ih(4)$ decomposes as a direct sum $S^{(3,1)}\oplus S^{(2,1,1)}$ by
the branching rule. But $K(4)$ is irreducible, $3$-dimensional.
Then, by the proof of Theorem $6$, there exists a basis of $K(4)$
whose vectors satisfy to the relations $(\bigtriangledown)$. Since
$S$ is a summand of $K(4)$ in $\W$, it is impossible to have a basis
of $S$, whose vectors satisfy to the same relations
$(\bigtriangledown)$. Thus, the uniqueness part in Theorem $13$
holds in the case $n=5$. \\\\
\textsc{The case $n=7$}: this case is special as there are four
inequivalent irreducible representations of $\ih(6)$ of degree $5$
and not only two. We have already studied two of them. It remains to
check wether $S$ could be isomorphic to $S^{(3,3)}$ or its conjugate
$S^{(2,2,2)}$. If so the restriction of $S$ to $\ih(5)$ would then
be isomorphic $S^{(3,2)}$ or to $S^{(2,2,1)}$. These two irreducible
representations are described in terms of matrix representations in
Fact $1$ of the thesis. We use the proof of Result $1$ to show that
it is impossible to have
$$S\downarrow_{_{\ih(5)}}\simeq
S^{(3,2)}\;\;\text{or}\;\;S\downarrow_{_{\ih(5)}}\simeq
S^{(2,2,1)}$$ First let $(w_1,\dots,w_5)$ be a basis of $S$ such
that the matrices of the left action of the $g_i$'s, $i=1,\dots,4$
in this basis are the $P_i$'s. The $w_i$'s now contain two more
nodes, namely nodes $6$ and $7$. Those don't impact our discussion
from the proof of Result $1$ which remains strictly identical and
leads to $l=r$, a contradiction with $l=-r^3$. Similarly, by
arguments already exposed before, it is impossible to have a basis
of vectors
$(v_1,\dots,v_5)$ of $S$ such that the matrices of the left action of the $g_i$'s, $i=1,\dots,4$ in this basis are the $Q_i$'s.\\\\
This ends the proof of Proposition $17$ and Theorem $13$. A
consequence of Theorem $13$ and its proof is the following result:
\begin{Prop}
Let $n$ be an integer with $n\geq 4$. Suppose $l=-r^3$ and
$r^{2n}=-1$. Then the unique irreducible $\dbw$-dimensional
invariant subspace of $\V^{(n)}$ is isomorphic to the Specht module
$S^{(n-2,1,1)}$.
\end{Prop}
\noin\textsc{Proof of the Proposition}: Name $\W$ the unique
irreducible $\dbw$-dimensional invariant subspace of Theorem $13$.
When $n=4$, the Proposition is proven by the proof of Theorem $6$.
When $n=5$, the only class of irreducible $\ih(5)$-module of degree
$6$ is $S^{(3,1,1)}$, so there is nothing to prove. Suppose $n\geq
6$ and suppose $\W$ is isomorphic to $S^{(3,1^{n-3})}$. Then, by the
branching rule the restriction module $\W\downarrow_{_{\ih(n-1)}}$
must be isomorphic to $S^{(2,1^{n-3})}\oplus S^{(3,1^{n-4})}$. We
recall from above that $\W$ decomposes as:
$$\W=K(n-1)\oplus S,$$ where $S$ is the unique irreducible $\ih(n-1)$-module of $\V^{(n)}$ of dimension $n-2$. Then $S$ must be isomorphic to $S^{(2,1^{n-3})}$. But we have seen that $S$ is isomorphic to $S^{(n-2,1)}$, not $S^{(2,1^{n-3})}$. Hence a contradiction and the fact that $\W$ is rather isomorphic to $S^{(n-2,1,1)}$.

\section{A Proof of the Main Theorem without Maple}

In this part, we proceed without using the results from previous
part. From part $8$, it suffices to show that the Main Theorem holds
for the small values $n\in\lb 3,4,5,6\rb$. In the case $n=3$, the
Main Theorem is proven with Theorem $4$ and Theorem $5$. Indeed, by
Theorem $4$, there exists a one-dimensional invariant subspace of
$\V$ if and only if $l\in\lb\unsur{r^3},-r^3\rb$. By Theorem $5$,
there exists an irreducible two-dimensional invariant subspace of
$\V$ if and only if $l\in\lb 1,-1\rb$. Similarly in the case $n=4$,
we know by Theorem $4$ (resp Theorem $5$, resp Corollary $4$) that
there exists an irreducible one (resp $3$, resp $2$)-dimensional
invariant subspace of $\V$ if and only if $l=\unsur{r^5}$ (resp
$l\in\lb \unsurr,-\unsurr,-r^3\rb$, resp $l=r$). Since the degrees
of the irreducible representations of $\ih(4)$ are $1$, $2$ or $3$,
this proves the Main Theorem in this case. Let's deal with the case
$n=5$. Suppose $\n^{(5)}$ is reducible and let $\W$ be an
irreducible invariant subspace of $\V$. Consider the $F$-vector
space $\W\cap\V_0$, with the same notations as before. If this
intersection is trivial, then the sum $\W+\V_0$ is direct and we
must have $\dw + \di\V_0\leq\di\V$, id est, $\dw\leq 4$. Then
$\dw\in\lb 1,4\rb$ as the degrees of the irreducible representations
of $\ih(5)$ are $1,4,5,6$. By Theorem $4$ and Theorem $5$, this
forces $l\in\lb\unsur{r^7},\unsur{r^2},-\unsur{r^2}\rb$. Suppose
$l\not\in\lb\unsur{r^7},\unsur{r^2},-\unsur{r^2}\rb$. Then we have
$0\subset\W\cap\V_0\subset\V_0$. By the case $n=4$ it follows that
$l\in\lb r,-r^3,\unsur{r^5},-\unsurr,\unsurr\rb$. Moreover, we have:
$$\begin{array}{ccccc}
\di(\W\cap\V_0)&=&\dw+\di\V_0-\di(\W+\V_0)&&\\
&\geq &\nts\nts\dw+\di\V_0-\di\V&=&\dw-4
\end{array}$$
Furthermore, by our assumption on $l$, we have $\dw\in\lb 5,6\rb$.
If $\dw=5$, then by the Result $1$, we know that it forces $l=r$.
From now on we suppose that $l\neq r$. So $\dw=6$ and
$\di(\W\cap\V_0)\geq 2$. Our assumption on $l$ is now $l\not\in\lb
r,\unsur{r^7},\unsur{r^2},-\unsur{r^2}\rb$ and $l\in\lb
-r^3,\unsur{r^5},\unsurr,-\unsurr\rb$. Our goal is to show that
$l=-r^3$. Since $\W\cap\V_0$ is a proper invariant subspace of
$\V_0$, it must be contained in $K(4)$. Therefore, we also have the
inequality $\di(\W\cap\V_0)\leq k(4)$, where we used the notations
of previous section. We will show that it is impossible to have
$l=\unsur{r^5}$ (unless $\unsur{r^5}=-r^3$) or $l=\unsurr$ or
$l=-\unsurr$. If $l=\unsur{r^5}$, then by Theorem $4$, there exists
a unique one-dimensional invariant subspace inside $\V_0$. In
particular $k(4)\neq 0$. Hence $k(4)\in\lb 1,2,3,4,5\rb$. At this
stage, we recall our assumption of semisimplicity for $\ih(5)$. When
$\ih(5)$ is semisimple, a fortiori $\ih(4)$ is semisimple. Suppose
$k(4)=2$. If $K(4)$ were not irreducible, it would contain a
one-dimensional submodule that has a one-dimensional summand by
semisimplicity of $\ih(4)$. This is impossible by uniqueness of the
one-dimensional invariant subspace of $\V_0$. Then $K(4)$ is an
irreducible two dimensional invariant subspace of $\V_0$. By
Corollary $4$, the existence of an irreducible two-dimensional
invariant subspace of $\V_0$ implies $l=r$, a contradiction with
$(r^2)^3\neq 1$. If $k(4)=3$, there are three possibilities. Either
$K(4)$ is irreducible and $l\in\lb -r^3,\unsurr,-\unsurr\rb$ by
Theorem $6$. Or $K(4)$ contains a one-dimensional submodule that has
a two dimensional summand by semisimplicity of $\ih(4)$. Then $l=r$
by uniqueness of the one-dimensional invariant subspace of $\V_0$
and Corollary $4$. Or $K(4)$ contains a two-dimensional submodule
that has a one-dimensional summand. Again, this forces $l=r$ by
uniqueness of the one-dimensional invariant subspace of $\V_0$ and
Corollary $4$. Gathering these results, the only possibility for $l$
that is compatible with $l=\unsur{r^5}$ and our assumption of
semisimplicity for $\ih(5)$ is to have $l=-r^3$. Suppose now
$k(4)=4$. Then $K(4)$ is reducible. If $K(4)$ contains a two-
dimensional invariant subspace, then this two-dimensional subspace
has a two-dimensional summand in $K(4)$. This is impossible by
uniqueness of any one-dimensional (resp two-dimensional) invariant
subspace inside $\V_0$. On the other hand, if $K(4)$ is a direct sum
of an irreducible $3$-dimensional submodule and a one-dimensional
submodule, we must have $l\in\lb -r^3,\unsurr,-\unsurr\rb$ and
$l=\unsur{r^5}$, which leaves the only possibility of
$l=-r^3=\unsur{r^5}$ for the parameters $l$ and $r$. Finally suppose
that $k(4)=5$. The story is similar. $K(4)$ is reducible and is
either a direct sum of a $4$-dimensional submodule and a
one-dimensional submodule or a direct sum of a $3$-dimensional
submodule and a $2$-dimensional submodule, with no further
decompositions allowed by uniqueness of an irreducible
two-dimensional submodule of $\V_0$ or a one-dimensional submodule
of $\V_0$ when these ones exist. Since there does not exist any
irreducible representations of $\ih(4)$ of degree $4$, the first
decomposition should still break, which is impossible. As for the
second decomposition it forces $l=r$, which contradicts our
assumption $l=\unsur{r^5}$. In summary, either $k(4)=1$ or $l=-r^3$.
To reach the goal it suffices to show that it is impossible to have
$k(4)=1$. When $k(4)=1$, the two inequalities above read:

$$\begin{array}{ccc}
\di(\W\cap\V_0)&\geq& 2\\
\di(\W\cap\V_0)&\leq& 1
\end{array},$$

\noindent a contradiction. To achieve our goal it remains to show
that it is impossible to have $l\in\lb\unsurr,-\unsurr\rb$. First we
show that for these values of $l$ we have $k(4)=3$. The scheme of
the proof is the same as in the case $l=\unsur{r^5}$. The existence
of a $3$-dimensional invariant subspace of $\V_0$ shows that
$k(4)\neq 0$. Then, we eliminate turn by turn the possibilities
$k(4)\in\lb 1,2,4,5\rb$. Immediately, if $k(4)=1$, then
$l=\unsur{r^5}$ by Theorem $4$, in contradiction with
$l\in\lb\unsurr,-\unsurr\rb$. Next, if $k(4)=2$, then $K(4)$ is
irreducible by Theorem $4$ and $l=r$ by Corollary $4$. This again
contradicts $l\in\lb\unsurr,-\unsurr\rb$ by semisimplicity of
$\ih(4)$. And in fact those two cases could right away be excluded
by a simple use of Theorem $6$. Indeed, we have seen in Theorem $6$
that when $l\in\lb\unsurr,-\unsurr\rb$, there exists a
$3$-dimensional invariant subspace inside $\V_0$. This implies that
$k(4)\geq 3$. If $k(4)=4$, by uniqueness in Theorem $4$ and in
Result $2$, $K(4)$ must be a direct sum of a $3$-dimensional
submodule and a one-dimensional submodule with no further
decomposition. This forces $l\in\lb -r^3,\unsurr,-\unsurr\rb$ by
Theorem $6$ and $l=\unsur{r^5}$ by Theorem $4$. Since it is
impossible to have $l=\unsur{r^5}$ and $l\in\lb\unsurr,-\unsurr\rb$,
we are led to conclude that $k(4)\neq 4$ when
$l\in\lb\unsurr,-\unsurr\rb$. Finally, if $k(4)=5$, we get $l=r$ by
the same arguments as those already described further above. This is
again impossible with $l\in\lb\unsurr,-\unsurr\rb$. Thus, we have
shown that when $l\in\lb\unsurr,-\unsurr\rb$, we have $k(4)=3$. An
immediate consequence is that $K(4)$ is irreducible. We also derive
from this fact that $\di(\W\cap\V_0)\leq 3$ by one of the two
inequalities of the beginning. Then, for these two values of $l$, we
have $2\leq\di(\W\cap\V_0)\leq 3$. Furthermore,
$l\in\lb\unsurr,-\unsurr\rb$ forces $\di(\W\cap\V_0)=3$. By equality
on the dimensions, we now get that $\W\cap\V_0=K(4)$ and we note on
the way that $\W\cap\V_0$ is irreducible. Suppose first that
$l=\unsurr$. Following Lemma $4$, equation $(40)$, it is an easy
verification using the table of the appendix or directly by hands
that $w_{14}-w_{23}\in K(4)$. Then $w_{14}-w_{23}$ also belongs to
$\W\cap\V_0$ and in particular to $\W$. Since $\W$ is a
$B(A_4)$-module, $e_4.(w_{14}-w_{23})$ must also belong to $\W$.
But,
$$e_4.(w_{14}-w_{23})=\unsur{r^2}\,x_{\al_4}$$
Then it comes $\W=\V$: contradiction.\\
Similarly for $l=-\unsurr$, we have:
$$x_{\al_1}+x_{\al_3}+x_{\al_1+\al_2+\al_3}-x_{\al_2}\in
K(4)=\W\cap\V_0$$ It follows that:
$$e_4.(x_{\al_1}+x_{\al_3}+x_{\al_1+\al_2+\al_3}-x_{\al_2})=\bigp 1+\unsur{r^2}\bigpd\,x_{\al_4}\in\W$$
Since $(r^2)^2\neq 1$ by semisimplicity of $\ih(5)$, this implies
that $x_{\al_4}\in\W$. But again $\W$ would then be the whole space
$\V$ by the arguments of \S$\;$8.1. We conclude that it is
impossible to have $l\in\lb\unsurr,-\unsurr\rb$. Thus, we have shown
that if $\n^{(5)}$ is reducible and $l\not\in\lb
r,\unsur{r^7},\unsur{r^2},-\unsur{r^2}\rb$, then $l=-r^3$. This says
exactly that if $\n^{(5)}$ is reducible, then $l\in\lb
r,-r^3,\unsur{r^7},\unsur{r^2},-\unsur{r^2}\rb$. Conversely, if
$l=\unsur{r^7}$, there exists a one-dimensional invariant subspace
of $\V$ by Theorem $4$, hence the representation is reducible; for
$l\in\lb\unsur{r^2},-\unsur{r^2}\rb$, there exists an irreducible
$4$-dimensional invariant subspace of $\V$ by Theorem $5$ and so the
representation is also reducible in that case. As for $l=r$ (resp
$l=-r^3$), it is a direct verification that the vector $\X$ (resp
$\Y$) of proposition $5$ of \S$\;$ 8.4 belongs to the proper
submodule $K(5)$ of $\V$. This achieves the proof of the Main
Theorem in the case $n=5$. It remains to do the case $n=6$. The
degrees of the irreducible representations of $\ih(6)$ are
$1,5,9,10,16$. The vector space $\V$ is $15$-dimensional. Hence, if
$\W$ is an irreducible submodule, then $\dw\in\lb 1,5,9,10\rb$. If
$\dw=1$ then $l=\unsur{r^9}$ by Theorem $4$. Also, if $\dw=5$ then
$l\in\lb\unsur{r^3},-\unsur{r^3}\rb$ by Theorem $5$. Suppose now
$l\not\in\lb\unsur{r^9},\unsur{r^3},-\unsur{r^3}\rb$. Then
$\W\cap\V_0\neq 0$. If $\W\cap\V_1=0$, then $\dw\leq 9$, which
forces, with the condition on $l$ above, $\dw=9$. Also we note that
$\W\ds\V_1=\V$ as $\dw+\di\V_1=9+6=15=\di\V$. Since $\W\subseteq
K(6)$, we must have $e_5.\W=0$. Also, by definition of $\V_1$, we
have $e_5.\V_1=0$. Then, it follows that $e_5.\V=0$, which is a
contradiction. Thus, if
$l\not\in\lb\unsur{r^9},\unsur{r^3},-\unsur{r^3}\rb$, then both
$\W\cap\V_0$ and $\W\cap\V_1$ are nonzero. By the case $n=4$ and the
case $n=5$ we now get:
$$\left\lbrace\begin{array}{l}
l\in\lb r,-r^3,\unsur{r^5},\unsurr,-\unsurr\rb\\
\qquad\qquad\&\\
l\in\lb r,-r^3,\unsur{r^7},\unsur{r^2},-\unsur{r^2}\rb
\end{array}\right.$$
which implies $l\in\lb r,-r^3\rb$. We conclude that if $\n^{(6)}$ is
reducible and $l\not\in\lb\unsur{r^9},\unsur{r^3},-\unsur{r^3}\rb$,
necessarily $l\in\lb r,-r^3\rb$. In other words, if $\n^{(6)}$ is
reducible then $l\in\lb
r,-r^3,\unsur{r^9},\unsur{r^3},-\unsur{r^3}\rb$. Conversely, for
$l=\unsur{r^9}$ (resp $l=\unsur{r^3}$ or $l=-\unsur{r^3}$), there
exists a one-dimensional (resp $5$-dimensional) invariant subspace
of $\V$ by Theorem $4$ (resp Theorem $5$). Thus, the representation
is reducible in these cases. As for $l=r$ (resp $l=-r^3$), we
already know from the case $n=5$ that $\X$ (resp $\Y$) belongs to
$K(5)$. And as has been seen in the proof of Proposition $5$, we
also have $e_6.\X=X_{46}.\X=X_{36}.\X=X_{26}.\X=X_{16}.\X=0$, with
the same equalities holding for $\Y$ when $l=-r^3$, so that $\X$
(resp $\Y$) belongs to $K(6)$. Which proves the reducibility of the
representation in these cases as well. The Main Theorem is thus
proven in the case $n=6$. And we have seen in Part $8$ that if
$(\P)_5$ and $(\P)_6$ hold, then $(\P)_n$ holds for all $n\geq 7$.
Conversely, for each of the values of $l$ in the set $\lb
\unsur{r^{2n-3}},\unsur{r^{n-3}},-\unsur{r^{n-3}},r,-r^3\rb$, the
representation $\n^{(n)}$ is reducible by Theorem $4$, Theorem $5$
and the proof of Proposition $5$.

\newpage
\appendix
\section{The program}
\begin{verbatim}
               SIZE := proc(n) binomial(n, 2) end proc


  G := proc(k, n)
local i, j, t, g;
    g;
    (l, s) -> array(1 .. SIZE(s), 1 .. SIZE(s));
    if n < 3 then ERROR(`you have entered an invalid value`)
    elif n = 3 then
        if k = 1 then
            g(k, n)[1, 1] := 1/l;
            g(k, n)[1, 2] := 1/r - r;
            g(k, n)[1, 3] := 0;
            g(k, n)[2, 1] := 0;
            g(k, n)[2, 2] := r - 1/r;
            g(k, n)[2, 3] := 1;
            g(k, n)[3, 1] := 0;
            g(k, n)[3, 2] := 1;
            g(k, n)[3, 3] := 0
        elif k = 2 then
            g(k, n)[1, 1] := 0;
            g(k, n)[1, 2] := 0;
            g(k, n)[1, 3] := 1;
            g(k, n)[2, 1] := 0;
            g(k, n)[2, 2] := 1/l;
            g(k, n)[2, 3] := (1/r - r)/l;
            g(k, n)[3, 1] := 1;
            g(k, n)[3, 2] := 0;
            g(k, n)[3, 3] := r - 1/r
        else ERROR(`you entered an invalid first coordinate`)
        end if
    else
        if k = n - 1 then
            for i to binomial(n - 2, 2) do for j to
                binomial(n, 2) do
                    if j <> i then g(k, n)[i, j] := 0
                    else g(k, n)[i, j] := r
                    end if
                end do
            end do;
            for i from 1 + binomial(n - 2, 2) to
            binomial(n - 1, 2) do for j to binomial(n, 2) do
                    if j <> i + n - 1 then g(k, n)[i, j] := 0
                    else g(k, n)[i, j] := 1
                    end if
                end do
            end do;
            for i from 2 + binomial(n - 1, 2) to
            binomial(n, 2) do
                for j from 1 + binomial(n - 1, 2) to
                binomial(n, 2) do
                    if j <> i then g(k, n)[i, j] := 0
                    else g(k, n)[i, j] := r - 1/r
                    end if
                end do;
                for j from 1 + binomial(n - 2, 2) to
                binomial(n - 1, 2) do
                    if j <> i - n + 1 then g(k, n)[i, j] := 0
                    else g(k, n)[i, j] := 1
                    end if
                end do
            end do;
            for j from 2 + binomial(n - 1, 2) to
            binomial(n, 2) do
                t := j - binomial(n - 1, 2) - 2;
                g(k, n)[1 + binomial(n - 1, 2), j] :=
                    (1/r - r)/(l*r^t)
            end do;
            g(k, n)[1 + binomial(n - 1, 2),
                1 + binomial(n - 1, 2)] := 1/l
        elif k < n - 1 then
            for i to binomial(n - 1, 2) do
                t := n - k - 2;
                for j from 1 + binomial(n - 1, 2) to
                binomial(n, 2) do
                    if j <> binomial(n - 1, 2) + n - k - 1
                    then g(k, n)[i, j] := 0
                    else
                        if i = binomial(k, 2) + 1 then
                            g(k, n)[i, j] :=
                            r^(t - 1) - r^(t + 1)
                        else g(k, n)[i, j] := 0
                        end if
                    end if
                end do
            end do;
            for i from 1 + binomial(n - 1, 2) to
            binomial(n, 2) do for j to binomial(n, 2) do
                    if
                    j = i and j <> binomial(n - 1, 2) + n - k
                     and j <> binomial(n - 1, 2) + n - k - 1
                    then g(k, n)[i, j] := r
                    elif i = binomial(n - 1, 2) + n - k - 1
                     and j = binomial(n - 1, 2) + n - k - 1
                    then g(k, n)[i, j] := r - 1/r
                    elif i = binomial(n - 1, 2) + n - k - 1
                     and j = binomial(n - 1, 2) + n - k then
                        g(k, n)[i, j] := 1
                    elif i = binomial(n - 1, 2) + n - k and
                    j = binomial(n - 1, 2) + n - k - 1 then
                        g(k, n)[i, j] := 1
                    else g(k, n)[i, j] := 0
                    end if
                end do
            end do;
            for i to binomial(n - 1, 2) do for j to
                binomial(n - 1, 2) do
                    g(k, n)[i, j] := G(k, n - 1)[i, j]
                end do
            end do
        else ERROR(`this element does not exist`)
        end if
    end if;
    Matrix(SIZE(n), g(k, n))
end proc


  IDENTITY := proc(n)
local identity, i;
    identity := Matrix(SIZE(n), SIZE(n));
    for i to SIZE(n) do identity[i, i] := 1 end do;
    Matrix(identity)
end proc


  e := proc(k, n)
    simplify(evalm(l*(
    `&*`(G(k, n), G(k, n)) + (1/r - r)*G(k, n) - IDENTITY(n))
    /(1/r - r)))
end proc


               E := proc(k, n) Matrix(e(k, n)) end proc


  ginv := proc(k, n)
    simplify(evalm(
    G(k, n) + (1/r - r)*IDENTITY(n) + (r - 1/r)*E(k, n)))
end proc


            GINV := proc(k, n) Matrix(ginv(k, n)) end proc


  T := proc(n)
local X, k, i, j;
    X;
    (s, t) -> Matrix(SIZE(n), SIZE(n));
    for i to n - 1 do for j from i + 1 to n do
            if j = i + 1 then X(i, j) := E(i, n)
            else X(i, j) := Matrix(simplify(evalm(`&*`(
                G(j - 1, n), X(i, j - 1), GINV(j - 1, n)))))
            end if
        end do
    end do;
    evalm(sum('sum('X(i, j)', 'j' = i + 1 .. n)',
        'i' = 1 .. n - 1))
end proc


    d := proc(n) LinearAlgebra:-Determinant(Matrix(T(n))) end proc


          NOTIRR := proc(n) solve({d(n) = 0}, {l}) end proc
\end{verbatim}
\newpage

And running the procedure NOTIRR for $n=3,4,5,6$ yields:

\begin{verbatim}

> NOTIRR(3);

                                                  3         1
     {l = -1}, {l = -1}, {l = 1}, {l = 1}, {l = -r }, {l = ----}
                                                             3
                                                            r

> NOTIRR(4);

  {l = r}, {l = r}, {l = 1/r}, {l = 1/r}, {l = 1/r}, {l = - 1/r},

                                        1            3          3
        {l = - 1/r}, {l = - 1/r}, {l = ----}, {l = -r }, {l = -r },
                                         5
                                        r

               3
        {l = -r }
\end{verbatim}
\newpage
\begin{verbatim}
> NOTIRR(5);

                                                     1
  {l = r}, {l = r}, {l = r}, {l = r}, {l = r}, {l = ----},
                                                      7
                                                     r

                1             1             1             1
        {l = - ----}, {l = - ----}, {l = - ----}, {l = - ----},
                 2             2             2             2
                r             r             r             r

              1           1           1           1            3
        {l = ----}, {l = ----}, {l = ----}, {l = ----}, {l = -r },
               2           2           2           2
              r           r           r           r

               3          3          3          3          3
        {l = -r }, {l = -r }, {l = -r }, {l = -r }, {l = -r }

> NOTIRR(6);

  {l = r}, {l = r}, {l = r}, {l = r}, {l = r}, {l = r}, {l = r},

                                1             1             1
        {l = r}, {l = r}, {l = ----}, {l = - ----}, {l = - ----},
                                 9             3             3
                                r             r             r

                1             1             1           1
        {l = - ----}, {l = - ----}, {l = - ----}, {l = ----},
                 3             3             3           3
                r             r             r           r

              1           1           1           1            3
        {l = ----}, {l = ----}, {l = ----}, {l = ----}, {l = -r },
               3           3           3           3
              r           r           r           r

               3          3          3          3          3
        {l = -r }, {l = -r }, {l = -r }, {l = -r }, {l = -r },

               3          3          3          3
        {l = -r }, {l = -r }, {l = -r }, {l = -r }
\end{verbatim}
\section{Table for $Sym(8)$} \hspace{-0.5in}
\begin{tabular}[c]{|p{.55in}|c|c|p{1in}|c|}
\hline Specht module $\;\;S^{\la}$ &
Ferrers diagram & conjugate Ferrers diagram & conjugate Specht module $\;\;S^{\la^{'}}$ & dimension\\
\hline $S^{(6,2)}$ & $\begin{array}{llll}\\
\begin{tabular}[c]{|p{.05in}|p{.05in}|p{.05in}|p{.05in}|p{.05in}|p{.05in}|}
\hline  &  &  &  &  & \\
\hline
\end{tabular}\\
\begin{tabular}[c]{|p{.05in}|p{.05in}|}
 &  \\ \hline\end{tabular}\\ \begin{array}{l}\end{array}
\end{array}$ & $\begin{array}{llll}\\
\begin{tabular}[c]{|p{.05in}|p{.05in}|p{.05in}|p{.05in}|}
\hline & \\
\hline & \\\hline\end{tabular}\\
\begin{tabular}{|p{.05in}|}
  \\
\hline\\
\hline\\
\hline\\
\hline\end{tabular}\\
\begin{array}{l}\\\end{array}\end{array}$ & $S^{(2,2,1,1,1,1)}$ & $20$\\ \hline
$S^{(6,1,1)}$ &
$\begin{array}{llll}\\\begin{tabular}[c]{|p{.05in}|p{.05in}|p{.05in}|p{.05in}|p{.05in}|p{.05in}|}
\hline & & & & & \\
\hline\end{tabular} \\ \begin{tabular}[c]{|p{.05in}|} \\\hline \\
\hline \end{tabular} \\ \begin{array}{l}\end{array}
\end{array}$ &
$\begin{array}{lll}\\\begin{tabular}[c]{|p{.05in}|p{.05in}|p{.05in}|}
\hline & & \\\hline \\ \cline{1-1} \\ \cline{1-1} \\ \cline{1-1} \\ \cline{1-1} \\
\cline{1-1} \end{tabular}\\ \begin{array}{l}\end{array}\end{array}$
& $S^{(3,1,1,1,1,1)}$ & $21$ \\ \hline $S^{(5,3)}$ & $\begin{array}{lll}\\
\begin{tabular}{|p{.05in}|p{.05in}|p{.05in}|p{.05in}|p{.05in}|}
\hline & & & & \\ \hline & & \\ \cline{1-3}
\end{tabular}\\\begin{array}{l}\end{array}\end{array}$ &
$\begin{array}{lll}\\ \begin{tabular}{|p{.05in}|p{.05in}|}\hline & \\ \hline & \\ \hline & \\
\hline
\\\cline{1-1} \\ \cline{1-1}
\end{tabular}\\\begin{array}{l}\end{array}\end{array}$ &
$S^{(2,2,2,1,1)}$ & $28$
\\\hline $S^{(5,2,1)}$ & $\begin{array}{lll}\\
\begin{tabular}{|p{.05in}|p{.05in}|p{.05in}|p{.05in}|p{.05in}|}
\hline & & & & \\ \hline & \\ \cline{1-2} \\\cline{1-1}
\end{tabular}\\\begin{array}{l}\end{array}\end{array}$ &
$\begin{array}{lll}\\ \begin{tabular}{|p{.05in}|p{.05in}|p{.05in}|}\hline & &\\\hline & \\ \cline{1-2} \\ \cline{1-1} \\
\cline{1-1} \\\cline{1-1}
\end{tabular}\\\begin{array}{l}\end{array}\end{array}$ &
$S^{(3,2,1,1,1)}$ & $64$
\\\hline $S^{(5,1,1,1)}$ & $\begin{array}{lll}\\
\begin{tabular}{|p{.05in}|p{.05in}|p{.05in}|p{.05in}|p{.05in}|}
\hline & & & & \\ \hline \\ \cline{1-1} \\\cline{1-1}\\\cline{1-1}
\end{tabular}\\\begin{array}{l}\end{array}\end{array}$ &
$\begin{array}{lll}\\ \begin{tabular}{|p{.05in}|p{.05in}|p{.05in}|p{.05in}|}\hline & & &\\\hline\\
\cline{1-1} \\
\cline{1-1} \\\cline{1-1}\\\cline{1-1}
\end{tabular}\\\begin{array}{l}\end{array}\end{array}$ &
$S^{(4,1,1,1,1)}$ & $35$
\\\hline
\end{tabular}
\newpage
\hspace{-0.35in}
\begin{tabular}[c]{|p{.55in}|c|c|p{1in}|c|} \hline
$S^{(4,4)}$ &
$\begin{array}{lll}\\\begin{tabular}{|p{.05in}|p{.05in}|p{.05in}|p{.05in}|}
\hline & & &
\\\hline & & & \\\hline\end{tabular}\\\begin{array}{l}\end{array}\end{array}$ &
$\begin{array}{lll}\\ \begin{tabular}{|p{.05in}|p{.05in}|} \hline &
\\\hline &\\\hline & \\\hline &
\\\hline\end{tabular}\\\begin{array}{l}\end{array}\end{array}$ &
$S^{(2,2,2,2)}$ & $14$ \\\hline $S^{(4,3,1)}$ &
$\begin{array}{lll}\\\begin{tabular}{|p{.05in}|p{.05in}|p{.05in}|p{.05in}|}\hline
& & & \\\hline & & \\\cline{1-3}
\\\cline{1-1}\end{tabular}\\\begin{array}{l}\end{array}\end{array}$ &
$\begin{array}{lll}\\\begin{tabular}{|p{.05in}|p{.05in}|p{.05in}|}\hline
& & \\\hline &\\\cline{1-2} & \\\cline{1-2}
\\\cline{1-1}\end{tabular}\\\begin{array}{l}\end{array}\end{array}$
& $S^{(3,2,2,1)}$ & $70$ \\\hline $S^{(4,2,2)}$ &
$\begin{array}{lll}\\\begin{tabular}{|p{.05in}|p{.05in}|p{.05in}|p{.05in}|}\hline
& & &\\\hline & \\\cline{1-2}
&\\\cline{1-2}\end{tabular}\\\begin{array}{l}\end{array}\end{array}$
&$\begin{array}{lll}\\\begin{tabular}{|p{.05in}|p{.05in}|p{.05in}|}
\hline & & \\\hline & & \\\hline \\\cline{1-1}
\\\cline{1-1}\end{tabular}\\\begin{array}{l}\end{array}\end{array}$
& $S^{(3,3,1,1)}$ & $56$\\\hline $S^{(4,2,1,1)}$ &
$\begin{array}{lll}\\\begin{tabular}{|p{.05in}|p{.05in}|p{.05in}|p{.05in}|}\hline
&&&\\\hline &\\\cline{1-2}
\\\cline{1-1}\\\cline{1-1}\end{tabular}\\\begin{array}{l}\end{array}\end{array}$
& self-conjugate & self-conjugate & $90$ \\\hline $S^{(3,3,2)}$ &
$\begin{array}{lll}\\\begin{tabular}{|p{.05in}|p{.05in}|p{.05in}|}\hline
5&4&2\\\hline 4&3&1\\\hline 2&1\\\cline{1-2}
\end{tabular}\\\begin{array}{l}\end{array}\end{array}$ &
self-conjugate & self-conjugate &
$\frac{8!}{5\ti4\ti3\ti2\ti4\ti2}=42$\\\hline
\end{tabular}
\newpage
\section{How the $X_{ij}$'s act} \hspace{-0.7in}
\begin{tabular}[t]{|c||c|c|c|} \hline
Reference & Multiplication on the tangles & index range & algebraic
result\\\hline $\text{(MR)}_k$ &
$\begin{array}{l}\\\epsfig{file=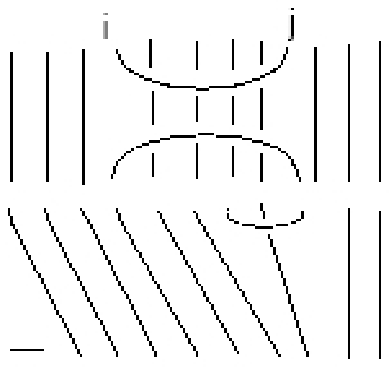,
height=4cm}\end{array}$ & $k=1,\dots,j-i-1$ &
$X_{ij}.\,w_{i+k,j}=l\,r^{k-1}\,w_{ij}$\\\hline $\text{(ML)}_k$ &
$\begin{array}{l}\\\epsfig{file=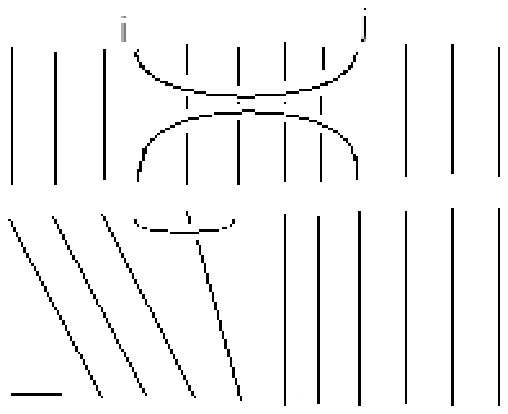, height=4cm}\end{array}$ &
$k=1,\dots,j-i-1$ &
$X_{ij}.\,w_{i,j-k}=\unsur{l\,r^{k-1}}\,w_{ij}$\\\hline
$\text{(TR)}_k$ & $\begin{array}{l}\\\epsfig{file=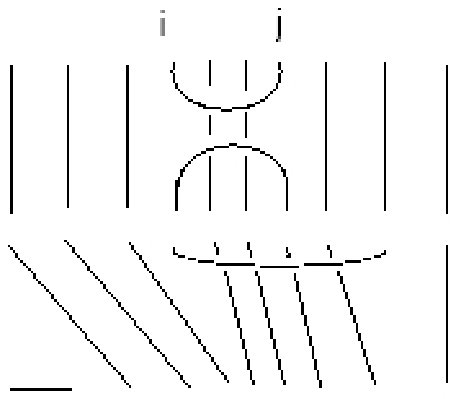,
height=4cm}\end{array}$ & $k=1,\dots,n-j$ &
$X_{ij}.\,w_{i,j+k}=l\,r^{k-1}\,w_{ij}$\\\hline $\text{(TL)}_k$ &
$\begin{array}{l}\\\epsfig{file=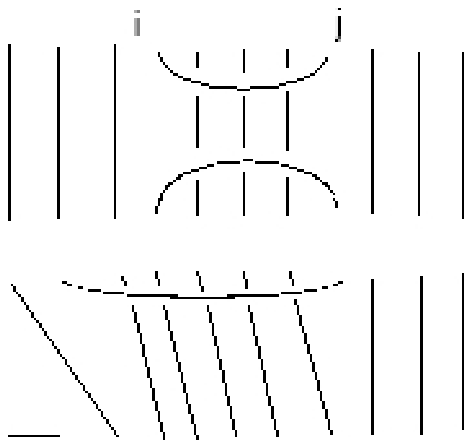, height=4cm}\end{array}$ &
$k=1,\dots,i-1$ &
$X_{ij}.\,w_{i-k,j}=\unsur{l\,r^{k-1}}\,w_{ij}$\\\hline
\end{tabular}
\newpage
\hspace{-1.87in}
\begin{tabular}[t]{|c||c|c|c|} \hline
Ref & Multiplication on the tangles & indices range & algebraic
result\\\hline $\text{(SR)}_k$ &$\begin{array}{l}\\\epsfig{file=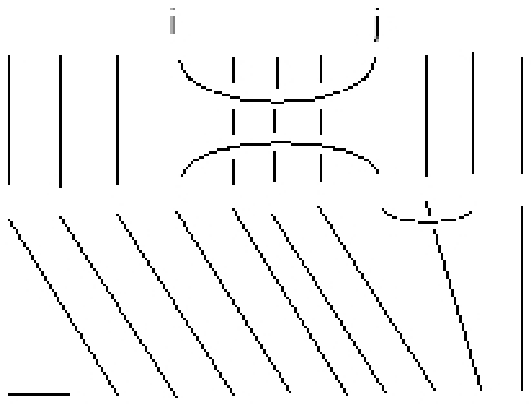,
height=4cm}\end{array}$ & $k=1,\dots,n-j$ &
$X_{ij}.\,w_{j,j+k}=r^{(k-1)+(j-i-1)}\,w_{ij}$\\\hline
$\text{(SL)}_{k}$ &$\begin{array}{l}\\\epsfig{file=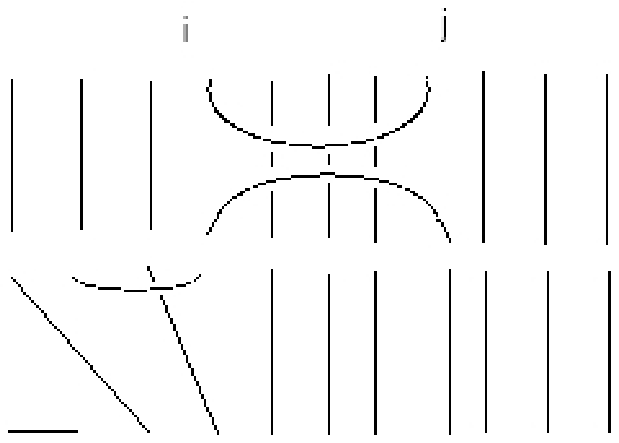,
height=4cm}\end{array}$ & $k=1,\dots,i-1$ &
$X_{ij}.\,w_{i-k,i}=\unsur{r^{(k-1)+(j-i-1)}}\,w_{ij}$\\\hline
$\text{(CR)}_{(s,t)}$ &$\begin{array}{l}\\\epsfig{file=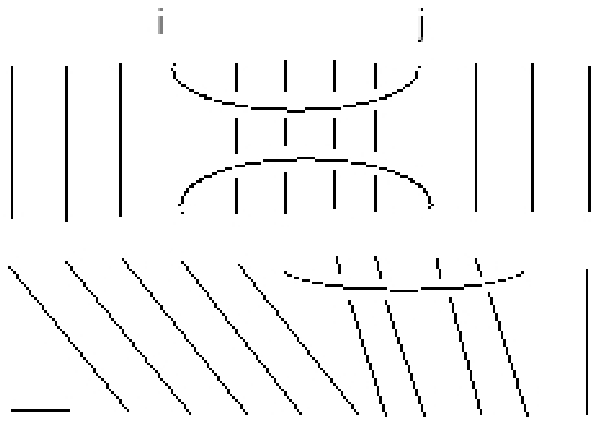,
height=4cm}\end{array}$ &
$\begin{array}{l}s=1,\dots,j-i-1\\t=1,\dots,n-j\end{array}$ &
$X_{ij}.\,w_{i+s,j+t}=(r^{t+s-1}-r^{t+s-3})(l-r)\,w_{ij}$\\\hline
$\text{(CL)}_{(s,t)}$ &$\begin{array}{l}\\\epsfig{file=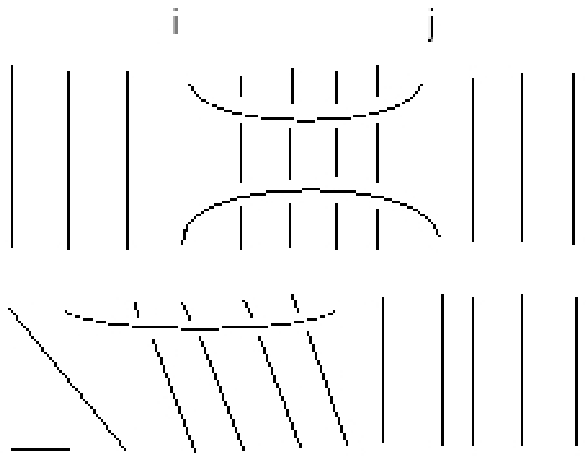,
height=4cm}\end{array}$ & $\begin{array}{l}
s=1,\dots,i-1\\t=1,\dots,j-i-1 \end{array}$ &
$X_{ij}.\,w_{i-s,j-t}=\bigp\unsur{r^{t+s-1}}-\unsur{r^{t+s-3}}\bigpd\bigp\unsur{l}-\unsurr\bigpd\,w_{ij}$\\\hline
\end{tabular}
\newpage
\section{Determinants of square submatrices of size $5$ of $T(5)$}
\begin{verbatim}

l:=r:N:=Matrix(T(5));

g:=proc () local M, d, i1, i2, i3, i4, i5, j1, j2, j3, j4, j5;
    for j1 to 10 do for j2 from j1+1 to 10 do for j3 from j2+1 to 10 do
         for j4 from j3+1 to 10 do for j5 from j4+1 to 10 do
             for i1 to 10 do for i2 from i1+1 to 10 do
                for i3 from i2+1 to 10 do for i4 from i3+1 to 10 do
                    for i5 from i4+1 to 10 do
                         with(linalg);
                         M := submatrix(N,[i1,i2, i3, i4, i5],
                             [j1, j2, j3, j4, j5]);
                         d := LinearAlgebra:-Determinant(
                              Matrix(M));
                         print(j1,j2,j3,j4,j5,i1,i2,i3,i4,i5,d)
                    end do
                  end do
                end do
              end do
            end do
          end do
        end do
      end do
    end do
  end do
end proc

\end{verbatim}
\newpage
\noin And some fragments of output (columns, then rows, then determinant):\\
$$\begin{array}{cccc}\textit{first line of output}&1,\,2,\,3,\,4,\,5,\,1,\,2,\,3,\,4,\,5,\,\mathbf{0}& &\\
&1,\,2,\,3,\,4,\,5,\,1,\,2,\,3,\,4,\,6,\,\mathbf{0}&&\\
&\vdots & \nts\nts\nts\nts\nts\rb\;\;\;\;\text{determinants are all zeros}&\\
\textbf{first}&1,\,2,\,3,\,4,\,7,\,1,\,2,\,3,\,4,\,5,\,\mathbf{0}&&\\
\textbf{$5$-tuple}&1,\,2,\,3,\,4,\,7,\,1,\,2,\,3,\,4,\,6,\,\mathbf{0}&&\\
\textbf{$[1,2,3,4,7]$}&1,\,2,\,3,\,4,\,7,\,1,\,2,\,3,\,4,\,7,\,\mathbf{\frac{2\,r^2+r^4+1}{r^2}}&\nts\nts\nts\leftarrow\,\textbf{first nonzero determinant}&\\
\textbf{is first}&1,\,2,\,3,\,4,\,7,\,1,\,2,\,3,\,4,\,8,\,\frac{2\,r^2+r^4+1}{r}&&\\
\textbf{"non-admissible"}&1,\,2,\,3,\,4,\,7,\,1,\,2,\,3,\,4,\,9,\,2\,r^2+r^4+1&&\\
\textbf{pattern}&1,\,2,\,3,\,4,\,7,\,1,\,2,\,3,\,4,\,10,\,r(2\,r^2+r^4+1)&&\\
\textbf{of columns}&1,\,2,\,3,\,4,\,7,\,1,\,2,\,3,\,5,\,6,\,0&
&\\&\vdots&\nts\nts\nts\nts\nts\rb\;\;\;\; \text{whatever
determinant}&\\
&5,\,6,\,7,\,8,\,10,\,5,\,6,\,7,\,9,\,10,\,0&&\\
&5,\,6,\,7,\,8,\,10,\,5,\,6,\,8,\,9,\,10,\,0&&\\
&5,\,6,\,7,\,8,\,10,\,5,\,7,\,8,\,9,\,10,\,-\frac{r^4+2\,r^2+1}{r^4}&&\\
&5,\,6,\,7,\,8,\,10,\,6,\,7,\,8,\,9,\,10,\,-\frac{r^4+2\,r^2+1}{r^3}&&\\
\textbf{last} &5,\,6,\,7,\,9,\,10,\,1,\,2,\,3,\,4,\,5,\,\mathbf{0}&&\\
\textbf{"admissible"}&5,\,6,\,7,\,9,\,10,\,1,\,2,\,3,\,4,\,6,\,\mathbf{0}
&&\\ \textbf{$5$-tuple}&\vdots&\nts\rb\;\;\;\;
\textbf{determinants are all zeros}&\\
\textbf{of columns}&5,\,6,\,7,\,9,\,10,\,6,\,7,\,8,\,9,\,10,\,\mathbf{0}&&\\
&\vdots &\nts\nts\nts\nts\nts\rb\;\;\;\;\text{whatever determinant}&\\
\textit{last line of
output}&6,\,7,\,8,\,9,\,10,\,6,\,7,\,8,\,9,\,10,\,\frac{r^4+2\,r^2+1}{r^2}\\
\end{array}$$
\begin{verbatim}




> with(linalg):l:=r:Determinant(Matrix(
                       submatrix(T(6),[1,2,3,4,7],[1,2,3,4,12])));

                                  0
\end{verbatim}
\newpage

\section{Reducibility of the L-K representation for $n=3$\\}
\begin{tabular}[c]{|c|c|c|c|}
\hline
$\begin{array}{l}\text{Ferrers}\\\text{diagrams}\\
\text{for $n=3$}\end{array}$
&$\begin{array}{l}\begin{array}{l}\end{array}\\\begin{tabular}[c]{|c|c|c|c|}\hline
& &
\\\hline\end{tabular}\end{array}$&$\begin{array}{l}\begin{array}{l}\end{array}\\\begin{tabular}[c]{|c|c|c|}\hline
& \\\hline
\\\cline{1-1}\end{tabular}\end{array}$& $\begin{array}{l}\begin{array}{l}\end{array}\\\begin{tabular}[c]{|c|}
\hline\\\hline\\\hline\\\hline
\end{tabular}\\\begin{array}{l}\end{array}\end{array}$\\
\hline $\begin{array}{l}\\ \text{Reducibility}\\
\;\;\;\text{of}\\\text{the L-K} \\\text{representation}\\\text{of
degree $3$}\\\begin{array}{l}\end{array}\end{array}$&
$l=\unsur{r^3}$&$l\in\lb 1,-1\rb$&$l=-r^3$\\\cline{1-4}
$\begin{array}{l}\\\text{Dimension}
\\\text{of the}\\
\text{corresponding}\\\text{invariant}\\\text{subspace}\\\text{of
$\V^{(3)}$}\\\text{ }\end{array}$&$1$&$2$&$1$\\\hline
$\begin{array}{l}\text{Spanning}\\\text{vectors}\end{array}$&
$\begin{array}{l}w_{12}\\\;\,+\\r\,w_{13}\\\;\,+\\r^2\,w_{23}\end{array}$ & $\begin{array}{l} l=1:\\
\left\lb\begin{array}{l}
\big(\unsurr-1\big)\,w_{12}+\big(w_{13}-\unsurr\,w_{23}\big)\\
\big(\unsurr-1\big)\,w_{23}+\big(w_{12}-\unsurr\,w_{13}\big)
\end{array}\right.\\
l=-1:\\
\left\lb\begin{array}{l}
\big(\unsurr+1\big)\,w_{12}+\big(w_{13}-\unsurr\,w_{23}\big)\\
\big(\unsurr+1)\,w_{23}-\big(w_{13}-\unsurr\,w_{13}\big)
\end{array}\right.
\end{array}$&
$\begin{array}{l}w_{12}\\\;\,-\\\unsurr\,w_{13}\\\;\,+\\\unsur{r^2}\,w_{23}\end{array}$\\
\hline
\end{tabular}
\newpage
\section{Reducibility of the L-K representation for $n=4$\\}
\hspace{-1.45in}
\begin{tabular}[c]{|c|c|c|c|c|c|}
\hline
$\begin{array}{l}\text{F.D.}\\
\text{for}\\n=4\end{array}$
&$\begin{array}{l}\begin{array}{l}\end{array}\\\begin{tabular}[c]{|c|c|c|c|c|}\hline
& & &
\\\hline\end{tabular}\end{array}$&$\begin{array}{l}\begin{array}{l}\end{array}\\\begin{tabular}[c]{|c|c|c|c|}\hline
& &\\\hline
\\\cline{1-1}\end{tabular}\end{array}$&  $\begin{array}{l}\begin{array}{l}\end{array}\\
\begin{tabular}[c]{|c|c|}\hline
&\\\hline&\\\hline\end{tabular}\end{array}$&
$\begin{array}{l}\begin{array}{l}\end{array}\\\begin{tabular}[c]{|c|c|}
\hline&\\\hline \\\cline{1-1}\\\cline{1-1}\end{tabular}\end{array}$
&$\begin{array}{l}\begin{array}{l}\end{array}\\\begin{tabular}[c]{|c|}
\hline\\\hline\\\hline\\\hline\\\hline
\end{tabular}\\\begin{array}{l}\end{array}\end{array}$\\
\hline $\begin{array}{l}\\ \text{Red.}\\
\;\;\;\text{of}\\\text{the L-K} \\\text{rep.}\\\text{of}\\\text{deg
$6$}\\\begin{array}{l}\end{array}\end{array}$&
$l=\unsur{r^5}$&$l\in\lb
\unsurr,-\unsurr\rb$&$l=r$&$l=-r^3$&$\begin{array}{l}\text{cannot}\\\text{occur}\end{array}$\\\hline
$\begin{array}{l}\\\text{Dim.}
\\\text{of the}\\
\text{inv.}\\\text{subspace}\\\text{of $\V^{(4)}$}\\\text{
}\end{array}$&$1$&$3$&$2$&$3$\\\cline{1-5}
$\begin{array}{l}\text{Spann.}\\\text{vect.}\end{array}$&
$\sum_{1\leq i<j\leq 4}r^{i+j}\,w_{ij}$&$\begin{array}{l}
\bullet\;(\unsurr-\e_l\,r)w_{12}\\
\;\;+(w_{13}-\unsurr\,w_{23})\\
\;\;+\e_l\,r(w_{14}-\unsurr\,w_{24})\\
\bullet\;(\unsurr-\e_l\,r)\,w_{23}\\
\;\;+(w_{24}-\e_l\unsurr\,w_{34})\\
\;\;+\e_l\,r(w_{12}-\unsurr\,w_{13})\\
\bullet\;(\unsurr-\e_l\,r)\,w_{34}\\
\;\;+(w_{13}-\unsurr\,w_{14})\\
\;\;+\e_l\,r(w_{23}-\unsurr\,w_{24})\\\\
\e_l=\begin{cases} 1 & \text{if $l=\unsurr$,}\\ -1 & \text{if
$l=-\unsurr$} \end{cases}
\end{array}$&$\begin{array}{l}\bullet\; (w_{13}-\unsurr\,w_{23})\\\;\;-\unsurr(w_{14}-\unsurr\,w_{24})\\
\bullet\;(w_{12}-\unsurr\,w_{13})\\
\;\;-\unsurr(w_{24}-\unsurr\,w_{34})\end{array}$&$\begin{array}{l}\bullet\;(\unsurr+\unsur{r^3})\,w_{34}\\
\;\;+(w_{13}-\unsurr\,w_{14})\\
\;\;+r(w_{23}-\unsurr\,w_{24}\\
\bullet\;(r+\unsurr)\,w_{14}\\\;\; -r\,w_{12}-r^2\,w_{13}\\
-\unsurr\,w_{34}-\unsur{r^2}\,w_{24}\\
\bullet\; (r+\unsur{r^3})\,w_{12}\\
\;\;+w_{24}+\unsurr\,w_{23}\\
\;\;-w_{13}-r\,w_{14}\end{array}$\\ \cline{1-5}
\end{tabular}
\section{Reducibility of the L-K representation for $n=5$}
\hspace{-1.4in}
\begin{tabular}[c]{|c|c|c|c|c|}
\hline
$\begin{array}{l}\text{F.D.}\\
\text{for}\\n=5\end{array}$
&$\begin{array}{l}\qquad\ccc\\\begin{array}{l}\end{array}\\\begin{tabular}[c]{|c|c|c|c|c|}\hline
& & & &
\\\hline\end{tabular}\end{array}$&$\begin{array}{l}\qquad\ccc\\\begin{array}{l}\end{array}\\\begin{tabular}[c]{|c|c|c|c|}
\hline & & &\\\hline
\\\cline{1-1}\end{tabular}\\\begin{array}{l}\end{array}\end{array}$&  $\begin{array}{l}\;\;\;\,\ccc\\\begin{array}{l}
\end{array}\\
\begin{tabular}[c]{|c|c|c|}\hline & &\\\hline
&\\\cline{1-2}\end{tabular}\\\begin{array}{l}\end{array}\end{array}$&
$\begin{array}{l}\begin{array}{l}\end{array}\\\begin{tabular}[c]{|c|c|c|}
\hline& &\\\hline
\\\cline{1-1}\\\cline{1-1}\end{tabular}\\\begin{array}{l}\end{array}\end{array}$
\\\hline $\begin{array}{l}\\ \text{Red.}\\
\;\;\;\text{of}\\\text{the L-K} \\\text{rep.}\\\text{of}\\\text{deg
$10$}\\\begin{array}{l}\end{array}\end{array}$&
$l=\unsur{r^7}$&$l\in\lb \unsur{r^2},-\unsur{r^2}\rb$&$l=r$&$l=-r^3$
\\\hline
$\begin{array}{l}\\\text{Dim.}
\\\text{of the}\\
\text{inv.}\\\text{subs.}\\\text{of $\V^{(5)}$}\\\text{
}\end{array}$&$1$&$4$&$5$&$6$\\\hline
$\begin{array}{l}\text{Spann.}\\\text{vect.}\end{array}$&
$\sum_{1\leq i<j\leq 5}r^{i+j}\,w_{ij}$& $\begin{array}{l}
\text{For $i=1,\dots,4$}\\
(\unsurr-\e_l\,r^2)\,w_{i,i+1}\\
+\sum_{k=i+2}^5\,r^{k-i-2}(w_{i,k}-\unsurr\,w_{i+1,k})\\
+\e_l\,\sum_{s=1}^{i-1}\,r^{s-i+3}\,(w_{s,i}-\unsurr\,w_{s,i+1})\\
\e_l=\begin{cases} 1 & \text{if $l=\unsur{r^2}$,}\\ -1 & \text{if
$l=-\unsur{r^2}$} \end{cases}\end{array}$& $\begin{array}{l}
\bullet\;
w_1^{(4)},\,w_2^{(4)}\\
\bullet\; w_1^{(5)},\,w_2^{(5)},\,w_3^{(5)}
\end{array}$ & $\begin{array}{l}
\bullet\;v_1^{(4)},v_2^{(4)},v_3^{(4)}\\
\bullet\;v_1^{(5)},v_2^{(5)},v_3^{(5)}\end{array}$
\\\hline
\end{tabular}
\\\\\\
$\bullet$\;\, A circled circ $\ccc$ on top of the Ferrers diagram
indicates that the conjugate irreducible representation cannot occur
in in the
Lawrence-Krammer representation. \\
$\bullet$ \;\, The vectors $w_1^{(4)}$, $w_2^{(4)}$ are the spanning
vectors of the unique irreducible $2$-dimensional invariant subspace
of $\V^{(4)}$.\\
\indent\hspace{-0.04in} The vectors $w_1^{(5)}$, $w_2^{(5)}$ and
$w_3^{(5)}$ are given by: $\left\lb\begin{array}{cc}
w_1^{(5)}:=&w_{15}-\unsurr\,w_{25}+r\,(w_{23}-r\,w_{13})\\
w_2^{(5)}:=& w_{25}-\unsurr\,w_{35}+r\,(w_{13}-r\,w_{12})\\
w_3^{(5)}:=& w_{35}-\unsurr\,w_{45}+r\,(w_{14}-r\,w_{13})\\
\end{array}\right.$\\\\
$\bullet$\;\, The vectors $v_1^{(4)},v_2^{(4)},v_3^{(4)}$ are the
spanning vectors of the irreducible $3$-dimensional invariant
subspace of $\V^{(4)}$ of the previous table.\\
\indent \hspace{-0.04in}The vectors $v_1^{(5)},v_2^{(5)},v_3^{(5)}$
are defined
by:\\
$v_k^{(5)}=w_{k+1,5}-r\,w_{k,5}+r^{5-k}\,w_{k,k+1}$, $\;k=1,\dots,3$

\section{Reducibility of the L-K representation for $n=6$}
\hspace{-1.53in}
\begin{tabular}[c]{|c|c|c|c|c|c|}
\hline
$\begin{array}{l}\text{F.D.}\\
\text{for}\\n=6\end{array}$
&$\begin{array}{l}\qquad\,\,\ccc\\\begin{array}{l}\end{array}\\\begin{tabular}[c]{|c|c|c|c|c|c|}\hline
& & & & &
\\\hline\end{tabular}\end{array}$&$\begin{array}{l}\qquad\,\ccc\\\begin{array}{l}\end{array}\\\begin{tabular}[c]{|c|c|c|c|c|}
\hline & & & &\\\hline
\\\cline{1-1}\end{tabular}\\\begin{array}{l}\end{array}\end{array}$&  $\begin{array}{l}\qquad\ccc\\\begin{array}{l}
\end{array}\\
\begin{tabular}[c]{|c|c|c|c|}\hline & & &\\\hline
&\\\cline{1-2}\end{tabular}\\\begin{array}{l}\end{array}\end{array}$&
$\begin{array}{l}\;\;\,\;\;\;\,\ccc\\\begin{array}{l}\end{array}\\\begin{tabular}[c]{|c|c|c|c|}
\hline& & &\\\hline
\\\cline{1-1}\\\cline{1-1}\end{tabular}\\\begin{array}{l}\end{array}\end{array}$&$\begin{array}{l}\;\;\;\;\,\ccc\\\begin{array}{l}\end{array}\\
\begin{tabular}[c]{|c|c|c|}\hline & &\\\hline & &\\\hline\end{tabular}\\\begin{array}{l}\end{array}\end{array}$
\\\hline $\begin{array}{l}\\ \text{Red.}\\
\;\;\;\text{of}\\\text{the L-K} \\\text{rep.}\\\text{of}\\\text{deg
$15$}\\\begin{array}{l}\end{array}\end{array}$&
$l=\unsur{r^9}$&$l\in\lb
\unsur{r^3},-\unsur{r^3}\rb$&$l=r$&$l=-r^3$&$\begin{array}{l}\text{cannot}\\\text{occur}\end{array}$
\\\hline
$\begin{array}{l}\\\text{Dim.}
\\\text{of the}\\
\text{inv.}\\\text{subs.}\\\text{of $\V^{(6)}$}\\\text{
}\end{array}$&$1$&$5$&$9$&$10$\\\cline{1-5}
$\begin{array}{l}\text{Spann.}\\\text{vect.}\end{array}$&
$\sum_{1\leq i<j\leq 6}r^{i+j}\,w_{ij}$& $\begin{array}{l}
\text{For $i=1,\dots,5$}\\
(\unsurr-\e_l\,r^3)\,w_{i,i+1}\\
+\sum_{k=i+2}^6\,r^{k-i-2}(w_{i,k}-\unsurr\,w_{i+1,k})\\
+\e_l\,\sum_{s=1}^{i-1}\,r^{s-i+4}\,(w_{s,i}-\unsurr\,w_{s,i+1})\\
\e_l=\begin{cases} 1 & \text{if $l=\unsur{r^3}$,}\\ -1 & \text{if
$l=-\unsur{r^3}$} \end{cases}\end{array}$& $\begin{array}{l} \cdot\;
\W^{(5)}\\
\cdot\;
w_1^{(6)}\\\cdot\;w_2^{(6)}\\\cdot\;w_3^{(6)}\\\cdot\;w_4^{(6)}
\end{array}$&
$\begin{array}{l}\begin{array}{l}\end{array}
\cdot\;v_s^{(t)},\left\lb\begin{array}{l}t\in\lb 4,5\rb\\1\leq s\leq 3\end{array}\right.\\
\cdot\;v_1^{(6)}\\
\cdot\;v_2^{(6)}\\
\cdot\;v_3^{(6)}\\
\cdot\;v_4^{(6)}
\end{array}$
\\\cline{1-5}
\end{tabular}
\\\\\\
$\bullet$\;\, A circled circ $\ccc$ on top of the Ferrers diagram
indicates that the conjugate irreducible representation cannot occur
in in the
Lawrence-Krammer representation. \\
$\bullet$ \;\, $\W^{(5)}=\lb
w_1^{(4)},\,w_2^{(4)},w_1^{(5)},w_2^{(5)},w_3^{(5)}\rb$, all the
vectors defined in Appendices $F$ and $G$.\\
$\bullet$\;\,$\begin{array}{cccccc} &w_{1}^{(6)}&=&w_{16}-\unsurr\,w_{26}&+&r^2\,(w_{23}-r\,w_{13})\\
\forall\;2\leq k\leq
4,&w_k^{(6)}&=&w_{k,n}-\unsurr\,w_{k+1,n}&+&r^2\,(w_{1,k+1}-r\,w_{1,k})\end{array}$\\
$\bullet$\;\, The vectors $v_s^{(t)}$ with $t\in\lb 4,5\rb$ and
$s\in\lb 1,2,3\rb$ are defined in Appendices $F$ and $G$.\\
$\bullet$\;\, The $v_k^{(6)}$, $k=1,\dots,4$ are given by:\\
\indent$v_k^{(6)}:=w_{k+1,6}-r\,w_{k,6}+r^{6-k}\,w_{k,k+1}$

\newpage
\begin{center}\Huge{Index}\end{center}

\textbf{\large{Index of the Theorems}}\\\\
Theorem $1$ \dotfill $15$\\
Theorem $2$ \dotfill $25$\\
Theorem $3$ \dotfill $44$\\
Theorem $4$ \dotfill $48$\\
Theorem $5$ \dotfill $52$\\
Theorem $6$ \dotfill $52$\\
Theorem $7$ \dotfill $86$\\
Theorem $8$ \dotfill $118$\\
Theorem $9$ \dotfill $133$\\
Theorem $10$ \dotfill $133$\\
Theorem $11$ \dotfill $142$\\
Theorem $12$ \dotfill $143$\\
Theorem $13$ \dotfill $145$\\

\textbf{\large{Index of the Propositions}}\\\\
Proposition $1$ \dotfill $23$\\
Proposition $2$ \dotfill $47$\\
Proposition $3$ \dotfill $86$\\
Proposition $4$ \dotfill $86$\\
Proposition $5$ \dotfill $95$\\
Proposition $6$ \dotfill $99$\\
Proposition $7$ \dotfill $107$\\
Proposition $8$ \dotfill $108$\\
Proposition $9$ \dotfill $108$\\
Proposition $10$ \dotfill $111$\\
Proposition $11$ \dotfill $119$\\
Proposition $12$ \dotfill $127$\\
Proposition $13$ \dotfill $130$\\
Proposition $14$ \dotfill $131$\\
Proposition $15$ \dotfill $131$\\
Proposition $16$ \dotfill $137$\\
Proposition $17$ \dotfill $145$\\
Proposition $18$ \dotfill $148$\\

\textbf{\large{Index of the Lemmas}}\\\\
Lemma $1$ \dotfill $10$\\
Lemma $2$ \dotfill $11$\\
Lemma $3$ \dotfill $15$\\
Lemma $4$ \dotfill $23$\\
Lemma $5$ \dotfill $54$\\
Lemma $6$ \dotfill $92$\\
Lemma $7$ \dotfill $96$\\
Lemma $8$ \dotfill $109$\\
Lemma $9$ \dotfill $126$\\
Lemma $10$ \dotfill $127$\\
Lemma $11$ \dotfill $136$\\
Lemma $12$ \dotfill $137$\\
Lemma $13$ \dotfill $139$\\

\textbf{\large{Index of the Corollaries}}\\\\
Corollary $1$ \dotfill $11$\\
Corollary $2$ \dotfill $13$\\
Corollary $3$ \dotfill $86$\\
Corollary $4$ \dotfill $103$\\
Corollary $5$ \dotfill $104$\\
Corollary $6$ \dotfill $117$\\
Corollary $7$ \dotfill $117$\\
Corollary $8$ \dotfill $119$\\
Corollary $9$ \dotfill $122$\\
Corollary $10$ \dotfill $123$\\
Corollary $11$ \dotfill $124$\\
Corollary $12$ \dotfill $125$\\
Corollary $13$ \dotfill $130$\\
Corollary $14$ \dotfill $133$\\
Corollary $15$ \dotfill $141$\\

\textbf{\large{Index of the Results}}\\\\
Result $1$ \dotfill $79$\\
Result $2$ \dotfill $101$\\

\textbf{\large{Index of the Claims}}\\\\
Claim $1$ \dotfill $82$\\
Claim $2$ \dotfill $109$\\
Claim $3$ \dotfill $120$\\

\textbf{\large{Index of the Facts}}\\\\
Fact $1$ \dotfill $77$\\
Fact $2$ \dotfill $116$\\

\end{document}